\catcode`\@=11
\font\tensmc=cmcsc10      
\def\smc{\tensmc}

\def\hcorrection#1{\advance\hoffset by #1 }
\def\vcorrection#1{\advance\voffset by #1 }
\def\wlog#1{}
\newif\iftitle@
\outer\def\title{\title@true\vglue 24\p@ plus 12\p@ minus 12\p@
   \bgroup\let\\=\cr\tabskip\centering
   \halign to \hsize\bgroup\tenbf\hfill\ignorespaces##\unskip\hfill\cr}
\def\endtitle{\cr\egroup\egroup\vglue 18\p@ plus 12\p@ minus 6\p@}
\outer\def\author{\iftitle@\vglue -18\p@ plus -12\p@ minus -6\p@\fi\vglue
    12\p@ plus 6\p@ minus 3\p@\bgroup\let\\=\cr\tabskip\centering
    \halign to \hsize\bgroup\smc\hfill\ignorespaces##\unskip\hfill\cr}
\def\endauthor{\cr\egroup\egroup\vglue 18\p@ plus 12\p@ minus 6\p@}
\outer\def\heading{\bigbreak\bgroup\let\\=\cr\tabskip\centering
    \halign to \hsize\bgroup\smc\hfill\ignorespaces##\unskip\hfill\cr}
\def\endheading{\cr\egroup\egroup\nobreak\medskip}
\outer\def\subheading#1{\medbreak\noindent{\tenbf\ignorespaces
      #1\unskip.\enspace}\ignorespaces}
\outer\def\proclaim#1{\medbreak\noindent\smc\ignorespaces
    #1\unskip.\enspace\sl\ignorespaces}
\outer\def\endproclaim{\par\ifdim\lastskip<\medskipamount\removelastskip
  \penalty 55 \fi\medskip\rm}
\outer\def\demo#1{\par\ifdim\lastskip<\smallskipamount\removelastskip
    \smallskip\fi\noindent{\smc\ignorespaces#1\unskip:\enspace}\rm
      \ignorespaces}
\outer\def\enddemo{\par\smallskip}
\newcount\footmarkcount@
\footmarkcount@=1
\def\makefootnote@#1#2{\insert\footins{\interlinepenalty=100
  \splittopskip=\ht\strutbox \splitmaxdepth=\dp\strutbox 
  \floatingpenalty=\@MM
  \leftskip=\z@\rightskip=\z@\spaceskip=\z@\xspaceskip=\z@
  \noindent{#1}\footstrut\rm\ignorespaces #2\strut}}
\def\footnote{\let\@sf=\empty\ifhmode\edef\@sf{\spacefactor
   =\the\spacefactor}\/\fi\futurelet\next\footnote@}
\def\footnote@{\ifx"\next\let\next\footnote@@\else
    \let\next\footnote@@@\fi\next}
\def\footnote@@"#1"#2{#1\@sf\relax\makefootnote@{#1}{#2}}
\def\footnote@@@#1{$^{\number\footmarkcount@}$\makefootnote@
   {$^{\number\footmarkcount@}$}{#1}\global\advance\footmarkcount@ by 1 }

\hyphenation{man-u-script man-u-scripts ap-pen-dix ap-pen-di-ces}
\hyphenation{data-base data-bases}
\ifx\amstexloaded@\relax\catcode`\@=13 
  \endinput\else\let\amstexloaded@=\relax\fi
\newlinechar=`\^^J
\def\eat@#1{}
\def\Space@.{\futurelet\Space@\relax}
\Space@. %
\newhelp\athelp@
{Only certain combinations beginning with @ make sense to me.^^J
Perhaps you wanted \string\@\space for a printed @?^^J
I've ignored the character or group after @.}
\def\futureletnextat@{\futurelet\next\at@}
{\catcode`\@=\active
\lccode`\Z=`\@ \lowercase
{\gdef@{\expandafter\csname futureletnextatZ\endcsname}
\expandafter\gdef\csname atZ\endcsname
   {\ifcat\noexpand\next a\def\next{\csname atZZ\endcsname}\else
   \ifcat\noexpand\next0\def\next{\csname atZZ\endcsname}\else
    \def\next{\csname atZZZ\endcsname}\fi\fi\next}
\expandafter\gdef\csname atZZ\endcsname#1{\expandafter
   \ifx\csname #1Zat\endcsname\relax\def\next
     {\errhelp\expandafter=\csname athelpZ\endcsname
      \errmessage{Invalid use of \string@}}\else
       \def\next{\csname #1Zat\endcsname}\fi\next}
\expandafter\gdef\csname atZZZ\endcsname#1{\errhelp
    \expandafter=\csname athelpZ\endcsname
      \errmessage{Invalid use of \string@}}}}
\def\atdef@#1{\expandafter\def\csname #1@at\endcsname}
\newhelp\defahelp@{If you typed \string\define\space cs instead of
\string\define\string\cs\space^^J
I've substituted an inaccessible control sequence so that your^^J
definition will be completed without mixing me up too badly.^^J
If you typed \string\define{\string\cs} the inaccessible control sequence^^J
was defined to be \string\cs, and the rest of your^^J
definition appears as input.}
\newhelp\defbhelp@{I've ignored your definition, because it might^^J
conflict with other uses that are important to me.}
\def\define{\futurelet\next\define@}
\def\define@{\ifcat\noexpand\next\relax
  \def\next{\define@@}%
  \else\errhelp=\defahelp@
  \errmessage{\string\define\space must be followed by a control 
     sequence}\def\next{\def\garbage@}\fi\next}
\def\undefined@{}
\def\preloaded@{}    
\def\define@@#1{\ifx#1\relax\errhelp=\defbhelp@
   \errmessage{\string#1\space is already defined}\def\next{\def\garbage@}%
   \else\expandafter\ifx\csname\expandafter\eat@\string
         #1@\endcsname\undefined@\errhelp=\defbhelp@
   \errmessage{\string#1\space can't be defined}\def\next{\def\garbage@}%
   \else\expandafter\ifx\csname\expandafter\eat@\string#1\endcsname\relax
     \def\next{\def#1}\else\errhelp=\defbhelp@
     \errmessage{\string#1\space is already defined}\def\next{\def\garbage@}%
      \fi\fi\fi\next}
\def\famzero{\fam\z@}

\def\lim{\mathop{\famzero lim}}

\def\textfont@#1#2{\def#1{\relax\ifmmode
    \errmessage{Use \string#1\space only in text}\else#2\fi}}
\textfont@\rm\tenrm
\textfont@\it\tenit
\textfont@\sl\tensl
\textfont@\bf\tenbf
\textfont@\smc\tensmc
\let\ic@=\/
\def\/{\unskip\ic@}
\def\textfonti{\the\textfont1 }
\def\t#1#2{{\edef\next{\the\font}\textfonti\accent"7F \next#1#2}}
\let\B=\=
\let\D=\.
\def~{\unskip\nobreak\ \ignorespaces}
{\catcode`\@=\active
\gdef\@{\char'100 }}
\atdef@-{\leavevmode\futurelet\next\athyph@}
\def\athyph@{\ifx\next-\let\next=\athyph@@
  \else\let\next=\athyph@@@\fi\next}
\def\athyph@@@{\hbox{-}}
\def\athyph@@#1{\futurelet\next\athyph@@@@}
\def\athyph@@@@{\if\next-\def\next##1{\hbox{---}}\else
    \def\next{\hbox{--}}\fi\next}
\def\.{.\spacefactor=\@m}
\atdef@.{\null.}
\atdef@,{\null,}
\atdef@;{\null;}
\atdef@:{\null:}
\atdef@?{\null?}
\atdef@!{\null!}   
\def\srdr@{\thinspace}                     
\def\drsr@{\kern.02778em}
\def\sldl@{\kern.02778em}
\def\dlsl@{\thinspace}
\atdef@"{\unskip\futurelet\next\atqq@}
\def\atqq@{\ifx\next\Space@\def\next. {\atqq@@}\else
         \def\next.{\atqq@@}\fi\next.}
\def\atqq@@{\futurelet\next\atqq@@@}
\def\atqq@@@{\ifx\next`\def\next`{\atqql@}\else\def\next'{\atqqr@}\fi\next}
\def\atqql@{\futurelet\next\atqql@@}
\def\atqql@@{\ifx\next`\def\next`{\sldl@``}\else\def\next{\dlsl@`}\fi\next}
\def\atqqr@{\futurelet\next\atqqr@@}
\def\atqqr@@{\ifx\next'\def\next'{\srdr@''}\else\def\next{\drsr@'}\fi\next}

\def\textfontii{\the\textfont2 }
\def\{{\relax\ifmmode\lbrace\else
    {\textfontii f}\spacefactor=\@m\fi}
\def\}{\relax\ifmmode\rbrace\else
    \let\@sf=\empty\ifhmode\edef\@sf{\spacefactor=\the\spacefactor}\fi
      {\textfontii g}\@sf\relax\fi}   
\def\nonhmodeerr@#1{\errmessage
     {\string#1\space allowed only within text}}
\def\linebreak{\relax\ifhmode\unskip\break\else
    \nonhmodeerr@\linebreak\fi}
\def\allowlinebreak{\relax
   \ifhmode\allowbreak\else\nonhmodeerr@\allowlinebreak\fi}
\newskip\saveskip@
\def\nolinebreak{\relax\ifhmode\saveskip@=\lastskip\unskip
  \nobreak\ifdim\saveskip@>\z@\hskip\saveskip@\fi
   \else\nonhmodeerr@\nolinebreak\fi}
\def\newline{\relax\ifhmode\null\hfil\break
    \else\nonhmodeerr@\newline\fi}
\def\nonmathaerr@#1{\errmessage
     {\string#1\space is not allowed in display math mode}}
\def\nonmathberr@#1{\errmessage{\string#1\space is allowed only in math mode}}
\def\mathbreak{\relax\ifmmode\ifinner\break\else
   \nonmathaerr@\mathbreak\fi\else\nonmathberr@\mathbreak\fi}
\def\nomathbreak{\relax\ifmmode\ifinner\nobreak\else
    \nonmathaerr@\nomathbreak\fi\else\nonmathberr@\nomathbreak\fi}
\def\allowmathbreak{\relax\ifmmode\ifinner\allowbreak\else
     \nonmathaerr@\allowmathbreak\fi\else\nonmathberr@\allowmathbreak\fi}
\def\pagebreak{\relax\ifmmode
   \ifinner\errmessage{\string\pagebreak\space
     not allowed in non-display math mode}\else\postdisplaypenalty-\@M\fi
   \else\ifvmode\penalty-\@M\else\edef\spacefactor@
       {\spacefactor=\the\spacefactor}\vadjust{\penalty-\@M}\spacefactor@
        \relax\fi\fi}
\def\nopagebreak{\relax\ifmmode
     \ifinner\errmessage{\string\nopagebreak\space
    not allowed in non-display math mode}\else\postdisplaypenalty\@M\fi
    \else\ifvmode\nobreak\else\edef\spacefactor@
        {\spacefactor=\the\spacefactor}\vadjust{\penalty\@M}\spacefactor@
         \relax\fi\fi}
\def\newpage{\relax\ifvmode\vfill\penalty-\@M\else\nonvmodeerr@\newpage\fi}
\def\nonvmodeerr@#1{\errmessage
    {\string#1\space is allowed only between paragraphs}}
\def\smallpagebreak{\relax\ifvmode\smallbreak
      \else\nonvmodeerr@\smallpagebreak\fi}
\def\medpagebreak{\relax\ifvmode\medbreak
       \else\nonvmodeerr@\medpagebreak\fi}
\def\bigpagebreak{\relax\ifvmode\bigbreak
      \else\nonvmodeerr@\bigpagebreak\fi}
\newdimen\captionwidth@
\captionwidth@=\hsize
\advance\captionwidth@ by -1.5in
\def\caption#1{}
\def\topspace#1{\gdef\thespace@{#1}\ifvmode\def\next
    {\futurelet\next\topspace@}\else\def\next{\nonvmodeerr@\topspace}\fi\next}
\def\topspace@{\ifx\next\Space@\def\next. {\futurelet\next\topspace@@}\else
     \def\next.{\futurelet\next\topspace@@}\fi\next.}
\def\topspace@@{\ifx\next\caption\let\next\topspace@@@\else
    \let\next\topspace@@@@\fi\next}
 \def\topspace@@@@{\topinsert\vbox to 
       \thespace@{}\endinsert}
\def\topspace@@@\caption#1{\topinsert\vbox to
    \thespace@{}\nobreak
      \smallskip
    \setbox\z@=\hbox{\noindent\ignorespaces#1\unskip}%
   \ifdim\wd\z@>\captionwidth@
   \centerline{\vbox{\hsize=\captionwidth@\noindent\ignorespaces#1\unskip}}%
   \else\centerline{\box\z@}\fi\endinsert}
\def\midspace#1{\gdef\thespace@{#1}\ifvmode\def\next
    {\futurelet\next\midspace@}\else\def\next{\nonvmodeerr@\midspace}\fi\next}
\def\midspace@{\ifx\next\Space@\def\next. {\futurelet\next\midspace@@}\else
     \def\next.{\futurelet\next\midspace@@}\fi\next.}
\def\midspace@@{\ifx\next\caption\let\next\midspace@@@\else
    \let\next\midspace@@@@\fi\next}
 \def\midspace@@@@{\midinsert\vbox to 
       \thespace@{}\endinsert}
\def\midspace@@@\caption#1{\midinsert\vbox to
    \thespace@{}\nobreak
      \smallskip
      \setbox\z@=\hbox{\noindent\ignorespaces#1\unskip}%
      \ifdim\wd\z@>\captionwidth@
    \centerline{\vbox{\hsize=\captionwidth@\noindent\ignorespaces#1\unskip}}%
    \else\centerline{\box\z@}\fi\endinsert}
\mathchardef\prime@="0230
\def\prime{{{}\prime@{}}}
\def\prim@s{\prime@\futurelet\next\pr@m@s}

\def\,{\relax\ifmmode\mskip\thinmuskip\else\thinspace\fi}
\def\!{\relax\ifmmode\mskip-\thinmuskip\else\negthinspace\fi}
\def\frac#1#2{{#1\over#2}}

\def\:{\nobreak\hskip.1111em{:}\hskip.3333em plus .0555em\relax}
\def\intic@{\mathchoice{\hskip5\p@}{\hskip4\p@}{\hskip4\p@}{\hskip4\p@}}
\def\negintic@
 {\mathchoice{\hskip-5\p@}{\hskip-4\p@}{\hskip-4\p@}{\hskip-4\p@}}
\def\intkern@{\mathchoice{\!\!\!}{\!\!}{\!\!}{\!\!}}
\def\intdots@{\mathchoice{\cdots}{{\cdotp}\mkern1.5mu
    {\cdotp}\mkern1.5mu{\cdotp}}{{\cdotp}\mkern1mu{\cdotp}\mkern1mu
      {\cdotp}}{{\cdotp}\mkern1mu{\cdotp}\mkern1mu{\cdotp}}}
\newcount\intno@             
\def\iint{\intno@=\tw@\futurelet\next\ints@} 
\def\iiint{\intno@=\thr@@\futurelet\next\ints@}
\def\iiiint{\intno@=4 \futurelet\next\ints@}
\def\idotsint{\intno@=\z@\futurelet\next\ints@}
\def\ints@{\findlimits@\ints@@}
\newif\iflimtoken@
\newif\iflimits@
\def\findlimits@{\limtoken@false\limits@false\ifx\next\limits
 \limtoken@true\limits@true\else\ifx\next\nolimits\limtoken@true\limits@false
    \fi\fi}
\def\multintlimits@{\intop\ifnum\intno@=\z@\intdots@
  \else\intkern@\fi
    \ifnum\intno@>\tw@\intop\intkern@\fi
     \ifnum\intno@>\thr@@\intop\intkern@\fi\intop}
\def\multint@{\int\ifnum\intno@=\z@\intdots@\else\intkern@\fi
   \ifnum\intno@>\tw@\int\intkern@\fi
    \ifnum\intno@>\thr@@\int\intkern@\fi\int}
\def\ints@@{\iflimtoken@\def\ints@@@{\iflimits@
   \negintic@\mathop{\intic@\multintlimits@}\limits\else
    \multint@\nolimits\fi\eat@}\else
     \def\ints@@@{\multint@\nolimits}\fi\ints@@@}
\def\Sb{_\bgroup\vspace@
        \baselineskip=\fontdimen10 \scriptfont\tw@
        \advance\baselineskip by \fontdimen12 \scriptfont\tw@
        \lineskip=\thr@@\fontdimen8 \scriptfont\thr@@
        \lineskiplimit=\thr@@\fontdimen8 \scriptfont\thr@@
        \Let@\vbox\bgroup\halign\bgroup \hfil$\scriptstyle
            {##}$\hfil\cr}
\def\endSb{\crcr\egroup\egroup\egroup}
\def\Sp{^\bgroup\vspace@
        \baselineskip=\fontdimen10 \scriptfont\tw@
        \advance\baselineskip by \fontdimen12 \scriptfont\tw@
        \lineskip=\thr@@\fontdimen8 \scriptfont\thr@@
        \lineskiplimit=\thr@@\fontdimen8 \scriptfont\thr@@
        \Let@\vbox\bgroup\halign\bgroup \hfil$\scriptstyle
            {##}$\hfil\cr}
\def\endSp{\crcr\egroup\egroup\egroup}
\def\Let@{\relax\iffalse{\fi\let\\=\cr\iffalse}\fi}
\def\vspace@{\def\vspace##1{\noalign{\vskip##1 }}}
\def\aligned{\,\vcenter\bgroup\vspace@\Let@\openup\jot\m@th\ialign
  \bgroup \strut\hfil$\displaystyle{##}$&$\displaystyle{{}##}$\hfil\crcr}
\def\endaligned{\crcr\egroup\egroup}
\def\matrix{\,\vcenter\bgroup\Let@\vspace@
    \normalbaselines
  \m@th\ialign\bgroup\hfil$##$\hfil&&\quad\hfil$##$\hfil\crcr
    \mathstrut\crcr\noalign{\kern-\baselineskip}}
\def\endmatrix{\crcr\mathstrut\crcr\noalign{\kern-\baselineskip}\egroup
                \egroup\,}
\newtoks\hashtoks@
\hashtoks@={#}
\def\format{\crcr\egroup\iffalse{\fi\ifnum`}=0 \fi\format@}
\def\format@#1\\{\def\preamble@{#1}%
  \def\c{\hfil$\the\hashtoks@$\hfil}%
  \def\r{\hfil$\the\hashtoks@$}%
  \def\l{$\the\hashtoks@$\hfil}%
  \setbox\z@=\hbox{\xdef\Preamble@{\preamble@}}\ifnum`{=0 \fi\iffalse}\fi
   \ialign\bgroup\span\Preamble@\crcr}

\def\cases{\left\{\,\vcenter\bgroup\vspace@
     \normalbaselines\openup\jot\m@th
       \Let@\ialign\bgroup$##$\hfil&\quad$##$\hfil\crcr
      \mathstrut\crcr\noalign{\kern-\baselineskip}}
\def\endcases{\endmatrix\right.}
\newif\iftagsleft@
\tagsleft@true
\def\TagsOnRight{\global\tagsleft@false}
\def\tag#1$${\iftagsleft@\leqno\else\eqno\fi
 \hbox{\def\pagebreak{\global\postdisplaypenalty-\@M}%
 \def\nopagebreak{\global\postdisplaypenalty\@M}\rm(#1\unskip)}%
  $$\postdisplaypenalty\z@\ignorespaces}
\interdisplaylinepenalty=\@M
\def\allowdisplaybreak@{\def\allowdisplaybreak{\noalign{\allowbreak}}}
\def\displaybreak@{\def\displaybreak{\noalign{\break}}}
\def\align#1\endalign{\def\tag{&}\vspace@\allowdisplaybreak@\displaybreak@
  \iftagsleft@\lalign@#1\endalign\else
   \ralign@#1\endalign\fi}
\def\ralign@#1\endalign{\displ@y\Let@\tabskip\centering\halign to\displaywidth
     {\hfil$\displaystyle{##}$\tabskip=\z@&$\displaystyle{{}##}$\hfil
       \tabskip=\centering&\llap{\hbox{(\rm##\unskip)}}\tabskip\z@\crcr
             #1\crcr}}
\def\lalign@
 #1\endalign{\displ@y\Let@\tabskip\centering\halign to \displaywidth
   {\hfil$\displaystyle{##}$\tabskip=\z@&$\displaystyle{{}##}$\hfil
   \tabskip=\centering&\kern-\displaywidth
        \rlap{\hbox{(\rm##\unskip)}}\tabskip=\displaywidth\crcr
               #1\crcr}}
\def\overrightarrow{\mathpalette\overrightarrow@}
\def\overrightarrow@#1#2{\vbox{\ialign{$##$\cr
    #1{-}\mkern-6mu\cleaders\hbox{$#1\mkern-2mu{-}\mkern-2mu$}\hfill
     \mkern-6mu{\to}\cr
     \noalign{\kern -1\p@\nointerlineskip}
     \hfil#1#2\hfil\cr}}}
\def\overleftarrow{\mathpalette\overleftarrow@}
\def\overleftarrow@#1#2{\vbox{\ialign{$##$\cr
     #1{\leftarrow}\mkern-6mu\cleaders\hbox{$#1\mkern-2mu{-}\mkern-2mu$}\hfill
      \mkern-6mu{-}\cr
     \noalign{\kern -1\p@\nointerlineskip}
     \hfil#1#2\hfil\cr}}}
\def\overleftrightarrow{\mathpalette\overleftrightarrow@}
\def\overleftrightarrow@#1#2{\vbox{\ialign{$##$\cr
     #1{\leftarrow}\mkern-6mu\cleaders\hbox{$#1\mkern-2mu{-}\mkern-2mu$}\hfill
       \mkern-6mu{\to}\cr
    \noalign{\kern -1\p@\nointerlineskip}
      \hfil#1#2\hfil\cr}}}
\def\underrightarrow{\mathpalette\underrightarrow@}
\def\underrightarrow@#1#2{\vtop{\ialign{$##$\cr
    \hfil#1#2\hfil\cr
     \noalign{\kern -1\p@\nointerlineskip}
    #1{-}\mkern-6mu\cleaders\hbox{$#1\mkern-2mu{-}\mkern-2mu$}\hfill
     \mkern-6mu{\to}\cr}}}
\def\underleftarrow{\mathpalette\underleftarrow@}
\def\underleftarrow@#1#2{\vtop{\ialign{$##$\cr
     \hfil#1#2\hfil\cr
     \noalign{\kern -1\p@\nointerlineskip}
     #1{\leftarrow}\mkern-6mu\cleaders\hbox{$#1\mkern-2mu{-}\mkern-2mu$}\hfill
      \mkern-6mu{-}\cr}}}
\def\underleftrightarrow{\mathpalette\underleftrightarrow@}
\def\underleftrightarrow@#1#2{\vtop{\ialign{$##$\cr
      \hfil#1#2\hfil\cr
    \noalign{\kern -1\p@\nointerlineskip}
     #1{\leftarrow}\mkern-6mu\cleaders\hbox{$#1\mkern-2mu{-}\mkern-2mu$}\hfill
       \mkern-6mu{\to}\cr}}}
\def\sqrt#1{\radical"270370 {#1}}
\def\dots{\relax\ifmmode\let\next=\ldots\else\let\next=\tdots@\fi\next}
\def\tdots@{\unskip\ \tdots@@}
\def\tdots@@{\futurelet\next\tdots@@@}
\def\tdots@@@{$\mathinner{\ldotp\ldotp\ldotp}\,
   \ifx\next,$\else
   \ifx\next.\,$\else
   \ifx\next;\,$\else
   \ifx\next:\,$\else
   \ifx\next?\,$\else
   \ifx\next!\,$\else
   $ \fi\fi\fi\fi\fi\fi}
\def\text{\relax\ifmmode\let\next=\text@\else\let\next=\text@@\fi\next}
\def\text@@#1{\hbox{#1}}
\def\text@#1{\mathchoice
 {\hbox{\everymath{\displaystyle}\def\textfonti{\the\textfont1 }%
    \def\textfontii{\the\textfont2 }\textdef@@ T#1}}
 {\hbox{\everymath{\textstyle}\def\textfonti{\the\textfont1 }%
    \def\textfontii{\the\textfont2 }\textdef@@ T#1}}
 {\hbox{\everymath{\scriptstyle}\def\textfonti{\the\scriptfont1 }%
   \def\textfontii{\the\scriptfont2 }\textdef@@ S\rm#1}}
 {\hbox{\everymath{\scriptscriptstyle}\def\textfonti{\the\scriptscriptfont1 }%
   \def\textfontii{\the\scriptscriptfont2 }\textdef@@ s\rm#1}}}
\def\textdef@@#1{\textdef@#1\rm \textdef@#1\bf
   \textdef@#1\sl \textdef@#1\it}

\def\textdef@#1#2{\def\next{\csname\expandafter\eat@\string#2fam\endcsname}%
\if S#1\edef#2{\the\scriptfont\next\relax}%
 \else\if s#1\edef#2{\the\scriptscriptfont\next\relax}%
 \else\edef#2{\the\textfont\next\relax}\fi\fi}
\scriptfont\itfam=\tenit \scriptscriptfont\itfam=\tenit
\scriptfont\slfam=\tensl \scriptscriptfont\slfam=\tensl
\mathcode`\0="0030
\mathcode`\1="0031
\mathcode`\2="0032
\mathcode`\3="0033
\mathcode`\4="0034
\mathcode`\5="0035
\mathcode`\6="0036
\mathcode`\7="0037
\mathcode`\8="0038
\mathcode`\9="0039
\def\Cal{\relax\ifmmode\let\next=\Cal@\else
     \def\next{\errmessage{Use \string\Cal\space only in math mode}}\fi\next}
\def\Cal@#1{{\fam2 #1}}
\def\bold{\relax\ifmmode\let\next=\bold@\else
   \def\next{\errmessage{Use \string\bold\space only in math
      mode}}\fi\next}\def\bold@#1{{\fam\bffam #1}}
\mathchardef\Gamma="0000
\mathchardef\Delta="0001
\mathchardef\Theta="0002
\mathchardef\Lambda="0003
\mathchardef\Xi="0004
\mathchardef\Pi="0005
\mathchardef\Sigma="0006
\mathchardef\Upsilon="0007
\mathchardef\Phi="0008
\mathchardef\Psi="0009
\mathchardef\Omega="000A
\mathchardef\varGamma="0100
\mathchardef\varDelta="0101
\mathchardef\varTheta="0102
\mathchardef\varLambda="0103
\mathchardef\varXi="0104
\mathchardef\varPi="0105
\mathchardef\varSigma="0106
\mathchardef\varUpsilon="0107
\mathchardef\varPhi="0108
\mathchardef\varPsi="0109
\mathchardef\varOmega="010A
\font\dummyft@=dummy
\fontdimen1 \dummyft@=\z@
\fontdimen2 \dummyft@=\z@
\fontdimen3 \dummyft@=\z@
\fontdimen4 \dummyft@=\z@
\fontdimen5 \dummyft@=\z@
\fontdimen6 \dummyft@=\z@
\fontdimen7 \dummyft@=\z@
\fontdimen8 \dummyft@=\z@
\fontdimen9 \dummyft@=\z@
\fontdimen10 \dummyft@=\z@
\fontdimen11 \dummyft@=\z@
\fontdimen12 \dummyft@=\z@
\fontdimen13 \dummyft@=\z@
\fontdimen14 \dummyft@=\z@
\fontdimen15 \dummyft@=\z@
\fontdimen16 \dummyft@=\z@
\fontdimen17 \dummyft@=\z@
\fontdimen18 \dummyft@=\z@
\fontdimen19 \dummyft@=\z@
\fontdimen20 \dummyft@=\z@
\fontdimen21 \dummyft@=\z@
\fontdimen22 \dummyft@=\z@
\def\fontlist@{\\{\tenrm}\\{\sevenrm}\\{\fiverm}\\{\teni}\\{\seveni}%
 \\{\fivei}\\{\tensy}\\{\sevensy}\\{\fivesy}\\{\tenex}\\{\tenbf}\\{\sevenbf}%
 \\{\fivebf}\\{\tensl}\\{\tenit}\\{\tensmc}}
\def\dodummy@{{\def\\##1{\global\let##1=\dummyft@}\fontlist@}}
\newif\ifsyntax@
\newcount\countxviii@
\def\newtoks@{\alloc@5\toks\toksdef\@cclvi}
\def\nopages@{\output={\setbox\z@=\box\@cclv \deadcycles=\z@}\newtoks@\output}
\def\syntax{\syntax@true\dodummy@\countxviii@=\count18
\loop \ifnum\countxviii@ > \z@ \textfont\countxviii@=\dummyft@
   \scriptfont\countxviii@=\dummyft@ \scriptscriptfont\countxviii@=\dummyft@
     \advance\countxviii@ by-\@ne\repeat
\dummyft@\tracinglostchars=\z@
  \nopages@\frenchspacing\hbadness=\@M}
\def\magstep#1{\ifcase#1 1000\or
 1200\or 1440\or 1728\or 2074\or 2488\or 
 \errmessage{\string\magstep\space only works up to 5}\fi\relax}
{\lccode`\2=`\p \lccode`\3=`\t 
 \lowercase{\gdef\tru@#123{#1truept}}}

\def\scaletype#1{\mag=#1\relax
 \hsize=\expandafter\tru@\the\hsize
 \vsize=\expandafter\tru@\the\vsize
 \dimen\footins=\expandafter\tru@\the\dimen\footins}

\def\scalefont#1#2\andcallit#3{\edef\font@{\the\font}#1\font#3=
  \fontname\font\space scaled #2\relax\font@}
\def\Mag@#1#2{\ifdim#1<1pt\multiply#1 #2\relax\divide#1 1000 \else
  \ifdim#1<10pt\divide#1 10 \multiply#1 #2\relax\divide#1 100\else
  \divide#1 100 \multiply#1 #2\relax\divide#1 10 \fi\fi}
\def\scalelinespacing#1{\Mag@\baselineskip{#1}\Mag@\lineskip{#1}%
  \Mag@\lineskiplimit{#1}}
\def\wlog#1{\immediate\write-1{#1}}
\catcode`\@=\active

\magnification=\magstep1
\baselineskip=5.5mm
\settabs 4\columns

\title
Local factorization and monomialization of morphisms
\endtitle
\author Steven Dale Cutkosky 
\footnote{partially supported by NSF}
\footnote" "{
} 
\endauthor

\subheading{Contents}
\vskip .2truein
\+ Introduction &&&1\cr

Statement of the main results

Geometry and valuations

Overview of the proof
\+ 1. Preliminaries&&& 19\cr

Valuations

Birational Transforms
\+ 2. Uniformizing Transforms&&& 24\cr
\+3. Rank 1&&& 27\cr

Perron Transforms

Monomialization in Rank 1
\+ 4. Monomialization&&& 112\cr
\+ 5. Factorization 1&&& 131\cr
\+ 6. Factorization 2&&& 141\cr

\heading Introduction \endheading
\subheading{Statement of the main results}

Suppose that we are given a  system of equations
$$
\align
x_1 &= f_1(y_1,y_2,\ldots,y_n)\tag *\\
&\vdots\\
x_n &= f_n(y_1,y_2,\ldots,y_n)
\endalign
$$
which is nondegenerate, in the sense that the Jacobian determinate  of the system is not (identically) zero.
This system is well understood in the special case that $f_1,\ldots, f_n$ are monomials in the variables
$y_1,\ldots, y_n$. For instance, by inverting the matrix $A$ of coefficients of the monomials, we can express   $y_1,\ldots, y_n$
 as rational functions of $d$-th roots of the variables $x_1,\ldots, x_n$, where $d$ is the determinant of $A$.  

When the $f_i$ are not monomials, it is not easy to analyze such a system.

Our main result shows that all solutions of a system (*) can be expressed in the following simple form.
There are  finitely many charts obtained from a composition of  monoidal tranforms in the variables $x$ and $y$
$$
\align
x_i &= \Phi_i(\overline x_1,\ldots,\overline x_n),\,\, 1\le i\le n\\
y_i &= \Psi_i(\overline y_1,\ldots,\overline y_n),\,\, 1\le i\le n
\endalign
$$
such that
the transform of the system (*) becomes a system of monomial equations 
$$
\align
\overline x_1 &= \overline y_1^{a_{11}} ..... \overline y_n^{a_{1n}}\\
\vdots&\\   
\overline x_n &= \overline y_1^{a_{n1}} ..... \overline y_n^{a_{nn}}
\endalign
$$
with $\text{det}(a_{ij})\ne0$.
A monoidal transform is a composition of
\item{1)} a change of variable
\item{2)} a transform  
$$
\align
x_1 &= x_1(1) x_2(1)\\
x_i &= x_i(1)\text{ if } i>1.
\endalign
$$

Our solution is constructive, as  it consists of  a series of algorithms. 

This result can be interpreted geometrically as follows. Suppose that $\phi:X\rightarrow Y$
is a generically finite morphism of varieties. Then it is possible to construct a finite number of
charts $X_i$ and $Y_i$ such that  $X_i\rightarrow Y_i$ are monomial mappings,
the mappings $X_i\rightarrow X$ and $Y_i\rightarrow Y$
are sequences of blowups of nonsingular subvarieties, 
and $X_i$ and $Y_i$ form complete systems, in the sense that they can be patched to obtain schemes which satisfy
the valuative criteria of properness.   
 
Our main result is  stated precisely in Theorem A.

\vfill\eject

 \proclaim{Theorem A} (Monomialization)
 Suppose that $R \subset S$ are  excellent regular local rings such that $\text{dim}(R)=\text{dim}(S)$,
 containing a field $k$ of
characteristic zero, such that the quotient field $K$ of $S$ is a finite extension of the quotient field $J$
of $R$.

Let $V$ be a valuation ring of $K$ which dominates $S$.  Then there
exist sequences of  monoidal transforms (blow ups of regular primes)
 $R \rightarrow R'$ and $S \rightarrow S'$
such that $V$ dominates $S'$, $S'$ dominates $R'$ and 
there are regular parameters $(x_1, .... ,x_n)$
in $R'$,  $(y_1, ... ,y_n)$ in $S'$, units $\delta_1,\ldots,\delta_n\in S'$ and a matrix $(a_{ij})$ of
nonnegative integers such that  $\text{Det}(a_{ij}) \ne 0$ and
$$
\align
x_1 &= y_1^{a_{11}} ..... y_n^{a_{1n}}\delta_1\tag **\\
\vdots&\\   
x_n &= y_1^{a_{n1}} ..... y_n^{a_{nn}}\delta_n.
\endalign
$$

\endproclaim

With the assumptions of Theorem A, An example of Abhyankar (Theorem 12 [Ab6])
shows that it is in general not possible to perform monoidal transforms along $V$ in $R$ and $S$ to obtain
 $R'\rightarrow S'$ such that $R'\rightarrow S'$ is (a localization of) a finite map.
As such, Theorem A is the strongest possible local result for generically finite maps.

A more geometric statement of Theorem A is given in Theorem B.
A complete variety over a field $k$ is an integral finite type $k$-scheme which satisfies the 
existence part of the valuative criterion for properness.  
Complete and separated is equivalent to proper.

\proclaim{Theorem B} Let $k$ be a field of characteristic zero, $\Phi:X\rightarrow Y$ a generically finite
morphism of integral nonsingular proper excellent $k$-schemes.   Then there are
  birational morphisms of nonsingular complete excellent $k$-schemes $\alpha:X_{1}\rightarrow X$
and $\beta: Y_{1}\rightarrow Y$, and a morphism $\Psi:X_1\rightarrow Y_1$ such that  the diagram
$$
\matrix
X_1&{\Psi}\atop{\rightarrow}&Y_1\\
\downarrow&&\downarrow\\
X&{\Phi}\atop{\rightarrow}&Y
\endmatrix
$$
commutes, $\alpha$ and $\beta$ are
 locally  products of blowups of nonsingular subvarieties, and $\Psi$ is locally a monomial mapping. 
That is, for every $z\in X_{1}$, there
exist affine neighborhoods $V_1$ of $z$, $V$ of $x=\alpha(z)$, such that $\alpha:V_1\rightarrow V$ 
is a  finite product of monoidal transforms,  
there
exist affine neighborhoods $W_1$ of $\Psi(z)$, $W$ of $y=\alpha(\Psi(z))$, such that $\beta:W_1\rightarrow W$ 
is a  finite product of monoidal transforms, and $\Psi:V_1\rightarrow W_1$ is a  mapping of the form (**) in
some uniformizing parameters of $V_1$ and $W_1$.
\endproclaim

Here a monoidal transform of a nonsingular $k$-scheme $S$ is the map $T\rightarrow S$ induced by an open subset $T$ of
$\text{Proj}(\oplus \Cal I^n)$, where $\Cal I$ is the ideal sheaf of a nonsingular subvariety of $S$.

In the special case of dimension two,  we can strengthen the conclusions of Theorem B.

\proclaim{Theorem B2} Let $k$ be a field of characteristic zero, $\Phi:S\rightarrow T$ a generically finite
morphism of integral nonsingular proper (projective) excellent $k$-surfaces.   Then there are
products of blowups of points (quadratic transforms) $\alpha:S_{1}\rightarrow S$
and $\beta: T_{1}\rightarrow T$, and a morphism $\Psi:S_1\rightarrow T_1$ such that  the diagram
$$
\matrix
S_1&{\Psi}\atop{\rightarrow}&T_1\\
\downarrow&&\downarrow\\
S&{\Phi}\atop{\rightarrow}&T
\endmatrix
$$
commutes, and $\Psi$ is locally a monomial mapping. 
That is, for every $z\in S_{1}$, there
exist affine neighborhoods $V_1$ of $z$ and  $W_1$ of $\Psi(z)$, such that 
 $\Psi:V_1\rightarrow W_1$ is a  mapping of the form (**) in
some uniformizing parameters of $V_1$ and $W_1$.
\endproclaim  

In the  case of complex surfaces,  a proof of Theorem B2
follows from results of Akbulut and King  (Chapter 7 of [AK]).

Stronger results hold for birational morphisms, morphisms which are an isomorphism on an open set.
A birational morphism of nonsingular  projective surfaces can be factored by a product of quadratic transforms.
This was proved by Zariski, over an algebraically closed field of arbitrary characteristic,
as a corollary to a local theorem on factorization (on page 589 of [Z3] and in section II.1 of [Z4]).
The most general form of this Theorem is due to Abhyankar, in Theorem 3 of his 1956 paper [Ab2].
Abhyankar proves that
an inclusion $R\subset S$ of regular local rings of dimension 2 with a common quotient field
can be factored by
a finite sequence of quadratic transforms (blowups of points). 

In higher dimensions, the simplest birational morphisms are the monoidal transforms. A monoidal transform is a
  blowup of a nonsingular subvariety.   Sally [S] and Shannon [Sh] have found  examples 
 of inclusions   $R\subset S$ of regular local rings of dimension 3 with a common quotient field
which cannot be factored by
a finite sequence of monoidal  transforms (blowups of points and nonsingular curves). 

In [C], we prove the following Theorem, which gives a positive answer to a conjecture of Abhyankar (page 237 [Ab5], [Ch]),
over fields of characteristic 0. In view of the counterexamples to a direct factorization, Theorem B is the best possible
local factorization result in dimension three.

\proclaim{Theorem C} (Theorem A  [C])
 Suppose that  $R\subset S$ are excellent regular local rings such that $\text{dim}(R)=\text{dim}(S)=3$,
containing a field $k$ of characteristic zero and  with a common quotient field $K$. 
  Let $V$ be a valuation ring of
$K$ which dominates $S$. Then there exists a regular local ring $T$, with quotient field $K$,
such that $T$ dominates $S$, $V$ dominates  $T$, and the inclusions $R\rightarrow T$ and $S\rightarrow T$
can be factored by sequences of monoidal transforms (blowups of regular primes).
$$
\matrix
&&V&&\\
&&\uparrow&&\\
&&T&&\\
&\nearrow&&\nwarrow&\\
R&&\longrightarrow && S
\endmatrix
$$
\endproclaim
\vskip .2truein

It is natural to ask if the generalization of this three dimensional factorization theorem is possible in all dimensions
by constructing a factorization by a sequence of   blowups and blowdowns with nonsingular centers along a valuation.
In this paper, we prove the following theorem which gives a positive answer to this question in all dimensions.

\proclaim{Theorem D}(Factorization 1)
 Suppose that $R \subset S$ are  excellent regular local rings of dimension $n\ge 3$, containing a field $k$ of
characteristic zero, with a common quotient field $K$.
Let $V$ be a valuation ring of $K$ which dominates $S$.
Then there 
exist sequences of  regular local rings contained in $K$

$$
\matrix
 &          &    R_1  &             &     &           &        &          &         &          & R_{n-2} &         &\\
 &  \nearrow &        &   \nwarrow  &     &  \nearrow & \cdots & \nwarrow &         & \nearrow &         &\nwarrow &\\
R&          &         &             & S_1 &           &        &          & S_{n-3} &          &         & &S_{n-2}=S\\
\endmatrix
$$
 such that each local ring is dominated by $V$
and each arrow is a sequence of  monoidal transforms (blow ups of regular primes).
Furthermore, we have inclusions
$R\subset S_i$ for all $i$.
\endproclaim

Theorem C follows from the special case $n=3$ of  Theorem D.

The proofs of the above theorems are essentially self contained in this paper. We only assume some
basic results on valuation theory (as can be found in [Ab3] and [ZS]) and  the basic resolution theorems of
Hironaka [H]. The Hironaka results are essentially only used in the case of a composite valuation, to establish the
existence of a nonsingular center of a composite valuation.

A long standing  conjecture in algebraic geometry is that one can factor a birational morphism
$X\rightarrow Y$ between nonsingular  projective varieties by a series of alternating blowups and blowdowns with
nonsingular centers (c.f. [P]). We will refer to this as the global factorization conjecture. 
In [P] an example is given of Hironaka, showing that it is not
possible in dimension $\ge 3$ to always factor birational morphisms of nonsingular varieties by blowups
with nonsingular centers.

Our Theorem D shows that there is no local obstruction to the  global factorization conjecture in any dimension.
We prove a local form of this conjecture.

\proclaim{Theorem E} Suppose that 
$X\rightarrow Y$ is a birational morphism of nonsingular projective 
$n$-dimensional  varieties, over a field of characteristic zero, and  
$\nu$ is a  valuation of the function field of $X$. Then there is
a sequence of projective birational morphisms of nonsingular varieties  
$$
\matrix
 &          &    X_1  &             &     &           &        &          &         &          & X_{n-1} &         &\\
 &  \swarrow &        &   \searrow  &     &  \swarrow & \cdots & \searrow &         & \swarrow &         &\searrow &\\
X&          &         &             & Y_1 &           &        &          & Y_{n-1} &          &         & &Y_{n-2}=Y\\
\endmatrix
$$
such that each morphism is a product of blowups of nonsingular subvarities in a Zariski neighbourhood of the
center of $\nu$.
\endproclaim

\proclaim{Theorem F} Let $k$ be a field of characteristic zero, $\phi:X\rightarrow Y$ a birational
morphism of integral nonsingular proper excellent $k$-schemes of dimension $n$.   Then there is
a sequence of  birational morphisms of nonsingular complete $k$-schemes $\alpha_i:X_{i+1}\rightarrow X_i$
and $\beta_i:X_{i+1}\rightarrow Y_{i+1}$ 
$$
\matrix
 &          &    X_1  &             &     &           &        &          &         &          & X_{n-1} &         &\\
 &  \swarrow &        &   \searrow  &     &  \swarrow & \cdots & \searrow &         & \swarrow &         &\searrow &\\
X&          &         &             & Y_1 &           &        &          & Y_{n-1} &          &         & &Y_{n-2}=Y\\
\endmatrix
$$
such that each morphism is locally a product of blowups of nonsingular subvarieties. 
That is, for every $z\in X_{i+1}$, there
exist affine neighborhoods $W$ of $z$, $U$ of $x=\alpha_{i}(z)$, $V$ of $y=\beta_i(z)$ such that $\alpha_i:W\rightarrow U$ and
$\beta_i:W\rightarrow V$ are finite products of monoidal transforms.
\endproclaim

Theorem F is proved in dimension $3$ by the author in [C]. The proof of Theorem F
is exactly the same, with the use of
Theorem D from this paper, which is valid in  all dimensions.
By Theorem D, for each valuation of the function field $K$ of $X$, there exist local rings for which the conclusions of Theorem
D hold. These local rings can be extended to affine varieties which are related by products of monoidal transforms.
By the quasi-compactness of the Zariski manifold (Theorem VI.17.40 [ZS]) all valuations of $K$ are centered at
finitely many of these affine constructions. We can then patch these affine varieties along the open sets where
they are isomorphic to get complete $k$-schemes as desired.

The Monomialization Theorem is a "resolution of singularities" type problem. Some of the difficulties which arise in it
are related to those which appear in the  problems of resolution of (char. 0) vector fields (c.f [Ca], [Se]),
 and in resolution of singularities in characteristic $p>0$ (c.f. [Ab4], [Co], [G], [L2]).
  Resolution of vector fields  is an open
problem (locally )in  dimension $\ge 4$ and is open (globally) in dimension $\ge 3$.
 Resolution of singularities in characteristic $p>0$
 is an open problem in dimension $\ge 4$.    

Some of the many important papers which are directly concerned with the global factorization problem are
Hironaka [H1], Danilov [D], Crauder [Cr], Pinkham [P].

An important special case where the global factorization problem has been solved is toric geometry.
The solution is in the series of papers Danilov [D2], Ewald [E], (dim 3) and Wlodarczyk [W], Morelli [M], 
Abramovich, [AMR] (dim n).

A birational morphism of nonsingular toric varieties can be thought of as a union of monomial mappings on affine
spaces. In toric geometry, the global factorization problem becomes more tractable than in the general
case of arbitrary polynomial mappings, since the problem can be translated into combinatorics.

Morelli's main result [M] is that a birational morphism of proper nonsingular  
toric varieties can be factored by one sequence of blowups (with nonsingular centers) followed
by one sequence of blowdowns (with nonsingular centers). [AMR] addresses some gaps and difficulties in the proof, 
and extends the result to toroidal morphisms.

Our main result, Theorem A - Monomialization, allows us to reduce the  factorization problem (locally) to
monomial mappings. If we then  make use of Morelli's result, which says (locally) that a birational monomial mapping can
be factored by one sequence of blowups, followed by one sequence of blowdowns, we obtain an even stronger local factorization 
theorem than  Theorem D.

Abhyankar has conjectured (page 237 [Ab5], [Ch]) that in all dimensions
it is possible to factor a birational mapping  along a valuation by a sequence of blowups followed by a sequence of blowdowns 
with nonsingular centers. 
This is the most optimistic possible local statement.

We prove the following Theorem, which proves Abhyankar's conjecture in all dimensions (over fields of characteristic 0).

\vfill\eject

\proclaim{Theorem G} (Factorization 2)
 Suppose that  $R\subset S$ are excellent regular local rings of dimension $n$,
containing a field $k$ of characteristic zero and  with a common quotient field $K$. 
  Let $V$ be a valuation ring of
$K$ which dominates $S$. Then there exists a regular local ring $T$, with quotient field $K$,
such that $T$ dominates $S$, $V$ dominates  $T$, and the inclusions $R\rightarrow T$ and $S\rightarrow T$
can be factored by sequences of monoidal transforms (blowups of regular primes).
$$
\matrix
&&V&&\\
&&\uparrow&&\\
&&T&&\\
&\nearrow&&\nwarrow&\\
R&&\longrightarrow && S
\endmatrix
$$
\endproclaim
\vskip .2truein

The solution  to Abhyankar's conjecture (as stated in [Ch]) is given  in Theorem H.

\proclaim{Theorem H} Suppose that $K$ is a field of algebraic functions over a
 field $k$ of characteristic zero, with 
$\text{trdeg}_kK=n$, $R$ and $S$ are regular local
rings, essentially of finite type over $k$, with quotient field $K$.  Let $V$ be a valuation ring of
$K$ which dominates $R$ and $S$. Then there exists a regular local ring $T$, essentially of finite type over $k$,
with quotient field $K$, dominated by $V$, containing $R$ and $S$,
such that   $R\rightarrow T$ and $S\rightarrow T$
can be factored by  products of monoidal transforms.
\endproclaim

In dimension 3, Theorems G and H have been proven  by the author in [C].
Theorem A, which shows that it is possible to monomialize a generically finite morphism along a valuation, is essential
in this proof.

Hironaka and Abhyankar (section 6 of chapter 0 [H] and page 254 [Ab5]) have conjectured that a birational
morphism of nonsingular projective varieties can be factored by a series of blowups followed by a series of blowdowns
with nonsingular centers. 

Our Theorem G shows that there is no local obstruction to this global factorization conjecture in any dimension.

We prove the following glocal analogue of Theorem G.

\proclaim{Theorem I} Let $k$ be a field of characteristic zero, $\phi:X\rightarrow Y$ a birational
morphism of integral nonsingular proper excellent $k$-schemes. Then there exists a nonsingular complete $k$-scheme
$Z$ and birational complete morphisms $\alpha:Z\rightarrow X$ and $\beta:Z\rightarrow Y$ making the diagram
$$
\matrix
&&Z&&\\
&\swarrow&&\searrow&\\
X&&\longrightarrow && Y
\endmatrix
$$
commute, such that $\alpha$ and $\beta$ are locally products of monoidal transforms. That is, for every $z\in Z$, there
exist affine neighborhoods $W$ of $z$, $U$ of $x=\alpha(z)$, $V$ of $y=\beta(z)$ such that $\alpha:W\rightarrow U$ and
$\beta:W\rightarrow V$ are finite products of monoidal transforms.
\endproclaim

Here a monoidal transform of a nonsingular $k$-scheme $S$ is the map $T\rightarrow S$ induced by an open subset $T$ of
$\text{Proj}(\oplus \Cal I^n)$, where $\Cal I$ is the ideal sheaf of a nonsingular subvariety of $S$.

Theorem I is proved in dimension $3$ by the author in [C]. The proof is exactly the same, with the use of
Theorem G from this paper, which is valid in  all dimensions.

\subheading{Geometry and valuations} 

A valuation ring of a field of algebraic functions $K$ will dominate some local ring of a projective model $V$ of $K$.
 This leads to the "valuative criterion for properness" (c.f. Theorem II.4.7 [Ha]).

The Zariski manifold $M$ of $K$ is a locally ringed space whose local rings are the valuations rings of $K$,
containing the ground field $k$ (c.f chapter VI, section 17 [ZS], [L1], section 6 of chaper 0 [H]).
$M$ satisfies the universal property that for any morphism of proper $k$-schemes $\Psi:X\rightarrow Y$ such that $X$ and $Y$ have
function fields (isomorphic to) $K$, there are projections $\pi_1:M\rightarrow X$ and 
$\pi_2:M\rightarrow Y$ making a commutative diagram
$$
\matrix
&&M&&\\
&\swarrow&&\searrow&\\
X&&\longrightarrow && Y.
\endmatrix
$$
When $K$ is a 1-dimensional function field, the only nontrivial valuation rings are the local rings of the
points on the nonsingular model of $K$. As such, a projective nonsingular curve can be identified with its
Zariski manifold (c.f. I.6 [Ha]). If $K$ has dimension $>1$, $K$ has many non-noetherian valuations, and
$M$ is far from being a $k$-scheme.

The main result of Zariski in [Z2] is his Theorem $U_1$, which states that for a valuation $B$ of a field of
algebraic functions $K$ over a ground field of characteristic 0, there is a projective model $V$ of $K$ on which the center of $B$
is at a nonsingular point of $V$.

Our Theorems A, D, G and H are direct analogues of Theorem $U_1$ for generically finite and birational
morphisms of varieties.

Zariski obtained a solution to "the classical problem of local uniformization" from his Theorem $U_1$. 
In the language of schemes (c.f. section 6 of chapter 0 [H]) Zariski's result shows that for any
integral proper $k$-scheme $X$ (where $k$ is a field of characteristic 0) there exists a complete nonsingular
integral $k$-scheme $Y$ and a birational morphism $Y\rightarrow X$. A complete variety over a field $k$
is an integral finite type $k$-scheme which satisfies the existence part of the valuative criterion for properness.

Our Theorems B, B2, F and I are analogous to Zariski's solution of "the classical problem of local uniformization".

\subheading{Overview of the proof}

The main thrust of the paper is to acheive monomialization.
Theorem A proves monomialization for generically finite extensions. 
The corollaries,  Theorems B through I are then easily obtained.

Theorem A is an immediate corollary of Theorem 4.3.
 Theorem  4.4 is a stronger version, valid for birational
extensions.

In fact, Theorems 4.3 and 4.4 prove more than monomialization. They produce a matrix of exponents $A=(a_{ij})$ which
has a very special form, depending on the rational rank of the rank 1 valuations composite with $V$.

Theorem 4.4  reduces the proof of
Theorem D (Local factorization) to the special case where $\text{dim } R=\text{ dim }S=n$ and $V$ has rank 1 and rational rank n. 
Factorization in the special case $n=3$ and $V$ has rational rank $n=3$ was solved by Christensen in [Ch].
 We generalize Christensen's
algorithm in Theorem 5.4 to prove factorization when $V$ has rational rank $n$. The proof of Theorem 5.4 uses only elementary
methods of linear algebra. Theorem D then follows from Theorem 4.4.

Now we will discuss the proof of Theorem 4.3.
The most difficult part of Theorem 4.3 is the case where $\nu$ has rank 1, which is proved in Theorem 4.1. Almost
the entirety of the paper (chapter 3) is devoted to the proof of Theorem 4.1.

Suppose that $\nu$ has rank 1 and rational rank $s$. Then it is not difficult to construct sequences of
monoidal transforms $R\rightarrow R(1)$ and $S\rightarrow S(1)$ such that $\nu$ dominates $S(1)$, $S(1)$ dominates $R(1)$,
$R(1)$ has regular parameters $(x_1(1),\ldots,x_n(1))$, $S(1)$ has regular parameters $(y_1(1),\ldots,y_n(1))$
such that
$$
\align
x_1(1)&=y_1(1)^{c_{11}(1)}\cdots y_s(1)^{c_{1s}(1)}\delta_1\\
&\vdots\\
x_s(1)&=y_1(1)^{c_{s1}(1)}\cdots y_s(1)^{c_{ss}(1)}\delta_s
\endalign
$$
where $\text{det}(c_{ij}(1))\ne 0$ and $\delta_i$ are units in $S(1)$.
This step is accomplished in the proof of Theorem 4.1.

The inductive step in the proof is Theorem 3.12, which starts with monoidal transform sequences (MTSs)
$R\rightarrow R(0)$ and $S\rightarrow S(0)$ such that $\nu$ dominates $S(0)$, $S(0)$ dominates $R(0)$,
$R(0)$ has regular parameters $(x_1(0),\ldots,x_n(0))$, $S(0)$ has regular parameters
$(y_1(0),\ldots,y_n(0))$ such that
$$
\align
x_1(0)&=y_1(0)^{c_{11}(0)}\cdots y_s(0)^{c_{1s}(0)}\delta_1\tag 1\\
&\vdots\\
x_s(0)&=y_1(0)^{c_{s1}(0)}\cdots y_s(0)^{c_{ss}(0)}\delta_s\\
x_{s+1}(0)&=y_{s+1}(0)\\
&\vdots\\
x_l(0)&=y_l(0)
\endalign
$$
where $\text{det}(c_{ij}(0))\ne 0$ and $\delta_i$ are units in $S(0)$, and
construct MTSs $R(0)\rightarrow R(t)$, $S(0)\rightarrow S(t)$ such that $\nu$ dominates $S(t)$, $S(t)$ dominates $R(t)$,
$R(t)$ has regular parameters $(x_1(t),\ldots,x_n(t))$, $S(t)$ has regular parameters $(y_1(t),\ldots,y_n(t))$ such that
$$
\align
x_1(t)&=y_1(t)^{c_{11}(t)}\cdots y_s(t)^{c_{1s}(t)}\delta_1\tag 2\\
&\vdots\\
x_s(t)&=y_1(t)^{c_{s1}(t)}\cdots y_s(t)^{c_{ss}(t)}\delta_s\\
x_{s+1}(t)&=y_{s+1}(t)\\
&\vdots\\
x_{l+1}(t)&=y_{l+1}(t)
\endalign
$$
where $\text{det}(c_{ij}(t))\ne 0$ and $\delta_i$ are units in $S(t)$.

To prove  Theorem 3.12, we make use of special sequences of monoidal transforms which are derived from the Cremona tranformations
constructed by Zariski in chapter B of [Z2], using an algorithm of Perron. We will call such transformations Perron transforms.
Our proof makes use of these transforms in local rings of etale extensions, giving the transforms the special form
$$
\align
x_1(i)&=x_1(i+1)^{a_{11}(i+1)}\cdots x_s(i+1)^{a_{1s}(i+1)}c_{i+1}^{a_{1,s+1}(i+1)}\tag 3\\
&\vdots\\
x_s(i)&=x_1(i+1)^{a_{s1}(i+1)}\cdots x_s(i+1)^{a_{ss}(i+1)}c_{i+1}^{a_{s,s+1}(i+1)}\\
x_r(i)&=x_1(i+1)^{a_{s+1,1}(i+1)}\cdots x_s(i+1)^{a_{s+1,s}(i+1)}(x_r(i+1)+1)c_{i+1}^{a_{s+1,s+1}(i+1)}
\endalign
$$
\vskip .2truein
$$
\align
y_1(i)&=y_1(i+1)^{b_{11}(i+1)}\cdots y_s(i+1)^{b_{1s}(i+1)}d_{i+1}^{b_{1,s+1}(i+1)}\\
&\vdots\\
y_s(i)&=y_1(i+1)^{b_{s1}(i+1)}\cdots y_s(i+1)^{b_{ss}(i+1)}d_{i+1}^{b_{s,s+1}(i+1)}\\
y_r(i)&=y_1(i+1)^{b_{s+1,1}(i+1)}\cdots y_s(i+1)^{b_{s+1,s}(i+1)}(y_r(i+1)+1)d_{i+1}^{b_{s+1,s+1}(i+1)}
\endalign
$$
where $\text{det}(a_{ij}(i+1))=\pm1$ and $\text{det}(b_{ij}(i+1))=\pm1$, $c_{i+1}$, $d_{i+1}$ are algebraic over $k$.

Zariski observes on page 343 of [Z1] that his Cremona transformations have "the same effect as the classical
Puiseux substitution $x=x_1^{\nu},y=x_1^{\mu}(c_1+y_1)$ used in the determination of the branches of the curve $\phi(x,y)=0$."
"The only difference - and advantage - is that our transformation does not lead to elements $x_1,y_1$ outside the field $k(x,y)$."

Our transforms (3) do induce  a field extension. They are the direct generalization of the classical Puiseux substitution
 to higher dimensions. We must
pay for the advantage of the simple form of the equations by introducing many difficulties arising from the need to make finite
etale extensions after each transform. We call a sequence of such transforms a uniformizing transform sequence (UTS).

Theorem 3.12 is proved by first constructing  UTSs such that (2) holds, and then using this partial solution
to construct sequences of monoidal transforms such that (2) holds. The UTSs are constructed in Theorems 3.8
and 3.12, and this is used to construct MTSs such that (2) holds in Theorems 3.9, 3.10, 3.11 and 3.12.

Underlying the whole proof is Zariski's algorithm for the reduction of the multiplicity of a polynomial along 
a rank 1 valuation, via Perron transforms. This algorithm is itself a a generalization of Newton's
algorithm to determine the branches of a curve singularity.

We will now give an outline of the proof of Theorem 3.12 (the inductive step). 
We will  give a formal construction, so that we need only consider UTSs, where the basic ideas are transparent.
  We will construct UTSs along $\nu$ (here the $U(i)$, $T(i)$ are complete local rings)
$$
\matrix
U(0)&\rightarrow&U(1)&\rightarrow&\cdots&\rightarrow&U(t)\\
\uparrow&&\uparrow&&&&\uparrow\\
T(0)&\rightarrow&T(1)&\rightarrow&\cdots&\rightarrow&T(t)
\endmatrix\tag 4
$$
such that $T(i)$ has regular parameters $(x_1(i),\ldots,x_n(i))$, $U(i)$ has regular parameters $(y_1(i),\ldots,y_n(i))$
such that
$$
\align
x_1(i)&=y_1(i)^{c_{11}(i)}\cdots y_s(i)^{c_{1s}(i)}\tag 5\\
&\vdots\\
x_s(i)&=y_1(i)^{c_{s1}(i)}\cdots y_s(i)^{c_{ss}(i)}\\
x_{s+1}(i)&=y_{s+1}(i)\\
&\vdots\\
x_{l}(i)&=y_{l}(i)
\endalign
$$
where $\text{det}(c_{ij}(i))\ne 0$  for $0\le i\le t$.

We will presume that $k$ is algebraically closed, and isomorphic to   the residue fields of $R$, $S$ and $V$.
We will also assume that
 various technical difficulties, such as the 
rank of $\nu$ increasing when $\nu$ is extended to the complete local ring  $U(i)$, do not occur.

In Theorem 3.8 it is shown that

\item{(6)} Given $f\in U(0)$, there exists (4) such that $f=y_1(t)^{d_1}\cdots y_s(t)^{d_s}\gamma$  
where $\gamma$ is a unit in $U(t)$.
\item{(7)} Given $f\in U(0)-k[[y_1,\ldots,y_l]]$, there exists (4) such that 
$$
f=P(y_1(t),\ldots,y_l(t))+y_1(t)^{d_1}\cdots y_s(t)^{d_s}y_{l+1}(t).
$$
where $P$ is a power series.
\vskip .2truein
Theorem 3.12 then shows that it is possible to construct a UTS (4) such that $x_{l+1}(t)=y_{l+1}(t)$, 
which allows us to conclude the
truth of the inductive step.

We will now give a more detailed analysis of these important steps. For simplicity, we will assume that 
$s=\text{rat rank}(\nu) =1$. This is the essential case.

Suppose that $R$ has regular parameters $(x_1,\ldots,x_n)$, and $2\le i\le n$.
Zariski constructed (in [Z1], and in a generalized form in [Z2]) a
MTS $R\rightarrow R(1)$ where $R(1)$ has regular parameters
$(x_1(1),x_2(1),\ldots,x_n(1))$ by the following method. Since $\nu$ has rational rank 1, we can identify the value group of $\nu$
 with a subgroup of $\bold R$.
$$
\nu\left(\frac{x_i}{x_1}\right)=\frac{a_{i1}(1)}{a_{11}(1)}
$$
where $a_{i1}(1),a_{11}(1)$ are relatively prime positive integers. We can then choose positive integers
$a_{i1}(1), a_{ii}(1)$ such that  $a_{11}(1)a_{ii}(1)-a_{1i}(1)a_{i1}(1)=1$. Then
$$
\align
&\nu(x_1^{a_{ii}(1)}x_i^{-a_{1i}(1)})>0\\
&\nu(x_1^{-a_{i1}(1)}x_i^{a_{11}(1)})=0
\endalign
$$  
There exits then a uniquely determined, nonzero $c_1\in k$ such that 
$$
\nu(x_1^{-a_{i1}(1)}x_i^{a_{11}(1)}-c_1)>0.
$$
We can define $x_j(1)$ for $1\le j\le n$ by
$$
\align
x_1&=x_1(1)^{a_{11}(1)}(x_i(1)+c_1)^{a_{1i}(1)}\tag 8\\
x_i& =x_1(1)^{a_{i1}(1)}(x_i(1)+c_1)^{a_{ii}(1)}\\
x_j&=x_j(1)\text{ if }j\ne 1\text{ or }i.
\endalign
$$
Set $R(1) = R[x_1(1),x_i(1)]_{(x_1(1),\ldots,x_n(1))}$.

Using such transformations, Zariski proves 
\proclaim{Theorem I1}(Zariski [Z1], [Z2]) Given $f\in R$, there exists a MTS along $\nu$
$$
R\rightarrow R(1)\rightarrow \cdots\rightarrow R(t)
$$
such that $f=x_1(t)^d\gamma$ where gamma is a unit in $R(t)$.
\endproclaim

In our analysis, we will consider the transformation (8) in formal coordinates. $\widehat R(1)$
has regular parameters $(\overline x_1(1),\overline x_2(1),\ldots,\overline x_n(1))$ defined by
$$
\align
\overline x_1(1)&= x_1(1)(x_i(1)+ c_1)^{\frac{a_{1i}(1)}{a_{11}(1)}}\\
\overline x_i(1)&=(x_i(1)+c_1)^{\frac{1}{a_{11}(1)}}-c_1^{\frac{1}{a_{11}(1)}}\\
\overline x_j(1)&=x_j(1)\text{ if }j\ne 1\text{ or }i.
\endalign
$$
In these coordinates, (8) becomes
$$
\align
x_1&=\overline x_1(1)^{a_{11}(1)}\\
x_i&=\overline x_1(1)^{a_{i1}(1)}(\overline x_i(1)+c_1^{\frac{1}{a_{11}(1)}})\\
x_j&=\overline x_j(1)\text{ if }j\ne 1\text{ or }i.
\endalign
$$

We have inclusions
$$
 T'(1)=R(1)\rightarrow  T''(1)\rightarrow T(1)=\widehat{ R(1)}
$$
where 
$$
 T''(1)= R(1)[(x_i(1)+ c_1)^{\frac{1}{a_{11}(1)}}]_{(\overline x_1(1),\ldots,\overline x_n(1))}
$$
is a localization of a finite etale extension of $R(1)$. We can extend $\nu$ to a valuation of the quotient
field of $\widehat{ R(1)}$ which dominates $\widehat{ R(1)}$. For simplicity, we will assume that this extension still
has rank 1.

We will construct  sequences of UTSs
$$
\matrix
 T'(1)&\rightarrow & T''(1)&\rightarrow &  T(1)=\widehat { T'(1)}&&&\\
&&\downarrow&&&&&\\
&& T'(2)&\rightarrow & T''(2)&\rightarrow &  T(2)=\widehat { T'(2)}&\\
&&&&\downarrow&&\\
&&&& T'(3)&\cdots&\\
&&&&&\vdots&
\endmatrix
$$
where each downward arrow is of the form (8). 

To prove the inductive step, we must construct  UTSs (4) starting with $R\rightarrow S$. 
By induction, we may assume that $R$ has regular parameters $(x_1,\ldots,x_n)$ and $S$ has regular parameters
$(y_1,\ldots,y_n)$ such that
$$
\align
x_1&=y_1^{t_0}\delta_1\\
x_2&=y_2\\
&\vdots\\
x_l&=y_l
\endalign
$$
where $\delta_1$ is a unit in $S$. (Recall that we are assuming that $s=1$.)
 By Hensel's lemma, $\delta_1^{\frac{1}{t_0}}$ is a unit in $\hat S$.
We then can start our sequence of UTSs by 
setting 
$ U''(0)=S[\delta_1^{\frac{1}{t_0}}]_{(\overline y_1,\ldots,\overline y_n)}$ and $ T''(0) = R$,
with regular parameters $(\overline x_1,\ldots,\overline x_n)$ and $(\overline y_1,\ldots,\overline y_n)$, 
such that 
$$
\align
\overline x_1&=\overline y_1^{t_0}\\
\overline x_2&=\overline y_2\\
&\vdots\\
\overline x_l&=\overline y_l
\endalign
$$
We will construct 2 types of  UTSs
$$
\matrix
 U''(0)&\rightarrow & U''(1)\\
\uparrow&&\uparrow\\
 T''(0)&\rightarrow & T''(1)
\endmatrix
$$
A transformation of Type I is defined  when $2\le i\le l$.  The equations defining the horizontal maps are then
$$
\align
\overline x_1&=\overline x_1(1)^{a_{11}(1)}\\
\overline x_i& = \overline x_1(1)^{a_{i1}(1)}(\overline x_i(1)+c_1)\\
\overline x_j&=\overline x_j(1)\text{ if }j\ne 1\text{ or }i,
\endalign
$$
$$
\align
\overline y_1&=\overline y_1(1)^{b_{11}(1)}\\
\overline y_i &= \overline y_1(1)^{b_{i1}(1)}(\overline y_i(1)+d_1)\\
\overline y_j&=\overline y_j(1)\text{ if }j\ne 1\text{ or }i.
\endalign
$$
$$
\frac{a_{i1}(1)}{a_{11}(1)} = \frac{\nu(\overline x_i)}{\nu(\overline x_1)}=\frac{\nu(\overline y_i)}{t_0\nu(\overline y_1)}
=\frac{b_{i1}(1)}{t_0b_{11}(1)}
$$
$(a_{11}(1),a_{i1}(1))=1$ implies $a_{11}(1)\mid t_0b_{11}(1)$.

We thus have
$$
\align
\overline x_1(1)&=\overline y_1(1)^{\frac{t_0b_{11}(1)}{a_{11}(1)}}=\overline y_1(1)^{t_1}\\
\overline x_2(1)&=\overline y_2(1)\\
&\vdots\\
\overline x_l(1)&=\overline y_l(1).
\endalign
$$
A transformation of Type II is defined  when $l<i$.  The equations defining the horizontal maps are then
$$
\overline x_j=\overline x_j(1)\text{ for }1\le j\le n
$$
$$
\align
\overline y_1&=\overline y_1(1)^{b_{11}(1)}\\
\overline y_i &= \overline y_1(1)^{b_{i1}(1)}(\overline y_i(1)+d_1)\\
\overline y_j&=\overline y_j(1)\text{ if }j\ne 1\text{ or }i.
\endalign
$$
In this case $ T''(0)= T''(1)$.

In this way, we can construct sequences of UTSs
$$
\matrix
S&\rightarrow &U(0)&\rightarrow &U(1)&\rightarrow&\cdots&\rightarrow &U(t)\\
\uparrow&&\uparrow&&\uparrow&&&&\uparrow\\
R&\rightarrow &T(0)&\rightarrow &T(1)&\rightarrow&\cdots&\rightarrow &T(t)
\endmatrix\tag A
$$
such that $T(k)$ has regular parameters $(\overline x_1(k),\ldots,\overline x_n(k))$
and $U(k)$ has regular parameters $(\overline y_1(k),\ldots,\overline y_n(k))$
related by
$$
\align
\overline x_1(k)&=\overline y_1(k)^{t_k}\tag B\\
\overline x_2(k)&=\overline y_2(k)\\
&\vdots\\
\overline x_l(k)&=\overline y_l(k)
\endalign
$$
for $0\le k\le t$. The transformations $T(k)\rightarrow T(k+1)$ and 
$U(k)\rightarrow U(k+1)$
are of type I or II, and we also allow changes of variables, replacing $\overline x_i(k)$
with $\overline x_i(k)-P(\overline x_1(k),\ldots,\overline x_{i-1}(k))$ and
replacing $\overline y_i(k)$
with $\overline y_i(k)-P(\overline x_1(k),\ldots,\overline x_{i-1}(k))$ if $2\le i\le l$,
for some power series $P$, and we
may replace $\overline y_i(k)$
with $\overline y_i(k)-P(\overline y_1(k),\ldots,\overline y_{i-1}(k))$ if $l<i$.

To prove the induction step,
we must prove Theorems I2 and I3 below.

\proclaim{Theorem I2}(Theorem 3.8 with $s=\text{ rat rank }\nu=1$)) 
\item{(1)} Given $f
\in U(0)$, there exists a sequence (A) such that  
$$
f=\gamma \overline y_1(t)^d
$$
with $\gamma$ a unit in $U(t)$. (If $f\in k[[\overline y_1,\ldots,\overline y_{\tau}]]$, the transformations of type
I and II in the sequence involve only the first $\tau$ variables.)
\item{(2)} Suppose that $f\in k[[\overline y_1,\ldots,\overline y_m]]-k[[\overline y_1,\ldots,\overline y_l]]$.
Then there exists a sequence (A) such that 
$$
f=P(\overline y_1(t),\ldots,\overline y_l(t))+\overline y_1(t)^{d_1}\overline y_m(t).
$$
\proclaim{Theorem I3}(Theorem 3.12)
there exists a sequence (A) such that $\overline x_{l+1}(t)=\overline y_{l+1}(t)$.
\endproclaim

\subheading{Outline of proof of (1) of Theorem I2} 
The proof is by induction on $\tau$. Suppose that (1) is true for $\tau-1$. We will assume that $\tau\le l$, which
is the essential case. Recall that
$$
\align
\overline x_1&=\overline y_1^{t_0}\\
\overline x_2&=\overline y_2\\
&\vdots\\
\overline x_l&=\overline y_l.
\endalign
$$
Let $\omega$ be a primitive $t_0^{th}$ root of unity. Set
$$
g(\overline x_1,\ldots,\overline x_{\tau}) 
= \prod_{i=0}^{t_0-1}f(\omega^i\overline y_1,\overline y_2,\ldots,\overline y_{\tau}).
$$
$f\mid g$ in $U(0)$. We will perform a UTS (A) to get $g=\overline x_1(t)^d\Lambda$ where $\Lambda$ is a unit in $T(t)$. Then
$f = \overline y_1(t)^{d_1'}\Lambda'$ where $\Lambda'$ is a  unit in $U(t)$.

To transform $g$ into the form $g=\overline x_1(t)^d\Lambda$, where $\Lambda$ is a unit,
 we will make use of an algorithm of Zariski ([Z1], [Z2])
to reduce the multiplicity of $g$.
Initially set $g= \overline x_1^bg_0$ where $\overline x_1$ does not divide $g_0$. Set 
$$
r=\text{mult}(g_0(0,\ldots,0,\overline x_{\tau})).
$$
$0\le r\le \infty$. If $r=0$, $g_0$ is a unit, and we are done. Suppose that $0<r$. We can write
$$
\align
g_0 &= \sum_{i=0}^{\infty}a_i(\overline x_1,\ldots,\overline x_{\tau-1})\overline x_{\tau}^i\\
& = \sum_{i=1}^da_{\alpha_i}\overline x_{\tau}^{\alpha_i}+\sum_ja_{\beta_j}\overline x_{\tau}^{\beta_j}
+\overline x_{\tau}^N\Omega
\endalign
$$
 where the terms $a_{\alpha_i}\overline x_{\tau}^{\alpha_i}$ have minimum value $\rho$, $N\nu(\overline x_{\tau})>\rho$,
and the $a_{\beta_j}\overline x_{\tau}^{\beta_j}$ terms are the finitely many remaining terms.
We must have $\alpha_i\le r$ for all $\alpha_i$, and $\alpha_i=r$ implies $\alpha_i$ is a unit.

By induction, we can perform UTSs in the first $\tau-1$ variables  to reduce to the case
$$
\align
a_{\alpha_i} &= \overline x_1^{\gamma_i}u_{\alpha_i}(\overline x_1,\ldots,\overline x_{\tau-1})\\
a_{\beta_j} &= \overline x_1^{\delta_j}u_{\beta_j}(\overline x_1,\ldots,\overline x_{\tau-1})\\
\endalign
$$

where $u_{\alpha_i}$ and $u_{\beta_j}$ are units in $T(0)$. 
Now we make a UTS
$$
\align
\overline x_1&=\overline x_1(1)^{a_{11}}\\
\overline x_{\tau}&=\overline x_1(1)^{a_{\tau 1}}(\overline x_{\tau}(1)+c_1).
\endalign
$$

$$
\align
g_0 &= \sum_{i=1}^d\overline x_1^{\gamma_ia_{11}+\alpha_ia_{\tau 1}}u_{\alpha_i}(\overline x_{\tau}(1)+c_1)^{\alpha_i}+\cdots\\
&= \overline x_1(1)^{\epsilon}(\sum_{i=1}^du_{\alpha_i}(\overline x_{\tau}(1)+c_1)^{\alpha_i}+\overline x_1\Omega).
\endalign
$$

Set
$$
g_1=\sum_{i=1}^du_{\alpha_i}(\overline x_{\tau}(1)+c_1)^{\alpha_i}+\overline x_1\Omega,
$$
$$
r_1 = \text{mult}(g_1(0,\ldots,0,\overline x_{\tau}(1)).
$$
$r_1<\infty$ and $r_1\le r$ since all $\alpha_i\le r$. Set
$$
\zeta(t) = g_1(0,\ldots,0,t-c_1) = \sum_{i=1}^du_{\alpha_i}(0,\ldots,0)t^{\alpha_i}.
$$
If we do not have a reduction in $r$, so that 
$r_1=r$,  
$$
\zeta(\overline x_{\tau}(1)+c_1) =e\overline x_{\tau}(1)^r
$$
for some nonzero $e\in k$. Thus
$\zeta(t) = e(t-c_1)^r$ has a nonzero $t^{r-1}$ term. We conclude that $\alpha_d=r$, $a_{\alpha_d}$ is a unit,
$\alpha_{d-1}=r-1$ and 
$$
\rho=\nu(a_{\alpha_d}\overline x_{\tau}^r) = \nu(a_{\alpha_{d-1}}\overline x_{\tau}^{r-1}).
$$
Thus 
$$
\nu(\overline x_{\tau}) = \nu(a_{\alpha_{d-1}}(\overline x_1,\ldots,\overline x_{\tau-1})).
$$
Since $\overline x_{\tau}^r$ is a minimum value term of $f$,
$$
\nu(\overline x_{\tau})\le \nu(\overline x_{\tau}^r)\le \nu(f).
$$
Now  make a change of variables in $ T(0)$, replacing
$\overline x_{\tau}$ with $\overline x_{\tau}'=\overline x_{\tau}-\lambda a_{\alpha_{d-1}}$
where $\lambda\in k$ is chosen to make $\nu(\overline x_{\tau})< \nu(\overline x_{\tau}')$.
Repeat the above procedure with these new variables. We eventually get a reduction in $r$. 
In fact if we didn't, we would have an infinite bounded sequence in $T(0)$
$$
\nu(\overline x_{\tau})<\nu(\overline x_{\tau}')<\cdots<\nu(f)
$$
which is impossible (by Lemma 1.3).

\subheading{Outline of proof of (2) of Theorem I2}
$$
f = \sum_{i=0}^{\infty}\sigma_i(\overline y_1,\ldots,\overline y_{m-1})\overline y_m^i.
$$
Set 
$$
Q = \sum_{i>0} \sigma_i(\overline y_1,\ldots,\overline y_{m-1})\overline y_m^i.
$$
After possibly permuting the variables $\overline y_{l+1},\ldots,\overline y_m$, we can assume that $Q\ne 0$.
 $Q = \overline y_1^{n_1}Q_0$. where $\overline y_1$ does not divide $Q_0$.
Set $r = \text{ mult }( Q_0(0,\ldots,0,\overline y_m))$. $1\le r\le \infty$. 

Suppose that $r>1$. Write
$$
Q_0 = \sum_{i=1}^d \sigma_{\alpha_i}(\overline y_1,\ldots,\overline y_{m-1})\overline y_m^{\alpha_i}+\cdots 
$$
where the $\sigma_{\alpha_i}\overline y_m^{\alpha_i}$ are the minimum value terms. By construction, all $\alpha_i>0$.
By (1) of Theorem I1, we can perform UTSs in the first m-1 variables to get $\sigma_i=u_{\alpha_i}\overline y_1^{\gamma_i}$
where $u_{\alpha_i}$ are units.
Then we can  perform a UTSs in $\overline y_1$ and $\overline y_m$ to get
$$
Q_0 = \overline y_1(1)^{\epsilon}(\sum_{i=1}^d\sigma_{\alpha_i}(\overline y_m(1)+c_1)^{\alpha_i}+\overline y_1(1)\Omega). 
$$
Set 
$$
Q_1 = \sum_{i=1}^d\sigma_{\alpha_i}(\overline y_m(1)+c_1)^{\alpha_i}+\overline y_1(1)\Omega-\sum_{i=1}^du_{\alpha_i}c_1^{\alpha_i}.
$$
Set $r_1 = \text{ mult }(Q_1(0,\ldots,0,\overline y_m(1))$.
$0<r_1<\infty$ and $r_1\le r$.
Suppose that we do not have a reduction in $r$, so that
$r_1=r$. Then as in the proof of (1) of Theorem I2,
  $\nu(\overline y_m)=\nu(\sigma_{\alpha_{d-1}}(\overline y_1,\ldots,\overline y_{m-1}))$. Now make a change of variables in $U(0)$,
replacing $\overline y_m$ with 
$$
\overline y_m'=\overline y_m-\lambda\sigma_{\alpha_{d-1}}(\overline y_1,\ldots,\overline y_{m-1}),
$$
where $\lambda\in k$ is chosen to make $\nu(\overline y_m)<\nu(\overline y_m')$. We have
$$
\nu(\overline y_m)\le \nu(\overline y_m^{r-1})\le \nu\left(\frac{\partial Q_0}{\partial \overline y_m}\right)
\le \nu\left(\frac{\partial f}{\partial \overline y_m}\right).
$$
Repeat the above algorithm with $\overline y_m$ replaced with $\overline y_m'$ in $U(0)$.
Since $\frac{\partial f}{\partial \overline y_m} = \frac{\partial f}{\partial \overline y_m'}$,
we  get a reduction in r after finitely many iterations, by an argument similar to that at the end of the proof
of (1) of Theorem I2.  

We can then repeat this procedure to eventually get an expression
$$
f = L(\overline y_1,\ldots,\overline y_{m-1})+\overline y_1^{n_0}\overline Q
$$
where $\text{ mult }\overline Q(0,\ldots,0,\overline y_m)=1$.

By induction, we can perform UTSs in the first $m-1$ variables only to get
$$
L = L'(\overline y_1,\ldots,\overline y_l)+\overline y_1^{n_1}\overline Q_1
$$
where $\text{ mult }(\overline Q_1(0,\ldots,0,\overline y_{m-1})=1$. Then $f$ is in the desired form.

\subheading{Outline of proof of Theorem I3}

By (2) of Theorem I2 we may assume that
$$
\align
\overline x_1&=\overline y_1^{t_0}\\
\overline x_2&=\overline y_2\\
&\vdots\\
\overline x_l&=\overline y_l\\
\overline x_{l+1}&=P(\overline y_1,\ldots,\overline y_l)+\overline y_1^{d_1}\overline y_{l+1}.
\endalign
$$
Let $\omega$ be a primitive $t_0^{th}$ root of unity. Set
$$
g(\overline x_1,\ldots,\overline x_{l+1}) = 
\prod_{i=0}^{t_0-1}(\overline x_{l+1}-P(\omega^i\overline y_1,\overline y_2,\ldots,\overline y_l)).
$$
$\overline y_{l+1}$ divides $g$ in $U(0)$. Set $r = \text{ mult }(g(0,\ldots,0,\overline x_{l+1}))$.
$1\le r<\infty$.

Suppose that $r=1$. Then in $T(0)$,
$$
\align
g&=\text{unit}(\overline x_{l+1}+\Phi(\overline x_1,\ldots,\overline x_l))\\
& = \text{unit}(P+\overline y_1^{d_1}\overline y_{l+1}+\Phi).
\endalign
$$
since $\overline y_{l+1}$ divides $g$, we must have $P=-\Phi$. We can then replace $\overline x_{l+1}$ with
$$
\overline x_{l+1}+\Phi = \overline y_1^{d_1}\overline y_{l+1},
$$
 which can be factored to achieve the
conclusions of  Theorem I3.

Now suppose that $r>1$. By (1) of Theorem I2, there exists a UTS in the first $l$ variables so that
$P = \overline y_1^{h_1}\overline P(\overline y_1,\ldots,\overline y_l)$, where $\overline P$ is a unit, and
$$
g= \sum_{i=1}^d\overline a_{\alpha_i}\overline x_1^{\gamma(i)}\overline x_{l+1}^{\alpha_i}+\cdots
$$
where the $\overline a_{\alpha_i}\overline x_{l+1}^{\alpha_i}$ 
are the minimum value terms and the $\overline a_{\alpha_i}(\overline x_1,\ldots,\overline x_l)$ are units.

\subheading{Case 1} Assume that $\nu(P)>\nu(\overline y_1^{d_1})$. Then
$\overline x_{l+1} = \overline y_1^{d_1}G$ where 
 $G = \overline y_1^{h_1-d_1}\overline P+\overline y_{l+1}$. Since $\text{ mult }(G(0,\ldots,0,\overline y_{l+1}) = 1$,
we can  factor this to get in the form of the conclusions of  Theorem I3.

\subheading{Case 2} Suppose that $\nu(P)\le \nu(\overline y_1^{d_1})$. Then
$\overline x_{l+1} = \overline y_1^{h_1}\overline P_1$
where $\overline P_1 = \overline P+\overline y_1^{d_1-h_1}\overline y_{l+1}$
is a unit. Perform a  UTS $T(0)\rightarrow T(1)$
defined by
$$
\align
\overline x_1 &= \overline x_1(1)^{a_{11}(1)}\\ 
\overline x_{l+1} &= \overline x_1(1)^{a_{l+1,1}(1)}(\overline x_{l+1}(1)+c_1).
\endalign
$$
We will show that this map factors through  $U(0)$. 
$$
\frac{\nu(\overline x_{l+1})}{\nu(\overline x_1)} = \frac{a_{l+1,1}(1)}{a_{11}(1)}
=\frac{h_1\nu(\overline y_{1})}{t_0\nu(\overline y_1)} 
= \frac{h_1}{t_0}.
$$
Thus $h_1 = t_1a_{l+1,1}(1)$ and $t_0=t_1a_{11}(1)$ for some positive integer $t_1$.
$$
\align
\overline x_1(1) &= \overline y_1^{t_1}\\
\overline x_{l+1}(1)+c_1&= \frac{\overline y_1^{h_1}}{\overline x_1(1)^{a_{l+1,1}}(1)}\overline P_1=\overline P_1
\endalign
$$
so that
$$
\overline x_{l+1} = P_1(\overline y_1,\ldots,\overline y_l)+\overline y_1^d\overline y_{l+1}.
$$
Set $g=\overline x_1(1)^eg_1$. $r_1 = \text{ mult } (g_1(0,\ldots,0,\overline x_{l+1}(1))\le r$.
If $r_1=r$, we can replace $\overline x_{l+1}$ with $\overline x_{l+1}-\sigma(\overline x_1,\ldots,\overline x_l)$
and repeat to eventually get $r_1<r$. 

We have $\overline y_{l+1}(1)=\overline y_{l+1}$ and $\overline y_{l+1}\mid g$, so that $\overline y_{l+1}(1)\mid g$.
$\overline x_1(1) = \overline y_1(1)^{t_1}$ implies $r_1>0$. Now we can repeat the above argument to eventually either reach
 $r=1$ or case 1.

\heading
Preliminaries
\endheading
\subheading{ Valuations}

\proclaim{Lemma 1.1} Suppose that $R$ is a  regular local ring, with quotient field $K$. Then $R=K\cap\hat R$ 
in the quotient field of $\hat R$.
\endproclaim
\demo{Proof}
Suppose that $f\in K\cap\hat R$. Then there exist $g,h\in R$ such that $f=\frac{g}{h}$, with $(g,h)=1$ in $R$. If $f\not\in R$, 
there exists an irreducible $s\in R$ such that $s\mid h$ but $s$ does not divide $g$. Let $s'\in \hat R$ be an 
irreducible such that $(s')\cap R=(s)$.  $hf=g$ in $\hat R$. $s'\mid h$ implies $s'\mid g$ in $\hat R$. hence 
$s\mid g$ in $R$, a contradiction.
\enddemo

\proclaim{Lemma 1.2} Let $R$ be a regular local ring with quotient field $K$, maximal ideal $m$.  
 Let $v$ be a rank 1 valuation of $K$ dominating $R$, with value group $\Gamma$, 
valuation ring $\Cal O_{\nu}$.
Let $\overline K$ be the completion of $K$ with respect to a metric  $|\cdot|$ associated to $v$. 
There exists a valuation $\overline v$ of $\overline K$ extending $v$, with valuation ring
 $\Cal O_{\overline v}$ such that
 $\Cal O_{\overline v}/m_{\overline v}=\Cal O_v/m_v$ and the 
value group of $\overline v$ is $\Gamma$ (c.f.  Theorems 1 and 2, Chapter 2 [Sch]). 

Then there exists a prime $p\in \hat R$ and an inclusion $\hat R/p\rightarrow \overline K$
which extends $R\rightarrow \overline K$.
\endproclaim

\demo{Proof} Let $\{a_n\}$ be a cauchy sequence in the $m$-adic topology of $R$. 
Let $v(m)=\rho>0$. Then $v(m^N)=N\rho$ implies $\{a_n\}$ is a fundamental sequence with respect to $|\cdot|$. 
Hence there is a natural map
$\phi:\hat R\rightarrow \overline K$ making
$$
\matrix
\hat R &  \phi\atop \rightarrow & \overline K\\
\uparrow & \nearrow &\\
R&&
\endmatrix
$$
commute. Let $P=\text{kernel } \phi$. 
\enddemo
\vskip .2truein

Lemma 1.3 extracts an argument from page 345 of [Z1].

\proclaim{Lemma 1.3} Let $R$ be a regular local ring containing a field of characteristic zero and $\nu$
 a rank 1  valuation of the quotient field of $R$ 
which has nonnegative value on $R$. 
 Suppose that $z_1,\ldots, z_n,\ldots$ is an infinite sequence of elements of $R$ such that
$$
\nu(z_1)<\nu(z_2)<\cdots<\nu(z_n)<\cdots
$$
is strictly increasing.
Then $\nu(z_n)$ has infinity for a limit. 
\endproclaim
\demo{Proof}
Let $m$ be the maximal ideal of $R$, $K$ the quotient field of $R$. Let $\{\overline t_i\}$ be a 
transcendence basis of $R/m$ over $k$. Lift $\overline t_i$ to $t_i\in R$. Let $L$ be the
field obtained by adjoining the $t_i$ to $k$. Then $L\subset R$. Let $L'$ be the algebraic closure of $L$ in $K$.
Then $L'\subset R$ since $R$ is normal. Let $\overline L$ be an algebraic closure of $L'$. 
$\overline K=K\otimes_{L'}\overline L$ is 
a field (c.f. Corollary 2, section 15, Chapter III [ZS]). Let $\overline \nu$ be an extension of 
$\nu$ to $\overline K$. $\overline \nu$
has rank 1 since $\overline K$ is algebraic over $K$. Let $\overline R$ be the localization of 
$R\otimes_{L'}\overline L$ at the center of $\overline \nu$. Then $\overline R$ is a regular local ring dominating $R$.
We can extend $\overline\nu$ to a valuation $\overline\nu$ dominating 
$\hat{\overline R}\cong \overline L[[x_1,\ldots,x_n]]$, a powerseries ring.

Let $\rho$ be a positive real number. Let $\sigma=\text{min}(\overline\nu(x_i))$. 
Let $n_{\rho}$ be the smallest integer such that $n_{\rho}\sigma>\rho$. Let 
$g(x_1,\cdots,x_n)\in \overline L[[x_1,\ldots,x_n]]$ be  
such that $\nu(g)\le \rho$. Write $g=g'+g''$ where $g'$ is a polynomial of degree $<n_{\rho}$, and $g''$ is a 
powerseries with terms of degree $\ge n_{\rho}$. Every form in $x,y,z$ of degree $m$ has value $\ge m\sigma$.
 Hence $\nu(g'')\ge n_{\rho}\sigma>\rho$. Since $\nu(g)\le \rho$, $\nu(g')=\nu(g)$. Thus if a powerseries has a value  
$\le\rho$, its value is the value of a  polynomial of degree $< n_{\rho}$. Hence,
 among the values assumed by 
elements of $\hat{\overline R}$, there is only a finite number of values which are less
 than or equal to a given fixed real number $\rho$.
\enddemo

\subheading{ Birational Transforms}

Suppose that $R$ is a regular local ring, with maximal ideal $m$, and that $x_1,\ldots, x_n\in R$ are such that 
$x_1,\ldots, x_n$ can be extended 
to a system of regular parameters $(x_1,\ldots,x_d)$ in $R$. Let $I$ be the ideal $I=(x_1,\ldots,x_n)$.

The blow up 
$$
\pi:\text{Proj}(\bigoplus_{n\ge0}I^n)\rightarrow\text{spec}(R)
$$
is called a monoidal transform of $\text{spec}(R)$. $\text{Proj}(\bigoplus_{n\ge0}I^n)$ is a regular scheme. let
$$
p\in\pi^{-1}(m)\subset\text{Proj}(\bigoplus_{n\ge0}I^n).
$$
 $p\in\text{spec}(R[\frac{x_1}{x_i},\cdots,\frac{x_n}{x_i}])$
for some $i$. Then
$$
R\rightarrow (R[\frac{x_1}{x_i},\cdots,\frac{x_n}{x_i}])_p
$$
is called a monoidal transform of $R$. If $n=d$, so that $I=m$,
$$
R\rightarrow (R[\frac{x_1}{x_i},\cdots,\frac{x_d}{x_i}])_p
$$
is called a quadratic transform.

In this section we state results of Abhyankar and Hironaka in a form which we will use.
The conclusions of Theorems 1.4 through 1.7 and Theorem 1.9 have been proved  by Hironaka [H] in equicharacteristic zero,
and have been proved by Abhyankar [Ab1], [Ab4] in positive characteristic,
 for varieties of dimension $\le 3$.

\proclaim{Definition 1.4} Let $R$ be a regular local ring. $f\in R$
is said to have simple normal crossings (SNCs), and be a SNC divisor,  if there exist
 regular parameters $(x_1,\ldots,x_n)$ in $R$ such that
$f=\text{unit }x_1^{a_1}\cdots x_n^{a_n}$ for some non-negative integers
$a_1,\ldots,a_n$.
\endproclaim

\proclaim{Theorem 1.5} Let $R$ be an excellent regular local ring,
containing a field  of characteristic zero. Let $X$ be a nonsingular
$R$-scheme, $f:X\rightarrow\text{spec}(R)$ a projective morphism, $h\in
R$. Then there exists a sequence of monoidal transforms $g:Y\rightarrow X$, such that $h$ has SNCs in $Y$. 
\endproclaim

\demo{Proof} Immediate from Main Theorem II(N) [H].

\proclaim{Theorem 1.6} Suppose that $R$, $S$ are excellent regular local rings containing a field $k$ of characteristic zero
such that $S$ dominates $R$. Let $\nu$ be a valuation of the quotient field $K$ of $S$ that dominates $S$,
$R\rightarrow R_1$ a monoidal transform such that $\nu$ dominates $R_1$. 
$R_1$ is a local ring on $X=\text{Proj}(\bigoplus_{n\ge0}p^n)$ for some prime $p\subset R$.
Let
$$
U=\{Q\in\text{spec}(S):pS_Q\text{ is invertible }\}
$$
an open subset of $\text{spec}(S)$. Then there exists a projective morphism  $f:Y\rightarrow \text{spec}(S)$
which is a product of monoidal transforms
such that if $S_1$ is the local ring of $Y$ dominated by $\nu$, then   $S_1$ dominates $R_1$, and
$(f)^{-1}(U)\rightarrow U$ is an isomorphism.
\endproclaim
\demo{Proof} 
Since $S$ is a UFD, we can write $pS=gI$, where $g\in S$, $I\subset S$ has height $\ge 2$. Then $U=\text{spec}(S)-V(I)$.
By Main Theorem II(N) [H], there exists a sequence of monoidal transforms $\pi:Y\rightarrow\text{spec}(S)$
such that $I\Cal O_Y$ is invertible, and $\pi^{-1}(U)\rightarrow U$ is an isomorphism. 
Let $S_1$ be the local ring of the center of $\nu$ on $Y$.
We have $pS_1 = hS_1$ for some $h\in p$. Hence $R[\frac{p}{h}]\subset S_1$, and since $\nu$ dominates $S_1$,
$R_1$ is the localization of $R[\frac{p}{h}]$ which is dominated by $S_1$.
\enddemo

\proclaim{Theorem 1.7} Suppose that $R$ is an excellent local domain containing a field of characteristic zero,
with quotient field $K$. Let $\nu$ be a valuation of $K$ dominating $R$. Suppose that $f\in K$ is such that $\nu(f)\ge0$.
Then there exists a MTS along $\nu$
$$
R\rightarrow R_1\rightarrow\cdots\rightarrow R_n
$$
such that $f\in R_n$.
\endproclaim

\demo{Proof} Write $f=\frac{a}{b}$ with $a,b\in R$. By  Main Theorem II(N) [H] applied to the ideal $I=(a,b)$ in $R$,
there exists a MTS along $\nu$, $R\rightarrow R_n$ such that $IR_n=\alpha R_n$ is a principal ideal. There exist
constants $c,d,u_1,u_2$ in $R_n$ such that $a=c\alpha,b=d\alpha,\alpha=u_1a+u_2b$. Then
$u_1c+u_2d=1$, so that $cR_n+dR_n=R_n$, and one of $c$ or $d$ is a unit in $R_n$. If $c$ is a unit, then
$0\le\nu(f)=\nu(\frac{c}{d})=\nu(c)-\nu(d)$ implies $\nu(d)=0$, and since $\nu$ dominates $R_n$,  $d$ is a unit and $f\in R_n$.
\enddemo

\proclaim{Theorem 1.8}(Abhyankar) Let $R$, $S$ be two dimensional regular 
local rings such that $R$ and $S$  have the same quotient 
field, and $S$ dominates $R$. Then there exists a unique finite sequence
$$
R_0\rightarrow R_1\rightarrow \cdots\rightarrow R_m
$$
of quadratic transforms such that $R_m=S$.
\endproclaim

\demo{Proof} This is Theorem 3 of [Ab2].

Suppose that $Y$ is an algebraic scheme, $X$, $D$ are subschemes of $Y$. Suppose that $g:Y'\rightarrow Y$,
$f:X'\rightarrow X$ are the monoidal transforms of $Y$ and $X$ with center $D$ and $D\cap X$ respectively. Then there
exists a unique isomorphism of $X'$ to a subscheme $X''$ of $Y'$ such that $g$ induces $f$ (c.f. chapter 0, section 2 [H]).
$X''$ is called the strict transform of $X$ be the monoidal transform $g$.

\proclaim{Theorem 1.9} 
Let $R$ be an excellent regular local ring, containing a field of characteristic zero. Let $W\subset \text{spec}(R)$ be an integral
subscheme, $V\subset \text{spec}(R)$ be the singular locus of $W$. Then there exists a sequence of monoidal transforms
$f:X\rightarrow\text{spec}(R)$ such that the strict transform of $W$ is nonsingular in $X$, and $f$ is an isomorphism over 
$\text{spec}(R)-V$.
\endproclaim
\demo{Proof} This is immediate from Theorem $I_2^{N,n}$ [H].
\enddemo

\proclaim{Theorem 1.10} Suppose that $R\subset S$ are $r$ dimensional local rings with a common quotient field $K$,
and respective maximal ideals $m$ and $n$ such that $S$ dominates $R$, $S/mS$ is a finite $R/m$ module, and
$R$ is normal and analytically irreducible. Then $R=S$.
\endproclaim
\demo{Proof} This is the version of Zariski's Main Theorem proved in Theorem 37.4 [N].
\enddemo

\proclaim{Theorem 1.11}(Theorem 1 [HHS]) Suppose that $R$ is an excellent regular local ring with quotient field $J$,
 $K$ is a finite
extension field of $J$ and $S$ is a regular local ring with quotient field $K$ such that $R\subset S$ and
$\text{dim}(R)=\text{dim}(S)$. Then $S$ is essentially of finite type over $R$.
\endproclaim

\demo{Proof} Let $(y_1,\ldots, y_n)$ be a system of regular parameters in $S$ and suppose that $K$ is generated by
$h_1,\ldots,h_r$ over $J$. Let $h_i=\frac{f_i}{g_i}$ for $1\le i\le r$ where $f_i,g_i\in S$. Let
$T$ be the normalization of   $R[y_1,\ldots,y_n,f_1,\ldots,f_r,g_1,\ldots,g_r]$, $q=m(S)\cap T$,
$U=T_q$. By Theorem 1.10 $T_q=U$. 
\enddemo

\proclaim{Theorem 1.12} Suppose that $R$ is an excellent regular local ring, with maximal ideal $m$,
$S$ is a regular local ring with maximal ideal $n$, such that $R\subset S$, $\text{dim}(R)=\text{dim}(S)$ 
and the quotient field of $S$ is a finite extension of the quotient field of $R$. Then there is an inclusion
$$
\hat R\subset \hat S
$$
where $\hat R$ is the $m$-adic completion of $R$, $\hat S$ is the $n$-adic completion of $S$.
\endproclaim

\demo{Proof} By Theorem 1.11 $S$ is essentially of finite type over $R$. Since $S$ is universally catenary, the
dimension formula (Theorem 15.6 [M]) holds.
$$
\text{dim} R+\text{trdeg}_RS=\text{dim}S+\text{restrdeg}_RS
$$
Since $R$ is analytically irreducible, $R$ is a subspace of $S$ by Theorem 10.13 [Ab4] ("A version of
Zariski's Subspace Theorem").
\enddemo

We will use the notation $m(R)$ to denote the maximal ideal of a local ring $R$, $k(R)$ to denote the residue field
$R/m(R)$. $\hat R$ or $R\,\,\hat{}\,$ will denote the $m(R)$-adic completion of $R$.

 Let $k$ be a field, $0\ne f(z_1,\ldots,z_n)\in k[[z_1,\ldots,z_n]]$. Let $m=(z_1,\ldots,z_n)$.
Define $\text{mult}(f)=r$ if $f\in m^r$, $f\not\in m^{r+1}$.

\heading
2 Uniformizing Transforms
\endheading

\proclaim{Definition 2.1}
Suppose that  
$R$ is a regular local ring. 
A monoidal transform sequence (MTS) is a sequence 
of ring homomorphisms 
$$
R=R_0\rightarrow R_1\rightarrow R_2\rightarrow \cdots\rightarrow R_n
$$
such that each map $R_i\rightarrow R_{i+1}$ is a finite product of monoidal transforms.
\endproclaim

\proclaim{Definition 2.2}
Suppose that 
$R$ is an excellent regular local ring  containing a field $k$ of characteristic zero,
 with quotient field $K$. 
A uniformizing transform sequence (UTS) is a sequence 
of ring homomorphisms
$$
\matrix
R&\rightarrow &\overline T_0''&\rightarrow &\overline T_0&&&&& &&&&& \\
 &            & \downarrow &               &             &\searrow &&&&&&&&&\\
 &            &\overline T_1'&\rightarrow  &\overline T_1''&\rightarrow &\overline T_1&&&&&&&&\\
 &              &           &               & \downarrow &               &             &\searrow &&&&&&&\\
 &             &               &            &\overline T_2'&\rightarrow  &\overline T_2''&\rightarrow &\overline T_2&&&&&&\\
&              &               &            &              &             &\downarrow     &             &            &\searrow&&&&&\\
&&&&&&\vdots&&&\vdots&&&&\\
&&&&&&&&&   &\downarrow   &             &               &    \searrow&\\           
&&&&&&&&&   &\overline T_n'&\rightarrow &\overline T_n''&\rightarrow &\overline T_n 
\endmatrix\tag 2.1
$$
such that 
 $\overline T_0=\hat R$, the completion of $R$ with respect to its maximal ideal, and for all $i$,
$\overline T_i$ is the completion 
with respect to its maximal ideal of a finite product of monoidal transforms $\overline T_i'$ of $\overline T_{i-1}''$.
For all $i$, 
 $\overline T_i''$ is a  a regular local ring  essentially of finite type over $\overline T_i'$ with quotient field $K_i$, 
such that $\overline T_i'\subset\overline T_i''\subset \overline T_i$ and 
$K_0$ is a finite extension of $K$,  $K_{i+1}$ is a finite extension of $K_{i}$ for all $i\ge0$. 
\endproclaim

To simplify notation, we will often denote the UTS (2.1) by $(R,\overline T_n'',\overline T_n)$ or by
$$
R\rightarrow \overline T_0\rightarrow \overline T_1\rightarrow\cdots\rightarrow \overline T_n.
$$
We will denote the UTS consisting of the maps
$$
\matrix
\overline T_{n-1}'&\rightarrow &\overline T_{n-1}''&\rightarrow&\overline T_{n-1}&&\\
                 &            & \downarrow        &            &                &\searrow&\\
                 &            & \overline T_n'    &\rightarrow  &\overline T_{n}''&\rightarrow&\overline T_n
\endmatrix
$$
by $\overline T_{n-1}\rightarrow \overline T_n$.

A UTS (2.1) is called a rational uniformizing sequence (RUTS) if there exists an associated MTS 
$$
R = R_0\rightarrow R_1\rightarrow\cdots\rightarrow R_n,
$$
maps $R_i\rightarrow \overline T_i'$ such that  $\hat R_i \cong \overline T_i$ for $1\le i\le n$, and all
squares in the resulting diagram
$$
\matrix
\overline T_0&\rightarrow&\overline T_1&\rightarrow&\cdots&\rightarrow&\overline T_n\\
\uparrow&    &           \uparrow     &           &      &           &\uparrow\\
R_0&\rightarrow&R_1&\rightarrow&\cdots&\rightarrow&R_n
\endmatrix\tag 2.2
$$
commute.

Suppose that $v$ is a  valuation of $K$ which dominates $R$ and
$$
R\rightarrow \overline T_0\rightarrow \overline T_1\rightarrow \cdots\rightarrow \overline T_n\tag 2.3
$$
is a UTS. 

Suppose that $v_0$ is an extension of $v$ to the quotient field of $\overline T_0$ such that $\nu_0$
dominates $\overline T_0$.  If $v_0$ dominates $\overline T_1'$ we can extend $v_0$ to a valuation $v_1$ of 
the quotient field of $\overline T_1$  which
dominates $\overline T_1$. 

Then if $v_1$ dominates $\overline T_2'$,
in the same manner
 we can extend $v_1$ to a valuation $v_2$ of the quotient field of $\overline T_2$ which
dominates $\overline T_2$. If we can inductively construct a sequence $v_1,\ldots, v_n$ of extensions of $v$ to the quotient
fields of $\overline T_i$ in this manner, (2.3) is called a UTS along $\nu$.  If there is no danger of confusion, we will denote the
extensions $v_i$ by $v$.

Suppose that (2.3) is a UTS along a rank 1 valuation $\nu$ of $K$. Let $\Gamma_{\nu}$ be the value goup of $\nu$. Suppose that $i$ is
such that $0\le i\le n$. Let $\Cal O_{\nu_i}$ be the valuation ring of $\nu_i$ and $\Gamma_{\nu_i}$ be the value group of $\nu_i$.
$\Gamma_{\nu}$ is a subgroup of $\Gamma_{\nu_i}$. Set 
$$
\overline \Gamma = \{\beta \in \Gamma_{\nu_i}\,|\,-\alpha\le\beta\le \alpha\text{ for some }\alpha \in \Gamma_{\nu}\}.
$$
$\overline \Gamma$ is an isolated subgroup of $\Gamma_{\nu_i}$, since $\Gamma_{\nu}$ is a subgroup,
 so there is a prime $a$ in $\Cal O_{\nu_i}$ (which could be 0) 
 such that  $\overline\Gamma$ is the isolated subgroup $\Gamma_a$ of $a$, 
(by Proposition 2.29 [Ab3] or Theorem 15, chapter VI, section 10 [ZS]). Set
$$
\overline p_i = a\cap \overline T_i = \{f\in \overline T_i\,|\,\nu(f)>\alpha\text{ for all }\alpha\in\Gamma_{\nu}\}.
$$
We will say that $\nu(f)=\infty$ if $f\in \overline p_i$.

For the rest of this chapter we will assume that $R\subset S$ are excellent regular local rings such that 
$\text{dim}(R)=\text{dim}(S)$,
containing a field $k$ of characteristic zero, such that the quotient field $K$ of $S$ is a
finite extension of the quotient field $J$ of $R$. We will fix a valuation $\nu$ of $K$ with
valuation ring $V$  such that $\nu$ dominates $S$. 

Note that the restriction of $\nu$ to $J$ has the same rank and rational rank that $\nu$ does
(Lemmas 1 and 2 of section 11, chapter VI [ZS]). 
Observe that
 $S$ is essentially of finite type over $R$ (Theorem 1.11) and
$\hat R\rightarrow \hat S$ is an inclusion (Theorem 1.12).

Suppose that $(R, \overline T_n'',\overline T_n)$ and $(S,\overline U_n'',\overline U_n)$ are UTSs. We will
say that $(R,\overline T_n'',\overline T_n)$ and $(S,\overline U_n'',\overline U_n)$ are compatible UTSs
(or a CUTS) if there are commutative diagrams of inclusions
$$
\matrix
\overline U_{i}'&\rightarrow &\overline U_i''&\rightarrow& \overline U_i\\
\uparrow      &            &\uparrow       &           &\uparrow\\
\overline T_{i}'&\rightarrow &\overline T_i''&\rightarrow& \overline T_i
\endmatrix\tag 2.4
$$
for $0\le i\le n$. In particular, the quotient field of $\overline U_i''$
is finite over the quotient field of $\overline T_i''$ for all $i$, and $\overline U_i''$ is essentially
of finite type over $\overline T_i''$ for all $i$.

We will say that UTSs along $\nu$ $(R,\overline T_n'',\overline T_n)$ and $(S,\overline U_n'',\overline U_n)$ are CUTS
 along $\nu$ if the  extensions of $\nu$ are compatible in (2.4).

If $(R,\overline T_n'',\overline T_n)$ and $(S,\overline U_n'',\overline U_n)$ are RUTSs and CUTSs, then we will say that
$(R,\overline T_n'',\overline T_n)$ and $(S,\overline U_n'',\overline U_n)$ are compatible RUTSs (or a CRUTS).

\proclaim{Lemma 2.3}
Suppose that the CRUTS $(R,\overline T_n'',\overline T_n)$ and $(S,\overline U_n'',\overline U_n)$
have respective associated MTSs
$$
R=R_0\rightarrow R_1\rightarrow\cdots\rightarrow R_n
$$
and
$$
S=S_0\rightarrow S_1\rightarrow\cdots\rightarrow S_n.
$$
Then there is a commutative diagram
$$
\matrix
R&\rightarrow& R_1&\rightarrow&\cdots&\rightarrow &   R_n\\
\uparrow&    &\uparrow&       &      &            &\uparrow\\
S&\rightarrow& S_1&\rightarrow&\cdots&\rightarrow& S_n
\endmatrix.\tag 2.5
$$
\endproclaim
\demo{Proof}
This follows from (2.4), (2.2) and Lemma 1.1, since then
$$
R_i\subset \hat R_i\cap J=\overline T_i\cap J\subset \overline U_i\cap K=\hat S_i\cap K=S_i
$$
for all $i$.
\enddemo

\heading 3. Rank 1 \endheading 
\subheading{Perron Transforms}

Throughout this chapter we will assume that
$R\subset S$ are excellent regular local rings such that 
$\text{dim}(R)=\text{dim}(S)$,
containing a field $k$ of characteristic zero, such that the quotient field $K$ of $S$ is a
finite extension of the quotient field $J$ of $R$. We will fix a valuation $\nu$ of $K$ with
valuation ring $V$  such that $\nu$ dominates $S$. 
We will further assume that
\item{1)} $\nu$ has rank 1 and arbitrary  rational rank  $s$ ($\le\text{dim}(S)$).
\item{2)}   $\text{dim}_R(\nu)=0$ and $\Cal O_{\nu}/m_{\nu}$ is algebraic over $k$.

Let $n=\text{dim}(R)=\text{dim}(S)$.
We will define
 2 types of UTSs.
Suppose that 
$(R,\overline T'',\overline T)$ is a UTS along $\nu$ and $\overline T''$ has regular parameters
$(\tilde x_1',\ldots,\tilde x_n')$ such that 
$$
\nu(\tilde x_1')=\tau_1,\ldots,\nu(\tilde x_s')=\tau_s
$$
 are rationally independent.
Let $\nu_0$ be an extension of $\nu$ to the quotient field of $\overline T''$ which dominates $\overline T''$.

We first define a UTS $\overline T\rightarrow \overline T(1)$ of type I along $\nu$.
 The MTS $\overline T''\rightarrow \overline T'(1)$
is defined as follows.
  $\overline T'(1) = T_h$ where $h$ is a positive integer  and $T_h$ is constructed as follows.

Set $\tau_i(0)=\tau_i$ for $1\le i\le s$. For each positive integer $h$ define $s$ positive, rationally independent
real numbers $\tau_1(h),\ldots,\tau_s(h)$ by the "Algorithm of Perron" (B.I of [Z2])
$$
\align
\tau_1(h-1)&=\tau_s(h)\\
 \tau_2(h-1)&=\tau_1(h)+a_2(h-1)\tau_s(h)\\
&\vdots\\
\tau_s(h-1)&=\tau_{s-1}(h-1)+a_s(h-1)\tau_s(h)
\endalign
$$ 
where 
$$
a_j(h-1)=\left[ \frac{\tau_j(h)}{\tau_1(h)}\right],
$$
the "greatest integer" in $ \frac{\tau_j(h)}{\tau_1(h)}$.
There are then nonnegative integers $A_i(h)$ such that 
$$
\tau_i = A_i(h)\tau_1(h)+A_i(h+1)\tau_2(h)+\cdots+A_i(h+s-1)\tau_s(h)
$$
for $1\le i\le s$.
$$
\text{det}\left(\matrix A_1(h)&\ldots&A_1(h+s-1)\\
\vdots&&\vdots\\
A_s(h)&\ldots&A_s(h+s-1)
\endmatrix\right) = (-1)^{h(s-1)}
$$
(See formula (4'), page 385 [Z2].) 
These numbers have the important property that 
$$
\text{Lim}_{h\rightarrow \infty}\frac{A_i(h)}{A_1(h)}=\frac{\tau_i}{\tau_1}\tag 3.1
$$
(See formula  (5), page 385 [Z2].)
Set $\tilde x_i(0)=\tilde x_i'$ for $1\le i\le n$.
Define $T_h$ by the sequence of MTSs along $\nu$
$$
\overline T''=\tilde T(0)\rightarrow \tilde T(1)\rightarrow\cdots\rightarrow \tilde T(h) = T_h =\overline T'(1)
$$
Where $\tilde T(i+1)=\tilde T(i)[\tilde x_1(i+1),\ldots,\tilde x_s(i+1)]_{(\tilde x_1(i+1),\ldots,\tilde x_n(i+1))}$
for $0\le i\le h-1$.  
$$
\align
\tilde x_1(i)&=\tilde x_s(i+1)\\
\tilde x_2(i)&=\tilde x_1(i+1)\tilde x_s(i+1)^{a_2(i)}\\
&\vdots\\
\tilde x_s(i)&=\tilde x_{s-1}(i+1)\tilde x_s(i+1)^{a_s(i)}
\endalign
$$
$\nu(\tilde x_j(i))=\tau_j(i)$ for $1\le j\le s$.
If we set $\overline x_i(1) = \tilde x_i(h)$, we then have regular parameters 
$(\overline x_1(1),\ldots,\overline x_n(1))$ in $\overline T'(1)$
satisfying
$$
\align
\tilde x_1' &= \overline x_1(1)^{A_1(h)}\cdots \overline x_s(1)^{A_1(h+s-1)}\tag 3.2\\
&\vdots\\
\tilde x_s' &= \overline x_1(1)^{A_s(h)}\cdots \overline x_s(1)^{A_s(h+s-1)}\\
\tilde x_{s+1}'&=\overline x_{s+1}(1)\\
&\vdots\\
\tilde x_n'&=\overline x_n(1)
\endalign
$$
Then $\overline T'(1)=\overline T''[\overline x_1(1),\ldots,\overline x_s(1)]
_{(\overline x_1(1),\ldots,\overline x_n(1))}$.
 Let $\overline T(1)$ be the completion of $\overline T'(1)$ at its maximal ideal.
 Set $\overline T''(1)=\overline T'(1)$.
Then for any extension $\nu_1$ of $\nu_0$ to the quotient field of $\overline T(1)$ which dominates $\overline T(1)$,
 $\overline T\rightarrow \overline T(1)$ is a UTS and
$(R,\overline T''(1),\overline T(1))$ is a UTS along $\nu$.  Note that 
$\nu(\overline x_1(1)),\ldots,\nu(\overline x_s(1))$ are rationally independent.

Now we define a UTS $\overline T\rightarrow \overline T(1)$ of type $\text{II}_r$
along $\nu$ (with the
restriction that $s+1\le r\le n$). The MTS
$\overline T''\rightarrow \overline T'(1)$ is constructed as follows.
Set $\nu(\tilde x_r') = \tau_r$. $\tau_r$ must be rationally dependent on $\tau_1,\ldots,\tau_s$
since $\nu$ has rational rank $s$. There are thus integers $\lambda,\lambda_1,\ldots,\lambda_s$ such that
$\lambda>0$, $(\lambda,\lambda_1,\ldots,\lambda_s)=1$ and 
$$
\lambda\tau_r=\lambda_1\tau_1+\cdots+\lambda_s\tau_s.
$$
First perform a MTS $\overline T''\rightarrow \tilde T(1)$ which is UTS along $\nu$ where $\tilde T(1)$ has regular parameters
$(\tilde x_1(1),\ldots,\tilde x_n(1))$ defined by
$\tilde x_i'=\tilde x_1(1)^{A_1(h)}\cdots\tilde x_s(1)^{A_s(h+s-1)}$ for $1\le i\le s$. Then
$\nu(\tilde x_i(1))= \tau_i(h)$ for $1\le i\le s$, $\nu(\tilde x_r(1))=\tau_r$. Set
$$
\lambda_i(h)=\lambda_1A_1(h+i-1)+\lambda_2A_2(h+i-1)+\cdots+\lambda_sA_s(h+i-1)
$$
for $1\le i\le s$. Then
$$
\lambda\tau_r=\lambda_1(h)\tau_1(h)+\cdots+\lambda_s(h)\tau_s(h).
$$
Take $h$ sufficiently large that all $\lambda_i(h)>0$. This is possible by (3.1), since
$\lambda_1\tau_1+\cdots+\lambda_s\tau_s>0$. We still have $(\lambda,\lambda_1(h),\ldots,\lambda_s(h))=1$ since
$\text{det}(A_i(h+j-1))=\pm1$.
After redindexing the $\tilde x_i(1)$, we may suppose that $\lambda_1(h)$ is not divisable by $\lambda$. Let
$\lambda_1(h)=\lambda\mu+\lambda'$, with $0<\lambda'<\lambda$. Now perform a MTS $\tilde T(1)\rightarrow\tilde T(2)$
along $\nu$ where $\tilde T(2)$ has regular parameters $(\tilde x_1(2),\ldots,\tilde x_n(2))$ defined by
$$
\align
\tilde x_1(1)&= \tilde x_r(2)\\
\tilde x_2(1)&= \tilde x_2(2)\\
&\vdots\\
\tilde x_s(1)&=\tilde x_s(2)\\
\tilde x_r(1)&=\tilde x_1(2)\tilde x_r(2)^{\mu}.
\endalign
$$
Set $\tau_i'=\nu(\tilde x_i(2))$ for all $i$. $\tau_1',\ldots,\tau_s',\tau_r'$ are positive and
$$
\lambda'\tau_r'=\lambda_1'\tau_1'+\cdots+\lambda_s'\tau_s'
$$
where $\lambda_1'=\lambda$, $\lambda_i'=-\lambda_i(h)$ for $2\le i\le s$. We have thus acheived a reduction in $\lambda$.
By repeating this procedure, we get a MTS $\overline T''\rightarrow \tilde T(\alpha)$ along $\nu$ where 
$\tilde T(\alpha)$ has regular parameters
$(\tilde x_1(\alpha),\ldots,\tilde x_n(\alpha))$
such that if $\overline \tau_i = \nu(\tilde x_i(\alpha))$, $\overline \tau_1,\ldots,\overline \tau_s$ are
rationally independent and
$$
\overline \tau_r=\overline \lambda_1\overline\tau_1+\cdots +\overline \lambda_s\overline\tau_s
$$
for some integers $\overline \lambda_i$. Now perform a MTS $\tilde  T(\alpha)\rightarrow \tilde T(\alpha+1)$
which is a UTS
of type I  along $\nu$ where 
$T(\alpha+1)$ has regular parameters
$(\tilde x_1(\alpha+1),\ldots,\tilde x_n(\alpha+1))$
such that if $\tau_i^* = \nu(\tilde x_i(\alpha+1))$, 
$$
\tau_r^*=\tilde \lambda_1\tau_1^*+\cdots +\tilde \lambda_s\tau_s^*
$$
for some positive integers $\tilde \lambda_i$. Finally perform a MTS $\tilde  T(\alpha+1)\rightarrow \tilde T(\alpha+2)$
where $\tilde T(\alpha+2)=\tilde T(\alpha+1)[N_r]_q$.
$$
N_r = \frac{\tilde x_r(\alpha+1)}{\tilde x_1(\alpha+1)^{\tilde\lambda_1}\cdots \tilde x_s(\alpha+1)^{\tilde\lambda_s}}
$$
and $q$ is the center of $\nu$ on $\tilde T(\alpha+1)[N_r]$. Set $\overline T'(1) = \tilde T(\alpha+2)$. Since
$\nu(N_r)=0$, $N_r$ has residue $c\ne0$ in $k(\overline T'(1))$.
Set $N_i=\tilde x_i(\alpha)$ for $1\le i\le s$. Then there exists a matrix $(a_{ij})$ such that
$$
\align
\tilde x_1'&=N_1^{a_{11}}\cdots N_r^{a_{1,s+1}}\\
&\cdots\\
\tilde x_s'&=N_1^{a_{s1}}\cdots N_r^{a_{s,s+1}}\\
\tilde x_r'&=N_1^{a_{s+1,1}}\cdots N_r^{a_{s+1,s+1}}
\endalign
$$
and $\text{det}(a_{ij})=\pm1$.
$\overline T'(1)$ is a localization of $\overline T''[N_1,\ldots,N_s,N_r]$.

Let $\overline T(1)$ be the completion of $\overline T'(1)$ at its maximal ideal. $\nu_0$ extends to a valuation 
 of the quotient field of $\overline T(1)$ which dominates $\overline T(1)$. Let $\nu_1$ be such an extension. 
$\overline T(1)$ has a regular system of parameters  $(x_1^*(1),\ldots,x_n^*(1))$
defined by

$$
\align
\tilde x_1' &= x_1^*(1)^{a_{11}}\cdots x_s^*(1)^{a_{1s}}
(x_r^*(1)+c)^{a_{1,s+1}}\tag 3.3\\
\vdots&\\
\tilde x_s' &= x_1^*(1)^{a_{s1}}\cdots x_s^*(1)^{a_{ss}}
(x_r^*(1)+c)^{a_{s,s+1}}\\
\tilde x_r' &= x_1^*(1)^{a_{s+1,1}}\cdots x_s^*(1)^{a_{s+1,s}}
(x_r^*(1)+c)^{a_{s+1,s+1}}.\\
\endalign
$$
$\text{Det}(a_{ij})= \pm 1$ and $\nu_1(x_1^*(1),\ldots,\nu_1(x_s^*(1))$ are rationally independent.
Set 
$$
\overline c = \text{Det}\left(\matrix a_{11}& \cdots & a_{1s}\\ \vdots&&\vdots\\ a_{s1}&\cdots& a_{ss}\endmatrix\right)
\text{Det}\left(\matrix a_{11}& \cdots & a_{1,s+1}\\ \vdots&&\vdots\\ a_{s+1,1}&\cdots& a_{s+1,s+1}\endmatrix\right)
\ne0
$$
since $\nu(\tilde x_1'),\ldots,\nu(\tilde x_s')$ are rationally independent.
Define rational numbers $\gamma_1,\ldots,\gamma_s$ by
$$
\left(\matrix \gamma_1\\ \vdots\\ \gamma_s\endmatrix\right) =
\left(\matrix a_{11}& \cdots & a_{1s}\\ \vdots&&\vdots\\ a_{s1}&\cdots& a_{ss}\endmatrix\right)^{-1}
\left(\matrix -a_{1, s+1}\\ \vdots \\ -a_{s,s+1}\endmatrix\right).
$$
$\gamma_i=\frac{m_i}{\overline c}$ for some $m_i\in\bold N$. By Cramer's rule,
$$
\left(\matrix a_{1,1}& \cdots & a_{1,s+1}\\ \vdots&&\vdots\\  a_{s,1}&\cdots& a_{s,s+1}\\
 a_{s+1,1}&\cdots& a_{s+1,s+1}\endmatrix\right)
\left(\matrix \gamma_1\\ \vdots\\ \gamma_s\\ 1 \endmatrix\right) =
\left(\matrix 0\\ \vdots \\ 0 \\ \frac{1}{\overline c}\endmatrix\right).
$$
$\gamma_1 a_{s+1,1}+ \cdots + \gamma_s a_{s_1,s} + a_{s+1,s+1} =  1/\overline c$.

Let $(b_{ij})=(a_{ij})^{-1}$. Then
$$
\align
N_1&=x_1^*(1)=(\tilde x_1')^{b_{1,1}}\cdots(\tilde x_s')^{b_{1,s}}(\tilde x_r')^{b_{1,s+1}}\\
&\vdots\\
N_s&=x_s^*(1)=(\tilde x_1')^{b_{s,1}}\cdots(\tilde x_s')^{b_{s,s}}(\tilde x_r')^{b_{s,s+1}}\\
N_r&=x_r^*(1)+c=(\tilde x_1')^{b_{s+1,1}}\cdots(\tilde x_s')^{b_{s+1,s}}(\tilde x_r')^{b_{s+1,s+1}}
\endalign
$$
$$
(x_1^*(1),\ldots,x_{r-1}^*(1),N_r',x_{r+1}^*(1),\ldots,x_n^*(1))
$$
are regular parameters in $\overline T'(1)$ where $N_r'=\prod(N_r-\sigma(c))$ where the product is over all conjugates
$\sigma(c)$ of $c$ over $k$ in an algebraic closure $\overline k$ of $k$. 
$$
k(\overline T'(1))\cong k(\overline T'')(c).
$$
We have $\left(\frac{x_r^*(1)}{c}+1\right)^{\frac{1}{\overline c}}\in\overline T(1)$ where 
$\left(\frac{x_r^*(1)}{c}+1\right)^{\frac{1}{\overline c}}$ is uniquely determined by the condition that it has
residue 1 in $k(\overline T(1))$.
Set 
$$
\overline T''(1)=\overline T''\left[c,\left(\frac{x_r^*(1)}{c}+1\right)^{\frac{1}{\overline c}}\right]
_{(x_1^*(1),\ldots,x_n^*(1))}.
$$
 $\overline T''(1)$ has regular parameters
$(\overline x_1(1),\ldots,\overline x_n(1))$ defined by
$$
\overline x_i(1)=\cases
x_i^*(1)\left(\frac{x_r^*(1)}{c}+1\right)^{-\gamma_i}& 1\le i\le s\\
\left(\frac{x_r^*(1)}{c}+1\right)^{\frac{1}{\overline c}}-1&i=r\\
x_i^*(1)&s<i,i\ne r
\endcases\tag 3.4
$$
We have
$$
\align
\tilde x_1'&=\overline x_1(1)^{a_{1,1}}\cdots\overline x_s(1)^{a_{1,s}}c^{a_{1,s+1}}\\
&\vdots\\
\tilde x_s'&=\overline x_1(1)^{a_{s,1}}\cdots\overline x_s(1)^{a_{s,s}}c^{a_{s,s+1}}\\
\tilde x_r'&=\overline x_1(1)^{a_{s+1,1}}\cdots\overline x_s(1)^{a_{s+1,s}}(\overline x_r(1)+1)c^{a_{s+1,s+1}}\\
\endalign
$$
Note that
$$
x_i^*(1)=\cases (\overline x_r(1)+1)^{\overline c \gamma_i}\overline x_i(1)&1\le i\le s\\
c[(\overline x_r(1)+1)^{\overline c}-1]&i=r\\
\overline x_i(1)&s<i,i\ne r
\endcases\tag 3.5
$$
Thus $\overline T\rightarrow \overline T(1)$
is a UTS and by our extension of $\nu$ to the quotient field of $\overline T(1)$,
 $(R,\overline T''(1),\overline T(1))$ is a UTS along $\nu$. 
We will call $\overline T\rightarrow \overline T(1)$ a UTS of type $\text{II}_r$.

\proclaim{Remark 3.1} In our constructions of UTSs  of types I and $\text{II}_r$, 
$\overline T''\rightarrow \overline T'(1)$ is a product of monoidal transforms 
$$
\overline T''=T_0\rightarrow T_1\rightarrow\cdots\rightarrow T_{t-1}\rightarrow T_t=\overline T'(1)
$$
where each $T_i\rightarrow T_{i+1}$ is a monoidal transform centered at a height 2 prime $a_i$ and
$a_i\overline T''(1) = (\overline x_1(1)^{d_1^i}\cdots \overline x_s(1)^{d_s^i})$ for some nonnegative integers $d_j^i$ for all $i$.
\endproclaim

\proclaim{Lemma 3.2}Suppose that $(R,\overline T'',\overline T)$ is a UTS along $\nu$,
$(x_1,\ldots,x_n)$ are regular parameters in $\overline T''$, and $\nu(x_1),\ldots,\nu(x_s)$ are
rationally independent.
\item{1)} Suppose that $M_1 = x_1^{d_1^1}\cdots x_s^{d_s^1},M_2 = x_1^{d_1^2}\cdots x_s^{d_s^2}\in \overline T''$
and $\nu(M_1)<\nu(M_2)$.
Then there exists a UTS of type I along $\nu$, $\overline T\rightarrow \overline T(1)$, such that $M_1\mid M_2$ in 
$\overline T'(1)$. 
\item{2)} Suppose that $M = x_1^{d_1}\cdots x_s^{d_s}$ is  such that the $d_i$ are integers and
$0<\nu(M)$. Then there exists a UTS of type I along $\nu$ $\overline T\rightarrow \overline T(1)$ such that $M\in \overline T(1)'$.
\endproclaim

\demo{Proof}
The proof of 1) is from Theorem 2 [Z2]. Consider the UTS with equations
 (3.2). In $\overline T'(1)$, 
$$
M_i=\overline x_1(1)^{d_1^iA_1(h)+\cdots+d_s^iA_s(h)} \cdots
\overline x_s(1)^{d_1^iA_1(h+s-1)+\cdots+d_s^iA_s(h+s-1)}
$$
for $i=1,2$. For $h>>0$
$$
d_1^2A_1(h+j-1)+\cdots+d_s^2A_s(h+j-1)>d_1^1A_1(h+j-1)+\cdots+d_s^1A_s(h+j-1)
$$
for $1\le j\le s$ by (3.1).

 To prove 2) just write $M = \frac{M_2}{M_1}$ where $M_i$ are monomials in $x_1,\ldots, x_s$.
Since $\nu(M_2)>\nu(M_1)$, 2) follows from 1).
\enddemo 

\proclaim{Lemma 3.3} Suppose that $(R,\overline T'',\overline T)$ and $(S,\overline U'',\overline U)$ is a CUTS
along $\nu$, $\overline T''$
 has regular parameters $(\overline x_1,\ldots, \overline x_n)$ and
$\overline U''$ has regular parameters $(\overline y_1,\ldots, \overline y_n)$, related by
$$
\align
\overline x_1 &= \overline y_1^{c_{11}}\cdots \overline y_s^{c_{1s}}\alpha_1\\
\vdots&\\
\overline x_s &= \overline y_1^{c_{s1}}\cdots \overline y_s^{c_{ss}}\alpha_s\\
\endalign
$$
such that $\alpha_1,\ldots,\alpha_s\in k(\overline U)$, 
$\nu(\overline x_1),\ldots,\nu(\overline x_s)$ are rationally independent and $\text{det}(c_{ij})\ne 0$.
 Suppose that $\overline T\rightarrow \overline T(1)$ is a UTS of type $I$ along $\nu$,
such that $\overline T'(1)=\overline T''(1)$ has regular parameters $(\overline x_1(1),\ldots,\overline x_n(1))$ with
$$
\align
\overline x_1 &= \overline x_1(1)^{a_{11}}\cdots \overline x_s(1)^{a_{1s}}\\
\vdots&\\
\overline x_s &= \overline x_1(1)^{a_{s1}}\cdots \overline x_s(1)^{a_{ss}}.
\endalign
$$
Then there exists a UTS of type I along $\nu$ $\overline U\rightarrow \overline U(1)$ 
such that $(R,\overline T''(1),\overline T(1))$ and $(S,\overline U''(1),\overline U(1))$ is a CUTS along $\nu$ and
$\overline U'(1)=\overline U''(1)$ has regular parameters $(\overline y_1(1),\ldots,\overline y_n(1))$ with 
$$
\align
\overline y_1 &= \overline y_1(1)^{b_{11}}\cdots \overline y_s(1)^{b_{1s}}\tag 3.6\\
\vdots&\\
\overline y_s &= \overline y_1(1)^{b_{s1}}\cdots \overline y_s(1)^{b_{ss}},
\endalign
$$
 and 
$$
\align
\overline x_1(1) &= \overline y_1(1)^{c_{11}(1)}\cdots \overline y_s(1)^{c_{1s}(1)}\alpha_1(1)\tag 3.7\\
\vdots&\\
\overline x_s(1) &= \overline y_1(1)^{c_{s1}(1)}\cdots \overline y_s(1)^{c_{ss}(1)}\alpha_s(1)
\endalign
$$
where $\alpha_1(1),\ldots,\alpha_s(1)\in k(\overline U(1))$,
$\nu(\overline x_1(1)),\ldots,\nu(\overline x_s(1))$ are rationally independent and $\text{det}(c_{ij}(1))\ne 0$.
\endproclaim

\demo{Proof} 
Let $(e_{ij})=(a_{ij})^{-1}$, $(d_{ij}) = (e_{ij})(c_{jk})$, an integral matrix. 
Let $\alpha_i(1) = \alpha_1^{e_{i1}}\cdots\alpha_s^{e_{is}}$ for $1\le i\le s$. Then
$$
\align
\overline x_1(1) &= \overline y_1^{d_{11}}\cdots \overline y_s^{d_{1s}}\alpha_1(1)\\
\vdots&\\
\overline x_s(1) &= \overline y_1^{d_{s1}}\cdots \overline y_s^{d_{ss}}\alpha_s(1)\\
\endalign
$$
By 2) of Lemma 3.2, we can construct a UTS (3.6) of type I $\overline U\rightarrow \overline U(1)$
such that  we have an inclusion $\overline T''(1)\subset \overline U(1)''$ and (3.7) holds. Then an
extension of $\nu$ from the quotient field of $\overline U$ which dominates $\overline U$ to a valuation of the
quotient field of $\overline U(1)$ which dominates $\overline U(1)$ restricts to an extension of $\nu$
 to the quotient field of $\overline T(1)$ which dominates $\overline T(1)$ so that
$(R,\overline T''(1),\overline T(1))$ and $(S,\overline U''(1),\overline U(1))$ is a CUTS along $\nu$.
\enddemo

\proclaim{Lemma 3.4} Suppose that $(R,\overline T'',\overline T)$, $(S,\overline U'',\overline U)$ is a CUTS along $\nu$,
$\overline T''$ has regular parameters $(\overline x_1,\ldots,\overline x_n)$ and
$\overline U''$ has regular parameters $(\overline y_1,\ldots,\overline y_n)$ such that
$$
\align
\overline x_1 &= \overline  y_1^{c_{11}}\cdots \overline y_s^{c_{1s}}\alpha_1\tag 3.8\\
\vdots&\\
\overline x_s &= \overline y_1^{c_{s1}}\cdots \overline y_s^{c_{ss}}\alpha_s\\
\overline x_{s+1} &= \overline y_{s+1}\\
\vdots&\\
\overline x_{l} &= \overline y_{l}
\endalign
$$
with $\alpha_1,\ldots,\alpha_s\in k(\overline U)$,
$\nu(\overline x_1),\ldots,\nu(\overline x_s)$  rationally independent, $\text{det}(c_{ij})\ne 0$.

Suppose that $\overline T\rightarrow \overline T(1)$ is a UTS of type $\text{II}_r$ along $\nu$, with $s+1\le r\le l$
such that $\overline T(1)''$ has regular parameters $(\overline x_1(1),\ldots,\overline x_n(1))$ with
$$
\align
\overline x_1 &= \overline x_1(1)^{a_{11}}\cdots \overline x_s(1)^{a_{1s}}c^{a_{1,s+1}}\tag 3.9\\
\vdots&\\
\overline x_s &= \overline x_1(1)^{a_{s1}}\cdots \overline x_s(1)^{a_{ss}}c^{a_{s,s+1}}\\
\overline x_r &= \overline x_1(1)^{a_{s+1,1}}\cdots \overline x_s(1)^{a_{s+1,s}}(\overline x_r(1)+1)c^{a_{s+1,s+1}}.
\endalign
$$
Then there exists a UTS of type $\text{II}_r$, (followed by a UTS of type I) 
$\overline U\rightarrow \overline U(1)$ along $\nu$
such that  $ \overline U(1)''$ has regular parameters $(\overline y_1(1),\ldots,\overline y_n(1))$ satisfying
$$
\align
\overline y_1 &= \overline y_1(1)^{b_{11}}\cdots \overline y_s(1)^{b_{1s}}d^{b_{1,s+1}}\tag 3.10\\
\vdots&\\
\overline y_s &= \overline y_1(1)^{b_{s1}}\cdots \overline y_s(1)^{b_{ss}}d^{b_{s,s+1}}\\
\overline y_r &= \overline y_1(1)^{b_{s+1,1}}\cdots \overline y_s(1)^{b_{s+1,s}}(\overline y_r(1)+1)d^{b_{s+1,s+1}},
\endalign
$$
$\overline T''(1)\subset \overline U''(1)$, and 
$$
\align
\overline x_1(1) &= \overline y_1(1)^{c_{11}(1)}\cdots \overline y_s(1)^{c_{1s}(1)}\alpha_1(1)\tag 3.11\\
\vdots&\\
\overline x_s(1) &= \overline y_1(1)^{c_{s1}(1)}\cdots \overline y_s(1)^{c_{ss}(1)}\alpha_s(1)\\
\overline x_{s+1}(1) &=\overline y_{s+1}(1)\\
\vdots&\\
\overline x_l(1) &=\overline y_l(1)
\endalign
$$
where $\alpha_1(1),\ldots,\alpha_s(1)\in k(\overline U(1))$,
$\nu(\overline x_1(1)),\ldots,\nu(\overline x_s(1))$ are rationally independent,\linebreak
 $\text{det}(c_{ij}(1))\ne 0$ and
$(R,\overline T''(1),\overline T(1))$, $(S,\overline U''(1),\overline U(1))$ is a CUTS along $\nu$.
\endproclaim

\demo{Proof}
Identify $\nu$ with our extension of $\nu$ to the quotient field of $\overline U$ which dominates $\overline U$.
Set $(g_{ij})=(a_{ij})^{-1}$,
$$
\align
A_1 &=\overline x_1^{g_{11}}\cdots\overline x_s^{g_{1s}}\overline x_r^{g_{1,s+1}}\\
&\vdots\\
A_s &=\overline x_1^{g_{s1}}\cdots\overline x_s^{g_{ss}}\overline x_r^{g_{s,s+1}}\\
A_r &=\overline x_1^{g_{s+1,1}}\cdots\overline x_s^{g_{s+1,s}}\overline x_r^{g_{s+1,s+1}}.
\endalign
$$
Then $\overline T(1)'$ is a localization of $\overline T''[A_1,\ldots,A_s,A_r]$. $\nu(A_i)>0$ for $1\le i\le s$ and $\nu(A_r)=0$.
We have
$$
\align
A_1 &=\overline y_1^{d_{11}}\cdots\overline y_s^{d_{1s}}\overline y_r^{d_{1,s+1}}\beta_1\\
&\vdots\\
A_s &=\overline y_1^{d_{s1}}\cdots\overline y_s^{d_{ss}}\overline y_r^{d_{s,s+1}}\beta_s\\
A_r &=\overline y_1^{d_{s+1,1}}\cdots\overline y_s^{d_{s+1,s}}\overline y_r^{d_{s+1,s+1}}\beta_r
\endalign
$$
where $\beta_i=\alpha_1^{g_{i1}}\cdots\alpha_s^{g_{is}}$ for $1\le i\le s$,
 $\beta_r=\alpha_1^{g_{s+1,1}}\cdots\alpha_s^{g_{s+1,s}}$ and
$$
(d_{ik})=(a_{ij})^{-1}\left(\matrix(c_{jk})&0\\0&1\endmatrix\right).\tag 3.12
$$
Define $B_i$ by
$$
\align
B_1&=\overline y_1^{h_{11}}\cdots\overline y_s^{h_{1s}}\overline y_r^{h_{1,s+1}}\\
&\vdots\\
B_s&=\overline y_1^{h_{s1}}\cdots\overline y_s^{h_{ss}}\overline y_r^{h_{s,s+1}}\\
B_r&=\overline y_1^{h_{s+1,1}}\cdots\overline y_s^{h_{s+1,s}}\overline y_r^{h_{s+1,s+1}}
\endalign
$$
where the matrix $(h_{ij})$ defines 
 a UTS of type $\text{II}_r$ $\overline U\rightarrow W$ along $\nu$ where $\overline W'$ is a localization of
$\overline U''[B_1,\ldots,B_s,B_r]$ with $\nu(B_i)>0$ for $1\le i\le s$ and $\nu(B_r)=0$.
Let $d$ be the residue of $B_r$ in $k(W)$. We have
$$
\align
A_1&=B_1^{e_{11}}\cdots B_s^{e_{1s}}B_r^{e_{1,s+1}}\beta_1\\
&\vdots\\
A_s&=B_1^{e_{s1}}\cdots B_s^{e_{ss}}B_r^{e_{s,s+1}}\beta_s\\
A_r&=B_1^{e_{s+1,1}}\cdots B_s^{e_{s+1,s}}B_r^{e_{s+1,s+1}}\beta_r
\endalign
$$
where $(e_{ij})=(d_{ij})(h_{ij})^{-1}$ is a matrix with integral coefficients. Since 
$\nu(A_r)=\nu(B_r)=0$ and $\nu(B_1),\cdots,\nu(B_s)$ are rationally independent, we have $e_{s+1,1}=\cdots=e_{s+1,s}=0$.
Then $\text{det}(e_{ij})\ne0$ implies $e_{s+1,s+1}\ne 0$. Since $\nu(A_1),\ldots,\nu(A_s)>0$,
by Lemma 3.2, we can perform a UTS of type
I $\overline W\rightarrow \overline U(1)$ along $\nu$ so that $\overline U'(1)$ is a localization of $\overline W'[C_1,\ldots,C_s]$
with 
$\nu(C_1),\ldots,\nu(C_s)$ rationally independent, and
$$
\align
B_1 &=C_1^{b_{11}''}\cdots C_s^{b_{1s}''}\\
&\vdots\\
B_s &=C_1^{b_{s1}''}\cdots C_s^{b_{ss}''}
\endalign
$$
to get
$$
\align
A_1&=C_1^{f_{11}}\cdots C_s^{f_{1s}}B_r^{f_{1,s+1}}\beta_1\tag 3.13\\
&\vdots\\
A_s&=C_1^{f_{s1}}\cdots C_s^{f_{ss}}B_r^{f_{s,s+1}}\beta_s\\
A_r&=C_1^{f_{s+1,1}}\cdots C_s^{f_{s+1,s}}B_r^{f_{s+1,s+1}}\beta_r
\endalign
$$
with $f_{ij}\ge 0$ for all $i,j$ and 
$$
(f_{ij}) = (e_{ij})\left(\matrix (b_{ij}'')&0\\0&1\endmatrix\right).
$$
Set
$$
(b_{ij})=(h_{ij})^{-1}\left(\matrix (b_{ij}'')&0\\0&1\endmatrix\right).
$$
$(3.13)$ implies $A_1,\ldots, A_r\in \overline U'(1)$. Thus $\overline U''\rightarrow \overline U'(1)$ is a MTS along $\nu$ and
$\overline T'(1)\subset \overline U'(1)$. Further, we have $f_{s+1,1}=\cdots=f_{s+1,s}=0$.

Extend $\nu$ from the quotient field of $\overline U$ to a valuation of the quotient field of $\overline U(1)$
which dominates $\overline U(1)$. $\overline U(1)$ has regular parameters $(y_1^*(1),\ldots, y_n^*(1))$
such that 
$$
\align
\overline y_1 &= y_1^*(1)^{b_{11}}\cdots y_s^*(1)^{b_{1s}}(y_r^*(1)+d)^{b_{1,s+1}}\tag 3.14\\
\vdots&\\
\overline y_s &= y_1^*(1)^{b_{s1}}\cdots y_s^*(1)^{b_{ss}}(y_r^*(1)+d)^{b_{s,s+1}}\\
\overline y_r &= y_1^*(1)^{b_{s+1,1}}\cdots y_s^*(1)^{b_{s+1,s}}(y_r^*(1)+d)^{b_{s+1,s+1}}.
\endalign
$$
 Set
$$
\overline U''(1) = \overline U'(1)\left[d,\left(\frac{y_r^*(1)}{d}+1\right)^{\frac{1}{\overline d}}\right]
_{(y_1^*(1),\ldots,y_n^*(1))}.
$$
We have a natural inclusion $\overline T(1)\subset\overline U(1)$.

Let $c\in k(\overline T'(1))$ be the residue of $A_r$. Then $c$ is the residue of $A_r$ in the residue field of
our extension of $\nu$ to the quotient field of $\overline U(1)$, since $\nu$ dominates $\overline U''(1)$. Set
$$
x_i^*(i)=\cases
A_i&1\le i\le s\\
\overline x_i&s<i,i\ne r\\
A_r-c&i=r\endcases
$$
Then $(x_1^*(1),\ldots,x_n^*(1))$ are regular parameters in $\overline T(1)$ such that
$$
\align
x_1^*(1) &= y_1^*(1)^{f_{11}}\cdots y_s^*(1)^{f_{1s}}(y_r^*(1)+d)^{f_{1,s+1}}\beta_1\tag 3.15\\
\vdots&\\
x_s^*(1) &= y_1^*(1)^{f_{s1}}\cdots y_s^*(1)^{f_{ss}}(y_r^*(1)+d)^{f_{s,s+1}}\beta_s\\
x_r^*(1)+c &= y_1^*(1)^{f_{s+1,1}}\cdots y_s^*(1)^{f_{s+1,s}}(y_r^*(1)+d)^{f_{s+1,s+1}}\beta_r.
\endalign
$$
They are related to the regular parameters $\overline x_i(1)$ in $\overline T(1)$ satisfying $(3.9)$ by
$$
x_i^*(1)=\cases
\overline x_i(1)(\overline x_r(1)+1)^{\overline c\gamma_i}&1\le i\le s\\
c[(\overline x_r(1)+1)^{\overline c}-1]&i=r
\endcases
$$
where
$$
\overline c=\text{Det}\left(\matrix a_{11}&\cdots& a_{1s}\\
\vdots&&\vdots\\
a_{s1}&\cdots& a_{ss}
\endmatrix \right)
\text{Det}\left(\matrix a_{11}&\cdots& a_{1,s+1}\\
\vdots&&\vdots\\
a_{s1}&\cdots& a_{s+1,s+1}
\endmatrix \right)
$$
$$
\left(\matrix \gamma_1\\ \vdots\\ \gamma_s \\ 1\endmatrix\right)
= (a_{ij})^{-1}\left(\matrix 0 \\ \vdots \\ 0 \\\frac{1}{\overline c}
\endmatrix\right).
$$
We have regular parameters $\overline y_i(1)$ in $\overline U''(1)$ satisfying $(3.10)$ with
$$
y_i^*(1)=\cases
\overline y_i(1)(\overline y_r(1)+1)^{\overline d\tau_i}&1\le i\le s\\
d[(\overline y_r(1)+1)^{\overline d}-1]&i=r
\endcases
$$
where
$$
\overline d=\text{Det}\left(\matrix b_{11}&\cdots& b_{1s}\\
\vdots&&\vdots\\
b_{s1}&\cdots& b_{ss}
\endmatrix \right)
\text{Det}\left(\matrix b_{11}&\cdots& b_{1,s+1}\\
\vdots&&\vdots\\
b_{s1}&\cdots& b_{s+1,s+1}
\endmatrix \right)
$$
$$
\left(\matrix \tau_1\\ \vdots\\ \tau_s \\ 1\endmatrix\right)
= (b_{ij})^{-1}\left(\matrix 0 \\ \vdots \\ 0 \\ \frac{1}{\overline d}
\endmatrix\right).
$$
$$
(f_{ij})\left(\matrix \tau_1\\ \vdots\\ \tau_s \\ 1\endmatrix\right)
=  (a_{ij})^{-1}(\left(\matrix (c_{ij}) & 0\\0&1\endmatrix\right)
\left(\matrix 0 \\ \vdots \\ 0 \\ \frac{1}{\overline d}
\endmatrix\right)
= (a_{ij})^{-1}\left(\matrix 0 \\ \vdots \\ 0 \\\frac{1}{\overline d}
\endmatrix\right)
=\frac{\overline c}{\overline d}\left(\matrix \gamma_1\\ \vdots\\ \gamma_s \\ 1\endmatrix\right).
$$
Substitute this in (3.15) to get 
$$
\align
\overline x_1(1)(\overline x_r(1)+1)^{\overline c\gamma_1}
& = \overline y_1(1)^{f_{11}}\cdots \overline y_s(1)^{f_{1s}}(\overline y_r(1)+1)
^{\overline c\gamma_1}d^{f_{1,s+1}}\beta_1\\
&\vdots\\
\overline x_s(1)(\overline x_r(1)+1)^{\overline c\gamma_s}
& = \overline y_1(1)^{f_{s1}}\cdots \overline y_s(1)^{f_{ss}}(\overline y_r(1)+1)
^{\overline c\gamma_s}d^{f_{s,s+1}}\beta_s\\
c(\overline x_r(1)+1)^{\overline c}
& = (\overline y_r(1)+1)^{\overline c}d^{f_{s+1,s+1}}\beta_r.
\endalign
$$
$\nu(\overline x_r(1))>0$ and $\nu(\overline y_r(1))>0$ imply 
$$
c=d^{f_{s+1,s+1}}\beta_r.\tag 3.16
$$
Our inclusion $\overline T(1)\subset \overline U(1)$ induces 
$$
\overline x_r(1)=\omega\overline y_r(1)+\omega-1
$$
in $\overline U(1)$ for some $\overline c$-th root of unity $\omega$. Since $\overline x_r(1)\in  m(\overline T(1))$ and
$\overline y_r(1)\in m(\overline U(1))$, we must have $\omega=1$.
We thus get
(3.11).
$$
\overline T''(1) = \overline T''[c,\overline x_r(1)]_{(\overline x_1(1),\ldots,\overline x_n(1))}
\subset \overline U''[d,\overline y_r(1)]_{(\overline y_1(1),\ldots,\overline y_n(1))} = \overline U''(1).
$$
An extension of $\nu_0$ to the quotient field of $\overline U(1)$ which dominates $\overline U(1)$ then makes 
$(R,\overline T''(1),\overline T(1))$, $(S,\overline U''(1),\overline U(1))$ a CUTS along $\nu$.
\enddemo

\proclaim{Lemma 3.5} Suppose that $(R,\overline T'',\overline T)$, $(S,\overline U'',\overline U)$ is a CUTS along $\nu$,
$\overline T''$ has regular parameters $(\overline x_1,\ldots,\overline x_n)$ and
$\overline U''$ has regular parameters $(\overline y_1,\ldots,\overline y_n)$ such that
$$
\align
\overline x_1 &= \overline  y_1^{c_{11}}\cdots \overline y_s^{c_{1s}}\alpha_1\tag 3.17\\
\vdots&\\
\overline x_s &= \overline y_1^{c_{s1}}\cdots \overline y_s^{c_{ss}}\alpha_s\\
\overline x_{s+1} &= \overline y_{s+1}\\
\vdots&\\
\overline x_{l} &= \overline y_{l}\\
\overline x_{l+1}&=\overline y_1^{d_1}\cdots\overline y_s^{d_s}\overline y_{l+1}
\endalign
$$
where $\alpha_1,\ldots,\alpha_s\in k(\overline U)$,
$\nu(\overline x_1),\ldots,\nu(\overline x_s)$ are rationally independent, $\text{det}(c_{ij})\ne 0$.

Suppose that $\overline T\rightarrow \overline T(1)$ is a UTS of type $\text{II}_{l+1}$ along $\nu$, 
such that $\overline T(1)''$ has regular parameters $(\overline x_1(1),\ldots,\overline x_n(1))$ with
$$
\align
\overline x_1 &= \overline x_1(1)^{a_{11}}\cdots \overline x_s(1)^{a_{1s}}c^{a_{1,s+1}}\tag 3.18\\
\vdots&\\
\overline x_s &= \overline x_1(1)^{a_{s1}}\cdots \overline x_s(1)^{a_{ss}}c^{a_{s,s+1}}\\
\overline x_{l+1} &= \overline x_1(1)^{a_{s+1,1}}\cdots \overline x_s(1)^{a_{s+1,s}}(\overline x_{l+1}(1)+1)c^{a_{s+1,s+1}}.
\endalign
$$
Then there exists a UTS of type $\text{II}_{l+1}$, (followed by a UTS of type I) 
$\overline U\rightarrow \overline U(1)$ along $\nu$ 
such that  $ \overline U(1)''$ has regular parameters $(\overline y_1(1),\ldots,\overline y_n(1))$ satisfying
$$
\align
\overline y_1 &= \overline y_1(1)^{b_{11}}\cdots \overline y_s(1)^{b_{1s}}d^{b_{1,s+1}}\tag 3.19\\
\vdots&\\
\overline y_s &= \overline y_1(1)^{b_{s1}}\cdots \overline y_s(1)^{b_{ss}}d^{b_{s,s+1}}\\
\overline y_{l+1} &= \overline y_1(1)^{b_{s+1,1}}\cdots \overline y_s(1)^{b_{s+1,s}}(\overline y_{l+1}(1)+1)d^{b_{s+1,s+1}},
\endalign
$$
 $\overline T''(1)\subset \overline U''(1)$, and 
$$
\align
\overline x_1(1) &= \overline y_1(1)^{c_{11}(1)}\cdots \overline y_s(1)^{c_{1s}(1)}\alpha_1(1)\tag 3.20\\
\vdots&\\
\overline x_s(1) &= \overline y_1(1)^{c_{s1}(1)}\cdots \overline y_s(1)^{c_{ss}(1)}\alpha_s(1)\\
\overline x_{s+1}(1) &=\overline y_{s+1}(1)\\
\vdots&\\
\overline x_l(1) &=\overline y_l(1)\\
\overline x_{l+1}(1) &=\overline y_{l+1}(1)
\endalign
$$
where $\alpha_1(1),\ldots,\alpha_s(1)\in k(\overline U(1))$,
$\nu(\overline x_1(1)),\ldots,\nu(\overline x_s(1))$ are rationally independent,\linebreak
 $\text{det}(c_{ij}(1))\ne 0$ and
$(R,\overline T''(1),\overline T(1))$, $(S,\overline U''(1),\overline U(1))$ is a CUTS along $\nu$.
\endproclaim

\demo{Proof} Change $r$ to $l+1$ in the proof of Lemma 3, and change $(d_{ik})$ to
$$
(d_{ik}) = (a_{ij})^{-1}
\left(\matrix c_{11} & \cdots & c_{1s}&0\\
\vdots&&\vdots&\vdots\\
c_{s1}&\cdots&c_{ss}&0\\
d_1&\cdots&d_s&1
\endmatrix\right)
$$
\enddemo

\proclaim{Lemma 3.6} Suppose that $(R,\overline T'',\overline T)$ and $(S,\overline U'',\overline U)$ is a  CUTS along $\nu$,
$\overline T''$ has regular parameters $(\overline x_1,\ldots,\overline x_n)$ and
$\overline U''$ has regular parameters $(\overline y_1,\ldots,\overline y_n)$ such that
$$
\align
\overline x_1 &= \overline  y_1^{c_{11}}\cdots \overline y_s^{c_{1s}}\alpha_1\\
\vdots&\\
\overline x_s &= \overline y_1^{c_{s1}}\cdots \overline y_s^{c_{ss}}\alpha_s\\
\overline x_{s+1} &=\overline y_{s+1}\\
\vdots&\\
\overline x_{l} &=\overline y_{l}\\
\overline x_{l+1} &= \overline y_1^{c_{s+1,1}}\cdots \overline y_s^{c_{s+1,s}}\delta
\endalign
$$
where $\alpha_1,\ldots,\alpha_s\in k(\overline U)$,
$\delta\in \overline U''$ is a unit, $\nu(\overline x_1),\ldots,\nu(\overline x_s)$ are rationally independent and
$\text{det}(c_{ij})\ne0$.

Suppose that $\overline T\rightarrow \overline T(1)$ is a UTS of type $\text{II}_{l+1}$ along $\nu$,
so that
 $\overline T(1)''$ has regular parameters $(\overline x_1(1),\ldots,\overline x_n(1))$ satisfying
$$
\align
\overline x_1 &= \overline x_1(1)^{a_{11}}\cdots \overline x_s(1)^{a_{1s}}c^{a_{1,s+1}}\\
\vdots&\\
\overline x_s &= \overline x_1(1)^{a_{s1}}\cdots \overline x_s(1)^{a_{ss}}c^{a_{s,s+1}}\\
\overline x_{l+1} &= \overline x_1(1)^{a_{s+1,1}}\cdots \overline x_s(1)^{a_{s+1,s}}(\overline x_{l+1}(1)+1)c^{a_{s+1,s+1}}.
\endalign
$$
Then there exists a UTS of type $\text{I}$ along $\nu$   $\overline U\rightarrow \overline U(1)$ 
such that $\overline U'(1)$ has regular parameters $(\hat y_1(1),\ldots,\hat y_n(1))$ with
$$
\align
\overline y_1 &= \hat y_1(1)^{b_{11}}\cdots \hat y_s(1)^{b_{1s}}\\
\vdots&\\
\overline y_s &= \hat y_1(1)^{b_{s1}}\cdots \hat y_s(1)^{b_{ss}}
\endalign
$$
and   $\overline U''(1)$ has regular parameters 
$(\overline y_1(1),\ldots,\overline y_n(1))$
such that $\overline y_i(1)=\epsilon_i\hat y_i(1)$ for $1\le i\le s$ for some units $\epsilon_i\in \overline U(1)''$,
$\overline T''(1)\subset \overline U''(1)$,
$$
\align
\overline x_1(1) &= \overline y_1(1)^{c_{11}(1)}\cdots \overline y_s(1)^{c_{1s}(1)}\alpha_1(1)\\
\vdots&\\
\overline x_s(1) &= \overline y_1(1)^{c_{s1}(1)}\cdots \overline y_s(1)^{c_{ss}(1)}\alpha_s(1)\\
\overline x_{s+1}(1)&=\overline y_{s+1}(1)\\
\vdots&\\
\overline x_{l}(1)&=\overline y_{l}(1)\\
\endalign
$$
where $\alpha_1(1),\ldots,\alpha_s(1)\in k(\overline U(1))$,
$\nu(\overline x_1(1)),\ldots,\nu(\overline x_s(1))$ are rationally independent,\linebreak
$\text{det}(c_{ij}(1))\ne0$, and $(R,\overline T''(1),\overline T(1))$ and $(S,\overline U''(1),\overline U(1))$ is a
CUTS along $\nu$.
\endproclaim

\demo{Proof}
Identify $\nu$ with our extension of $\nu$ to the quotient field of $\overline U$ which dominates $\overline U$. 
Set $(g_{ij})=(a_{ij})^{-1}$,
$$
\align
A_1&= \overline x_1^{g_{11}}\cdots\overline x_s^{g_{1s}}\overline x_{l+1}^{g_{1,s+1}}\\
&\vdots\\
A_s&= \overline x_1^{g_{s1}}\cdots\overline x_s^{g_{ss}}\overline x_{l+1}^{g_{s,s+1}}\\
A_{l+1}&= \overline x_1^{g_{s+1,1}}\cdots\overline x_s^{g_{s+1,s}}\overline x_{l+1}^{g_{s+1,s+1}}.
\endalign
$$
$\overline T(1)'$ is a localization of $\overline T''[A_1,\cdots,A_s,A_{l+1}]$. $\nu(A_i)>0$ for
$1\le i\le s$ and $\nu(A_{l+1})=0$.

$$
\align
A_1 &= \overline y_1^{d_{11}}\cdots \overline y_s^{d_{1s}}\delta_1e_1\\
\vdots&\\
A_s &= \overline y_1^{d_{s1}}\cdots \overline y_s^{d_{ss}}\delta_se_s\\
A_{l+1} &= \overline y_1^{d_{s+1,1}}\cdots \overline y_s^{d_{s+1,s}}\delta_{l+1}\\
\endalign
$$
where $e_i\in k(\overline U)$,
$(d_{ik}) = (a_{ij})^{-1}(c_{jk})$ and $\delta_i$ are units in $\overline U''$
such that $\delta_i$ has residue 1 in $k(\overline U)$
for $1\le i\le s$.
$\nu(A_{l+1})=0$ and $\nu(\overline y_1),\ldots,\nu(\overline y_s)$ rationally independent implies 
$$
d_{s+1,1} = \cdots = d_{s+1,s}=0.
$$
Since $\nu(A_i)>0$ for  $1\le i\le s$, by Lemma 3.2 we can 
perform a UTS of type I along $\nu$ $\overline U\rightarrow \overline U(1)$ where $\overline U'(1)$ has regular parameters
$(\hat y_1(1),\ldots,\hat y_n(1))$ satisfying 
$$
\align
\overline y_1 &= \hat y_1(1)^{b_{11}(1)}\cdots \hat y_s(1)^{b_{1s}(1)}\\
\vdots&\\
\overline y_s &= \hat y_1(1)^{b_{s1}(1)}\cdots \hat y_s(1)^{b_{ss}(1)}
\endalign
$$
to get 
$$
\align
A_1 &= \hat y_1(1)^{c_{11}(1)}\cdots \hat y_s(1)^{c_{1s}(1)}\delta_1e_1\\
\vdots&\\
A_s &= \hat y_1(1)^{c_{s1}(1)}\cdots \hat y_s(1)^{c_{ss}(1)}\delta_se_s\\
A_{l+1} &= \delta_{l+1}
\endalign
$$
where all $c_{ij}(1)\ge 0$. Thus 
$$
\overline T''[A_1,\ldots,A_s,A_{l+1}]\subset \overline U'(1)=
\overline U''[\hat y_1(1),\ldots,\hat y_s(1)]_{(\hat y_1(1),\ldots,\hat y_n(1))},
$$
and since $\nu$ dominates $\overline U'(1)$ and $\overline T'(1)$, $\overline U'(1)$ dominates $\overline T'(1)$.

Now extend $\nu$ from the quotient field of $\overline U$ to a valuation of the quotient field of 
$\overline U(1)$ which dominates $\overline U(1)$. 
$\overline T(1)''$ has regular parameters $(x_1^*(1),\ldots, x_n^*(1))$ with
$$
\align
\overline x_1 &= x_1^*(1)^{a_{11}}\cdots x_s^*(1)^{a_{1s}}(x_{l+1}^*(1)+c)^{a_{1,s+1}}\\
\vdots&\\
\overline x_s &= x_1^*(1)^{a_{s1}}\cdots x_s^*(1)^{a_{ss}}(x_{l+1}^*(1)+c)^{a_{s,s+1}}\\
\overline x_{l+1} &= x_1^*(1)^{a_{s+1,1}}\cdots x_s^*(1)^{a_{s+1,s}}(x_{l+1}^*(1)+c)^{a_{s+1,s+1}}.
\endalign
$$
$$
\overline T''(1) = \overline T''\left[c,\left(\frac{x_{l+1}^*(1)}{c}+1\right)^{\frac{1}{\overline c}}\right]
_{(x_1^*(1),\ldots,x_n^*(1))}.
$$
$(\overline x_1(1),\ldots,\overline x_n(1))$ are regular parameters in $\overline T''(1)$ which satisfy
$$
\align
\overline x_1(1) & = x_1^*(1)\left(\frac{x_{l+1}^*(1)}{c}+1\right)^{-\gamma_1} = \hat y_1(1)^{c_{11}(1)}\cdots \hat y_s(1)^{c_{1s}(1)}
\delta_1e_1\left(\frac{\delta_{l+1}}{c}\right)^{-\gamma_1}\\
\vdots&\\
\overline x_s(1) &= x_s^*(1)\left(\frac{x_{l+1}^*(1)}{c}+1\right)^{-\gamma_s} =\hat y_1(1)^{c_{s1}(1)}\cdots \hat y_s(1)^{c_{ss}(1)}
\delta_se_s\left(\frac{\delta_{l+1}}{c}\right)^{-\gamma_s}\\
\overline x_{l+1}(1) &= \left(\frac{x_{l+1}^*(1)}{c}+1\right)^{\frac{1}{\overline c}}-1=\left(\frac{\delta_{l+1}}{c}\right)
^{\frac{1}{\overline c}}-1
\endalign
$$
Set $(e_{ij})=(c_{ij}(1))^{-1}$, 
$$
\epsilon_i=\delta_1^{e_{i1}}\cdots\delta_s^{e_{is}}\left(\frac{\delta_{l+1}}{c}\right)^{-\gamma_1e_{i1}-\cdots-\gamma_se_{is}}
$$
 for
$1\le i\le s$. Define 
$$
\overline y_i(1) =
\cases \epsilon_i \hat y_i(1)& 1\le i\le s\\
\hat y_i(1) & s+1\le i.
\endcases
$$
Then the conclusions of Lemma 3.6 hold with
$$
\overline U''(1)= \overline U'(1)[c,\left(\frac{\delta_{l+1}}{c}\right)
^{\frac{1}{\overline c}},\epsilon_1,\ldots,\epsilon_s]_{(\overline y_1(1),\ldots,\overline y_n(1))}.
$$ 
\enddemo

\subheading{Monomialization in rank 1}

\proclaim{Theorem 3.8}
Suppose that $(R,\overline T'',\overline T)$ and $(S,\overline U'',\overline U)$  is a CUTS along $\nu$
such that $\overline T''$ contains the subfield  $k(c_0)$ for some $c_0\in \overline T''$ and $\overline U''$
contains a subfield isomorphic to $k(\overline U'')$,
$\overline T''$ has
regular parameters $(\overline z_1,\ldots, \overline z_n)$ and $\overline U''$ has regular parameters 
$(\overline w_1,\ldots, \overline w_n)$
such that 
$$
\align
\overline z_1 &= \overline w_1^{c_{11}}\cdots \overline w_s^{c_{1s}}\phi_1\\
\vdots&\\
\overline z_s &= \overline w_1^{c_{s1}}\cdots \overline w_s^{c_{ss}}\phi_s\\
\overline z_{s+1} &= \overline w_{s+1}\\
\vdots&\\
\overline z_l&=\overline w_l.
\endalign
$$
where $\phi_1,\ldots,\phi_s\in k(\overline U'')$,
$\nu(\overline z_1),\ldots,\nu(\overline z_s)$ are rationally independent, $\text{det}(c_{ij})\ne0$.

Suppose that one of the following three conditions hold.
\item{1)} $f\in k(\overline U)[[\overline w_1,\ldots, \overline w_m]]$ for some $m$ such that $s\le m\le n$ 
with $\nu(f)<\infty$. 
\item{2)} $f\in k(\overline U)[[\overline w_1,\ldots, \overline w_m]]$ for some $m$ such that $s< m\le n$ 
with $\nu(f)=\infty$ and $A>0$ is given.
\item{3)} 
$$
f \in \left(k(\overline U)[[\overline w_1,\ldots,\overline w_m]]
 - k(\overline U)[[\overline w_1,\ldots,\overline w_l]]\right)\cap \overline U''
$$
\item{}for some $m$ such that $l<m\le n$.
Then there exists a CUTS along $\nu$ $(R,\overline T''(t),\overline T(t)$ and $(S,\overline U''(t),\overline U(t))$
$$
\matrix
\overline U = &\overline U(0) & \rightarrow &\overline U(1) &\rightarrow &\cdots &\rightarrow & \overline U(t)\\
&\uparrow &&\uparrow&&&&\uparrow\\
\overline T= &\overline T(0) &\rightarrow &\overline T(1)&\rightarrow &\cdots&\rightarrow& \overline T(t)
\endmatrix\tag 3.21
$$
such that
$\overline T''(i)$ has regular parameters $(\overline z_1(i),\ldots, \overline z_n(i))$,
 $\overline U''(i)$ has regular parameters $(\overline w_1(i),\ldots, \overline w_n(i))$ satisfying
$$
\align
\overline z_1(i) &= \overline w_1(i)^{c_{11}(i)}\cdots \overline w_s(i)^{c_{1s}(i)}\phi_1(i)\\
\vdots&\\
\overline z_s(i) &= \overline w_1(i)^{c_{s1}(i)}\cdots \overline w_s(i)^{c_{ss}(i)}\phi_s(i)\\
\overline z_{s+1}(i) &= \overline w_{s+1}(i)\\
\vdots&\\
\overline z_l(i)&=\overline w_l(i)
\endalign
$$
$\overline T''(i)$ contains a subfield  $k(c_0,\ldots,c_i)$ and $\overline U''(i)$ contains a subfield
isomorphic to $k(\overline U(i))$. 
$\phi_1(i),\ldots,\phi_s(i)\in k(\overline U(i))$,
$\nu(\overline z_1(i)),\ldots,\nu(\overline z_s(i))$ are rationally independent, $\text{det}(c_{ij}(i))\ne0$
for $0\le i\le t$.
In case 1) we have
$$
f = \overline w_1(t)^{d_1}\cdots\overline w_s(t)^{d_s}u(\overline w_1(t),\ldots,\overline w_m(t))
$$
where $u\in k(\overline U(t))[[\overline w_1(t),\ldots,\overline w_m(t)]]$ is a unit power series.

In case 2) we have
$$
f = \overline w_1(t)^{d_1}\cdots\overline w_s(t)^{d_s}\Sigma(\overline w_1(t),\ldots,\overline w_m(t))
$$
where $\Sigma\in k(\overline U(t))[[\overline w_1(t),\ldots,\overline w_m(t)]]$,
$\nu(\overline w_1(t)^{d_1}\cdots\overline w_s(t)^{d_s})>A$.

In case 3) we have
$$
f = P(\overline w_1(t),\ldots,\overline w_l(t))+
   \overline w_1(t)^{d_1}\cdots\overline w_s(t)^{d_s}H
$$
for some powerseries $P\in k(\overline U(t))[[\overline w_1(t),\ldots,\overline w_l(t)]]$,
$$
H= u(\overline w_{m}(t)+\overline w_1(t)^{g_1}\cdots\overline w_s(t)^{g_s}\Sigma)
$$
where $u\in k(\overline U(t))[[\overline w_1(t),\ldots,\overline w_m(t)]]$ is a unit,
$\Sigma\in k(\overline U(t))[[\overline w_1(t),\ldots,\overline w_{m-1}(t)]]$
and $\nu(\overline w_m(t))\le \nu(\overline w_1(t)^{g_1}\cdots\overline w_s(t)^{g_s})$.

$(3.21)$ will be
such that $\overline T''(\alpha)$ has regular parameters
$$
(\overline z_1(\alpha),\ldots, \overline z_n(\alpha))\text{ and } 
(\tilde{\overline  z}'_1(\alpha),\ldots, \tilde{\overline z}'_n(\alpha)),
$$
$\overline U''(\alpha)$ has regular parameters
$$
(\overline w_1(\alpha),\ldots, \overline w_n(\alpha))\text{ and }
(\tilde{\overline  w}'_1(\alpha),\ldots, \tilde{\overline w}'_n(\alpha))
$$
where $\overline z_i(0)=\overline z_i$ and  $\overline w_i(0)=\overline w_i$
for $1\le i\le n$.
$(3.21)$ will consist of three  types of CUTS.
\item{M1)} $\overline T(\alpha)\rightarrow \overline T(\alpha+1)$ and $\overline U(\alpha)\rightarrow 
\overline U(\alpha+1)$ are of type I.
\item{M2)} $\overline T(\alpha)\rightarrow \overline T(\alpha+1)$ is of type $\text{II}_r$,  $s+1\le r\le l$,
and $\overline U(\alpha)\rightarrow \overline U(\alpha+1)$ 
is a transformation  of type $\text{II}_r$,  followed by a transformation of type I.
\item{M3)} $\overline T(\alpha)= \overline T(\alpha+1)$ and 
$\overline U(\alpha)\rightarrow \overline U(\alpha+1)$ is of type $\text{II}_r$ ($l+1\le r\le m$).

We will find polynomials $P_{i,\alpha}$ so that the variables will be related by:

$$
\tilde{\overline z}'_i(\alpha) = 
\cases
\overline z_i(\alpha)-P_{i,\alpha}(\overline z_1(\alpha),\ldots,
\overline z_{i-1}(\alpha))  &\text{ if }s+1\le i\le l\\  
\overline z_i(\alpha) & \text{ otherwise}
\endcases
$$

$$
\tilde{\overline w}'_i(\alpha) = 
\cases
\tilde{\overline z}'_i(\alpha)&\text{ if } s+1\le i\le l\\
\overline w_i(\alpha)-P_{i,\alpha}(\overline w_1(\alpha),\ldots,
\overline w_{i-1}(\alpha))  &\text{ if }l+1\le i\le m\\  
\overline w_i(\alpha) & \text{ otherwise}
\endcases
$$
The coefficients of $P_{i,\alpha}$ will be in $k(c_0,\ldots,c_{\alpha})$ if $i\le l$, and will be in 
$k(\overline U(\alpha))$ if $i>l$.
For all $\alpha$ we will have
$$
\align
\overline z_1(\alpha) &= \overline w_1(\alpha)^{c_{11}(\alpha)}\cdots \overline w_s(\alpha)^{c_{1s}(\alpha)}\phi_1(\alpha)\tag 3.22\\
\vdots&\\
\overline z_s(\alpha) &= \overline w_1(\alpha)^{c_{s1}(\alpha)}\cdots \overline w_s(\alpha)^{c_{ss}(\alpha)}\phi_s(\alpha)\\
\overline z_{s+1}(\alpha) &= \overline w_{s+1}(\alpha)\\
\vdots&\\
\overline z_l(\alpha)&=\overline w_l(\alpha)
\endalign
$$
and
$$
\align
\tilde{\overline z}'_1(\alpha) &= \tilde{\overline w}'_1(\alpha)^{c_{11}(\alpha)}\cdots 
\tilde{\overline w}'_s(\alpha)^{c_{1s}(\alpha)}\phi_1(\alpha)\tag 3.23\\
\vdots&\\
\tilde{\overline z}'_s(\alpha) &= \tilde{\overline w}'_1(\alpha)^{c_{s1}(\alpha)}\cdots 
\tilde{\overline w}'_s(\alpha)^{c_{ss}(\alpha)}\phi_s(\alpha)\\
\tilde{\overline z}'_{s+1}(\alpha) &= \tilde{\overline w}'_{s+1}(\alpha)\\
\vdots&\\
\tilde{\overline z}_l(\alpha)&=\tilde{\overline w}'_l(\alpha)
\endalign
$$
where $\phi_1(\alpha),\ldots,\phi_s(\alpha)\in k(\overline U(\alpha))$. $\overline T''(\alpha)$
contains a subfield $k(c_0,\ldots,c_{\alpha})$  and 
$\overline U''(\alpha)$ contains a subfield isomorphic to $k(\overline U(\alpha))$.

In a transformation $\overline T(\alpha)\rightarrow \overline T(\alpha+1)$ of type I $\overline T''(\alpha+1)$ will have
regular parameters $(\overline z_1(\alpha+1),\ldots,\overline z_n(\alpha+1))$ 
$$
\align
\tilde{\overline z}'_1(\alpha) &= \overline z_1(\alpha+1)^{a_{11}(\alpha+1)}\cdots \overline z_s(\alpha+1)^{a_{1s}(\alpha+1)}\tag 3.24\\
\vdots&\\
\tilde{\overline  z}'_s(\alpha) &= \overline z_1(\alpha+1)^{a_{s1}(\alpha+1)}\cdots \overline z_s(\alpha+1)^{a_{ss}(\alpha+1)}
\endalign
$$
and $c_{\alpha+1}$ is defined to be $1$.
In a transformation $\overline T(\alpha)\rightarrow \overline T(\alpha+1)$ of type $\text{II}_r$  ($s+1\le r\le l$) 
$\overline T''(\alpha+1)$ will have regular parameters $(\overline z_1(\alpha+1),\ldots,\overline z_n(\alpha+1))$ 
$$
\align
\tilde{\overline  z}'_1(\alpha) &= \overline z_1(\alpha+1)^{a_{11}(\alpha+1)}\cdots \overline z_s(\alpha+1)^{a_{1s}(\alpha+1)}
c_{\alpha+1}^{a_{1,s+1}(\alpha+1)}\tag 3.25\\
\vdots&\\
\tilde{\overline z}'_s(\alpha) &=
 \overline z_1(\alpha+1)^{a_{s1}(\alpha+1)}\cdots \overline z_s(\alpha+1)^{a_{ss}(\alpha+1)}c_{\alpha+1}^{a_{s,s+1}(\alpha+1)}\\
\tilde{\overline z}'_r(\alpha) &= \overline z_1(\alpha+1)^{a_{s+1,1}(\alpha+1)}\cdots \overline z_s(\alpha+1)^{a_{s+1,s}(\alpha+1)}
(\overline z_r(\alpha+1)+1)c_{\alpha+1}^{a_{s+1,s+1}(\alpha+1)}
\endalign
$$
In a transformation $\overline U(\alpha)\rightarrow \overline U(\alpha+1)$ of type I
$\overline U''(\alpha+1)$ will have regular parameters $(\overline w_1(\alpha+1),\ldots,\overline w_n(\alpha+1))$
$$
\align
\tilde{\overline w}'_1(\alpha) &= \overline w_1(\alpha+1)^{b_{11}(\alpha+1)}\cdots \overline w_s(\alpha+1)^{b_{1s}(\alpha+1)}\tag 3.26\\
\vdots&\\
\tilde{\overline  w}'_s(\alpha) &= \overline w_1(\alpha+1)^{b_{s1}(\alpha+1)}\cdots \overline w_s(\alpha+1)^{b_{ss}(\alpha+1)}
\endalign
$$
and $d_{\alpha+1}$ is defined to be 1.
In a transformation $\overline U(\alpha)\rightarrow \overline U(\alpha+1)$ of type $\text{II}_r$  ($s+1\le r\le m$)
$\overline U''(\alpha+1)$ will have regular parameters $(\overline w_1(\alpha+1),\ldots,\overline w_n(\alpha+1))$
$$
\align
\tilde{\overline  w}'_1(\alpha) &= \overline w_1(\alpha+1)^{b_{11}(\alpha+1)}\cdots \overline w_s(\alpha+1)^{b_{1s}(\alpha+1)}
d_{\alpha+1}^{b_{1,s+1}(\alpha+1)}\tag 3.27\\
\vdots&\\
\tilde{\overline w}'_s(\alpha) &=
 \overline w_1(\alpha+1)^{b_{s1}(\alpha+1)}\cdots \overline w_s(\alpha+1)^{b_{ss}(\alpha+1)}d_{\alpha+1}^{b_{s,s+1}(\alpha+1)}\\
\tilde{\overline w}'_r(\alpha) &= \overline w_1(\alpha+1)^{b_{s+1,1}(\alpha+1)}\cdots \overline w_s(\alpha+1)^{b_{s+1,s}(\alpha+1)}
(\overline w_r(\alpha+1)+1)d_{\alpha+1}^{b_{s+1,s+1}(\alpha+1)}
\endalign
$$
In a transformation of type M2) $c_{\alpha+1}$ is related to $d_{\alpha+1}$ by (3.16) of the proof of Lemma 3.4.

We will call a UTS (CUTS) as in (3.21) a UTS (CUTS) in the first m variables.
\endproclaim

\demo{Proof} 
We will first  show that it is possible to construct a UTS along $\nu$
$$
\overline T \rightarrow \overline T(1)\rightarrow\cdots\rightarrow\overline T(t)\tag 3.28
$$
so that the conditions 1)', 2)' and 4) below hold.

\item{4)} Suppose that $s\le m\le l$. Then there exists a UTS (3.28) in the first $m$ variables 
 such that 
$$
p_m(i) = \{f\in k(\overline T(i))[[\overline z_1(i),\ldots,\overline z_m(i)]]\,|\, \nu(f)=\infty\}
$$
has the form 
$$
p_m(t) = (\overline z_{r(1)}(t)-Q_{r(1)}(\overline z_1(t),\cdots,\overline z_{r(1)-1}),\ldots,
\overline z_{r(\tilde m)}(t)-Q_{r(\tilde m)}(\overline z_1(t),\cdots,\overline z_{r(\tilde m)-1}))\tag P(m)
$$
for some $0\le \tilde m\le m-s$ and $s<r(1)< r(2)<\cdots<r(\tilde m)\le m$, where $Q_{r(i)}$ are power series with coefficients in 
$k(c_0,\ldots,c_t)$.

\item{1')} Suppose that  $h\in k(\overline T)[[\overline z_1,\ldots, \overline z_m]]$
 for some $m$ with $s\le m\le n$ 
and $\nu(h)<\infty$. Then there exists a UTS (3.28), in  the first $m$ variables 
 such that $P(m)$ holds in $\overline T(t)$ and 
$$
h = \overline z_1(t)^{d_1}\cdots\overline z_s(t)^{d_s}u(\overline z_1(t),\ldots,\overline z_m(t))
$$
where $u$ is a unit power series with coefficients in $k(\overline T(t))$. If 
$h\in k(c_0)[[\overline z_1,\ldots,\overline z_m]]$ then $u$ has coefficients in $k(c_0,\ldots, c_t)$.

\item{2')} Suppose that 
$h\in k(\overline T)[[\overline z_1,\ldots, \overline z_m]]$ for some $m$ with $s< m\le n$,
and $\nu(h)=\infty$ and $A>0$ is given. Then there exists a UTS (**), in the first $m$ variables 
 such that $P(m)$ holds in $\overline T(t)$ and 
$$
h = \overline z_1(t)^{d_1}\cdots\overline z_s(t)^{d_s}\Sigma(\overline z_1(t),\ldots,\overline z_m(t))
$$
where $\nu(\overline z_1(t)^{d_1}\cdots\overline z_s(t)^{d_s})>A$, $\Sigma$ is a power series with coefficients in 
$k(\overline T(t))$. If $h\in k(c_0)[[\overline z_1,\ldots,\overline z_m]]$, then $\Sigma$ has coefficients in 
$k(c_0,\ldots,c_t)$.

We will establish 1'), 2') and 4) by proving the following inductive statements.

\item{$A(\overline m)$:} 1'), 2') and 4) for $m<\overline  m$ imply 4) for $m=\overline m$.

\item{$B(\overline m)$:} 1'), 2') for $m<\overline m$ and 4) for $m=\overline m$ imply 1') and 2') for $m=\overline m$.

\subheading{Proof of $A(\overline m)$} 
By assumption there exists a UTS $\overline T\rightarrow \overline T(t)$ satisfying 4) for $\overline m-1$. After replacing
$\overline T''(0)$ with $\overline T''(t)$ and replacing $c_0$ with a primitive element of $k(c_0,\ldots,c_t)$ over $k$,
we may assume that
$$
p_{\overline m-1} = (\overline z_{r(1)}-Q_{r(1)}(\overline z_1,\cdots,\overline z_{r(1)-1}),\ldots,
\overline z_{r(\tilde{ m-1})}-Q_{r(\tilde{ m-1})}(\overline z_1,\cdots,\overline z_{r(\tilde{ m-1})-1})).
$$
where $Q_{r(i)}$ are power series with coefficients in $k(c_0)$.
If $p_{\overline m-1}k(\overline T)[[\overline z_1,\ldots, \overline z_m]]
=p_{\overline m}$ we are done. So suppose that there
exists $\overline f\in p_{\overline m}-p_{\overline m-1}k(\overline T)[[\overline z_1,\ldots, \overline z_{\overline m}]]$.
Let $L$ be a Galois closure of $k(\overline T)$ over $k(c_0)$, $G$ be the Galois group of $L$ over $k(c_0)$. Set
$$
f=\prod_{\sigma\in G}\sigma(\overline f)\in k(c_0)[[\overline z_1,\ldots,\overline z_{\overline m}]].
$$
 $f\in p_{\overline m}\cap k(c_0)[[\overline z_1,\ldots,\overline z_{\overline m}]]$ and $\nu(f)=\infty$ since $\overline f|f$
in $k(\overline T)[[\overline z_1,\ldots,\overline z_m]]$. Suppose 
$$
f\in p_{\overline m-1}k(\overline T)[[\overline z_1,\ldots,\overline z_{\overline m}]].
$$
Then $f\in p_{\overline m-1}L[[\overline z_1,\ldots,\overline z_{\overline m}]]$ which is a prime ideal, and 
$\sigma(\overline f)\in p_{\overline m-1}L[[\overline z_1,\ldots,\overline z_{\overline m}]]$ for some $\sigma\in G$.
But
$$
\sigma\left( p_{\overline m-1}L[[\overline z_1,\ldots,\overline z_{\overline m}]]\right) 
=  p_{\overline m-1}L[[\overline z_1,\ldots,\overline z_{\overline m}]]
$$
for all $\sigma\in G$. Thus 
$$
\overline f\in \left(p_{\overline m-1}L[[\overline z_1,\ldots,\overline z_{\overline m}]]\right)\cap 
k(\overline T)[[\overline z_1,\ldots,\overline z_{\overline m}]]
=p_{\overline m-1}k(\overline T)[[\overline z_1,\ldots,\overline z_{\overline m}]]
$$
a contradiction. Thus $f\not\in p_{\overline m-1}k(\overline T)[[\overline z_1,\ldots,\overline z_m]]$.
$$
f=\sum_{i=0}^{\infty}a_i(\overline z_1,\ldots,\overline z_{\overline m-1})\overline z_{\overline m}^i.
$$
where the $a_i$ have coefficients in $k(c_0)$.
By assumption $\nu(a_i)<\infty$ for some $i$. Set $r=\text{mult}(f(0,\ldots,0,\overline z_m))$.
$$
f = \sum_{i=1}^d a_i\overline z_{\overline m}^{f_i} + \sum_ja_j\overline z_{\overline m}^{f_j}
+ \sum_ka_k\overline z_{\overline m}^{f_k}+ \overline z_{\overline m}^{\alpha}\Omega
$$
where the first sum consists of the terms of minimal value $\rho=\nu(a_i\overline z_{\overline m}^{f_i})$,
$1\le i\le d$,
 $\alpha\nu(\overline z_{\overline m})>\rho$,
the second sum is a finite sum of terms of finite value 
$\nu(a_j\overline z_{\overline m}^{f_j})>\rho$ and the third sum is a finite sum of terms $a_k\overline z_{\overline m}^{f_k}$ of
infinite value.

Set 
$$
\overline R = k(c_0)[[\overline z_1,\ldots,\overline z_{\overline m-1}]]/\left(p_{\overline m-1}\cap
k(c_0)[[\overline z_1,\ldots,\overline z_{\overline m-1}]]\right).
$$
  $\nu$ induces a rank 1 valuation on the quotient field of $\overline R$.

Given $d$ in the value group $\Gamma_{\nu}\subset \bold R$ of $\nu$, let
$$
I_d=\{f\in \overline R|\nu(f)\ge d\}.
$$
There is a  set of real numbers
$$
d_1<d_2<\cdots<d_i<\tag 3.29
$$
with $\text{Lim}_{i\rightarrow\infty}d_i=\infty$ such that $d_i$
 are the possible finite values of elements of $\overline R$ and $\cap_{i=1}^{\infty}I_{d_i}=0$.
 Thus there is a function $\gamma(i)$ such that
$I_{d_i}\subset m(\overline R)^{\gamma(i)}$ and $\gamma(i)\rightarrow \infty$ as $i\rightarrow \infty$.

By  assumption, we can construct a UTS in the first $\overline m-1$ variables along $\nu$
so that for all $i,j,k$
$$
\align
a_i&=\overline z_1(t_1)^{e_1(i)}\cdots \overline z_s(t_1)^{e_s(i)}\overline a_i\\
a_j&=\overline z_1(t_1)^{f_1(j)}\cdots \overline z_s(t_1)^{f_s(j)}\overline a_j\\
a_k&=\overline z_1(t_1)^{g_1(k)}\cdots \overline z_s(t_1)^{g_s(k)}\Sigma_k
\endalign
$$
in $k(c_0,\ldots,c_{t_1})[[\overline z_1(t_1),\ldots,\overline z_{\overline m-1}(t_1)]]$
where $\overline a_i$, $\overline a_j$ are units and 
$$
\nu(\overline z_1(t_1)^{g_1(k)}\cdots \overline z_s(t_1)^{g_s(k)})>\rho.
$$

Now perform a UTS of type $\text{II}_{\overline m}$ and a UTS of type I along $\nu$ to get
$$
f=\overline z_1(t_2)^{d_1}\cdots\overline z_s(t_2)^{d_s}f_1.
$$
where $f_1\in k(c_0,\ldots,c_{t_2})[[\overline z_1(t_2),\ldots,\overline z_{\overline m}(t_2)]]$.
Set
$$
\lambda_i=c_{t_2}^{a_{1,s+1}(t_2)e_1(i)+\cdots+a_{s,s+1}(t_2)e_s(i)+a_{s+1,s+1}(t_2)f_i}
$$
for $1\le i\le d$. Then 
$$
f_1 = \sum_{i=1}^d\lambda_i\overline a_i(\overline z_{\overline m}(t_2)+1)^{f_i}
+\overline z_1(t_2)\cdots\overline z_s(t_2)\Lambda,
$$
for some series $\Lambda\in \overline T(t_2)$.
Set $r_1=\text{mult}(f_1(0,\ldots,0,\overline z_{\overline m}(t_2))<\infty$. $r_1\le f_d\le r$.

The residue of $\overline a_i$ in 
$$
\overline T(t_2)/(\overline z_1(t_2),\ldots,\overline z_{\overline m-1}(t_2),\overline z_{\overline m+1}(t_2),
\cdots,\overline z_n(t_2))\cong k(\overline T(t_2))[[\overline z_{\overline m}(t_2)]]
$$
is a nonzero constant $\tilde a_i\in k(\overline T(t_2))$ for $1\le i\le d$. Set
$$
\zeta(t) = f_1(0,\ldots,0,t-1) = \sum_{i=1}^d\lambda_i\tilde a_it^{f_i}.
$$
Suppose that $r_1=r$. Then $f_d=r$ and 
$$
\zeta(\overline z_{\overline m}(t_2)+1)=
\sum_{i=1}^d\lambda_i\tilde a_i(\overline z_{\overline m}(t_2)+1)^{f_i} 
= \lambda_d\tilde a_d\overline z_{\overline m}(t_2)^{f_d}\tag 3.30
$$
Thus $\zeta(t)=\lambda_d\tilde a_d(t-1)^r$ has a nonzero $t^{r-1}$ term, so that $f_{d-1}=r-1$ and $\tilde a_{d-1}\ne 0$.
Therefore $a_d=\overline a_d$ and  $\nu(\overline z_{\overline m})=\nu(a_{d-1})$. Define $\tau(0)$ by
$\nu(\overline z_{\overline m})=\nu(a_{d-1})=d_{\tau(0)}$ in (3.29). Then $a_{d-1}=h+\phi$ with 
$h\in p_{\overline m-1}\cap k(c_0)[[\overline z_1,\ldots,\overline z_{\overline m-1}]]$
and 
$$
\phi\in m(
k(c_0)[[\overline z_1,\ldots,\overline z_{\overline m-1}]])^{\gamma(\tau(0))}.
$$
Let $\alpha=\tilde a_d\in k(c_0)$ be the constant term of the power series 
$a_d\in k(c_0)[[\overline z_1,\ldots,\overline z_{\overline m-1}]]$.
Expanding out the LHS of $(3.30)$, we have
$$
\lambda_dr\tilde a_d+\lambda_{d-1}\tilde a_{d-1}=0.
$$
$$
\align
\frac{\overline z_{\overline m}}{\phi}&=\frac {\overline z_{\overline m}^r}{(a_{d-1}-h)\overline z_{\overline m}^{r-1}}\\
&\\
&=\frac{(\overline z_{\overline m}(t_2)+1)^r\lambda_d}
{\left(\lambda_{d-1}\overline a_{d-1}-\frac{hc_{t_2}^{a_{s+1,s+1}(t_2)(r-1)}}{\overline z_1(t_2)^{a_{s+1,1}(t_2)}\cdots
\overline z_s(t_2)^{a_{s+1,s}(t_2)}}\right)
(\overline z_{\overline m}(t_2)+1)^{r-1}}\\
&\\
&=\frac{(\overline z_{\overline m}(t_2)+1)\lambda_d}
{\lambda_{d-1}\overline a_{d-1}-\frac{hc_{t_2}^{a_{s+1,s+1}(t_2)(r-1)}}
{\overline z_1(t_2)^{a_{s+1,1}(t_2)}\cdots\overline z_s(t_2)^{a_{s+1,s}(t_2)}}}\\
&\\
&=\frac{ \frac{\overline z_{\overline m}(t_2)\lambda_d}{\lambda_{d-1}\overline a_{d-1}}
+\frac{\lambda_d}{\lambda_{d-1}\overline a_{d-1}}}
{1-\frac{hc_{t_2}^{a_{s+1,s+1}(t_2)(r-1)}}
{\lambda_{d-1}\overline a_{d-1}\overline z_1(t_2)^{a_{s+1,1}(t_2)}\cdots\overline z_s(t_2)^{a_{s+1,s}(t_2)}}}
\endalign
$$
has residue $-\frac{1}{r\alpha}$ in $\Cal O_{\nu}/m_{\nu}$. (Recall that $\nu(h)=\infty$). Thus
$\nu(\overline z_{\overline m}+\frac{1}{r\alpha}\phi)>\nu(\overline z_{\overline m})$. Since $\frac{1}{r\alpha}\in k(c_0)$,
$\frac{1}{r\alpha}\phi\in m(k(c_0)[[\overline z_1,\ldots,\overline z_{\overline m-1}]])^{\gamma(\tau(0))}$ and
there exists 
$$
A_1\in m(k(c_0)[\overline z_1,\ldots,\overline z_{\overline m-1}])^{\gamma(\tau(0))}
$$
such that $\nu(\overline z_{\overline m}-A_1)>\nu(\overline z_{\overline m})$.
Set $\overline z_{\overline m}^{(1)}=\overline z_{\overline m}-A_1$.

 Repeat the above
algorithm, with $\overline z_{\overline m}$ replaced by $\overline z_{\overline m}^{(1)}$. If we do not acheive a reduction $r_1<r$,
we can make an infinite sequence of change of variables 
$$
\overline z_{\overline m}^{(i)}=\overline z_{\overline m}^{(i-1)}-A_{i}
$$
such that $A_i\in k(c_0)[\overline z_1,\ldots,\overline z_{\overline m-1}]$,  
$\nu(A_i)=\nu(\overline z_{\overline m}^{(i-1)})$, $\nu(A_i)=d_{\tau(i-1)}$, 
$$
A_i\in m(k(c_0)[\overline z_1,\ldots,\overline z_{\overline m-1}])^{\gamma(\tau(i-1))}
$$
and
$$
\tau(0)<\tau(1)<\cdots<\tau(i)<\cdots
$$
Then
$$
\overline z_{\overline m}^{(i)}-\overline z_{\overline m}^{(j)}\in 
m(k(c_0)[\overline z_1,\ldots,\overline z_{\overline m-1}])^{\text{min}\{\gamma(\tau(i-1)),\gamma(\tau(j-1))\}}.
$$
Thus $\{\overline z_{\overline m}^{(i)}\}$ is a Cauchy sequence, and there exists 
a series 
$$
\overline A(\overline z_1,\ldots,\overline z_{\overline m-1})\in 
k(c_0)[[\overline z_1,\ldots,\overline z_{\overline m-1}]]
$$
 such that
$$
\overline z_{\overline m}^{\infty} = \text{Lim}_{i\rightarrow \infty}\overline z_{\overline m}^{(i)}
=\overline z_{\overline m}-\overline A
$$
and $\nu(\overline z_{\overline m}^{\infty}) = \infty$. Thus $\overline z_{\overline m}^{\infty}\in p_{\overline m}$.

Set $\lambda(\alpha) = \gamma(\tau(\alpha))$. For all $\alpha$  there are series
$a_i,a_j,a_k$ and exponents $f_i,f_j,f_k$ such that we can write
$$
f = [a_1(\overline z_{\overline m}^{(\alpha)})^{f_1}+\cdots+a_r(\overline z_{\overline m}^{(\alpha)})^r]
+\Sigma a_j(\overline z_{\overline m}^{(\alpha)})^{f_j}+\Sigma a_k(\overline z_{\overline m}^{(\alpha)})^{f_k}
+(\overline z_{\overline m}^{(\alpha)})^{r+1}\Omega
$$
where the terms in the first sum satisfy
$$
\nu(a_i(\overline z_{\overline m}^{(\alpha)})^{f_i})=r\nu(\overline z_{\overline m}^{(\alpha)})=rd_{\tau(\alpha)},
$$
$a_r$ is a unit, the terms in the second sum satisfy
$\nu(a_j(\overline z_{\overline m}^{(\alpha)})^{f_j})>rd_{\tau(\alpha)}$, and the terms in the third sum satisfy $\nu(a_k)=\infty$.
Set $m=m(k(c_0)[[\overline z_1,\ldots,\overline z_{\overline m-1}]])$. Since
$\nu(a_i)\ge d_{\tau(\alpha)}$ implies $a_i\in p_{\overline m-1}
+m^{\lambda(\alpha)}$,
we have 
$$
f\equiv a_r(\overline z_{\overline m}^{(\alpha)})^r\text{ mod }(m^{\lambda(\alpha)}+p_{\overline m-1}\overline T
+(\overline z_{\overline m}^{(\alpha)})^{r+1})
$$
so that
$$
f\in (\overline z_{\overline m}^{(\alpha)})^r+m^{\lambda(\alpha)}+p_{\overline m-1}\overline T
= (\overline z_{\overline m}^{(\alpha)})^r+m^{\lambda(\alpha)}+p_{\overline m-1}\overline T.
$$
Thus 
$$
f\in \cap_{\alpha=1}^{\infty}( (\overline z_{\overline m}^{(\alpha)})^r+m^{\lambda(\alpha)}+p_{\overline m-1}\overline T)
=  (\overline z_{\overline m}^{(\infty)})^r+p_{\overline m-1}\overline T.
$$
Since the $a_r$ are units, we have
$$
f = u(\overline z_m-\overline A(\overline z_1,\ldots,\overline z_{\overline m-1}))^r+h\tag 3.31
$$
where $u$ is a unit power series, $h\in p_{\overline m-1}\overline T$.

Suppose that we reach a reduction $r_1<r$ after a finite number of iterations.  We can repeat
the whole algorithm with $f$ replaced with $f_1$, $r$ replaced with $r_1$,
$c_0$ with a primitive element of $k(c_0,\ldots,c_{t_2})$ over $k$, $\overline T''$ with $\overline T''(t_2)$.
 (Recall that $k(c_0,\ldots,c_{t_2})\subset \overline T''(t_2)$)
 We have $\nu(f_1)=\infty$, so that the algorithm cannot
terminate with $r=0$, and we must produce $\overline z_{\overline m}^{\infty}(t)$ such that
$$
\overline z_{\overline m}^{\infty}(t)=\overline z_{\overline m}(t)
-\overline A(\overline z_1(t),\ldots,\overline z_{\overline m-1}(t)),
$$
with $\overline A\in k(c_0,\ldots,c_t)[[\overline z_1(t),\ldots,\overline z_{\overline m-1}(t)]]$
and  $\nu(\overline z_{\overline m}^{\infty}(t))=\infty$. 
In particular, the algorithm produces
 $\overline z_{\overline m}(t)-Q_{\overline m}(\overline z_1(t),\ldots,\overline z_{\overline m-1}(t))$
 of infinite value.

By 4) for $m=\overline m-1$, we can now construct a further UTS in the first $\overline m-1$ variables 
along $\nu$, so that 
$$
\align
&p_{\overline m - 1}(t)\\
&  = (\overline z_{r(1)}(t)-Q_{r(1)}(\overline z_1(t),\cdots,\overline z_{r(1)-1}(t)),\ldots,
\overline z_{r(\tilde {m-1})}(t)-Q_{r(\tilde {m-1})}(\overline z_1(t),\cdots,\overline z_{r(\tilde{ m-1})-1})(t)).
\endalign
$$
Now suppose $g\in p_{\overline m}(t)$ and $\nu(g)=\infty$. Then there exists
$g_0\in k(\overline T)[[\overline z_1(t),\ldots,\overline z_{\overline m-1}(t)]]$
and $g_1\in\overline T$
such that $g = g_0+(\overline z_{\overline m}(t)-Q_{\overline m})g_1$ and $\nu(g_0)=\infty$. Thus $g_0\in p_{\overline m-1}(t)$,
showing that 
$$
p_{\overline m}=
(\overline z_{r(1)}(t)-Q_{r(1)}(\overline z_1(t),\cdots,\overline z_{r(1)-1}),\ldots,
\overline z_{\overline m}-Q_{\overline m}(\overline z_1,\ldots,\overline z_{\overline m-1})).
$$

\subheading{Proof of $B(\overline m)$}
\subheading{Case 1) $\nu(h)<\infty$}
There exists a UTS $\overline T\rightarrow \overline T(t)$ satisfying 4) for $m=\overline m$. After replacing $\overline T''(0)$
 with $\overline T''(t)$ and replacing $c_0$ with a primitive element of $k(c_0,\ldots,c_t)$ over $k$,
we may assume that 
$$
p_{\overline m}=
(\overline z_{r(1)}-Q_{r(1)}(\overline z_1,\cdots,\overline z_{r(1)-1}),\ldots,
\overline z_{r(\tilde m)}-Q_{r(\tilde m)}(\overline z_1,\ldots,\overline z_{r(\tilde m)-1})).
$$
where the $Q_{r(i)}$ are power series with coefficients in $k(c_0)$. Let $L$ be a Galois closure of $k(\overline T)$
over $k(c_0)$ and
G be the Galois group of $L$ over $k(c_0)$. Set
$$
g=\prod_{\sigma\in G}\sigma(h)\in k(c_0)[[\overline z_1,\ldots,\overline z_m]].
$$
$g\in k(c_0)[[\overline z_1,\ldots,\overline z_m]]$.
Suppose $\nu(g)=\infty$. Then $g\in p_{\overline m}L[[\overline z_1,\ldots,\overline z_{\overline m}]]$
which is a prime ideal, invariant under $G$. Thus 
$$
h\in \left(p_{\overline m}L[[\overline z_1,\ldots,\overline z_{\overline m}]]\right)\cap 
k(\overline T)[[\overline z_1,\ldots,\overline z_{\overline m}]]=p_{\overline m}
$$
which implies $\nu(h)=\infty$, a contradiction. Thus $\nu(g)\ne\infty$. We will construct a UTS so that
$$
g=u\overline z_1(t)^{e_1}\cdots\overline z_s(t)^{e_s}
$$
where $u$ is a unit power series in $k(c_0,\ldots,c_t)[[\overline z_1(t),\ldots,\overline z_{\overline m}(t)]]$
and $h\in k(\overline T)[[\overline z_1(t),\ldots,\overline z_{\overline m}(t)]]$. Since $h\mid g$ in 
$k(\overline T)[[\overline z_1(t),\ldots,\overline z_{\overline m}(t)]]$, we will then have $h$ in the desired form in 
$\overline T(t)$.

First suppose that $\overline m=s$ (note that $p_s=0$). 
Set $\tau_i = \nu(\overline z_i)$ for $1\le i\le s$. $g$ has an expansion
$$
g=\sum_{i\ge 1} a_i\overline z_1^{b_1(i)}\cdots\overline z_s^{b_s(i)}
$$
where the $a_i\in k(c_0)$ and the terms have increasing value. 
Set 
$$
c= \frac{3}{\text{min}\{\frac{\tau_i}{\tau_1}\}}(b_1(1)+b_2(1)\frac{\tau_2}{\tau_1}+\cdots+b_s(1)\frac{\tau_s}{\tau_1}).
$$
We can perform a UTS of type I where $\overline T\rightarrow \overline T(1)$ is such that
$$
\mid \frac{\tau_i}{\tau_1}-\frac{a_{ij}(1)}{a_{1j}(1)}\mid < \frac{\tau_i}{2\tau_1}
$$
for $1\le i\le s$.
$$
\overline z_1^{b_1(i)}\cdots\overline z_s^{b_s(i)}
= \overline z_1(1)^{b_1(i)a_{11}(1)+\cdots+b_s(i)a_{s1}(1)}
\cdots
\overline z_s(1)^{b_1(i)a_{1s}(1)+\cdots+b_s(i)a_{ss}(1)}.
$$
Suppose that $i$ is such that $b_1(i)+\cdots +b_s(i)>c$. Then for all  $1\le j\le s$ we have
$$
\align
b_1(i)a_{1j}(1)+\cdots+b_s(i)a_{sj}(1) &=
a_{1j}(1)\left(b_1(i) + b_2(i)\frac{a_{2j}(1)}{a_{1j}(1)}+\cdots+b_s(i)\frac{a_{sj}(1)}{a_{1j}(1)}\right)\\
&\ge \frac{a_{1j}(1)}{2}\left(b_1(i) + b_2(i)\frac{\tau_2}{\tau_1}+\cdots+b_s(i)\frac{\tau_s}{\tau_1}\right)\\
&\ge a_{1j}(1)(b_1(i)+\cdots+b_s(i))\frac{\text{min}(\frac{\tau_i}{\tau_1})}{2}\\
&>\frac{3a_{1j}(1)}{2}\left(b_1(1) + b_2(1)\frac{\tau_2}{\tau_1}+\cdots+b_s(1)\frac{\tau_s}{\tau_1}\right)\\
&>a_{1j}(1)\left(b_1(1) + b_2(1)\frac{a_{2j}(1)}{a_{1j}(1)}+\cdots+b_s(1)\frac{a_{sj}(1)}{a_{1j}(1)}\right)\\
&=b_1(1)a_{1j}(1)+\cdots+b_s(1)a_{sj}(1).
\endalign
$$
By Lemma 3.2 we may choose the $a_{ij}(1)$ so that the inequality 
$$
b_1(i)a_{1j}(1)+\cdots+b_s(i)a_{sj}(1) > b_1(1)a_{1j}(1)+\cdots+b_s(1)a_{sj}(1)
$$
also holds for $1\le j\le s$ for the finitely many $i$ such that $b_1(i)+\cdots +b_s(i)\le c$. Then $g$ has the desired form  in 
$\overline T(1)$.

 Now assume that  $s< \overline m$.

Set $g=\overline z_1^{d_1}\cdots\overline z_s^{d_s}g_0$
 where $\overline z_i$ does not divide $g_0$
for $1\le i\le s$. Set 
$$
r = \text{mult}(g_0(0,\ldots,0,\overline z_{\overline m}).
$$
 $0\le r\le \infty$. 
We will also have an induction on 
$r$. If $r=0$ we are done, so suppose that $r>0$.
$$
\align
g_0
&= \sum_{i=1}^d \sigma_i(\overline z_1,\cdots,\overline z_{\overline m-1})\overline z_{\overline m}^{a_i}
+ \sum_j \sigma_j(\overline z_1,\cdots,\overline z_{\overline m-1})\overline z_{\overline m}^{a_j}\\
&+ \sum_k \sigma_k(\overline z_1,\cdots,\overline z_{\overline m-1})\overline z_{\overline m}^{a_k}
+\overline z_{\overline m}^a\Psi
\endalign
$$
where 
the coefficients of $\sigma_i$, $\sigma_j$, $\sigma_k$ and $\Psi$ are in $k(c_0)$,
$\Psi$ is a power series in $\overline z_1,\cdots,\overline z_{\overline m}$,
the first sum is over terms of minimum value $\rho$, $a$ satisfies $a\nu(\overline z_{\overline m})>\rho$, and the (finitely many)
remaining terms of finite value are in the second sum, the (finitely many) remaining terms of infinite value are in the third sum.

By 1'), 2') for $m<\overline m$ there is a UTS $\overline T\rightarrow \overline T(\alpha)$ in the first $\overline m-1$
variables along $\nu$ 
 such that
$$
\sigma_i=\overline z_1(\alpha)^{c_1(i)}\cdots \overline z_s(\alpha)^{c_s(i)}
\overline u_i
$$
 for all $i$,
$$
\sigma_j=\overline z_1(\alpha)^{c_1(j)}\cdots \overline z_s(\alpha)^{c_s(j)}
\overline u_j
$$
 for all $j$ and
$$
\sigma_k=\overline z_1(\alpha)^{c_1(k)}\cdots \overline z_s(\alpha)^{c_s(k)}
\overline u_k
$$
 for all $k$
where $\overline u_i,\overline u_j, \overline u_k\in 
k(c_0,\ldots,c_{\alpha})[[\overline z_1(\alpha),\ldots,\overline z_{\overline m-1}(\alpha)]]$, 
 $\overline u_i,\overline u_j$
are units and 
$$
\nu(\overline z_1(\alpha)^{c_1(k)}\cdots \overline z_s(\alpha)^{c_s(k)})>\rho.
$$
 Now perform a UTS 
 of type $\text{II}_{\overline m}$ $\overline T(\alpha)\rightarrow \overline T(\alpha+1)$ along $\nu$ to get 
$$
\align
g_0= & \overline z_1(\alpha+1)^{e_1}\cdots\overline z_s(\alpha+1)^{e_s}
(\sum_i\lambda_i\overline u_i(\overline z_{\overline m}(\alpha+1)+1)^{a_i})
\\
&+\sum_j\overline z_1(\alpha+1)^{e_1^j}\cdots\overline z_s(\alpha+1)^{e_s^j}\lambda_j\overline u_j
(\overline z_{\overline m}(\alpha+1)+1)^{a_j}\\
&+\sum_k\overline z_1(\alpha+1)^{e_1^k}\cdots\overline z_s(\alpha+1)^{e_s^k}\lambda_k\overline u_k
(\overline z_{\overline m}(\alpha+1)+1)^{a_k}\\
&+(\overline z_1(\alpha+1)^{a_{s+1,1}}\cdots\overline z_s(\alpha+1)^{a_{s+1,s}})^a
\Psi'.
\endalign
$$
where
$$
\lambda_i=c_{\alpha+1}^{a_{1,s+1}(\alpha+1)c_1(i)+\cdots+a_{s,s+1}(\alpha+1)c_s(i)+a_{s+1,s+1}(\alpha+1)a_i}.
$$
Then perform a UTS of type I $\overline T(\alpha+1)\rightarrow \overline T(\alpha+2)$ along $\nu$ to get 
$$
g_0 = \overline z_1(\alpha+2)^{d_1(\alpha+2)}\cdots\overline z_s(\alpha+2)^{d_s(\alpha+2)}g_1
$$
where 
$$
g_1=\sum_{i=1}^d\lambda_i\overline u_i(\overline z_{\overline m}
(\alpha+2)+1)^{a_i}+\overline z_1(\alpha+2)\cdots\overline z_s(\alpha+2)\Psi_{1}
$$
is a power series with coefficients in $k(c_0,\ldots,c_{\alpha+1})$,
$\Psi_{1}$ a power series in \linebreak
$\overline z_1(\alpha+2),\ldots,\overline z_{\overline m}(\alpha+2)$.
Set $r_1=\text{mult}(g_1(0,\ldots,0,\overline z_{\overline m}(\alpha+2))$. $r_1<\infty$ and $r_1\le r$.

Suppose that $r_1=r$. Then as in (3.30) in the proof of $A(\overline m)$, $\overline z_{\overline m}^r$
 is a minimal value term in $g_0$, so that $a_d=r$,  $a_{d-1}=r-1$, $\sigma_{d-1}\ne 0$, and
 $\nu(\sigma_{d-1})=\nu(\overline z_{\overline m})$. 

As in the proof of $A(\overline m)$, there exists $A_1\in k(c_0)[\overline z_1,\ldots,\overline z_{\overline m-1}]$ such that
we can
 make a change of variable, replacing $\overline z_{\overline m}$ with $\overline z'_{\overline m}
=\overline z_{\overline m}-A_1$ to get $\nu(\overline z_{\overline m}-A_1)>\nu(\overline z_{\overline m})$.
 We have
$$
\nu(\overline z_{\overline m})\le \nu(\overline z_{\overline m}^r)\le \nu(g_0)
$$
since $\overline z_{\overline m}^r$ is a minimal value term in $g_0$. Now repeat this procedure.
 If we do not achieve $r_1<r$ after a
finite number of iterations, we get an infinite sequence
$$
\nu(\overline z_{\overline m})<\nu(\overline z'_{\overline m})<\cdots<\nu(\overline z^{(i)}_{\overline m})<\cdots
$$
such that $\nu(\overline z^{(i)}_{\overline m})\le \nu(g_0)$ for all $i$. By Lemma 1.3, this is impossible.

Thus after replacing $\overline z_{\overline m}$ with 
$$
\tilde{\overline z}_{\overline m}' =\overline z_{\overline m} - P_{\overline m,0}(\overline z_1,\ldots,\overline z_{\overline m-1})
$$
 for some appropriate polynomial
$P_{\overline m,0}\in k(c_0)[\overline z_1,\ldots,\overline z_{\overline m-1}]$,
 we achieve a reduction $r_1<r$  in $\overline T''(\alpha+2)$. By induction on $r$,  we can  construct a UTS $\overline T
\rightarrow \overline T(t)$ along $\nu$ such that
$$
g=\overline z_1(t)^{d_1}\cdots\overline z_s(t)^{d_s}\overline u(\overline z_1(t),\ldots,\overline z_{\overline m}(t))
$$
where $\overline u$ is a unit power series with coefficients in $k(c_0,\ldots,c_t)$.

By 4) for $m=\overline m$ we can perform a further UTS to get
$$
p_{\overline m}(t) = (\overline z_{r(1)}(t)-Q_{r(1)}(\overline z_1(t),\cdots,\overline z_{r(1)-1}),\ldots,
\overline z_{r(\tilde m)}(t)-Q_{r(\tilde m)}(\overline z_1(t),\cdots,\overline z_{r(\tilde m)-1}))
$$
while preserving 
$$
g=\overline z_1(t)^{d_1}\cdots\overline z_s(t)^{d_s}\overline u(\overline z_1(t),\ldots,\overline z_{\overline m}(t))
$$
where $\overline u$ is a unit.

\subheading{Case 2) $\nu(h)=\infty$}

By 4) for $m=\overline m$, we can assume that 
$$
p_{\overline m}=
(\overline z_{r(1)}(t)-Q_{r(1)}(\overline z_1(t),\cdots,\overline z_{r(1)-1}),\ldots,
\overline z_{r(\tilde m)}-Q_{r(\tilde m)}(\overline z_1,\ldots,\overline z_{r(\tilde m)-1})).
$$
where the $Q_{r(i)}$ are series with coefficients in $k(c_0)$.
Then 
$$
h=\sum_{i=1}^{\tilde m}\sigma_i(\overline z_{r(i)}-Q_{r(i)})
$$
for some $\sigma_i\in k(\overline T)[[\overline z_1,\ldots,\overline z_m]]$.
Choose $b$ so that $b\nu(m(\overline T))>A$. There are polynomials $P_{r(i)}(\overline z_1,\ldots,z_{r(i)-1})$
in $k(c_0)[\overline z_1,\ldots,\overline z_{r(i)-1}]$ such that
$Q_{r(i)}-P_{r(i)}\in m(\overline T)^b$ and 
$$
\nu(\overline z_{r(i)}-P_{r(i)})>A.
$$
 Make a change of variables replacing 
$\overline z_{r(i)}$ with $\overline z_{r(i)}-P_{r(i)}$ for $1\le i\le \tilde m$. Then construct the UTS
$\overline T\rightarrow \overline T(t_1)$ which is a sequence of UTSs of type $\text{II}_r$ for $s+1\le r\le\overline m$,
followed by a MTS of type I to get
$$
g=\overline z_1(t_1)^{d_1}\cdots\overline z_s(t_1)^{d_s}\Sigma
$$
with $\nu(\overline z_1(t_1)^{d_1}\cdots\overline z_s(t_s)^{d_s})>A$. By 4) for $m=\overline m$, we can perform a UTS along $\nu$
in the first $\overline m$ variables  to get
$$
p_{\overline m}(t) = (\overline z_{r(1)}(t)-Q_{r(1)}(\overline z_1(t),\cdots,\overline z_{r(1)-1}(t)),\ldots,
\overline z_{r(\tilde m)}(t)-Q_{r(\tilde m)}(\overline z_1(t),\cdots,\overline z_{r(\tilde m)-1}(t)))
$$
while preserving 
$$
g=\overline z_1(t)^{d_1}\cdots\overline z_s(t)^{d_s}\Sigma
$$
with $\nu(\overline z_1(t)^{d_1}\cdots\overline z_s(t)^{d_s})>A$.

\subheading{Proof when 1) or 2) holds}
 \subheading{Case 1} Suppose that $s\le m\le l$. 
After performing a CUTS in the first $m$ variables, we may assume that 
$$
p_m = (\overline z_{r(1)}-Q_{r(1)}(\overline z_1,\cdots,\overline z_{r(1)-1}),\ldots,
\overline z_{r(\tilde m)}-Q_{r(\tilde m)}(\overline z_1,\cdots,\overline z_{r(\tilde m)-1})).
$$
where the coefficients of $Q_{r(i)}$ are in $k(c_0)$.
$f\in k(\overline U)[[\overline w_1,\ldots,\overline w_m]]$

First suppose that $\nu(f)<\infty$.
Let $\overline T_i=k(c_0)[[\overline z_1,\ldots,\overline z_i]]$,  
$\overline U_i=k(\overline U)[[\overline w_1,\ldots,\overline w_i]]$ for $s\le i\le m$.
Let $d=\text{det}(c_{ij})$, $(d_{ij})$ be the adjoint matrix of $(c_{ij})$. Then
$$
\align
\overline w_1 & = \overline z_1^{\frac{d_{11}}{d}}\cdots\overline z_s^{\frac{d_{1s}}{d}}\lambda_1\\
&\vdots\\
\overline w_s & = \overline z_1^{\frac{d_{s1}}{d}}\cdots\overline z_s^{\frac{d_{ss}}{d}}\lambda_s\\
\overline w_{s+1}&=\overline z_{s+1}\\
&\vdots\\
\overline w_m&=\overline z_m
\endalign
$$
where
$$
\lambda_i=\phi_1^{-\frac{d_{i1}}{d}}\cdots\phi_s^{-\frac{d_{is}}{d}}\text{ for }1\le i\le s
$$
Given a CUTS (3.21), set
$\sigma(i)$ to be the largest possible $\alpha$ such that after possibly permuting the parameters $\overline z_{s+1}(i),\ldots
,\overline z_{ m}(i)$, $\nu$ induces a rank 1 valuation on the quotient field of 
$k(\overline T(i))[[\overline z_1(i),\ldots,
 \overline z_{\alpha}(i)]]$. (Since $\nu(\overline z_1(i),\ldots,\nu(\overline z_s(i)))$ are rationally independent, $\sigma(i)\ge s$.)

If $\sigma(i)$ drops during the course of the proof,
 we can start the corresponding algorithm again with 
this smaller value of $\sigma(i)$. Eventually $\sigma(i)$ must stabilize, so we may assume that 
$\sigma(i)$ is constant throughout the
proof.

$\nu(\overline z_1^{d_{i1}}\cdots\overline z_s^{d_{is}})>0$ for $1\le i\le s$. By Lemma 3.6, there exists a CUTS
of type M1) $\overline T\rightarrow \overline T(1)$, $\overline U\rightarrow \overline U(1)$ such that
$\overline z_1^{d_{i1}}\cdots\overline z_s^{d_{is}}\in \overline T_m(1) 
= k(\overline T)[[\overline z_1(1),\ldots,\overline z_m(1)]]$
for all $i$. Since $p_m\overline T_m(1)$ is a prime and we may assume that $\sigma(1)=\sigma(0)$, we have
$p_m\overline T_m(1)=p_m(1)$. 

 Let $\omega$
be a primitive $d$th root of unity.
Let $L$ be a Galois closure of $k(\overline U)(\omega,\lambda_1,\ldots,\lambda_s)$ over $k(c_0)$ with Galois group $G$.
Set 
$\overline W= L[[\overline z_1(1)^{\frac{1}{d}},\ldots,\overline z_s(1)^{\frac{1}{d}},\overline z_{s+1}(1),\ldots,\overline z_m(1)]]$.
Given $i_1,\ldots,i_s\in \bold N$, Define a $k$-automorphism
 $\sigma_{i_1\cdots i_s}:\overline W\rightarrow \overline W$
by $\sigma_{i_1\cdots i_s}(\overline z_j(1)^{\frac{1}{d}}) = \omega^{i_j}\overline z_j(1)^{\frac{1}{d}}$ for
$1\le j\le s$. 

Our extension of
$\nu$ to the quotient field $F$ of $\overline U(1)$  extends to a valuation of the finite field extension
generated by $L$ and
$F(\overline z_1(1)^{\frac{1}{d}},\ldots,\overline z_s(1)^{\frac{1}{d}})$, which induces  valuations 
 on the quotient fields of $\overline T_m$, $\overline U_m$
and $\overline W$ which are compatible with the inclusions $\overline T_m\subset \overline U_m\subset \overline W$.
$\overline T_m(1)\rightarrow \overline W$ is finite, $p_m\overline W$ is prime implies
$$
p_m\overline W_m=\{h\in \overline W_m\,|\, \nu(h)=\infty\}.
$$
Thus 
$$
p_m\overline U_m=\{h\in \overline U_m\,|\, \nu(h)=\infty\}.
$$

Set $\overline g=\prod\sigma_{i_1\cdots i_s}(f), g=\prod_{\tau\in G}\tau(\overline g)\in
k(c_0)[[\overline z_1(1),\ldots,\overline z_m(1)]]\subset \overline T_m(1)$.
 Suppose $\nu(g)=\infty$.
Then $g\in p_m\overline T_m(1)$ implies $g\in p_m\overline W$ which implies $\tau\sigma_{i_1\cdots i_s}(f)\in p_m\overline W$
for some $\tau,\sigma_{i_1\cdots i_s}$ since $p_m\overline W$ is prime. But $\tau\sigma_{i_1\cdots i_s}(p_m\overline W)
=p_m\overline W$ implies $f\in p_m\overline W\cap \overline U_m= p_m\overline U_m$ so that $\nu(f) = \infty$. This is
a contradiction.
Thus $\nu(g)<\infty$.

By 1') (and Lemmas 3.3 and 3.4)
we can construct a CUTS $(R,\overline T''(t),\overline T(t))$ and $(S,\overline U''(t),\overline U(t))$
in the first $m$ variables to 
transform $g$ into the form 
$$
g =\overline z_1(t)^{\tilde d_1}\cdots\overline z_s(t)^{\tilde d_s}\overline u\tag 3.32
$$
in $\overline T(t)$ where $\overline u(\overline z_1(t),\ldots,\overline z_m(t))$ is a unit power series with coefficients in
$k(c_0,\ldots,c_t)$.
Then $f\mid g$ in $\overline U(t)$ implies
$$
f=\overline w_1(t)^{ d_1}\cdots\overline w_s(t)^{ d_s} u
$$
where $ u$ is a unit in $\overline U(t)$.
But $f$ is a series in $\overline w_1(t),\ldots,\overline w_m(t)$
with coefficients in $k(\overline U(t))$.  Thus
$$
f=\overline w_1(t)^{ d_1}\cdots\overline w_s(t)^{ d_s} u(\overline w_1(t),
\ldots,\overline w_m(t)).
$$
where the coefficients of $u$ are in $k(\overline U(t))$.
 
Now suppose that $\nu(f)=\infty$.  
$p_m\overline U_m$ is the set of elements of $\overline U_m$ of infinite value. 
Otherwise, as argued above, we can perform a UTS $\overline T\rightarrow \overline T(1)$ to get $\sigma(1)<\sigma(0)$. 
Thus it suffices by 4) to prove the
theorem when $f=\overline z_{r(i)}-Q_{r(i)}$ is a generator of $p_m$. This follows from 2').

\subheading{Case 2} Suppose that $m>l$.
 The proof is
by induction on $m-l$, assuming that it is true for smaller differences $m-l$.
 
First suppose that $\nu(f)<\infty$. 
Set 
$$
f=\overline w_1^{d_1}\cdots\overline w_s^{d_s}f_0
$$
 where $\overline w_i$ does not divide $f_0$
for $1\le i\le s$. Set $r = \text{mult}(f_0(0,\ldots,0,\overline w_m)$. $0\le r\le \infty$. 
We will also have an induction on 
$r$. If $r=0$ we are done, so suppose that $r>0$.
$$
\align
f_0
&= \sum_i \sigma_i(\overline w_1,\cdots,\overline w_{m-1})\overline w_m^{a_i}
+ \sum_j \sigma_j(\overline w_1,\cdots,\overline w_{m-1})\overline w_m^{a_j}\\
&+ \sum_k \sigma_k(\overline w_1,\cdots,\overline w_{m-1})\overline w_m^{a_k}
+\overline w_m^a\Psi
\endalign
$$
where $\sigma_i,\sigma_j,\sigma_k$ are power series with coefficients in $k(\overline U)$,
$\Psi$ is a power series in $\overline w_1,\cdots,\overline w_{m}$ with coefficients in $k(\overline U)$,
the first sum is over terms of minimum value $\rho$, $a$ satisfies $a\nu(\overline w_m)>\rho$,  the (finitely many)
 terms  in the second sum have finite value and the (finitely many)
 terms  in the second sum have infinite value.

By induction there is a  CUTS $(R,\overline T''(\alpha),\overline T(\alpha))$ and $(S,\overline U''(\alpha),\overline U(\alpha))$
in the first $m-1$ variables  such that
$$
\sigma_i=\overline w_1(\alpha)^{c_1(i)}\cdots \overline w_s(\alpha)^{c_s(i)}
\overline u_i
$$
 for all $i$,
$$
\sigma_j= \overline w_1(\alpha)^{c_1(j)}\cdots \overline w_s(\alpha)^{c_s(j)}
\overline u_j
$$
 for all $j$ and 
$$
\sigma_k= \overline w_1(\alpha)^{c_1(k)}\cdots \overline w_s(\alpha)^{c_s(k)}
\overline u_k
$$
 for all $k$
where $\overline u_i,\overline u_j, \overline u_k
\in k(\overline U(\alpha))[[\overline w_1(\alpha),\ldots,\overline w_{m-1}(\alpha)]]$,
 $\overline u_i,\overline u_j$
are units for all $i,j$ and 
$$
\nu(\overline w_1(\alpha)^{c_1(k)}\cdots \overline w_s(\alpha)^{c_s(k)})>\rho
$$
 for all $k$. Now perform a 
CUTS of type M3) where $\overline U(\alpha)\rightarrow \overline U(\alpha+1)$ is of type $\text{II}_m$ to get 
$$
\align
f_0= & \overline w_1(\alpha+1)^{e_1}\cdots\overline w_s(\alpha+1)^{e_s}
(\sum_i\lambda_i\overline u_i(\overline w_m(\alpha+1)+1)^{a_i})
\\
&+\sum_j\overline w_1(\alpha_1)^{e_1^j}\cdots\overline w_s(\alpha+1)^{e_s^j}\lambda_j
\overline u_j(\overline w_m(\alpha+1)+1)^{a_j}\\
&+\sum_k\overline w_1(\alpha+1)^{e_1^k}\cdots\overline w_s(\alpha+1)^{e_s^k}\lambda_k
\overline u_k(\overline w_m(\alpha+1)+1)^{a_k}\\
&+\left(\overline w_1(\alpha+1)^{b_{s+1,1}(\alpha+1)}\cdots
\overline w_s(\alpha+1)^{b_{s+1,s}(\alpha+1)}
\right)^a\Psi'
\endalign
$$
$$
\lambda_i=d_{\alpha+1}^{c_1(i)b_{1,s+1}(\alpha+1)+\cdots+c_s(i)b_{s,s+1}(\alpha+1)+b_{s+1,s+1}(\alpha+1)a_i}.
$$
Now perform a CUTS of type M1) 
$\overline T(\alpha+1)\rightarrow \overline T(\alpha+2)$, $\overline U(\alpha+1)\rightarrow \overline U(\alpha+2)$ to get 
$$
f_0 = \overline w_1(\alpha+2)^{d_1(\alpha+2)}\cdots\overline w_s(\alpha+2)^{d_s(\alpha+2)}f_1 
$$
where 
$$
f_1=\sum_{i=1}^d\lambda_i\overline u_i(\overline w_m(\alpha+2)+1)^{a_i}+\overline w_1(\alpha+2)\cdots\overline w_s(\alpha+2)\Psi_{1},
$$
$\Psi_{1}$ a power series in $\overline w_1(\alpha+2),\ldots,\overline w_m(\alpha+2)$
with coefficients in $k(\overline U(\alpha+2))$.
Set $r_1=\text{mult}(f_1(0,\ldots,0,\overline w_m(\alpha+2))$. $r_1<\infty$ and $r_1\le r$.

As in the proof of Case 1) of $B(\overline m)$, there is a polynomial $P_{m,0}
\in k(\overline U)[\overline w_1,\ldots,\overline w_{m-1}]$
 such that if we replace $\overline w_m$ with
$$
\tilde{\overline w}_m' =\overline w_m-P_{m,0}
$$ 
we get a reduction $r_1<r$ in $\overline U(\alpha+2)$. By induction, we can construct a CUTS as desired. 

Suppose that $\nu(f)=\infty$. Given a CUTS (3.21) and $i$ such that $s\le i\le n$,
set 
$$
a_i(t)=\{h\in k(\overline U(t))[[\overline w_1(t),\ldots,\overline w_i(t)]]|\nu(h)=\infty\}.
$$
$$
f=\sum_{i=0}^{\infty}\sigma_i(\overline w_1,\ldots,\overline w_{m-1})\overline w_m^i.
$$
If $\nu(\sigma_i)=\infty$ for all $i$, we can put $f$ in the desired form by induction on $m$ applied to a finite
set of generators of the ideal generated by the $\sigma_i$.

Suppose some $\nu(\sigma_i)<\infty$ for some $i$.
As in the proof of $A(\overline m)$, we can perform a UTS in the first $m$ variables to get
$$
f=\overline w_1(t_1)^{d_1}\cdots\overline w_s(t_1)^{d_s}f_1
$$
such that as in (3.31), there is a series $\overline A(\overline w_1(t_1),\ldots,\overline w_{m-1}(t_1))$ 
with coefficients in $k(\overline U(t_1))$ such that
$\nu(\overline w_m(t_1)-\overline A)=\infty$ and $f_1=u(\overline w_m(t_1)-\overline A)^r+h$
where $u\in k(\overline U(t_1))[[\overline w_1(t_1),\ldots,\overline w_m(t_1)]]$
 is a unit power series, $h\in a_{m-1}(t_1)$ and $r>0$. By induction on $m$, we are reduced to the case
$$
f=\overline w_m-\overline A(\overline w_1,\ldots,\overline w_{m-1}).
$$
We can then put $f$ in the desired form using the argument of Case 2) of the proof of $B(\overline m)$.

\subheading{Proof when 3) holds}
 Suppose that  $f$ is as in 3) of the statement of the
theorem. 

$$
f=\sum_{i=0}^{\infty}b_i(\overline w_1,\ldots,\overline w_{m-1})\overline w_m^i.
$$
Set $Q=\sum_{i=1}^{\infty}b_i\overline w_m^i$. After reindexing the $\overline w_i$,
 $l+1\le i\le m$, we may
assume that $Q\ne 0$.  $Q = \overline w_1^{n_1(0)}\cdots\overline w_s^{n_s(0)}Q_0$ where
$\overline w_i$ does not divide $Q_0$ for $1\le i\le s$. Set $r = \text{mult}(Q_0(0,\ldots,0,\overline w_m)$.
$1\le r\le \infty$.
The proof will be by induction on $r$. Suppose that $r>1$. 

$\nu(\frac{\partial f}{\partial \overline w_m})<\infty$ since
$\frac{\partial f}{\partial \overline w_m}\in\overline U''$
and $\nu$ restricts to a rank 1 valuation of the quotient field of $\overline U''$. 
Thus there must be some $i>0$ such that $\nu(b_i)<\infty$.

$$
\align
Q_0 &= \sum_{i=1}^d\sigma_i(\overline w_1,\ldots,\overline w_{m-1})\overline w_m^{a_i}
+\sum_j \sigma_j(\overline w_1,\ldots,\overline w_{m-1})\overline w_m^{a_j}\\
&+\sum_k \sigma_k(\overline w_1,\ldots,\overline w_{m-1})\overline w_m^{a_k}
+\overline w_m^a\Psi
\endalign
$$
where the first sum is over terms of minimum value $\rho$, $a$ satisfies $a\nu(\overline w_m)>\rho$,
the (finitely many)  terms  in the second sum have finite value and
the (finitely many)  terms  in the third sum have infinite value.

By 1) of the Theorem there is a CUTS $(R,\overline T''(\alpha),\overline T(\alpha))$ and 
$(S,\overline U''(\alpha),\overline U(\alpha))$ in the first $m-1$ variables such that
$$
\sigma_i = \overline w_1(\alpha)^{c_1^i(\alpha)}\cdots
\overline w_s(\alpha)^{c_s^i(\alpha)}\overline u_i
$$
for all $i$,
$$
\sigma_j = \overline w_1(\alpha)^{c_1^j(\alpha)}\cdots
\overline w_s(\alpha)^{c_s^j(\alpha)}\overline u_j
$$
for all $j$ and
$$
\sigma_k = \overline w_1(\alpha)^{c_1^i(\alpha)}\cdots
\overline w_s(\alpha)^{c_s^k(\alpha)}\overline u_k
$$
for all $k$

 where $\overline u_i,\overline u_j,\overline u_k\in 
k(\overline U(\alpha))[[\overline w_1(\alpha),\ldots,\overline w_{m-1}(\alpha)]]$,
$\overline u_i,\overline u_j$ are units and 
$$
\nu(\overline w_1(\alpha)^{c_1^i(\alpha)}\cdots
\overline w_s(\alpha)^{c_s^k(\alpha)})>\rho.
$$
 Perform a UTS $\overline U(\alpha)\rightarrow  \overline U(\alpha+1)$ of type $\text{II}_{m}$ to get 
$$
\align
Q_0 &=  \overline w_1(\alpha+1)^{c_1^i(\alpha+1)}\cdots\overline w_s(\alpha+1)^{c_s^i(\alpha+1)}
(\sum_i \lambda_i\overline u_i(\overline w_m(\alpha+1)+1)^{a_i})\\
&+\sum_j\overline w_1(\alpha+1)^{c_1^j(\alpha+1)}\cdots
\overline w_s(\alpha+1)^{c_s^j(\alpha+1)}\lambda_j\overline u_j(\overline w_m(\alpha+1)+1)^{a_j}\\
&+\sum_k\overline w_1(\alpha+1)^{c_1^k(\alpha+1)}\cdots
\overline w_s(\alpha+1)^{c_s^k(\alpha+1)} \lambda_k \overline u_k(\overline w_m(\alpha+1)+1)^{a_k}\\
&+\left(\overline w_1(\alpha+1)^{b_{s+1,1}(\alpha+1)}\cdots
\overline w_s(\alpha+1)^{b_{s+1,s}(\alpha+1)}\right)^a
\Psi'.
\endalign
$$
$$
\lambda_i= d_{\alpha+1}^{c_1^i(\alpha)b_{1,s+1}(\alpha+1)+\cdots+c^i_s(\alpha)b_{s,s+1}(\alpha+1)+b_{s+1,s+1}(\alpha+1)a_i}.
$$
Now perform a CUTS of type MI) $\overline T(\alpha+1)\rightarrow \overline T(\alpha+2)$,
$\overline U(\alpha+1)\rightarrow \overline U(\alpha+2)$ to get
$$
\align
Q_0 
=  \overline w_1(\alpha+2)^{c_1^i(\alpha+2)}\cdots\overline w_s(\alpha+2)^{c_s^i(\alpha+2)}
(&\sum_{i=1}^d \lambda_i\overline u_i(\overline w_m(\alpha+2)+1)^{a_i}\\
&+\overline w_1(\alpha+2)\cdots\overline w_s(\alpha+2)\Omega).
\endalign
$$ 
Set 
$$
Q_1 = \sum_{i=1}^d \lambda_i\overline u_i(\overline w_m(\alpha+2)+1)^{a_i}
+\overline w_1(\alpha+2)\cdots\overline w_s(\alpha+2)\Omega - \sum_{i=1}^d \lambda_i\overline u_i.
$$
Set $r_1 = \text{mult}(Q_1(0,\ldots,0,\overline w_m(\alpha+2))$. $0<r_1<\infty$ and $r_1\le r$.

Suppose $r_1=r$. Then as in (3.30) in the proof of $A(\overline m)$, $a_d=r$,
 $\sigma_{r-1}\ne 0$
and 
$\nu(\sigma_{r-1})=\nu(\overline w_m)$. 

As in the argument  of the proof of $A(\overline m)$, there is a polynomial 
$$
A_1\in k(\overline U)[\overline w_1,\ldots,\overline w_{m-1}]\subset \overline U''
$$
 such that
we can make a change of variables, replacing $\overline w_m$
 with $\overline w'_m=\overline w_m-A_1$,  to get
 $\nu(\overline w_m')>\nu(\overline w_m)$. We have
$$
\nu(\overline w_m)\le\nu(\overline w_m^{r-1})\le \nu(\frac{\partial Q_0}{\partial \overline w_m})
$$ 
since $\overline w_m^{r-1}$ is a minimum value term of $\frac{\partial Q_0}{\partial \overline w_m}$.
$\frac{\partial f}{\partial \overline w_m} \in \overline U''$ and
$$
\frac{\partial f}{\partial \overline w_m} 
= \overline w_1^{n_1(0)}\ldots\overline w_s^{n_s(0)}\frac{\partial Q_0}{\partial \overline w_m}
$$
implies $\nu(\frac{\partial Q_0}{\partial \overline w_m})<\infty$. Now repeat the above procedure. Since
$$
\frac{\partial Q_0}{\partial \overline w'_m} = \frac{\partial Q_0}{\partial \overline w_m}
$$
we will achieve a reduction in $r$ after a finite number of iterations by Lemma 1.3.

Thus after replacing $\overline w_m$ with 
$$
\tilde{\overline w}'_m = \overline w_m-P_{m,0}
$$
 for some 
$P_{m,0}\in k(\overline U)[\overline w_1,\ldots,\overline w_{m-1}]$, we achieve a reduction $r_1<r$ 
in $\overline U(\alpha+2)$. 

Thus we can construct a CUTS $(R,\overline T''(\beta),\overline T(\beta))$ and $(S,\overline U''(\beta),\overline U(\beta))$
 such that
$$
f = L(\overline w_1(\beta),\ldots,\overline w_{m-1}(\beta)) + \overline w_1(\beta)^{\alpha_1}
\cdots,\overline w_{s}(\beta)^{\alpha_s}\overline Q
$$
where $\text{mult}(\overline Q(0,\ldots,0,\overline w_m(\beta)) = 1$.
Set
$$
\tau = \nu(\overline w_1(\beta)^{\alpha_1}
\cdots,\overline w_{s}(\beta)^{\alpha_s}).
$$
Suppose that $L$ is not in $k(\overline U(\beta))[[\overline w_1(\beta),\ldots,\overline w_l(\beta)]]$.
Set 
$$
A = k(\overline U(\beta))[[\overline w_1(\beta),\ldots,\overline w_{m-1}(\beta)]].
$$
 We can write
$L=f_1+H$, with 
$$
f_1\in k(\overline U(\beta))[\overline w_1(\beta),\ldots,\overline w_{m-1}(\beta)]
\subset \overline U''(\beta),
$$
 $H\in m(A)^{\alpha}$ where
$\nu(m(A)^{\alpha})>\tau$.
After permuting the variables $\{\overline w_{l+1}(\beta),\ldots,\overline w_{m-1}(\beta)\}$ we may assume that
$\frac{\partial f_1}{\partial\overline w_{m-1}(\beta)}\ne 0$. Thus 
$\nu(\frac{\partial f_1}{\partial\overline w_{m-1}(\beta)})<\infty$
since $\frac{\partial f_1}{\partial\overline w_{m-1}(\beta)}\in \overline U''(\beta)$.
 By induction on $m$, we can perform a CUTS in the first $m-1$ variables
 to get 
$$
f_1=L'(\overline w_1(\gamma),\ldots,\overline w_l(\gamma))+\overline w_1(\gamma)^{\beta_1}\cdots\overline w_s(\gamma)^{\beta_s}
\overline Q_1
$$
so that
$$
\overline Q_1 = \overline u(\overline w_{m-1}(\gamma)
+\overline w_1(\gamma)^{\overline g_1}\cdots\overline w_s(\gamma)^{\overline g_s}\overline\Sigma)
$$
where $\overline u\in k(\overline U(\gamma))[[\overline w_1(\gamma),\ldots,\overline w_{m-1}(\gamma)]]$
is a unit,  $\Sigma\in k(\overline U(\gamma))[[\overline w_1(\gamma),\ldots,\overline w_{m-2}(\gamma)]]$,
 and
$$
\nu(\overline w_{m-1}(\gamma))
\le \nu(\overline w_1(\gamma)^{\overline g_1}\cdots\overline w_s(\gamma)^{\overline g_s}).
$$

Now perform a CUTS consisting of CUTSs of type M2), $s+1\le r\le m-1$ (with $P_{r,t}=0$ for $\gamma+1\le t\le \delta-1$)
  and a CUTS of type M1)
to get 
$$
H = \overline w_1(\delta)^{d_1(\delta)}\cdots\overline w_s(\delta)^{d_s(\delta)}\Psi
$$
with
$$
\nu(\overline w_1(\delta)^{d_1(\delta)}\cdots\overline w_s(\delta)^{d_s(\delta)})>
\tau.
$$
$$
\overline Q_1 = \overline w_1(\delta)^{e_1(\delta)}\cdots\overline w_s(\delta)^{e_s(\delta)}
\overline u'(\overline w_{m-1}(\delta) + \phi')
$$
for some  $\phi'\in k(\overline U(\delta))[[\overline w_1(\delta),\ldots,\overline w_{m-2}(\delta)]]$,
and unit $\overline u'\in k(\overline U(\delta))[[\overline w_1(\delta),\ldots,\overline w_{m-1}(\delta)]]$.
After possibly interchanging $\overline w_{m-1}(\delta)$ and $\overline w_{m}(\delta)$ 
and performing a CUTS of type M1), we have $f$ in the  form
$$
f=L(\overline w_1(\delta),\ldots,\overline w_l(\delta))+\overline w_1(\delta)^{\alpha_1}\cdots
\overline w_s(\delta)^{\alpha_s}Q
$$
where $\text{mult } Q(0,\ldots,0,\overline w_m(\delta))=1$. Thus $Q= u(\overline w_m(\delta)+\Omega)$ where
$$
u\in k(\overline U(\delta))[[\overline w_1(\delta),\ldots,\overline w_m(\delta)]]
$$
 is a unit and
$\Omega\in k(\overline U(\delta))[[\overline w_1(\delta),\ldots,\overline w_{m-1}(\delta)]]$.
After replacing $\overline w_m(\delta)$ with $\overline w_m(\delta)+\Psi$, where $\Psi\in 
 k(\overline U(\delta))[\overline w_1(\delta),\ldots,\overline w_{m-1}(\delta)]\subset \overline U''(\delta)$,
 we can assume
that 
$$
\Omega\in (\overline w_1(\delta),\ldots,\overline w_{m-1}(\delta))^B
$$
 where $B$ is arbitrarilly large.

If $\nu(Q)<\infty$, we can choose $B$ so large that $\nu(Q)=\nu(\overline w_m(\delta))<\nu(\Omega)$. Then by
the conclusions of 1) and 2) of the Theorem, we can perform a CUTS in the first $m-1$ variables
to get 
$$
\Omega = \overline w_1(\epsilon)^{g_1}\cdots\overline w_s(\epsilon)^{g_s}\Sigma
$$ 
with $\nu(\overline w_1(\epsilon)^{g_1}\cdots\overline w_s(\epsilon)^{g_s})>\nu(\overline w_m(\epsilon))$.

If $\nu(Q)=\infty$, we must have $\nu(\Omega)=\nu(\overline w_m(\delta))<\infty$. Then by 1) of the
Theorem, we can perform a CUTS in the first $m-1$ variables to get
$$
\Omega = \overline w_1(\epsilon)^{g_1}\cdots\overline w_s(\epsilon)^{g_s}\Sigma
$$ 
with $\nu(\overline w_1(\epsilon)^{g_1}\cdots\overline w_s(\epsilon)^{g_s})=\nu(\overline w_m(\epsilon))$.
\enddemo

\proclaim{Theorem 3.9}
Suppose that $ T''(0)\subset \hat R$ is a regular local ring essentially of finite type over $R$
such that the quotient field of $ T''(0)$ is finite over $J$, $U''(0)\subset \hat S$ is a
regular local ring essentially of finite type over $S$
such that the quotient field of $ U''(0)$ is finite over $K$,  $ T''(0)\subset  U''(0)$,
$ T''(0)$ contains a subfield isomorphic to $k(c_0)$ for some $c_0\in k(T''(0))$ and $ U''(0)$ contains a subfield 
isomorphic to $ k(U''(0))$.
Suppose that  $R$ has regular parameters $(x_1,\ldots, x_n)$, $S$ has regular parameters $(y_1,\ldots,y_n)$,
 $T''(0)$ has
regular parameters $(\tilde{\overline x}_1,\ldots, \tilde{\overline x}_n)$ and $\ U''(0)$ has regular parameters 
$(\tilde{\overline y}_1,\ldots, \tilde{\overline y}_n)$
such that 
$$
\align
\tilde{\overline x}_1 &= \tilde{\overline y}_1^{c_{11}}\cdots \tilde{\overline y}_s^{c_{1s}}\phi_1\\
\vdots&\\
\tilde{\overline x}_s &= \tilde{\overline y}_1^{c_{s1}}\cdots \tilde{\overline y}_s^{c_{ss}}\phi_s\\
\tilde{\overline x}_{s+1} &= \tilde{\overline y}_{s+1}\\
\vdots&\\
\tilde{\overline x}_l&=\tilde{\overline y}_l
\endalign
$$
where $\phi_1,\ldots,\phi_s\in k(U''(0))$, $\nu(\tilde{\overline x}_1),\ldots,\nu(\tilde{\overline x}_s)$
are rationally independent, $\text{det}(c_{ij})\ne 0$. Suppose that there exists a regular local ring $\tilde R\subset R$
such that $(x_1,\ldots, x_l)$ are regular parameters in $\tilde R$ and $k(\tilde R)\cong k(c_0)$.
For $1\le i\le l$, there exists   $\gamma_i\in k(c_0)[[x_1,\ldots, x_l]]\cap T''(0)$
such that $\gamma_i\equiv 1\text{ mod }(x_1,\ldots,x_l)$ and
$$
x_i=\cases
\gamma_i\tilde{\overline x}_i& 1\le i\le l\\
\tilde{\overline x}_i& l<i\le n.
\endcases
$$
In particular 
$k(c_0)[[x_1,\ldots, x_l]] = k(c_0)[[\tilde{\overline x}_1,\ldots, \tilde{\overline x}_l]]$.
 There exists $\gamma_i^y\in U''(0)$ such that 
$y_i=\gamma_i^y\tilde{\overline y}_i$, $\gamma_i^y\equiv 1\text{ mod } m(U''(0))$ for $1\le i\le n$. 

Suppose that one of the following three conditions holds.

1) $f\in k(U''(0))[[\tilde{\overline y}_1,\ldots, \tilde{\overline y}_m]]$ for some $m$ with $l\le m\le n$ and $\nu(f)<\infty$. 

2) $f\in k(U''(0))[[\tilde{\overline y}_1,\ldots, \tilde{\overline y}_m]]$ for some $m$ with $l< m\le n$, $\nu(f)=\infty$,
and $A\in \bold N$ is given. 

3) $f \in  U''(0) -  k(U''(0))[[\tilde{\overline y}_1,\ldots,\tilde{\overline y}_l]]$.

Then there exists a positive integer $N_0$ such that for $N\ge N_0$, we can construct  a CRUTS along $\nu$ 
$(R, T''(t), T(t))$ and $(S,U''(t),U(t))$ with associated MTSs
$$
\matrix
S & \rightarrow & S(t)\\
\uparrow &&\uparrow\\
R &\rightarrow & R(t)
\endmatrix
$$
such that the following holds.
$T''(t)$ contains a subfield $k(c_0,\ldots,c_t)$, $U''(t)$ contains a subfield isomorphic to $k(U(t))$,
$R(t)$ has regular parameters $(x_1(t),\ldots,x_n(t))$,
$T''(t)$ has regular parameters $(\tilde{\overline x}_1(t),\ldots, \tilde{\overline x}_n(t))$, $S(t)$ has regular parameters
$(y_1(1),\ldots, y_n(t))$,
 $U''(t)$ has regular parameters $(\tilde{\overline y}_1(t),\ldots, \tilde{\overline y}_n(t))$
such that
$$
x_i(t)=\cases \gamma_i(t)\tilde{\overline x}_i(t)& 1\le i\le l\\
\tilde{\overline x}_i(t)& l<i\le n
\endcases
$$ 
where $\gamma_i(t)\in k(c_0,\ldots,c_t)[[x_1(t),\ldots, x_l(t)]]$  units such that 
$\gamma_i(t)\equiv 1\text{ mod } (x_1(t),\ldots, x_l(t))$.
In particular,
$$
k(c_0,\ldots, c_t)[[x_1(t),\ldots, x_l(t)]]=
k(c_0,\ldots,c_t)[[\tilde{\overline x}_1(t),\ldots, \tilde{\overline x}_l(t)]].
$$
For $1\le i\le n$ there exists $\gamma_i^y(t)\in U''(t)$ such that $y_i(t)=\gamma_i^y(t)\tilde{\overline y}_i(t)$,
$\gamma_i^y(t)\equiv 1\text{ mod }m(U''(t))$.
$$
\align
\tilde{\overline x}_1(t) &= \tilde{\overline y}_1(t)^{c_{11}(t)}\cdots \tilde{\overline y}_s(t)^{c_{1s}(t)}\phi_1(t)\tag 3.33\\
\vdots&\\
\tilde{\overline x}_s(t) &= \tilde{\overline y}_1(t)^{c_{s1}(t)}\cdots \tilde{\overline y}_s(t)^{c_{ss}(t)}\phi_s(t)\\
\tilde{\overline x}_{s+1}(t) &= \tilde{\overline y}_{s+1}(t)\\
\vdots&\\
\tilde{\overline x}_l(t)&=\tilde{\overline y}_l(t)
\endalign
$$
$\phi_1(t),\ldots,\phi_s(t)\in  k(U(t))$, $\nu(\tilde{\overline x}_1(t)),\ldots,\nu(\tilde{\overline x}_s(t))$
are rationally independent, $\text{det}(c_{ij}(t))\ne 0$ and there exists a regular local ring $\tilde R(t)\subset R(t)$
such that $(x_1(t),\ldots, x_l(t))$ are regular parameters in $\tilde R(t)$ and $k(\tilde R(t))\cong k(c_0,\ldots,c_t)$.
Furthermore,
$x_i(t)=x_i$ for $l+1\le i\le n$, $y_i(t)=y_i$ for $m+1\le i\le n$, so that the CRUTS is in the first $m$ variables where $m=n$
in case 3).
Set $n_{t,l} = m\left(k(U(t))
[[\tilde{\overline y}_1(t),\ldots,\tilde{\overline y}_l(t)]]\right)$.

In case 1) we have
$$
f \equiv \tilde{\overline y}_1(t)^{d_1}\cdots\tilde{\overline y}_s(t)^{d_s}u(\tilde{\overline y}_1(t),\ldots,\tilde{\overline y}_m(t))
\text{ mod }m(U(t))^N\tag 3.34
$$
where $u$ is a unit power series. Further if $f\in k(\overline U)[[\tilde{\overline y}_1,\ldots,\tilde{\overline y}_l]]$,
$$
f \equiv \tilde{\overline y}_1(t)^{d_1}\cdots\tilde{\overline y}_s(t)^{d_s}u(\tilde{\overline y}_1(t),\ldots,\tilde{\overline y}_l(t))
\text{ mod }n_{t,l}^N.
$$
In case 2) we have
$$
f \equiv \tilde{\overline y}_1(t)^{d_1}\cdots\tilde{\overline y}_s(t)^{d_s}
\Sigma(\tilde{\overline y}_1(t),\ldots,\tilde{\overline y}_m(t))
\text{ mod }m(U(t))^N\tag 3.35
$$
with $\nu(\tilde{\overline y}_1(t)^{d_1}\cdots\tilde{\overline y}_s(t)^{d_s})>A$.
 Further if $f\in k(\overline U)[[\tilde{\overline y}_1,\ldots,\tilde{\overline y}_l]]$,
$$
f \equiv \tilde{\overline y}_1(t)^{d_1}\cdots\tilde{\overline y}_s(t)^{d_s}u(\tilde{\overline y}_1(t),\ldots,\tilde{\overline y}_l(t))
\text{ mod }n_{t,l}^N.
$$
In case 3) we have
$$
f \equiv P(\tilde{\overline y}_1(t),\ldots,\tilde{\overline y}_l(t))+
   \tilde{\overline y}_1(t)^{d_1}\cdots\tilde{\overline y}_s(t)^{d_s}H
\text{ mod }m(U(t))^N\tag 3.36
$$
where $P$ is a series with coefficients in $k(U(t))$ and 
$$
H= u(\tilde{\overline y}_{l+1}(t)+ \tilde{\overline y}_1(t)^{g_1}\cdots\tilde{\overline y}_s(t)^{g_s}\Sigma)
$$
where $u\in U(t)$ is a unit, $\Sigma \in k(U(t))[[\tilde{\overline y}_1(t),\ldots,\tilde{\overline y}_l(t),
\tilde{\overline y}_{1+2}(t),\ldots,\tilde{\overline y}_n(t)]]$ and
$\nu(\tilde{\overline y}_{l+1}(t))\le \nu(\tilde{\overline y}_1(t)^{g_1}\cdots\tilde{\overline y}_s(t)^{g_s})$.
\endproclaim

\demo{Proof} 
Set $\overline T = \hat R$, $\overline U = \hat S$, $\overline T'' = T''(0)$, $\overline U'' = U''(0)$.
 Set $\overline z_i = \tilde{ \overline x}_i$,
 $\overline w_i =\tilde{ \overline y}_i$
for $1\le i\le n$. In case 3) set $m=n$. By Theorem 3.8 there is a CUTS along $\nu$ $(R,\overline T''(t),\overline T(t))$ and
$(S,\overline U''(t),\overline U(t))$
$$
\matrix
\overline U(0) & \rightarrow & \overline U(t)\\
\uparrow &&\uparrow\\
\overline T(0) &\rightarrow & \overline T(t)
\endmatrix
$$
so that in the notation  of Theorem 3.8 and its proof, for $0\le \alpha\le t$,
$\overline T''(\alpha)$ has regular parameters
$$
(\overline z_1(\alpha),\ldots, \overline z_n(\alpha))\text{ and } 
(\tilde{\overline z}_1'(\alpha),\ldots, \tilde{\overline z}_n'(\alpha)),
$$ 
$\overline U''(\alpha)$ has parameters
$$
(\overline w_1(\alpha),\ldots, \overline w_n(\alpha))\text{ and }
(\tilde{\overline w}_1'(\alpha),\ldots, \tilde{\overline w}_n'(\alpha))
$$
such that 
in case 1) we have
$$
f = \overline w_1(t)^{d_1}\cdots\overline w_s(t)^{d_s}u(\overline w_1(t),\ldots,\overline w_m(t))
$$
where $u$ is a unit power series.
In case 2) we have
$$
f = \overline w_1(t)^{d_1}\cdots\overline w_s(t)^{d_s}\Sigma(\overline w_1(t),\ldots,\overline w_m(t))
$$
where $\nu(\overline w_1(t)^{d_1}\cdots\overline w_s(t)^{d_s})>A$.
In case 3) we have
$$
f = P(\overline w_1(t),\ldots,\overline w_l(t))+
   \overline w_1(t)^{d_1}\cdots\overline w_s(t)^{d_s}H.
$$
for some powerseries $P\in k(\overline U(t))[[\overline w_1(t),\ldots,\overline w_l(t)]]$,
$$
H=u(\overline w_n(t)+\overline w_1(t)^{g_1}\cdots\overline w_s(t)^{g_s}\Sigma)
$$
 where $u\in\overline U(t)$ is a unit,
$\Sigma\in k(\overline U(t))[[\overline w_1(t),\ldots,\overline w_{n-1}(t)]]$ and
$$
\nu(\overline w_n(t))\le \nu(\overline w_1(t)^{g_1}\cdots\overline w_s(t)^{g_s}).
$$

\subheading{Step 1} Fix $N>0$. To begin with, we will construct commutative diagrams of inclusions of
regular local rings
$$
\matrix
U'(\alpha)&\rightarrow&U''(\alpha)&\rightarrow &U(\alpha)\\
\uparrow&             &\uparrow    &           &\uparrow\\ 
T'(\alpha)&\rightarrow&T''(\alpha)&\rightarrow &T(\alpha)
\endmatrix\tag 3.37
$$
for $1\le \alpha\le t$ such that
$T(\alpha) = T'(\alpha)\,\,\hat{}\,$, $U(\alpha)= U'(\alpha)\,\,\hat{}\,$ for all $\alpha$,
$T''(\alpha)$ has regular parameters 
$$
(\overline x_1(\alpha),\ldots, \overline x_n(\alpha)), (\tilde{\overline  x}_1(\alpha),\ldots, \tilde{\overline x}_n(\alpha)),
(\tilde{\overline  x}'_1(\alpha),\ldots, \tilde{\overline x}'_n(\alpha)).
$$
$U''(\alpha)$ has regular parameters
$$
(\overline y_1(\alpha),\ldots, \overline y_n(\alpha)),
(\tilde{\overline  y}_1(\alpha),\ldots, \tilde{\overline y}_n(\alpha)),
(\tilde{\overline  y}'_1(\alpha),\ldots, \tilde{\overline y}'_n(\alpha)).
$$
where $\overline x_i(0)=\tilde{\overline x}_i$ and  
$\overline y_i(0)=\tilde{\overline y}_i$
for $1\le i\le n$. We will have  isomorphisms
$$
\align
&\eta_T^{\alpha}:k(T(\alpha))\rightarrow k(\overline T(\alpha))\text{ and }\tag 3.38\\
&\eta_U^{\alpha}:k(U(\alpha))\rightarrow k(\overline U(\alpha))
\endalign
$$
such that the diagrams

$$
\matrix
k(\overline T(\alpha))&\rightarrow         &k(\overline T(\alpha+1))\\
\uparrow&&\uparrow\\
k(T(\alpha))&\rightarrow         &k(T(\alpha+1))
\endmatrix\tag 3.39
$$
 and 
$$
\matrix
k(\overline U(\alpha)) & \rightarrow & k(\overline U(\alpha+1))\\
\uparrow&&\uparrow\\
k(U(\alpha))&\rightarrow         &k(U(\alpha+1))\\
\endmatrix
$$
commute for $0\le \alpha\le t-1$. 
For all $\alpha$ we will have
$$
\align
\overline x_1(\alpha) &= \overline y_1(\alpha)^{c_{11}(\alpha)}\cdots \overline y_s(\alpha)^{c_{1s}(\alpha)}\phi_1(\alpha)
\tag 3.40\\
\vdots&\\
\overline x_s(\alpha) &= \overline y_1(\alpha)^{c_{s1}(\alpha)}\cdots \overline y_s(\alpha)^{c_{ss}(\alpha)}\phi_s(\alpha)\\
\overline x_{s+1}(\alpha) &= \overline y_{s+1}(\alpha)\\
\vdots&\\
\overline x_l(\alpha)&=\overline y_l(\alpha)
\endalign
$$
with $\phi_1(\alpha),\ldots,\phi_s(\alpha)\in k(U(\alpha))$ the coefficients of (3.22) of Theorem 3.8, $(c_{ij}(\alpha))$
the exponents of (3.22) of Theorem 3.8. 

We will construct (3.37) inductively. Suppose that (3.37) has been constructed out to $T(\alpha)\rightarrow U(\alpha)$ and
regular parameters $(\overline x_1(\alpha),\ldots,\overline x_n(\alpha))$ in $T''(\alpha)$ and
$(\overline y_1(\alpha),\ldots,\overline y_n(\alpha))$ in $U''(\alpha)$ have been defined so that (3.38) and (3.40) hold.

If we identity $k(T(\alpha))$ with $k(\overline T(\alpha))$ and $k(U(\alpha))$ with $k(\overline U(\alpha))$ we
have isomorphisms $T(\alpha)\cong k(\overline T(\alpha))[[\overline x_1(\alpha),\ldots,\overline x_n(\alpha)]]$
and $U(\alpha)\cong k(\overline U(\alpha))[[\overline y_1(\alpha),\ldots,\overline y_n(\alpha)]]$.

We can choose $\Lambda_{\alpha}$ and $\Omega_{i,\alpha}$
 arbitrarily subject to the following conditions, to define regular parameters in $T''(\alpha)$ by
$$
\tilde{\overline x}_i(\alpha) = 
\cases
\overline x_r(\alpha)+\Lambda_{\alpha}(\overline x_1(\alpha),\ldots,
\overline x_{l}(\alpha))  &\text{ if }\overline T(\alpha-1)\rightarrow \overline T(\alpha)
\text{ is of type II}_r\text{ and }i=r\\  
\overline x_i(\alpha) & \text{ otherwise}
\endcases\tag 3.41
$$

with $\Lambda_{\alpha}\in k(c_0,\ldots, c_{\alpha})[[\overline x_1(\alpha),\ldots,\overline x_l(\alpha)]]\cap T''(\alpha)$
and $\text{mult}(\Lambda_{\alpha})\ge N$. 
We will take $\Lambda_0 = 0$.

Recall that the $P_{i,\alpha}$ constructed in Theorem 3.8 are polynomials with coefficients in $k(c_0,\ldots,c_{\alpha})$
if $i\le l$.
Define
$$
\tilde{\overline x}'_i(\alpha) = 
\cases
\tilde{\overline x}_i(\alpha)-P_{i,\alpha}(\tilde{\overline x}_1(\alpha),\ldots,
\tilde{\overline x}_{i-1}(\alpha))+\Omega_{i,\alpha}(\tilde{\overline x}_1(\alpha),\ldots,
\tilde{\overline x}_{l}(\alpha)),  &\text{ if }s+1\le i\le l\\  
\tilde{\overline x}_i(\alpha) & \text{ otherwise}
\endcases\tag 3.42
$$
with $\Omega_{i,\alpha}\in k(c_0,\ldots,c_{\alpha})[[\overline x_1(\alpha),\ldots,\overline x_l(\alpha)]]\cap T''(\alpha)$
and $\text{mult}(\Omega_{i,\alpha})\ge N$. 
(If $P_{i,\alpha}=0$, or if $1\le i\le s$ we will have 
$\overline x_i'(\alpha)=\tilde{\overline x}_i(\alpha)$.)
We then have 
$$
\align
&k(c_0,\ldots,c_{\alpha})[[\overline x_1(\alpha),\ldots, \overline x_l(\alpha)]]
=k(c_0,\ldots,c_{\alpha})[[\tilde{\overline  x}_1(\alpha),\ldots, \tilde{\overline x}_l(\alpha)]]\\
&=k(c_0,\ldots,c_{\alpha})[[\tilde{\overline  x}'_1(\alpha),\ldots, \tilde{\overline x}'_l(\alpha)]]
\endalign
$$
and 
$$
\overline x_i(\alpha) 
=\tilde{\overline x}_i(\alpha) = \tilde{\overline x}_i'(\alpha)=x_i
$$
for $l<i\le n$. Define
$$
\tilde{\overline y}_i(\alpha) = 
\cases
\tilde{\overline x}_i(\alpha)&\text{ if } s+1\le i\le l\\
\overline y_i(\alpha)+\Lambda_{\alpha}(\overline y_1(\alpha),\ldots,
\overline y_{n}(\alpha)),  &\text{ if }S(\alpha-1)\rightarrow S(\alpha)\text{ is of type  II}_r\\
&\text{ \,\,and }i=r\ge l+1\\  
\overline y_i(\alpha) & \text{ otherwise}
\endcases\tag 3.43
$$
with $\Lambda_{\alpha}\in U''(\alpha)$
and $\text{mult}(\Lambda_{\alpha})\ge N$. We will take $\Lambda_0=0$.

 Recall that the $P_{i,\alpha}$ constructed in Theorem 3.8 are polynomials with coefficents in $k(\overline U(\alpha))$
for $l+1\le i$.
Define
$$
\tilde{\overline y}'_i(\alpha) = 
\cases
\tilde{\overline x}'_i(\alpha)&\text{ if } s+1\le i\le l\\
\tilde{\overline y}_i(\alpha)-P_{i,\alpha}(\tilde{\overline y}_1(\alpha),\ldots,
\tilde{\overline y}_{i-1}(\alpha))+\Omega_{i,\alpha}(\tilde{\overline y}_1(\alpha),\ldots,
\tilde{\overline y}_{n}(\alpha))  &\text{ if }l+1\le i\le m\\  
\tilde{\overline y}_i(\alpha) & \text{ otherwise}
\endcases\tag 3.44
$$
with $\Omega_{i,\alpha}\in U''(\alpha)$
and $\text{mult}(\Omega_{i,\alpha})\ge N$.

These variables are such that for all $\alpha$,
$$
k(U(\alpha))[[\overline y_1(\alpha),\ldots, \overline y_l(\alpha)]]
=k(U(\alpha))[[\tilde{\overline  y}_1(\alpha),\ldots, \tilde{\overline y}_l(\alpha)]]
=k(U(\alpha))[[\tilde{\overline  y}'_1(\alpha),\ldots, \tilde{\overline y}'_l(\alpha)]]
$$
and  $\overline y_i(\alpha)=\tilde{\overline y}_i(\alpha)=\tilde{\overline y}'_i(\alpha)=y_i$
for $m<i\le n$.

If $\overline T(\alpha)\rightarrow \overline T(\alpha+1)$ is of type $I$, defined by (3.24) of Theorem 3.8,
 $T(\alpha)\rightarrow T(\alpha+1)$ will be the UTS  of type I such that $T''(\alpha+1)$ has regular parameters
$(\overline x_1(\alpha+1),\ldots, \overline x_n(\alpha+1))$ defined by
$$
\align
\tilde{\overline x}'_1(\alpha) &= \overline x_1(\alpha+1)^{a_{11}(\alpha+1)}\cdots \overline x_s(\alpha+1)^{a_{1s}(\alpha+1)}
\tag 3.45\\
\vdots&\\
\tilde{\overline  x}'_s(\alpha) &= \overline x_1(\alpha+1)^{a_{s1}(\alpha+1)}\cdots \overline x_s(\alpha+1)^{a_{ss}(\alpha+1)}.
\endalign
$$

Suppose that $\overline T(\alpha)\rightarrow \overline T(\alpha+1)$ is of type $\text{II}_r$, defined by (3.25)
of Theorem 3.8. Set
$(e_{ij})=(a_{ij}(\alpha+1))^{-1}$,
$$
\align
M_1&=\tilde{\overline x}_1'(\alpha)^{e_{11}}\cdots\tilde{\overline x}_s'(\alpha)^{e_{1s}}
\tilde{\overline x}_r'(\alpha)^{e_{1,s+1}}\tag 3.46\\
&\vdots\\
M_s& =\tilde{\overline x}_1'(\alpha)^{e_{s1}}\cdots\tilde{\overline x}_s'(\alpha)^{e_{ss}}
\tilde{\overline x}_r'(\alpha)^{e_{s,s+1}}\\
M_r& =\tilde{\overline x}_1'(\alpha)^{e_{s+1,1}}\cdots\tilde{\overline x}_s'(\alpha)^{e_{s+1,s}}
\tilde{\overline x}_r'(\alpha)^{e_{s+1,s+1}}.
\endalign
$$
Let $k_1$ be the integral closure of $k$ in $T(\alpha)$. Set 
$$
A = \left( T''(\alpha)[M_1,\ldots, M_s,M_r]\otimes_{k_1}k(\overline T(\alpha+1))\right)_a
$$
where 
$$
a=(M_1,\ldots,M_s,\tilde{\overline x}_{s+1}'(\alpha),\ldots,\tilde{\overline x}_{r-1}'(\alpha),
M_r-c_{\alpha+1},\tilde{\overline x}_{r+1}'(\alpha),\ldots,\tilde{\overline x}_{n}'(\alpha)).
$$
Set $q_1=T''(\alpha)[M_1,\ldots,M_s,M_r]\cap a$ where the intersection is in $A$. Define
$$
T'(\alpha+1) = T''(\alpha)[M_1,\ldots,M_s,M_r]_{q_1}.
$$
 Define $T(\alpha+1) = T'(\alpha+1)\,\,\hat{}\,$.
Our inclusion $T'(\alpha+1)\subset A$ induces an isomorphism
$\eta_T^{\alpha+1}:k(T(\alpha+1))\rightarrow k(\overline T(\alpha+1))$. We can thus identify
$c_{\alpha+1}$ with $(\eta_T^{\alpha+1})^{-1}(c_{\alpha+1})$.
$T(\alpha+1)$ has regular parameters $(\hat{\overline x}_1(\alpha+1),\ldots,\hat{\overline x}_n(\alpha+1))$
defined by
$$
\hat{\overline x}_i(\alpha+1)=\cases M_i&1\le i\le s\\
M_r-c_{\alpha+1}&i=r\\
\tilde{\overline x}_i'(\alpha)& s<i,i\ne r
\endcases\tag 3.47
$$
Set 
$$
T''(\alpha+1) = T''(\alpha)\left[c_{\alpha+1},\left(\frac{\hat{\overline x}_r(\alpha+1)}{c_{\alpha+1}}+1\right)
^{\frac{1}{\overline c_{\alpha+1}}}\right]_{(\hat{\overline x}_1(\alpha+1),\ldots,\hat{\overline x}_n(\alpha+1))}\tag 3.48
$$
where $\left(\frac{\hat{\overline x}_r(\alpha+1)}{c_{\alpha+1}}+1\right)
^{\frac{1}{\overline c_{\alpha+1}}}$ has residue 1 in $k(T(\alpha+1))$. $T''(\alpha+1)$ has regular parameters
$(\overline x_1(\alpha+1),\ldots,\overline x_n(\alpha+1))$ defined by
$$
\overline x_i(\alpha+1) =\cases
\hat{\overline x}_i(\alpha+1) \left(\frac{\hat{\overline x}_r(\alpha+1)}{ c_{\alpha+1}}+1\right)
^{-\gamma_i(\alpha+1)}&1\le i\le s\\
\left(\frac{\hat{\overline x}_r(\alpha+1)}{c_{\alpha+1}}+1\right)
^{\frac{1}{\overline c_{\alpha+1}}}-1&i=r\\
\hat{\overline x}_i(\alpha+1)&s<i,i\ne r
\endcases\tag 3.49
$$
Then $T(\alpha)\rightarrow T(\alpha+1)$ is a UTS of type $\text{II}_r$ with
$$
\align
\tilde{\overline  x}'_1(\alpha) &= \overline x_1(\alpha+1)^{a_{11}(\alpha+1)}\cdots \overline x_s(\alpha+1)^{a_{1s}(\alpha+1)}
c_{\alpha+1}^{a_{1,s+1}(\alpha+1)}\tag 3.50\\
\vdots&\\
\tilde{\overline x}'_s(\alpha) &=
 \overline x_1(\alpha+1)^{a_{s1}(\alpha+1)}\cdots \overline x_s(\alpha+1)^{a_{ss}(\alpha+1)}
c_{\alpha+1}^{a_{s,s+1}(\alpha+1)}\\
\tilde{\overline x}'_r(\alpha) &= \overline x_1(\alpha+1)^{a_{s+1,1}(\alpha+1)}\cdots \overline x_s(\alpha+1)^{a_{s+1,s}(\alpha+1)}
(\overline x_r(\alpha+1)+1)c_{\alpha+1}^{a_{s+1,s+1}(\alpha+1)}
\endalign
$$

If $\overline U(\alpha)\rightarrow \overline U(\alpha+1)$ is of type I, defined by (3.26) of Theorem 3.8, 
 $U(\alpha)\rightarrow U(\alpha+1)$ will be the UTS of type I such that $U''(\alpha+1)$ has regular parameters
$(\overline y_1(\alpha+1),\ldots,\overline y_s(\alpha+1))$ defined by
$$
\align
\tilde{\overline y}'_1(\alpha) &= \overline y_1(\alpha+1)^{b_{11}(\alpha+1)}\cdots \overline y_s(\alpha+1)^{b_{1s}(\alpha+1)}
\tag 3.51\\
\vdots&\\
\tilde{\overline  y}'_s(\alpha) &= \overline y_1(\alpha+1)^{b_{s1}(\alpha+1)}\cdots \overline y_s(\alpha+1)^{b_{ss}(\alpha+1)}.
\endalign
$$

Suppose that $\overline U(\alpha)\rightarrow \overline U(\alpha+1)$ is of type $\text{II}_r$, defined by (3.27)
of Theorem 3.8.
 Set
$(f_{ij})=(b_{ij}(\alpha+1))^{-1}$,
$$
\align
N_1&=\tilde{\overline y}_1'(\alpha)^{f_{11}}\cdots\tilde{\overline y}_s'(\alpha)^{f_{1s}}
\tilde{\overline y}_r'(\alpha)^{f_{1,s+1}}\tag 3.52\\
&\vdots\\
N_s& =\tilde{\overline y}_1'(\alpha)^{f_{s1}}\cdots\tilde{\overline y}_s'(\alpha)^{e_{ss}}
\tilde{\overline y}_r'(\alpha)^{f_{s,s+1}}\\
N_r& =\tilde{\overline y}_1'(\alpha)^{f_{s+1,1}}\cdots\tilde{\overline y}_s'(\alpha)^{f_{s+1,s}}
\tilde{\overline y}_r'(\alpha)^{f_{s+1,s+1}}.
\endalign
$$
Let $k_2$ be the integral closure of $k$ in $U(\alpha)$. Set 
$$
B = \left( U''(\alpha)[N_1,\ldots, N_s,N_r]\otimes_{k_2}k(\overline U(\alpha+1))\right)_{b}\tag 3.53
$$
where 
$$
b=(N_1,\ldots,N_s,\tilde{\overline y}_{s+1}'(\alpha),\ldots,\tilde{\overline y}_{r-1}'(\alpha),
N_r-d_{\alpha+1},\tilde{\overline y}_{r+1}'(\alpha),\ldots,\tilde{\overline y}_{s+1}'(\alpha)).
$$
Set $q_2=U''(\alpha)[N_1,\ldots,N_s,N_r]\cap b$ where the intersection is in $B$. Define
$$
U'(\alpha+1) = U''(\alpha)[N_1,\ldots,N_s,N_r]_{q_2}.\tag 3.54
$$
Define $U(\alpha+1) = U'(\alpha+1)\,\,\hat{}\,$.
Our inclusion $U'(\alpha+1)\subset B$ induces an isomorphism
$\eta_U^{\alpha+1}:k(U(\alpha+1))\rightarrow k(\overline U(\alpha+1))$. We can thus identify
$d_{\alpha+1}$ with $(\eta_U^{\alpha+1})^{-1}(d_{\alpha+1})$. 
$U(\alpha+1)$ has regular parameters $(\hat{\overline y}_1(\alpha+1),\ldots,\hat{\overline y}_n(\alpha+1))$
defined by
$$
\hat{\overline y}_i(\alpha+1)=\cases N_i&1\le i\le s\\
N_r-d_{\alpha+1}&i=r\\
\tilde{\overline y}_i'(\alpha)& s<i,i\ne r
\endcases\tag 3.55
$$
Set 
$$
U''(\alpha+1) = U''(\alpha)\left[d_{\alpha+1},\left(\frac{\hat{\overline y}_r(\alpha+1)}{d_{\alpha+1}}+1\right)
^{\frac{1}{\overline d_{\alpha+1}}}\right]_{(\hat{\overline y}_1(\alpha+1),\ldots,\hat{\overline y}_n(\alpha+1))}\tag 3.56
$$
where $\left(\frac{\hat{\overline y}_r(\alpha+1)}{d_{\alpha+1}}+1\right)
^{\frac{1}{\overline d_{\alpha+1}}}$ has residue 1 in $k(U(\alpha+1))$. $U''(\alpha+1)$ has regular parameters
$(\overline y_1(\alpha+1),\ldots,\overline y_n(\alpha+1))$ defined by
$$
\overline y_i(\alpha+1) =\cases
\hat{\overline y}_i(\alpha+1) \left(\frac{\hat{\overline y}_r(\alpha+1)}{d_{\alpha+1}}+1\right)
^{-\tau_i(\alpha+1)}&1\le i\le s\\
\left(\frac{\hat{\overline y}_r(\alpha+1)}{d_{\alpha+1}}+1\right)
^{\frac{1}{\overline d_{\alpha+1}}}-1&i=r\\
\hat{\overline y}_i(\alpha+1)&s<i,i\ne r
\endcases\tag 3.57
$$
Then $U(\alpha)\rightarrow U(\alpha+1)$ is a UTS of type $\text{II}_r$ with
$$
\align
\tilde{\overline  y}'_1(\alpha) &= \overline y_1(\alpha+1)^{b_{11}(\alpha+1)}\cdots \overline y_s(\alpha+1)^{b_{1s}(\alpha+1)}
d_{\alpha+1}^{b_{1,s+1}(\alpha+1)}\tag 3.58\\
\vdots&\\
\tilde{\overline y}'_s(\alpha) &=
 \overline y_1(\alpha+1)^{b_{s1}(\alpha+1)}\cdots \overline y_s(\alpha+1)^{b_{ss}(\alpha+1)}d_{\alpha+1}^{b_{s,s+1}
(\alpha+1)}\\
\tilde{\overline y}'_r(\alpha) &= \overline y_1(\alpha+1)^{b_{s+1,1}(\alpha+1)}\cdots \overline y_s(\alpha+1)^{b_{s+1,s}(\alpha+1)}
(\overline y_r(\alpha+1)+1)d_{\alpha+1}^{b_{s+1,s+1}(\alpha+1)}.
\endalign
$$

We will now prove that (3.37), (3.39) and (3.40) hold for $\alpha+1$. The essential case is when 
$\overline T(\alpha)\rightarrow \overline T(\alpha+1)$ is of type $\text{II}_r$ with $s+1\le r\le l$.

By (3.13) of Lemma 3.4 in the construction of $\overline T(\alpha)\rightarrow \overline T(\alpha+1)$ and
$\overline U(\alpha)\rightarrow \overline U(\alpha+1)$,
$$
\align
M_1&=N_1^{g_{11}}\cdots N_s^{g_{1s}}N_r^{g_{1,s+1}}\beta_1\tag 3.59\\
&\vdots\\
M_s&=N_1^{g_{s1}}\cdots N_s^{g_{ss}}N_r^{g_{s,s+1}}\beta_s\\
M_r&=N_r^{g_{s,s+1}}\beta_r
\endalign
$$
$\beta_i=\phi_1(\alpha)^{e_{i1}}\cdots\phi_s(\alpha)^{e_{is}}\in k(U(\alpha))\subset U''(\alpha)$ for $1\le i\le s$,
$$
\beta_r=\phi_1(\alpha)^{e_{s+1,1}}\cdots\phi_s(\alpha)^{e_{s+1,s}}.
$$
$$
(g_{ij}) = (a_{ij}(\alpha+1))^{-1}
\left(\matrix (c_{ij}(\alpha))&0\\ 0&1\endmatrix\right)(b_{ij}(\alpha+1))
$$
$g_{s+1,1}=\cdots= g_{s+1,s}=0$ and $g_{ij}\ge 0$ for all $i,j$. Thus
$T''(\alpha)[M_1,\ldots,M_s,M_r]\subset U''(\alpha)[N_1,\ldots,N_s,N_r]$. Our inclusion
$k(\overline T(\alpha+1))\rightarrow k(\overline U(\alpha+1))$ induces  an identification
$c_{\alpha+1}=d_{\alpha+1}^{g_{s+1,s+1}}\beta_r$. Then by (3.16) of Lemma 3.4,
$$
\align
M_r-c_{\alpha+1}&=N_r^{g_{s+1,s+1}}\beta_r-c_{\alpha+1}=(N_r^{g_{s+1,s+1}}-d_{\alpha+1}^{g_{s+1,s+1}})\beta_r\\
&=\prod_{i=1}^{g_{s+1,s+1}}(N_r-\omega^id_{\alpha+1})\beta_r
\endalign
$$
where $\omega$ is a primitive $g_{s+1,s+1}$-th root of unity (in an algebraic closure of $k(U(\alpha+1))$).
Thus $N_r-d_{\alpha+1}$ divides $M_r-c_{\alpha+1}$ in $U''(\alpha)[N_1,\ldots,N_s,N_r]$ and
we have an inclusion $A\subset B$ which induces $T'(\alpha+1)\subset U'(\alpha+1)$
and $T(\alpha+1)\subset U(\alpha+1)$. Thus (3.39) holds for $\alpha+1$.

By the argument of Lemma 3.4 in the construction of $\overline T(\alpha+1)\rightarrow \overline U(\alpha+1)$ and 
(3.47)-(3.49), (3.55)-(3.57),
we have that $T''(\alpha+1)\subset U''(\alpha+1)$ and (3.40) holds for $\alpha+1$.

\subheading{Step 2} 

Suppose that $T(t)\rightarrow U(t)$ is constructed as in Step 1, and $f$ satisfies 1), 2) or 3) in the statement of
Theorem 3.9. We will show that $f$ satisfies the respective equation (3.34), (3.35) or (3.36) in $U(t)$.
It suffices to prove the following  statement.

Suppose that $0\le j\le t$ and 
$\overline f_j$ 
is defined by $f(\overline w_1,\cdots,\overline w_n)=
\overline f_j(\overline w_1(j),\cdots,\overline  w_n(j))$ in $\overline U(j)$.
Then 
$$
f(\tilde {\overline y}_1,\cdots,\tilde {\overline y}_n)\equiv\overline f_j(\overline y_1(j),\cdots,\overline y_n(j))
\text{ mod }m(U(j))^{N}.\tag 3.60
$$
The statement (3.60) will be proved by induction on $j$. By induction, suppose that
$$
f(\tilde {\overline y}_1,\cdots, \tilde {\overline y}_n)
\equiv \overline f_j(\overline y_1(j),\cdots,\overline y_n(j))\text{ mod } m(U(j))^{N}.
$$
We have  $\overline f_j(\overline w_1(j),\cdots,\overline w_n(j))
=\overline f_{j+1}(\overline w_1(j+1),\cdots,\overline w_n(j+1))$ in $\overline U(j+1)$.

There are series $\overline P_{i,j}$ with coefficients in $k(\overline U(j))$ such that
$$
\overline P_{i,j}(\tilde{\overline w}'_1(j),\cdots,\tilde{\overline w}'_{i-1}(j))=\cases
P_{i,j}(\overline z_1(j),\cdots,\overline z_{i-1}(j))& s+1\le i\le l\\
P_{i,j}(\overline w_1(j),\cdots,\overline w_{i-1}(j)) &l+1<i
\endcases 
$$
We have
$$
\align
f(\tilde {\overline y}_1,\cdots,\tilde {\overline y}_n)
\equiv&\overline f_j(\overline y_1(j),\cdots,\overline y_n(j))\text{ mod }m(U(j+1))^{N}\\
\equiv& \overline f_j(\tilde{\overline y}'_1(j),\cdots,\tilde{\overline y}'_s(j),\\
&\tilde{\overline y}'_{s+1}(j)+\overline P_{s+1,j}(\tilde{\overline y}'_1(j),\cdots,\tilde{\overline y}'_{s}(j)),
\cdots,\\
&\tilde{\overline y}'_{n}(j)+\overline P_{n,j}(\tilde{\overline y}'_1(j),\cdots,
\tilde{\overline y}'_{n-1}(j)))
\text{ mod }m(U(j+1))^N\\
\equiv &\overline f_{j+1}(\overline y_1(j+1),\cdots,\overline y_n(j+1))\text{ mod }m(U(j+1))^N
\endalign
$$
Set $n_{\alpha,l}=m\left(k(U(\alpha))[[\overline y_1(\alpha),\ldots,\overline y_l(\alpha)]]\right)$ for $1\le \alpha\le t$.
In the case of $f(\tilde{\overline y}_1,\ldots,\tilde{\overline y}_l)\in 
k(U)[[\tilde{\overline y}_1,\ldots,\tilde{\overline y}_l]]$,
 the above argument is valid with $n$ replaced by $l$
 and $m(U(j+1))$ replaced by $n_{j+1,l}$, since $\overline U(0)\rightarrow \overline U(t)$ is then a UTS in the
first $l$ variables.

\subheading{Step 3} Now we will construct, with suitable choice of the series $\Lambda_{\alpha}$ and $\Omega_{i,\alpha}$,
 in (3.41)-(3.44)
 of Step 1, and our fixed $N$,
 a CRUTS $(R,T''(t),T(t))$ and $(S,U''(t),U(t))$ 
 with associated MTS $(3.61)$

$$
\matrix
S=&S(0) & \rightarrow &S(1) &\rightarrow &\cdots &\rightarrow & S(t)\\
&\uparrow &&\uparrow&&&&\uparrow\\
R=&R(0) &\rightarrow &R(1)&\rightarrow &\cdots&\rightarrow& R(t)
\endmatrix\tag 3.61
$$
such that $R(\alpha)$ has regular parameters
$$
(x_1(\alpha),\ldots, x_n(\alpha)), (\tilde x_1(\alpha),\ldots, \tilde x_n(\alpha)),
(\tilde x_1'(\alpha),\ldots, \tilde x_n'(\alpha)).
$$
$S(\alpha)$ has regular parameters
$$
(y_1(\alpha),\ldots, y_n(\alpha)), (\tilde y_1(\alpha),\ldots, \tilde y_n(\alpha)),
 (\tilde y_1'(\alpha),\ldots, \tilde y_n'(\alpha)).
$$
$(3.61)$ will consist of three  types of MTSs.
\item{M1)} $R(\alpha)\rightarrow R(\alpha+1)$ and $S(\alpha)\rightarrow S(\alpha+1)$ are of type I.
\item{M2)} $R(\alpha)\rightarrow R(\alpha+1)$ is of type $\text{II}_r$,  $s+1\le r\le l$,
and $S(\alpha)\rightarrow S(\alpha+1)$ 
is a MTS  of type $\text{II}_r$,  followed by a MTS of type I.
\item{M3)} $R(\alpha)= R(\alpha+1)$ and $S(\alpha)\rightarrow S(\alpha+1)$ is of type $\text{II}_r$ ($l+1\le r\le m$).

There exists for all $\alpha$ a regular local ring $\tilde R(\alpha)\subset R(\alpha)$ such that $\tilde R(\alpha)$ has
regular parameters $(x_1(\alpha),\ldots,x_l(\alpha))$ and 
${\tilde R}(\alpha)\,\,\hat{}\, \cong k(c_0,\ldots,c_{\alpha})[[x_1(\alpha),\ldots,x_l(\alpha)]]$.

 The series $\Lambda_{\alpha}$ in (3.41) is chosen so that
$$
x_i(\alpha)=\gamma_i(\alpha)\tilde{\overline x}_i(\alpha)
$$
 for $1\le i\le l$ where 
$\gamma_i(\alpha)\in k(c_0,\ldots,c_{\alpha})[[x_1(\alpha),\ldots x_l(\alpha)]]\cap T''(\alpha)$ are units
such that $\gamma_i(\alpha)\equiv 1\text{ mod }(x_1(\alpha),\ldots,x_l(\alpha))$. In fact,
in conjunction with an appropriate choice of $\lambda_{\alpha-1}$ in (2) below,
 we will have
$x_r(\alpha) = \gamma_r(\alpha)\overline x_r(\alpha)+\psi_{\alpha}$
where 
$$
\psi_{\alpha}\in k(c_0,\ldots,c_{\alpha})[[\overline x_1(\alpha),\ldots,\overline x_l(\alpha)]],
$$ 
 and $\text{mult}(\psi_{\alpha})\ge N$ if $R(\alpha-1)\rightarrow R(\alpha)$ is of type $\text{II}_r$.

 The series  $\Omega_{i,\alpha}$ in (3.42)  is chosen so that
 we can define regular parameters $\tilde x_i(\alpha)$
in $R(\alpha)$ by 
$$
\tilde x_i(\alpha)=\cases
\gamma_i(\alpha)\tilde{\overline x}'_i(\alpha)&1\le i\le l\\
\tilde{\overline x}'_i(\alpha)&l<i
\endcases \tag 3.62
$$
(If $P_{i,\alpha}=0$, or if $1\le i\le s$ we will have 
 $\tilde x_i(\alpha)=x_i(\alpha)$.)
Define regular parameters $\tilde x'_i(\alpha)$
in $R(\alpha)$ by 
$$
\tilde x'_i(\alpha) = 
\cases
\lambda_{\alpha}\tilde x_r(\alpha), 
& \text{ if }R(\alpha)\rightarrow R(\alpha+1)\text{ is of type II}_r\text{ and } i=r\\
\tilde x_i(\alpha) & \text{ otherwise}
\endcases\tag 3.63
$$
We will have $\lambda_{\alpha}\in \tilde R(\alpha)\subset k(c_0,\ldots,c_{\alpha})[[x_1(\alpha),\ldots, x_l(\alpha)]]$,
$\lambda_{\alpha}
\equiv 1\text{ mod }m(\tilde R(\alpha))^N$.

These variables are such that for all $\alpha$, 
$$
\align
&k(c_0,\ldots,c_{\alpha})[[x_1(\alpha),\ldots, x_l(\alpha)]]
=k(c_0,\ldots,c_{\alpha})[[\tilde x_1(\alpha),\ldots, \tilde x_l(\alpha)]]\\
&=k(c_0,\ldots,c_{\alpha})[[\tilde x_1'(\alpha),\ldots, \tilde x_l'(\alpha)]]\\
&=k(c_0,\ldots,c_{\alpha})[[\overline x_1(\alpha),\ldots, \overline x_l(\alpha)]]
=k(c_0,\ldots,c_{\alpha})[[\tilde{\overline  x}_1(\alpha),\ldots, \tilde{\overline x}_l(\alpha)]]\\
&=k(c_0,\ldots,c_{\alpha})[[\tilde{\overline  x}'_1(\alpha),\ldots, \tilde{\overline x}'_l(\alpha)]]
\endalign
$$
and 
$$
x_i(\alpha)=\tilde x_i(\alpha) =\tilde x_i'(\alpha) =\overline x_i(\alpha) 
=\tilde{\overline x}_i(\alpha) = \tilde{\overline x}_i'(\alpha)
$$
for $l<i\le n$.

 $k(c_0,\ldots,c_{\alpha})\subset k(R(\alpha))$ and 
$k(\tilde R(\alpha))\cong k(c_0,\ldots,c_{\alpha})$
for all $i$.

 The series $\Lambda_{\alpha}$ in (3.43) is chosen so that
$$
y_i(\alpha)=\gamma_i^y(\alpha)\tilde{\overline y}_i(\alpha)
$$
  where 
$\gamma_i^y(\alpha)\in U''(\alpha)$ is a  unit such that $\gamma_i^y(\alpha)\equiv 1\text{ mod }m(U''(\alpha))$
for $1\le i\le m$. In fact,
we will have 
$y_r(\alpha) = \gamma_r^y(\alpha)\overline y_r(\alpha)+\psi_{\alpha}$
where $\psi_{\alpha}\in  S(\alpha)\,\,\hat{}\,$, 
  and $\text{mult}(\psi_{\alpha})\ge N$ if $S(\alpha-1)\rightarrow S(\alpha)$ is of type $\text{II}_r$,
with $l+1\le r$.

The series  $\Omega_{i,\alpha}$ in (3.44) is chosen so that we can define regular parameters $\tilde y_i(\alpha)$
in $S(\alpha)$ by 
$$
\tilde y_i(\alpha) =  \gamma_i^y(\alpha)\tilde{\overline y}'_i(\alpha).\tag 3.64
$$
which satisfy
$$
\tilde y_i(\alpha)=
\cases 
y_i(\alpha)&\text{ if }1\le i\le s\text{ or }m<i\le n\\
\tilde  x_i(\alpha)&\text{ if }s+1\le i\le l\\
\gamma_i^y(\alpha)\tilde{\overline y}'_i(\alpha),&\text{ if } l+1\le i\le m.\
\endcases
$$
We will have $\gamma_i^y(\alpha) = \gamma_i(\alpha)$ if $s+1\le i\le l$. Define
regular parameters $\tilde y'_i(\alpha)$
in $S(\alpha)$ by 
$$
\tilde y'_i(\alpha) = 
\cases
\tilde x_i'(\alpha)&\text{ if }s+1\le i\le l\\
\lambda_{\alpha}\tilde y_r(\alpha)
& \text{ if }S(\alpha)\rightarrow S(\alpha+1)\text{ is of type II}_r\text{ with }r \ge l+1, i=r\\
\tilde y_i(\alpha) & \text{ otherwise}
\endcases\tag 3.65
$$

We will have $\lambda_{\alpha}\in  S(\alpha),\lambda_{\alpha}
\equiv 1\text{ mod }m(S(\alpha))^N$.

These variables are such that for all $\alpha$,
$y_i(\alpha)=\tilde y_i(\alpha)=y_i$
for $m<i\le n$.

Suppose that we have constructed the CRUTS  out to $(R,T''(\alpha),T(\alpha))$ and \linebreak
$(S,U''(\alpha),U(\alpha))$, the MTS out to
$$
\matrix
S & \rightarrow & S(\alpha)\\
\uparrow &&\uparrow\\
R &\rightarrow & R(\alpha)
\endmatrix,
$$
we have constructed $\tilde R(\alpha)\subset R(\alpha)$
and have defined regular parameters $(x_1(\alpha),\ldots,x_n(\alpha))$ in $R(\alpha)$, 
$(\overline x_1(\alpha),\ldots,\overline x_n(\alpha))$,
 $(\tilde{\overline x}_1(\alpha),\ldots,\tilde{\overline x}_n(\alpha))$ in $T''(\alpha)$, 
$(y_1(\alpha),\ldots,y_n(\alpha))$ in $S(\alpha)$,
$(\overline y_1(\alpha),\ldots,\overline y_n(\alpha))$ and
 $(\tilde{\overline y}_1(\alpha),\ldots,\tilde{\overline y}_n(\alpha))$ in $U''(\alpha)$.

\subheading{Case 1}
Suppose that both $\overline T(\alpha)\rightarrow \overline T(\alpha+1)$ and 
$\overline U(\alpha)\rightarrow \overline U(\alpha+1)$ are of type I. By assumption
$x_i(\alpha)=\gamma_i(\alpha)\tilde{\overline x}_i(\alpha)$ for $1\le i\le l$. 
$$
-\gamma_i(\alpha)P_{i,\alpha}(\tilde{\overline x}_1(\alpha),\ldots,
\tilde{\overline x}_{i-1}(\alpha))\in k(c_0,\ldots,c_{\alpha})[[\tilde{\overline x}_1(\alpha),\ldots,
\tilde{\overline x}_l(\alpha)]]={\tilde R}(\alpha)\,\,\hat{}\,,
$$
 the completion of $\tilde R(\alpha)$ for $s+1\le i\le l$.
Thus there exists $A\in\tilde R(\alpha)\subset R(\alpha)$, $\overline \Omega_{i,\alpha}\in 
m\left(k(c_0,\ldots,c_{\alpha})
[[\tilde{\overline x}_1(\alpha),\ldots,\tilde{\overline x}_l(\alpha)]]\right)^N$ such that
$A-\overline \Omega_{i,\alpha}=-\gamma_i(\alpha)P_{i,\alpha}$. 
Set $\Omega_{i,\alpha}=\gamma_i(\alpha)^{-1}\overline \Omega_{i,\alpha}$.
$$
\align
&\gamma_i(\alpha)[\tilde{\overline x}_i(\alpha)-P_{i,\alpha}+\Omega_{i,\alpha}]\\
&=\gamma_i(\alpha)\tilde{\overline x}_i(\alpha)-\gamma_i(\alpha)P_{i,\alpha}+\overline \Omega_{i,\alpha}\\
&=x_i(\alpha)+A\in \tilde R(\alpha).
\endalign
$$
Thus by suitable choice of the  $\Omega_{i,\alpha}$, we have regular parameters $\tilde x'_i(\alpha)$ 
in $R(\alpha)$ and 
regular parameters $\tilde{\overline x}_i'(\alpha)$ in $T''(\alpha)$ satisfying (3.62) and (3.63).
We can also define $\Omega_{i,\alpha}$ for $l+1\le i\le m$ to get
 regular parameters $\tilde y'_i(\alpha)$ in $S(\alpha)$, and regular parameters
$\tilde{\overline y}_i'(\alpha)$ in $U''(\alpha)$ satisfying (3.64) and (3.65).
 Define 
$R(\alpha)\rightarrow R(\alpha+1)$ and  $S(\alpha)\rightarrow S(\alpha+1)$ by
$$
\align
\tilde x_1'(\alpha) &= x_1(\alpha+1)^{a_{11}(\alpha+1)}\cdots x_s(\alpha+1)^{a_{1s}(\alpha+1)}\\
\vdots&\\
\tilde x_s'(\alpha) &= x_1(\alpha+1)^{a_{s1}(\alpha+1)}\cdots x_s(\alpha+1)^{a_{ss}(\alpha+1)}
\endalign
$$
and
$$
\align
\tilde y_1'(\alpha) &= y_1(\alpha+1)^{b_{11}(\alpha+1)}\cdots y_s(\alpha+1)^{b_{1s}(\alpha+1)}\\
\vdots&\\
\tilde y_s'(\alpha) &= y_1(\alpha+1)^{b_{s1}(\alpha+1)}\cdots y_s(\alpha+1)^{b_{ss}(\alpha+1)}.
\endalign
$$
Then $T(\alpha+1)= R(\alpha+1)\,\,\hat{}\,$ and $U(\alpha+1)= S(\alpha+1)\,\,\hat{}\,$.
 Set $\tilde{\overline x}_i(\alpha+1)=\overline x_i(\alpha+1)$,
$\tilde{\overline y}_i(\alpha+1)=\overline y_i(\alpha+1)$ for all $i$.

Set $(e_{ij})=(a_{ij}(\alpha+1))^{-1}$, set
$$
\gamma_i(\alpha+1)=\cases
\gamma_1(\alpha)^{e_{i1}}\cdots\gamma_s(\alpha)^{e_{is}}& 1\le i\le s\\
\gamma_i(\alpha)&s<i\le l
\endcases
$$
$\gamma_i(\alpha+1)\in
k(c_0,\ldots,c_{\alpha+1})[[ x_1(\alpha+1),\ldots, x_l(\alpha+1)]]\cap T''(\alpha+1)$
for $1\le i\le l$. 
$$
x_i(\alpha+1)=\cases
\gamma_i(\alpha+1)\tilde{\overline x}_i(\alpha+1)& 1\le i\le l\\
\tilde{\overline x}_i(\alpha+1)&l+1<i\le n
\endcases
$$
Set $(f_{ij})=(b_{ij}(\alpha+1))^{-1}$, set
$$
\gamma_i^y(\alpha+1)=\cases
\gamma_1^y(\alpha)^{f_{i1}}\cdots\gamma_s^y(\alpha)^{f_{is}}& 1\le i\le s\\
\gamma_i^y(\alpha)&s<i\le n
\endcases
$$
 Then
$y_i(\alpha+1)=
\gamma_i^y(\alpha+1)\tilde{\overline y}_i(\alpha+1)\text{ for }1\le i\le n$.
Set 
$$
\tilde R(\alpha+1)=\tilde R(\alpha)[x_1(\alpha+1),\ldots,x_s(\alpha+1)]_{(x_1(\alpha+1),\ldots,x_l(\alpha+1))}.
$$

\subheading{Case 2}
Suppose that both $\overline T(\alpha)\rightarrow \overline T(\alpha+1)$ and 
$\overline U(\alpha)\rightarrow \overline U(\alpha+1)$ are of type $\text{II}_r$
with $s+1\le r\le l$.
By suitable choice of the  $\Omega_{i,\alpha}$ as in Case 1, we have regular parameters $\tilde x_i(\alpha)$ in $R(\alpha)$
and regular parameters $\tilde y_i(\alpha)$ in $S(\alpha)$ satisfying (3.62) and (3.64). 
Set $(e_{ij}) = (a_{ij}(\alpha+1))^{-1}$.
Choose $\lambda_{\alpha}\in \tilde R(\alpha)\subset R(\alpha)$ in (3.63) so that 
$$
\align
& \gamma_1(\alpha)^{e_{s+1,1}}\cdots\gamma_s(\alpha)^{e_{s+1,s}}
\gamma_r(\alpha)^{e_{s+1,s+1}}\lambda_{\alpha}^{e_{s+1,s+1}}\tag 3.66\\
&\equiv 1\text{ mod }(x_1(\alpha),\ldots,x_l(\alpha))^Nk(c_0,\ldots,c_{\alpha})[[x_1(\alpha),\ldots,x_l(\alpha)]].
\endalign
$$

Set 
$$
\align
A_1&=\tilde x_1'(\alpha)^{e_{11}}\cdots\tilde x_s'(\alpha)^{e_{1s}}\tilde x_r'(\alpha)^{e_{1,s+1}}\tag 3.67\\
&\vdots\\
A_s&=\tilde x_1'(\alpha)^{e_{s1}}\cdots\tilde x_s'(\alpha)^{e_{ss}}\tilde x_r'(\alpha)^{e_{s,s+1}}\\
A_r&=\tilde x_1'(\alpha)^{e_{s+1,1}}\cdots\tilde x_s'(\alpha)^{e_{s+1,s}}\tilde x_r'(\alpha)^{e_{s+1,s+1}}
\endalign
$$
where $(e_{ij})=(a_{ij}(\alpha+1))^{-1}$.

Let $k_3$ be the integral closure of $k$ in $R(\alpha)$. We have 
$$
k_3\rightarrow k(R(\alpha))\rightarrow k(T(\alpha))\rightarrow k(T(\alpha+1)).
$$
Set 
$$
C=\left( R(\alpha)[A_1,\ldots,A_s,A_r]\otimes_{k_3}k(T(\alpha+1))\right)_c
$$
where $c=(A_1,\ldots,A_s,\tilde x_{s+1}'(\alpha),\ldots,\tilde x_{r-1}'(\alpha),
A_r-c_{\alpha+1},\tilde x_{r+1}'(\alpha),\ldots,\tilde x_n'(\alpha))$.
Set $q_3=R(\alpha)[A_1,\ldots,A_s,A_r]\cap c$ where the intersection is in $C$. Define
$$
R(\alpha +1)=R(\alpha)[A_1,\ldots,A_s,A_r]_{q_3}.
$$
 Our construction gives an isomorphism
$k(R(\alpha+1))\cong k(T(\alpha+1))$. Define regular parameters $(x_1^*(\alpha+1),\ldots,x_n^*(\alpha+1))$ in
$R(\alpha+1)\,\,\hat{}\,$ by
$$
x_i^*(\alpha+1)=\cases
A_i&1\le i\le s\\
A_r-c_{\alpha+1}&i=r\\
\tilde x_i'(\alpha)&s<i,i\ne r
\endcases\tag 3.68
$$

$$
\align
\tilde x_1'(\alpha) &= x_1^*(\alpha+1)^{a_{11}(\alpha+1)}\cdots x_s^*(\alpha+1)^{a_{1s}(\alpha+1)}
(x_r^*(\alpha+1)+c_{\alpha+1})^{a_{1,s+1}(\alpha+1)}\\
\vdots&\\
\tilde x_s'(\alpha) &= x_1^*(\alpha+1)^{a_{s1}(\alpha+1}\cdots x_s^*(\alpha+1)^{a_{ss}(\alpha+1)}
(x_r^*(\alpha+1)+c_{\alpha+1})^{a_{s,s+1}(\alpha+1)}\\
\tilde x_r'(\alpha) &= x_1^*(\alpha+1)^{a_{s+1,1}(\alpha+1)}\cdots x_s^*(\alpha+1)^{a_{s+1,s}(\alpha+1)}
(x_r^*(\alpha+1)+c_{\alpha+1})^{a_{s+1,s+1}(\alpha+1)}
\endalign
$$
Set $\tilde R(\alpha+1)=\tilde R(\alpha)[A_1,\ldots,A_s,A_r]_q$ where $q=\tilde R(\alpha)[A_1,\ldots,A_s,A_r]\cap c$.
$$
k(\tilde R(\alpha+1))\cong k(\tilde R(\alpha))(c_{\alpha+1})\cong k(c_0,\ldots,c_{\alpha+1})
$$
and 
$$
\tilde R(\alpha+1)\,\,\hat{}\,\,\cong k(c_0,\ldots,c_{\alpha+1})[[x_1^*(\alpha+1),\ldots,x_l^*(\alpha+1)]].
$$
Set
$$
\align
B_1&=\tilde y_1'(\alpha)^{f_{11}}\cdots\tilde y_s'(\alpha)^{f_{1s}}\tilde y_r'(\alpha)^{f_{1,s+1}}\tag 3.69\\
&\vdots\\
B_s&=\tilde y_1'(\alpha)^{f_{s1}}\cdots\tilde y_s'(\alpha)^{f_{ss}}\tilde y_r'(\alpha)^{f_{s+1,s+1}}\\
B_r&=\tilde y_1'(\alpha)^{f_{s+1,1}}\cdots\tilde y_s'(\alpha)^{f_{s+1,s}}\tilde y_r'(\alpha)^{f_{s+1,s+1}}
\endalign
$$
where $(f_{ij})=(b_{ij}(\alpha+1))^{-1}$. Let $k_4$ be the integral closure of $k$ in $S(\alpha)$. We have
$$
k_4\rightarrow k(S(\alpha))\rightarrow k(U(\alpha))\rightarrow k(U(\alpha+1)).
$$
Set 
$$
D=\left(S(\alpha)[B_1,\ldots,B_s,B_r]\otimes_{k_4}k(U(\alpha+1))\right)_d
$$
where $d=(B_1,\ldots,B_s,\tilde y_{s+1}'(\alpha),\ldots,\tilde y_{r-1}'(\alpha),B_r-d_{\alpha+1},
\tilde y_{r+1}'(\alpha),\ldots,\tilde y_n'(\alpha))$. Set $q_4=S(\alpha)[B_1,\ldots,B_s,B_r]\cap d$ where
the intersection is in $D$. Define 
$$
S(\alpha+1)=S(\alpha)[B_1,\ldots,B_s,B_r]_{q_4}.
$$
Our construction gives an isomorphism $k(S(\alpha+1))\cong k(U(\alpha+1))$.
Define regular parameters $(y_1^*(\alpha+1),\ldots,y_n^*(\alpha+1))$ in $S(\alpha+1)\,\,\hat{}\,$ by
$$
y_i^*(\alpha+1)=\cases
B_i&1\le i\le s\\
B_r-d_{\alpha+1}&i=r\\
\tilde y_i'(\alpha)&s<i,i\ne r
\endcases\tag 3.70
$$

$$
\align
\tilde y_1'(\alpha) &= y_1^*(\alpha+1)^{b_{11}(\alpha+1)}\cdots y_s^*(\alpha+1)^{b_{1s}(\alpha+1)}
(y_r^*(\alpha+1)+d_{\alpha+1})^{b_{1,s+1}(\alpha+1)}\\
\vdots&\\
\tilde y_s'(\alpha) &= y_1^*(\alpha+1)^{b_{s1}(\alpha+1}\cdots y_s^*(\alpha+1)^{b_{ss}(\alpha+1)}
(y_r^*(\alpha+1)+d_{\alpha+1})^{b_{s,s+1}(\alpha+1)}\\
\tilde y_r'(\alpha) &= y_1^*(\alpha+1)^{b_{s+1,1}(\alpha+1)}\cdots y_s^*(\alpha+1)^{b_{s+1,s}(\alpha+1)}
(y_r^*(\alpha+1)+d_{\alpha+1})^{b_{s+1,s+1}(\alpha+1)}.
\endalign
$$

Set $\sigma_i = \gamma_1(\alpha)^{e_{i1}}\cdots\gamma_s(\alpha)^{e_{is}}
\gamma_r(\alpha)^{e_{i,s+1}}\lambda_{\alpha}^{e_{i,s+1}}$ for $1\le i\le s$, and set
$$
\sigma_r = \gamma_1(\alpha)^{e_{s+1,1}}\cdots\gamma_s(\alpha)^{e_{s+1,s}}
\gamma_r(\alpha)^{e_{s+1,s+1}}\lambda_{\alpha}^{e_{s+1,s+1}}.
$$
By (3.62), (3.63), (3.46) and (3.67)
$$
\align
A_1&=\sigma_1M_1\tag 3.71\\
&\vdots\\
A_s&=\sigma_sM_s\\
A_r&=\sigma_rM_r.
\endalign
$$
Thus $R(\alpha)[A_1,\ldots,A_s,A_r]\subset T''(\alpha)[M_1,\ldots,M_s,M_r]$.
We then have an inclusion $C\rightarrow  A$ which induces an inclusion $R(\alpha+1)\rightarrow T'(\alpha+1)$.
By  (3.68), (3.71), (3.47) and (3.49),
$$
x_i^*(\alpha+1)
=\sigma_i(\overline x_r(\alpha+1)+1)^{\overline c_{\alpha+1}\gamma_i(\alpha+1)}\overline x_i(\alpha+1)
$$
for $1\le i\le s$ and
$$
\align
x_r^*(\alpha+1)+c_{\alpha+1}
&=\sigma_r (\hat{\overline x}_r(\alpha+1)+c_{\alpha+1})\\
&=\sigma_r\left(c_{\alpha+1}(\overline x_r(\alpha+1)+1)^{\overline c_{\alpha+1}}\right)\\
&=\sigma_r(\overline c_{\alpha+1}c_{\alpha+1}u\overline x_r(\alpha+1)+c_{\alpha+1}]
\endalign
$$
where $u\in\bold Q[\overline x_r(\alpha+1)]$ is a polynomial with constant term 1.
$$
\frac{x_r^*(\alpha+1)}{\overline c_{\alpha+1}c_{\alpha+1}}=\sigma_ru\overline x_r(\alpha+1)
+\frac{\sigma_r-1}{\overline c_{\alpha+1}}
$$
By (3.66),  $T(\alpha+1)=R(\alpha+1)\,\,\hat{}\,$ since $k(R(\alpha+1))\cong k(T(\alpha+1))$,
$(x_1^*(\alpha+1),\ldots,x_n^*(\alpha+1))$ are regular parameters in $T(\alpha+1)$.

Thus there exists 
$$
\Omega\in (\overline x_1(\alpha+1),\ldots,\overline x_l(\alpha+1))^Nk(c_0,\ldots,c_{\alpha+1})
[[\overline x_1(\alpha+1),\ldots,\overline x_l(\alpha+1)]]
$$
such that
$$
\sigma_ru\overline x_r(\alpha+1)+\Omega\in\tilde R(\alpha+1)\subset R(\alpha+1).
$$
Set
$$
x_i(\alpha+1)=\cases
x_i^*(\alpha+1)&i\ne r\\
\sigma_ru\overline x_r(\alpha+1)+\Omega&i=r
\endcases
$$
$$
\gamma_i(\alpha+1)=\cases
\sigma_i(\overline x_r(\alpha+1)+1)^{\overline c_{\alpha+1}\gamma_i(\alpha+1)}&1\le i\le s\\
\sigma_ru&i=r\\
\gamma_i(\alpha)&s<i\le l, i\ne r
\endcases
$$
Set $\Lambda_{\alpha+1}=\sigma_r^{-1}u^{-1}\Omega$. By definition
$$
\tilde{\overline x}_i(\alpha+1)=\cases
\overline x_r(\alpha+1)+\Lambda_{\alpha+1}& i=r\\
\overline x_i(\alpha+1)&i\ne r
\endcases
$$
Then $(x_1(\alpha+1),\ldots,x_n(\alpha+1))$ are regular parameters in $R(\alpha+1)$ and 
$$
x_i(\alpha+1)=\gamma_i(\alpha+1)\tilde{\overline x}_i(\alpha+1)
$$
for $1\le i\le l$.
Since $\sigma_i\in T''(\alpha)$ for all $i$, $c_{\alpha+1}\in T''(\alpha+1)$, $\overline c_{\alpha+1}\gamma_i(\alpha+1)$
are integers and
$(\overline x_1(\alpha+1),\ldots,\overline x_n(\alpha+1))$ are regular parameters in $T''(\alpha+1)$,
 $\gamma_i(\alpha+1)\in T''(\alpha+1)$ for all $i$ and
$(\tilde{\overline x}_1(\alpha+1),\ldots,\tilde{\overline x}_n(\alpha+1))$
are regular parameters in $T''(\alpha+1)$.

Set $\sigma_i^y = \gamma_1^y(\alpha)^{f_{i1}}\cdots\gamma_s^y(\alpha)^{f_{is}}\gamma_r^y(\alpha)^{f_{i,s+1}}
\lambda_{\alpha}^{f_{i,s+1}}$
for $1\le i\le s$ and set
$$
\sigma_r^y = \gamma_1^y(\alpha)^{f_{s+1,1}}\cdots\gamma_s^y(\alpha)^{f_{s+1,s}}\gamma_r^y(\alpha)^{f_{s+1,s+1}}
\lambda_{\alpha}^{f_{s+1,s+1}}.
$$
By (3.64), (3.65), (3.52) and (3.69)
$$
\align
B_1 &=\sigma_1^yN_1\tag 3.72\\
&\vdots\\
B_s&=\sigma_s^yN_s\\
B_r&=\sigma_r^yN_r.
\endalign
$$
Thus $S(\alpha)[B_1,\ldots,B_s,B_r]\subset U''(\alpha)[N_1,\ldots,N_s,N_r]$. We then have an inclusion $D\subset B$ 
which induces an inclusion $S(\alpha+1)\rightarrow U'(\alpha+1)$. By (3.70), (3.72), (3.55) and (3.57)

$$
\align
y_i^*(\alpha+1)&= 
\sigma_i^y \hat{\overline y}_i(\alpha+1)\\
&=\sigma_i^y(\overline y_r(\alpha+1)+1)^{\overline d_{\alpha+1}\tau_i(\alpha+1)}\overline y_i(\alpha+1)
\endalign
$$
 for $1\le i\le s$
so that $y_i^*(\alpha+1)=\text{unit }\overline y_i(\alpha+1)$ in $U(\alpha+1)$ for $1\le i\le s$, and
$$
\align
y_r^*(\alpha+1)+d_{\alpha+1} 
 &=\sigma_r^y (\hat{\overline y}_r(\alpha+1)+d_{\alpha+1})\\
&=\sigma_r^yd_{\alpha+1}(\overline y_r(\alpha+1)+1)^{\overline d_{\alpha+1}}\\
&=\sigma_r^y[\overline d_{\alpha+1}d_{\alpha+1}u\overline y_r(\alpha+1)+d_{\alpha+1}]
\endalign
$$
where $u\in\bold Q[\overline y_r(\alpha+1)]$ is a polynomial with constant term 1. 
$\sigma_r^y\equiv 1\text{ mod }m(U(\alpha))$.
Thus  
$$
\sigma_r^y\equiv 1\text{ mod } (\overline y_1(\alpha+1),\ldots,\overline y_{r-1}(\alpha+1),
\overline y_{r+1}(\alpha+1),\ldots, \overline y_n(\alpha+1)).
$$
$$
\align
\frac{y_r^*(\alpha+1)}{\overline d_{\alpha+1}d_{\alpha+1}}&=\sigma_r^yu\overline y_r(\alpha+1)
+\frac{\sigma_r^y-1}{\overline d_{\alpha+1}}\\
&\equiv u\overline y_r(\alpha+1)\text{ mod } (\overline y_1(\alpha +1),\ldots,\overline y_{r-1}(\alpha+1),
\overline y_{r+1}(\alpha+1),\ldots,\overline y_n(\alpha+1)).
\endalign
$$
Thus  $U(\alpha+1) = S(\alpha+1)\,\,\hat{}\,$ since $(y_1^*(\alpha+1),\ldots,y_n^*(\alpha+1))$ are regular 
parameters in $U(\alpha+1)$ and $k(S(\alpha+1))\cong k(U(\alpha+1))$.
By Lemma 2.3, $T(\alpha+1)\rightarrow U(\alpha+1)$ induces a map $R(\alpha+1)\rightarrow S(\alpha+1)$.
Set 
$$
y_i(\alpha+1) = \cases
y_i^*(\alpha+1)&i\ne r\\
x_r(\alpha+1)&i=r
\endcases
$$
$$
\gamma_i^y(\alpha+1)=\cases
\sigma_i^y(\overline y_r(\alpha+1)+1)^{\overline d_{\alpha+1}\tau_i(\alpha+1)}& 1\le i\le s\\
\gamma_i(\alpha+1)&i=r\\
\gamma_i^y(\alpha)&s<i\le n, i\ne r
\endcases
$$
By definition
$$
\tilde{\overline y}_i(\alpha+1)=\cases
\overline y_i(\alpha+1)&i\ne r\\
\tilde{\overline x}_i(\alpha+1)&i=r
\endcases
$$
Then $(y_1(\alpha+1),\ldots,y_n(\alpha+1))$ are regular parameters in $S(\alpha+1)$ and 
$$
y_i(\alpha+1)=\gamma_i^y(\alpha+1)\tilde{\overline y}_i(\alpha+1)
$$
for $1\le i\le n$,
 $\gamma_i^y(\alpha+1)\in U''(\alpha+1)$ and 
$(\tilde{\overline y}_1(\alpha+1),\ldots,\tilde{\overline y}_n(\alpha+1))$ are regular parameters in $U''(\alpha+1)$.

\subheading{Case 3} Suppose that $\overline T(\alpha)=\overline T(\alpha+1)$ and 
$\overline U(\alpha)\rightarrow \overline U(\alpha+1)$ is of type $\text{II}_r$
with $l+1\le r\le m$.
By suitable choice of the  $\Omega_{i,\alpha}$ as in Case 1, we have regular parameters 
 $\tilde y_i(\alpha)$ in $S(\alpha)$ satisfying  (3.64). 
Set $(f_{ij}) = (b_{ij}(\alpha+1))^{-1}$.
Choose $\lambda_{\alpha}$ in (3.65) so that 
$$
 \gamma_1^y(\alpha)^{f_{s+1,1}}\cdots\gamma_s^y(\alpha)^{f_{s+1,s}}
\gamma_r^y(\alpha)^{f_{s+1,s+1}}\lambda_{\alpha}^{f_{s+1,s+1}}
\equiv 1\text{ mod } m(U(\alpha))^N.
$$
As in the argument of Case 2, we can define, by (3.69),
$$
S(\alpha+1) = S(\alpha)[B_1,\ldots,B_s,B_r]_{q_4}
$$
 so that $S(\alpha+1)\,\,\hat{}\,$ has regular parameters
$(y_1^*(\alpha+1),\ldots,y_n^*(\alpha+1))$ defined by (3.70).
Set $\sigma_i^y = \gamma_1^y(\alpha)^{f_{i1}}\cdots\gamma_s^y(\alpha)^{f_{is}}
\gamma_r^y(\alpha)^{f_{i,s+1}}\lambda_{\alpha}^{f_{i,s+1}}$ for $1\le i\le s$ and
set$$
\sigma_r = \gamma_1^y(\alpha)^{f_{s+1,1}}\cdots\gamma_s^y(\alpha)^{f_{s+1,s}}
\gamma_r^y(\alpha)^{f_{s+1,s+1}}\lambda_{\alpha}^{f_{s+1,s+1}}.
$$

Then equations (3.72) hold, and we have an inclusion $S(\alpha+1)\rightarrow U'(\alpha+1)$ as in the argument of Case 2.

By (3.70), (3.72), (3.55) and (3.57),
$$
\align
y_i^*(\alpha+1)&= 
\sigma_i \hat{\overline y}_i(\alpha+1)\\
&=\sigma_i(\overline y_r(\alpha+1)+1)^{\overline d_{\alpha+1}\tau_i}\overline y_i(\alpha+1)
\endalign
$$
 for $1\le i\le s$ and

$$
\align
y^*_r(\alpha+1)+d_{\alpha+1} 
&=\sigma_r (\hat{\overline y}_r(\alpha+1)+d_{\alpha+1})\\
&=\sigma_r^yd_{\alpha+1}(\overline y_r(\alpha+1)+1)^{\overline d_{\alpha+1}}\\
&=\sigma_r^y[\overline d_{\alpha+1}d_{\alpha+1}u\overline y_r(\alpha+1)+d_{\alpha+1}]
\endalign
$$
where $u\in\bold Q[\overline y_r(\alpha+1)]$ is a polynomial with constant term 1.
$$
\frac{y_r^*(\alpha+1)}{\overline d_{\alpha+1}d_{\alpha+1}}=\sigma_r^yu\overline y_r(\alpha+1)
+\frac{\sigma_r^y-1}{\overline d_{\alpha+1}}
$$
Recall that 
$\sigma_r^y\equiv 1\text{ mod }m(U(\alpha))^N$. Thus $U(\alpha+1)=S(\alpha+1)\,\,\hat{}\,$ since\linebreak
$(y_1^*(\alpha+1),\ldots,y_n^*(\alpha+1))$ are regular parameters in $U(\alpha+1)$ and $k(S(\alpha+1))\cong k(U(\alpha+1))$.

Thus there exists 
$$
\Omega\in (\overline y_1(\alpha+1),\ldots,\overline y_n(\alpha+1))^NU(\alpha+1)
 = (\overline y_1(\alpha+1),\ldots,\overline y_n(\alpha+1))^N
S(\alpha+1)\,\,\hat{}
$$
 such that 
$$
\sigma_ru\overline y_r(\alpha+1)+\Omega\in S(\alpha+1).
$$
Set
$$
y_i(\alpha+1)=\cases
y_i^*(\alpha+1) &i\ne r\\
\sigma_ru\overline y_r(\alpha+1)+\Omega&i=r
\endcases
$$
\vskip .2truein
$$
\gamma_i^y(\alpha+1) =\cases
\sigma_i(\overline y_r(\alpha+1)+1)^{\overline d_{\alpha+1}\tau_i}&1\le i\le s\\
\sigma_ru&i=r\\
\gamma_i^y(\alpha)&s<i\le n, i\ne r
\endcases
$$
Set $\Lambda_{\alpha+1}=\sigma_r^{-1}u^{-1}\Omega$. By definition
$$
\tilde{\overline y}_i(\alpha+1)=\cases
\overline y_i(\alpha+1)&i\ne r\\
\overline y_r(\alpha+1)+\Lambda_{\alpha+1}&i=r
\endcases
$$
Then $( y_1(\alpha+1),\ldots,y_n(\alpha+1))$ are regular parameters in $S(\alpha+1)$ and
$$
y_i(\alpha+1)=\gamma_i^y(\alpha+1)\tilde{\overline y}_i(\alpha+1)
$$
 for  $1\le i\le n$. $\gamma_i^y(\alpha+1)\in U''(\alpha+1)$ for all $i$ and
 $(\tilde{\overline y}_1(\alpha+1),\ldots,\tilde{\overline y}_n(\alpha+1))$
are regular  parameters in $U''(\alpha+1)$.

\subheading{Step 4}
It remains to show that the CRUTS (3.37), (3.61); $(R,T''(t),T(t))$ and $(S,U''(t),U(t))$,
constructed in step 3 is a CRUTS along $\nu$ if $N$ is sufficiently large.

We have an extension of $\nu$ to the quotient field of $\overline U(t)$ which dominates $\overline U(t)$.
Define 
$\tilde U(t)=
\overline U(t)/B(t)$ where $B(t)$ is the prime ideal of elements of $\overline U(t)$ of infinite value.
Let $G(t)$ be the quotient field of $\tilde U(t)$.
Let 
$\overline K$ be the completion of $K$ with respect to a metric associated to $\nu$ (c.f Lemma 1.2),
$\overline G(t)$ be the completion of $G(t)$ with respect to a metric  associated to $\nu$.
We have a natural inclusion of complete fields $\overline K\rightarrow
\overline G(t)$.  
Suppose that for some $0\le \beta\le t$
$$
\matrix
U(0)&\rightarrow &U(1)&\rightarrow &\cdots &\rightarrow &U(\beta)\\
\uparrow&&\uparrow   &            &        &           &\uparrow\\
T(0)&\rightarrow& T(1)&\rightarrow &\cdots&\rightarrow &T(\beta)
\endmatrix
$$
is a CRUTS along $\nu$.
Then by Lemma 1.2, we have natural maps:
$$
U(i)=\hat S(i)\rightarrow \overline K\text{ for $0\le i\le \beta$}.
$$
Let $\Cal O_w$ be the valuation ring of the natural  extension
$w$ of $\nu$ to $\overline G(t)$,  $m_w$ be the maximal ideal of $\Cal O_w$ and $\Gamma_w$ be the value group of $w$.
We have an inclusion $k(\overline U(t))\rightarrow \Cal O_w/m_w$.

We will prove the following inductive statement on $\alpha$ with $0\le\alpha\le t$.
 Given a positive element $\lambda_{\alpha}'\in \Gamma_w$, such that 
$$
\lambda_{\alpha}'>\text{max}\{w(\overline w_1(\alpha)),\ldots,w(\overline w_n(\alpha))\},\tag 3.73
$$
 there exists a positive element $N_{\alpha}$ such that if $N\ge N_{\alpha}$, and 
$$
\matrix
U(0)&\rightarrow &U(1)&\rightarrow &\cdots &\rightarrow &U(\alpha)\\
\uparrow&&\uparrow   &            &        &           &\uparrow\\
T(0)&\rightarrow& T(1)&\rightarrow &\cdots&\rightarrow &T(\alpha)
\endmatrix\tag 3.74
$$
is  a CRUTS (3.37) as constructed in Step 3,  then
\item{A1)}   $w(\overline y_i(\alpha))>0$ for $1\le i\le n$.
\item{A2)}
$$
\align
\overline y_i(\alpha)&= \overline w_i(\alpha)\text{ for } 1\le i\le s\\
w(\overline w_i(\alpha)-\overline y_i(\alpha))&>\lambda_{\alpha}'\text{ for } s+1\le i\le n
\endalign
$$
\item{A3)} $U(\alpha)\rightarrow \Cal O_w$ and there is a commutative diagram
$$
\matrix k(\overline U(\alpha)) &\rightarrow &\Cal O_w/m_w\\
\uparrow\eta_U^{\alpha}&\nearrow&\\
k(U(\alpha))&&
\endmatrix
$$

Since  $\overline y_i(0)=\overline w_i(0)$, for  $1\le i\le n$, and $U''(0)=\overline U''(0)$,  the statement is true for $m=0$.

Suppose the inductive statement is true for CRUTS of length $\alpha$, and for any given $\lambda_{\alpha}'\in \Gamma_w$
satisfying (3.73). 
We will prove it for sequences of length $\alpha+1$, and any given $\lambda'_{\alpha+1}$ such that 
$$
\lambda_{\alpha+1}'>\text{max}\{w(\overline w_1(\alpha+1)),\ldots,w(\overline w_n(\alpha+1))\}.
$$

Choose 
$\lambda'=\lambda_{\alpha+1}'$ if $\overline U(\alpha)\rightarrow \overline U(\alpha+1)$ is of type I,
$$
\lambda'>\lambda'_{\alpha+1}+b_{s+1,1}(\alpha+1)w(\overline w_1(\alpha+1))+\ldots+b_{s+1,s}(\alpha+1)w(\overline w_s(\alpha+1))
$$
if $\overline U(\alpha)\rightarrow\overline U(\alpha+1)$ is of type $\text{II}_r$.
By induction, there exists $N''$ such that $N\ge N''$ implies
$w(\overline w_i(\alpha)-\overline y_i(\alpha))>\lambda'$ for all $i$, and
 we can further choose $N''$ so large that 
$Nw(m(\overline U(\alpha)))>\lambda'$. Then $\nu(\overline y_i(\alpha))=\nu(\overline w_i(\alpha))$ for all $i$
and $Nw(m(U(\alpha)))>\lambda'$.

$\nu(\Lambda_{\alpha})>\lambda'$ implies  $w(\overline w_i(\alpha)-\tilde{\overline y}_i(\alpha))>\lambda'$ for $s+1\le i\le n$.
$\nu(\Omega_{i,\alpha})>\lambda'$ implies  $w(\tilde{\overline w}'_i(\alpha)-\tilde{\overline y}'_i(\alpha))
>\lambda'$ for $s+1\le i\le n$. 
 Thus there exists $\sigma_i\in\Cal O_w$ with $w(\sigma_i)>\lambda'$
such that 
$\tilde{\overline y}_i'(\alpha)=\tilde{\overline w}'_i(\alpha)+\sigma_i$ for $s+1\le i\le n$ and
$\tilde{\overline y}'_i(\alpha)=\tilde{\overline w}'_i(\alpha)$ for $1\le i\le s$.

Suppose that $\overline U(\alpha)\rightarrow \overline U(\alpha+1)$ is of type I.
$$
\tilde{\overline y}_i'(\alpha) =\prod_{j=1}^s\overline y_j(\alpha+1)^{b_{ij}(\alpha+1)}
$$
and
$$
\tilde{\overline w}_i'(\alpha) =\prod_{j=1}^s\overline w_j(\alpha+1)^{b_{ij}(\alpha+1)}
$$
for $1\le i\le s$.
Thus $\overline y_j(\alpha+1) = \overline w_{j+1}(\alpha+1)$ for $1\le j\le s$
and A1), A2) and A3) hold for $\alpha+1$.

Suppose that $\overline U(\alpha)\rightarrow \overline U(\alpha+1)$ is of type $\text{II}_r$. Set $b_{ij}=b_{ij}(\alpha+1)$.
Set $(f_{ij}) = (b_{ij}(\alpha+1))^{-1}$. $\overline U(\alpha+1)$ has regular parameters $(\hat{\overline w}_1(\alpha+1).
\ldots,\hat{\overline w}_n(\alpha+1))$ such that
$$
\align
\hat{\overline w}_1(\alpha+1) &= \tilde{\overline w}'_1(\alpha)^{f_{11}}\cdots\tilde{\overline w}'_s(\alpha)^{f_{1s}}
\tilde{\overline w}'_r(\alpha)^{f_{1,s+1}}\tag 3.75\\
&\vdots\\
\hat{\overline w}_s(\alpha+1) &= \tilde{\overline w}'_1(\alpha)^{f_{s1}}\cdots\tilde{\overline w}'_s(\alpha)^{f_{ss}}
\tilde{\overline w}'_r(\alpha)^{f_{s,s+1}}\\
\hat{\overline w}_r(\alpha+1)+d_{\alpha+1} 
&= \tilde{\overline w}'_1(\alpha)^{f_{s+1,1}}\cdots\tilde{\overline w}'_s(\alpha)^{f_{s+1,s}}
\tilde{\overline w}'_r(\alpha)^{f_{s+1,s+1}}
\endalign
$$
Recall equation (3.52).
$$
\align
N_1 &= \tilde{\overline y}'_1(\alpha)^{f_{11}}\cdots\tilde{\overline y}'_s(\alpha)^{f_{1s}}
\tilde{\overline y}'_r(\alpha)^{f_{1,s+1}}\tag 3.76\\
&\vdots\\
N_s &= \tilde{\overline y}'_1(\alpha)^{f_{s1}}\cdots\tilde{\overline y}'_s(\alpha)^{f_{ss}}
\tilde{\overline y}'_r(\alpha)^{f_{s,s+1}}\\
N_r 
&= \tilde{\overline y}'_1(\alpha)^{f_{s+1,1}}\cdots\tilde{\overline y}'_s(\alpha)^{f_{s+1,s}}
\tilde{\overline y}'_r(\alpha)^{f_{s+1,s+1}}
\endalign
$$
There is a natural map
$$
U''(\alpha)[N_1,\ldots,N_s,N_r]\rightarrow \overline K\rightarrow G(t).
$$
From (3.75) and (3.76) we have
$$
\frac{N_i}{\hat{\overline w}_i(\alpha+1)}=\left(\frac{\tilde{\overline y}_r'(\alpha)}{\tilde{\overline w}_r'(\alpha)}\right)
^{f_{i,s+1}}\in\Cal O_w\tag 3.77
$$
for $1\le i\le s$ and 
$$
\frac{N_r}{\hat{\overline w}_r(\alpha+1)+d_{\alpha+1}}
=\left(\frac{\tilde{\overline y}_r'(\alpha)}{\tilde{\overline w}_r'(\alpha)}\right)
^{f_{s+1,s+1}}\in\Cal O_w\tag 3.78
$$
All of these ratios have residue 1 in $\Cal O_w/m_w$. Thus $U''(\alpha)[N_1,\ldots,N_s,N_r]\rightarrow\Cal O_w$.
Then since $k(\overline U(\alpha+1))\subset\Cal O_w$, and we have an inclusion 
$$
k(U(\alpha))\atop{\eta_{\alpha}}{\rightarrow}
k(\overline U(\alpha))\rightarrow k(\overline U(\alpha+1))\rightarrow \Cal O_w/m_w.
$$
There is a natural map
$$
U''(\alpha)[N_1,\ldots,N_s,N_r]\otimes_{k_2}k(\overline U(\alpha+1))\rightarrow \Cal O_w
$$
where $k_2$ is the integral closure of $k$ in $U(\alpha)$.

by (3.77) and (3.78), $\nu(N_i)=\nu(\hat{\overline w}_i(\alpha+1))>0$ for $1\le i\le s$ and
$N_r \mapsto u(\hat{\overline w}_r(\alpha+1)+d_{\alpha+1})$ where $u\in\Cal O_w$ is a unit, $u\equiv 1\text{ mod }m_w$. 
Thus the residue of $N_r$ in   $\Cal O_w/m_w$ is $d_{\alpha+1}$.
Thus $w(N_r-d_{\alpha+1})>0$ and we have a map $U'(\alpha+1)\rightarrow \Cal O_w$ by (3.53) and (3.54), which induces a map
$U(\alpha+1)\rightarrow \Cal O_w$ such that
$$
\matrix
k(\overline U(\alpha+1))&\rightarrow &k(\Cal O_w)\\
\uparrow \eta^U_{\alpha+1}&\nearrow&\\
k(U(\alpha+1))&&
\endmatrix
$$
commutes, verifying A3) for $\alpha+1$.

The regular parameters $(\hat{\overline y}_1(\alpha+1),\ldots,\hat{\overline y}_n(\alpha+1))$ in $U(\alpha+1)$
defined by (3.55) satisfy
$$
\hat{\overline y}_i(\alpha+1) =\cases
N_i&1\le i\le s \\
N_r-d_{\alpha+1}&i=r\\
\tilde{\overline y}_i'(\alpha)&s<i,i\ne r
\endcases\tag 3.79
$$

$\overline U(\alpha+1)$ has regular parameters $(\overline w_1(\alpha+1),\ldots,\overline w_n(\alpha+1))$
defined by
 $$
\overline w_i(\alpha+1)=\cases
\hat{\overline w}_i(\alpha+1)\left(\frac{\hat{\overline w}_r(\alpha+1)}{d_{\alpha+1}}+1\right)^{-\tau_i(\alpha+1)}&1\le i\le s\\
 \left(\frac{\hat{\overline w}_r(\alpha+1)}{d_{\alpha+1}}+1\right)^{\frac{1}{\overline d_{\alpha+1}}}-1&i=r\\
\hat{\overline w}_i(\alpha+1)&s<i,i\ne r
\endcases\tag 3.80
$$                                            
 The regular parameters $(\overline y_1(\alpha+1),\ldots,\overline y_n(\alpha+1))$ of $U(\alpha+1)$ are defined by (3.57),
$$
\overline y_i(\alpha+1)=\cases
\hat{\overline y}_i(\alpha+1)\left(\frac{\hat{\overline y}_r(\alpha+1)}{d_{\alpha+1}}+1\right)^{-\tau_i(\alpha+1)}&1\le i\le s\\
 \left(\frac{\hat{\overline y}_r(\alpha+1)}{d_{\alpha+1}}+1\right)^{\frac{1}{\overline d_{\alpha+1}}}-1&i=r\\
\hat{\overline y}_i(\alpha+1)&s<i,i\ne r
\endcases\tag 3.81
$$

$$
\align
\tilde{\overline y}_1'(\alpha)&= \overline y_1(\alpha+1)^{b_{11}(\alpha+1)}\cdots\overline y_s(\alpha+1)^{b_{1s}(\alpha+1)}
d_{\alpha+1}^{b_{1,s+1}(\alpha+1)}\tag 3.82\\
&\vdots\\
\tilde{\overline y}_s'(\alpha)&= \overline y_1(\alpha+1)^{b_{s1}(\alpha+1)}\cdots\overline y_s(\alpha+1)^{b_{ss}(\alpha+1)}
d_{\alpha+1}^{b_{s,s+1}(\alpha+1)}\\
\tilde{\overline y}_r'(\alpha)&= \overline y_1(\alpha+1)^{b_{s+1,1}(\alpha+1)}\cdots\overline y_s(\alpha+1)^{b_{s+1,s}(\alpha+1)}
(\overline y_r(\alpha+1)+1)d_{\alpha+1}^{b_{s+1,s+1}(\alpha+1)}
\endalign
$$
\vskip .2truein
$$
\align
\tilde{\overline w}_1'(\alpha)&= \overline w_1(\alpha+1)^{b_{11}(\alpha+1)}\cdots\overline w_s(\alpha+1)^{b_{1s}(\alpha+1)}
d_{\alpha+1}^{b_{1,s+1}(\alpha+1)}\tag 3.83\\
&\vdots\\
\tilde{\overline w}_s'(\alpha)&= \overline w_1(\alpha+1)^{b_{s1}(\alpha+1)}\cdots\overline w_s(\alpha+1)^{b_{ss}(\alpha+1)}
d_{\alpha+1}^{b_{s,s+1}(\alpha+1)}\\
\tilde{\overline w}_r'(\alpha)&= \overline w_1(\alpha+1)^{b_{s+1,1}(\alpha+1)}\cdots\overline w_s(\alpha+1)^{b_{s+1,s}(\alpha+1)}
(\overline w_r(\alpha+1)+1)d_{\alpha+1}^{b_{s+1,s+1}(\alpha+1)}
\endalign
$$
Comparing (3.82) and (3.83), we see that for $1\le i\le s$
$\overline y_i(\alpha+1)=\lambda_i\overline w_i(\alpha+1)$ for some $\overline d_{\alpha+1}$-th
 roots of unity $\lambda_i$. By our construction of UTSs (at the beginning of Chapter 2) 
$\left(\frac{\hat{\overline y}_r(\alpha+1)}{d_{\alpha+1}}+1\right)^{\frac{1}{\overline d_{\alpha+1}}}$ and
$\left(\frac{\hat{\overline w}_r(\alpha+1)}{d_{\alpha+1}}+1\right)^{\frac{1}{\overline d_{\alpha+1}}}$ have residue 1
in $\Cal O_w/m_w$. By (3.77), (3.79), (3.80) and (3.81) we see that $\frac{\overline y_i(\alpha+1)}{\overline w_i(\alpha+1)}$
 has residue 1 in $\Cal O_w/m_w$ for $1\le i\le s$.
Thus  we have $\overline y_i(\alpha+1)=\overline w_i(\alpha+1)$ for $1\le i\le s$, proving the first half of 
A2) for $\alpha+1$.
$$
\align
\overline w_r(\alpha+1)-\overline y_r(\alpha+1)&=
\frac{\tilde{\overline w}_r'(\alpha)}{d_{\alpha+1}^{b_{s+1,s+1}}
\overline w_1(\alpha+1)^{b_{s+1,1}}\cdots\overline w_s(\alpha+1)^{b_{s+1,s}}}\\
&-\frac{\tilde{\overline y}_r'(\alpha)}{d_{\alpha+1}^{b_{s+1,s+1}}
\overline y_1(\alpha+1)^{b_{s+1,1}}\cdots\overline y_s(\alpha+1)^{b_{s+1,s}}}\\
&=\frac{\tilde{\overline w}_r'(\alpha)-\tilde{\overline y}_r'(\alpha)}{d_{\alpha+1}^{b_{s+1,s+1}}
\overline w_1(\alpha+1)^{b_{s+1,1}}\cdots\overline w_s(\alpha+1)^{b_{s+1,s}}}
\endalign
$$
Thus
$$
\align
w(\overline w_r(\alpha+1)-\overline y_r(\alpha+1))&=w(\tilde{\overline w}_r'(\alpha)-\tilde{\overline y}_r'(\alpha))
-b_{s+1,1}w(\overline w_1(\alpha+1))-\cdots-b_{s+1,s}w(\overline w_s(\alpha+1))\\
&>\lambda'-b_{s+1,1}w(\overline w_1(\alpha+1))-\cdots-b_{s+1,s}w(\overline w_s(\alpha+1))\\
&>\lambda_{\alpha+1}'
\endalign
$$
verifying the second half of A2), and A1) for $\alpha+1$.
\enddemo

\proclaim{Theorem 3.10}
Suppose that $ T''(0)\subset \hat R$ is a regular local ring essentially of finite type over $R$
such that the quotient field of $T''(0)$ is finite over $J$, $U''(0)\subset \hat S$ is a
regular local ring essentially of finite type over $S$
such that the quotient field of $U''(0)$ is finite over $K$, $T''(0)\subset  U''(0)$,
$T''(0)$ contains a subfield isomorphic to $k(c_0)$ for some $c_0\in k(T''(0)$ and $U''(0)$
contains a subfield isomorphic to $k(U''(0))$. 
Suppose that  $R$ has regular parameters $(x_1,\ldots, x_n)$, $S$ has regular parameters $(y_1,\ldots, y_n)$,
$T''(0)$ has regular parameters $(\overline x_1,\ldots,\overline x_n)$
 and $ U''(0)$ has regular parameters $(\overline y_1,\ldots, \overline y_n)$
such that
$$
\align
\overline x_1 &= \overline y_1^{c_{11}}\cdots \overline y_s^{c_{1s}}\phi_1\\
\vdots&\\
 \overline x_s &= \overline y_1^{c_{s1}}\cdots \overline y_s^{c_{ss}}\phi_s\\
\overline  x_{s+1} &= \overline y_{s+1}\\
\vdots&\\
\overline  x_l&=\overline y_l
\endalign
$$
where $\phi_1,\ldots,\phi_s\in k(U''(0))$, $\nu(\overline x_1),\ldots,\nu(\overline x_s)$
are rationally independent, $\text{det}(c_{ij})\ne 0$. Suppose that there exists a regular local ring
$\tilde R\subset R$ such that $(x_1,\ldots,x_l)$ are regular parameters in $\tilde R$,  $k(\tilde R)\cong k(c_0)$ 
and 
$$
x_i=\cases
\gamma_i\overline x_i&1\le i\le l\\
\overline x_i&l<i\le n
\endcases
$$
 with 
$\gamma_i\in k(c_0)[[x_1,\ldots,x_l]]\cap T''(0)$ for $1\le i\le l$ and
$\gamma_i\equiv 1\text{ mod } ( x_1,\ldots, x_l)$, there exist
$\gamma_i^y\in U''(0)$ such that $y_i=\gamma_i^y\overline y_i$, $\gamma_i^y\equiv 1\text{ mod }m(U''(0))$ for $1\le i\le n$.
Further suppose that 
$$
x_{l+1} = \overline P(\overline y_1,\ldots,\overline y_l)+
   \overline y_1^{d_1}\cdots\overline y_s^{d_s}\overline H+\Omega
$$
where $\overline P$ is a power series with coefficients in $k(U''(0))$, 
$$
\overline H = \overline u(\overline y_{l+1}+\overline y_1^{\overline g_1}\cdots \overline y_s^{\overline g_s}\overline\Sigma)
$$
where $\overline u\in  U''(0)\,\,\hat{}\,$ is a unit, 
$\overline\Sigma\in k(U''(0))[[\overline y_1,\ldots,\overline y_l,\overline y_{l+2},\ldots,\overline y_n]]$,
 $\nu(\overline y_{l+1})\le \nu(\overline y_1^{\overline g_1}\cdots\overline y_s^{\overline g_s})$ and
$\Omega\in m(U(0))^N$ with $N\nu(m(U(0)))>\nu(\overline y_1^{d_1}\cdots\overline y_s^{d_s}\overline y_{l+1})$.
Then there exists a CRUTS along $\nu$ $(R, T''(t'), T(t'))$ and $(S, U''(t'), U(t'))$ with associated MTS  
$$
\matrix
S & \rightarrow & S(t')\\
\uparrow &&\uparrow\\
R &\rightarrow & R(t')
\endmatrix
$$
such that $T''(t')$ contains a subfield isomorphic to $k(c_0,\ldots,c_{t'})$,
$U''(t')$ contains a subfield isomorphic to $k(U(t'))$,
$R(t')$ has regular parameters $(x_1(t'),\ldots,x_n(t'))$, $ T''(t')$ has regular parameters
$(\tilde{\overline x}_1(t'),\ldots \tilde{\overline x}_n(t'))$,
$S(t')$ has regular parameters $(y_1(t'),\ldots, y_n(t'))$,
 $ U''(t')$ has regular parameters $(\tilde{\overline y}_1(t'),\ldots, \tilde{\overline y}_n(t'))$
where
$$
\align
\tilde{\overline x}_1(t') &= \tilde{\overline y}_1(t')^{c_{11}(t')}\cdots \tilde{\overline y}_s(t')^{c_{1s}(t')}\phi_1(t')\\
\vdots&\\
 \tilde{\overline x}_s(t') &= \tilde{\overline y}_1(t')^{c_{s1}(t')}\cdots \tilde{\overline y}_s(t')^{c_{ss}(t')}\phi_s(t')\\
\tilde{\overline x}_{s+1}(t') &= \tilde{\overline y}_{s+1}(t')\\
\vdots&\\
\tilde{\overline x}_l(t')&=\tilde{\overline y}_l(t')\\
 x_{l+1}(t')=\tilde{\overline x}_{l+1}(t') &= P(\tilde{\overline y}_1(t'),\ldots,\tilde{\overline y}_l(t'))+
   \tilde{\overline y}_1(t')^{d_1(t')}\cdots\tilde{\overline y}_s(t')^{d_s(t')}H
\endalign
$$
where $P,H$ are power series with coefficients in $k(U(t'))$, 
with 
$$
\text{mult}(H(0,\ldots,0,\tilde{\overline y}_{l+1}(t'),0,\ldots,0)=1,
$$
$\phi_1(t'),\ldots,\phi_s(t')\in k(U(t'))$, $\nu(\tilde{\overline x}_1(t')),\ldots,\nu(\tilde{\overline x}_s(t'))$
are rationally independent,\linebreak
 $\text{det}(c_{ij}(t'))\ne 0$. There exists a regular local ring $\tilde R(t')\subset R(t')$ 
such that $(x_1(t'),\ldots,x_l(t'))$ are regular parameters in $\tilde R(t')$ and $k(\tilde R(t'))\cong k(c_0,\ldots,
c_{t'})$.
$$
x_i(t')=\cases
\gamma_i(t')\tilde{\overline x}_i(t')&1\le i\le l\\
\tilde{\overline x}_i(t')&l<i\le n
\endcases
$$
 with 
$\gamma_i(t')\in k(c_0,\ldots,c_{t'})[[x_1(t'),\ldots,x_l(t')]]\cap T''(t')$ units 
 for $1\le i\le l$, such that
$\gamma_i(t')\equiv 1\text{ mod }(x_1(t'),\ldots,x_l(t'))$ and for $1\le i\le n$
there exists $\gamma_i^y(t')\in U''(t')$ such that $y_i(t')=\gamma_i^y(t')\tilde{\overline y}_i(t')$,
$\gamma_i^y(t')\equiv 1\text{ mod }m(U''(t'))$.
\endproclaim

\demo{Proof}
Perform the following sequence of CUTSs of type M2) for $0\le r\le l-s-1$
$$
\align
\overline x_1(r) &= \overline x_1(r+1)^{a_{11}(r+1)}\cdots\overline x_s(r+1)^{a_{1s}(r+1)}
c_{r+1}^{a_{1,s+1}(r+1)}\\
\vdots&\\
\overline x_s(r) &= \overline x_1(r+1)^{a_{s1}(r+1)}\cdots\overline x_s(r+1)^{a_{ss}(r+1)}
c_{r+1}^{a_{s,s+1}(r+1)}\\
\overline x_{s+r+1}(r) &= \overline x_1(r+1)^{a_{s+1,1}(r+1)}\cdots\overline x_s(r+1)^{a_{s+1,s}(r+1)}
(\overline x_{s+r+1}(r+1)+1)c_{r+1}^{a_{s+1,s+1}(r+1)}
\endalign
$$
$$
\align
\overline y_1(r) &= \overline y_1(r+1)^{b_{11}(r+1)}\cdots\overline y_s(r+1)^{b_{1s}(r+1)}
d_{r+1}^{b_{1,s+1}(r+1)}\\
\vdots&\\
\overline y_s(r) &= \overline y_1(r+1)^{b_{s1}(r+1)}\cdots\overline y_s(r+1)^{b_{ss}(r+1)}
d_{r+1}^{b_{s,s+1}(r+1)}\\
\overline y_{s+r+1}(r) &= \overline y_1(r+1)^{b_{s+1,1}(r+1)}\cdots\overline y_s(r+1)^{b_{s+1,s}(r+1)}
(\overline y_{s+r+1}(r+1)+1)d_{r+1}^{b_{s+1,s+1}(r+1)},
\endalign
$$
the sequence of CUTSs of type M3) for $l-s\le r\le n-s-1$
$$
\align
\overline y_1(r) &= \overline y_1(r+1)^{b_{11}(r+1)}\cdots\overline y_s(r+1)^{b_{1s}(r+1)}
d_{r+1}^{b_{1,s+1}(r+1)}\\
\vdots&\\
\overline y_s(r) &= \overline y_1(r+1)^{b_{s1}(r+1)}\cdots\overline y_s(r+1)^{b_{ss}(r+1)}
d_{r+1}^{b_{s,s+1}(r+1)}\\
\overline y_{s+r+1}(r) &= \overline y_1(r+1)^{b_{s+1,1}(r+1)}\cdots\overline y_s(r+1)^{b_{s+1,s}(r+1)}
(\overline y_{s+r+1}(r+1)+1)d_{r+1}^{b_{s+1,s+1}(r+1)},
\endalign
$$
followed by a CUTS of type M1), with $t'=n-s+1$
$$
\align
\overline x_1(t'-1) &= \overline x_1(t')^{a_{11}(t')}\cdots\overline x_s(t')^{a_{1s}(t')}\\
\vdots&\\
\overline x_s(t'-1) &= \overline x_1(t')^{a_{s1}(t')}\cdots\overline x_s(t')^{a_{ss}(t')}\\
\endalign
$$
$$
\align
\overline y_1(t'-1) &= \overline y_1(t')^{b_{11}(t')}\cdots\overline y_s(t')^{b_{1s}(t')}\\
\vdots&\\
\overline y_s(t'-1) &= \overline y_1(t')^{b_{s1}(t')}\cdots\overline y_s(t')^{b_{ss}(t')},
\endalign
$$
so that  
$$
\align
x_{l+1} &= P(\overline y_1(t'),\ldots,\overline y_l(t'))\\
&+\overline y_1(t')^{d_1(t')}\cdots\overline y_s(t')^{d_s(t')}\overline u
(\overline y_{l+1}(t')+1)\lambda+\Sigma)\\
&+\overline y_1(t')^{d_1(t')+1}\cdots\overline y_s(t')^{d_s(t')+1}\Psi
\endalign
$$
where 
$$
\Sigma\in k(U(t'))[[\overline y_1(t'),\ldots,\overline y_l(t'),\overline y_{l+2}(t'),\ldots,\overline y_n(t')]],
$$
$P(\overline y_1(t'),\ldots,\overline y_l(t'))= \overline P(\overline y_1,\ldots,\overline y_l)$,
$\overline y_1(t')\cdots\overline y_s(t')\mid (\overline u-\overline u(0,\ldots,0))$,
and
$$
\Omega=\overline y_1(t')^{d_1(t')+1}\cdots\overline y_s(t')^{d_s(t')+1}\Psi
$$
with $\lambda=d_{l-s+1}^{b_{s+1,s+1}(l-s+1)}\in k(U(t'))$, $\Psi\in U(t')$.
Then after  replacing  $P$ with
$$
P+\lambda\overline u(0,\ldots,0)\overline y_1(t')^{d_1(t')}\cdots\overline y_s(t')^{d_s(t')},
$$
 we can put $\overline x_{l+1}(t')= x_{l+1}$ in the desired form with 
$$
H = \overline u \,\overline y_{l+1}(t')\lambda+(\overline u-\overline u(0,\ldots,0))\lambda+
\overline u\Sigma+\overline y_1(t')\cdots \overline y_s(t')\Psi.
$$
The proof that
$$
\matrix U(0)&\rightarrow&U(t')\\
\uparrow&&\uparrow\\
T(0)&\rightarrow&T(t')
\endmatrix
$$
is a CRUTS is a simplification of the argument of step 3 in the proof of Theorem 3.9. We will give
an outline of the proof.

We can define MTSs $R\rightarrow R(t')$ and $S\rightarrow S(t')$ such that $R(r)\,\,\hat{}\,$, $S(r)\,\,\hat{}\,$
have respective regular parameters 
$(x_1^*(r),\ldots,x_n^*(r))$ and $(y_1^*(r),\ldots,y_n^*(r))$ for $0\le r\le n$.  For $0\le r\le l-s-1$
$$
\align
 x_1^*(r) &=  x_1^*(r+1)^{a_{11}(r+1)}\cdots x_s^*(r+1)^{a_{1s}(r+1)}(x_{s+r+1}^*(r+1)+c_{r+1})^{a_{1,s+1}(r+1)}\\
\vdots&\\
 x_s^*(r) &=  x_1^*(r+1)^{a_{s1}(r+1)}\cdots x_s^*(r+1)^{a_{ss}(r+1)}(x_{s+r+1}^*(r+1)+c_{r+1})^{a_{s,s+1}(r+1)}\\
 x_{s+r+1}^*(r) &=  x_1^*(r+1)^{a_{s+1,1}(r+1)}\cdots x_s^*(r+1)^{a_{s+1,s}(r+1)}(x_{s+r+1}^*(r+1)+c_{r+1})^{a_{s+1,s+1}(r+1)}
\endalign
$$
$$
\align
 y_1^*(r) &=  y_1^*(r+1)^{b_{11}(r+1)}\cdots y_s^*(r+1)^{b_{1s}(r+1)}
(y_{s+r+1}^*(r+1)+d_{r+1})^{b_{1,s+1}(r+1)}\\
\vdots&\\
 y_s^*(r) &=  y_1^*(r+1)^{b_{s1}(r+1)}\cdots y_s^*(r+1)^{b_{ss}(r+1)}
(y_{s+r+1}^*(r+1)+d_{r+1})^{b_{s,s+1}(r+1)}\\
 y_{s+r+1}^*(r) &=  y_1^*(r+1)^{b_{s+1,1}(r+1)}\cdots y_s^*(r+1)^{b_{s+1,s}(r+1)}
(y_{s+r+1}^*(r+1)+d_{r+1})^{b_{s+1,s+1}(r+1)},
\endalign
$$
For $l-s\le r\le n-s-1$
$$
\align
 y_1^*(r) &=  y_1^*(r+1)^{b_{11}(r+1)}\cdots y_s^*(r+1)^{b_{1s}(r+1)}
(y_{s+r+1}^*(r+1)+d_{r+1})^{b_{1,s+1}(r+1)}\\
\vdots&\\
 y_s^*(r) &=  y_1^*(r+1)^{b_{s1}(r+1)}\cdots y_s^*(r+1)^{b_{ss}(r+1)}
(y_{s+r+1}^*(r+1)+d_{r+1})^{b_{s,s+1}(r+1)}\\
 y_{s+r+1}^*(r) &=  y_1^*(r+1)^{b_{s+1,1}(r+1)}\cdots y_s^*(r+1)^{b_{s+1,s}(r+1)}
(y_{s+r+1}^*(r+1)+d_{r+1})^{b_{s+1,s+1}(r+1)},
\endalign
$$
followed by a MTS of type M1), with $t'=n-s+1$
$$
\align
 x_1^*(t'-1) &=  x_1^*(t')^{a_{11}(t')}\cdots x_s^*(t')^{a_{1s}(t')}\\
\vdots&\\
 x_s^*(t'-1) &=  x_1^*(t')^{a_{s1}(t')}\cdots x_s^*(t')^{a_{ss}(t')}\\
\endalign
$$
$$
\align
 y_1^*(t'-1) &=  y_1^*(t')^{b_{11}(t')}\cdots y_s^*(t')^{b_{1s}(t')}\\
\vdots&\\
 y_s^*(t'-1) &=  y_1^*(t')^{b_{s1}(t')}\cdots y_s^*(t')^{b_{ss}(t')}.\\
\endalign
$$
For $1\le r\le l-s$ we have
$$
x_i^*(r)=\cases
\lambda_i(r)\overline x_i(r)& 1\le i\le s\\
\lambda_{s+r}(r)\overline x_{s+r}(r)+\Phi_{s+r}&i=s+r\\
x_i^*(r-1)&s<i,i\ne s+r
\endcases
$$
where $\lambda_i(r),\Phi_{s+r}\in k(c_0,\ldots,c_r)[[\overline x_1(r),\ldots,\overline x_l(r)]]$,
the $\lambda_i(r)$ are units and 
$$
\Phi_{s+r}\in (\overline x_1(r),\ldots,\overline x_{s+r-1}(r),\overline x_{s+r+1}(r),\ldots,\overline x_l(r)).
$$
For $1\le r\le n-s$ we have
$$
y_i^*(r) =\cases \lambda_i^y(r)\overline y_i(r)&1\le i\le s\\
\lambda_{s+r}^y(r)\overline y_{s+r}(r)+\Phi^y_{s+r}&i=s+r\\
y_i^*(r-1)&s<i,i\ne s+r
\endcases
$$
where $\lambda_i^y(r),\Phi_{s+r}^y\in k(U(r))[[\overline y_1(r),\ldots,\overline y_{s+r-1}(1),\overline y_{s+r+1}(r),
\ldots,\overline y_n(r)]]$, the $\lambda_i^y(r)$ are units.
 
$R(t')$ has regular parameters $(x_1(t'),\ldots,x_n(t'))$ where
$$
x_i(t')=\cases
x_i^*(t')&1\le i\le s,l+1\le i\le n\\
\prod(x_i^*(t')+c_{i-s}-\sigma(c_{i-s}))&s+1\le i\le l
\endcases
$$
where the product is over the distinct conjugates $\sigma(c_{i-s})$ of $c_{i-s}$ in 
an algebraic closure $\overline k$ of $k$ over $k$.
For $0\le r\le l-s-1$ define
$$
\tilde R(r+1) = \tilde R(r)[x_1^*(r+1),\ldots,x_s^*(r+1),\tau(r+1)]_{q_{r+1}}
$$
where
$$
\tau(r+1) = \prod(x^*_{s+r+1}(r+1)+c_{r+1}-\sigma(c_{r+1}))
$$
is the product over the distinct conjugates $\sigma(c_{r+1})$ of $c_{r+1}$ over $k$,
$$
q_{r+1} = m(R(r+1))\cap \tilde R(r)[x_1^*(r+1),\ldots,x_s^*(r+1),\tau(r+1)].
$$
For $l-s\le r\le n-s+1$ define $\tilde R(r+1)=\tilde R(r)$.
Define 
$$
\tilde R(t')=\tilde R(n-s)[x_1^*(t'),\ldots,x_s^*(t')]_{q_{t'}}
$$
where
$$
q_{t'}=m(R(t'))\cap \tilde R(n-s)[x_1^*(t'),\ldots,x_s^*(t')].
$$

$S(t')$ has regular parameters $(y_1(t'),\ldots,y_n(t'))$ where
$$
y_i(t')=\cases
y_i^*(t')&1\le i\le s\\
x_i(t')&s+1\le i\le l\\
\prod(y_i^*(t')+d_{i-s}-\sigma(d_{i-s}))&l+1\le i\le n
\endcases
$$
where the product is over the distinct conjugates $\sigma(d_{i-s})$ of $d_{i-s}$ in an algebraic
closure $\overline k$ of $k$ over $k$.

Now set
$$
\tilde{\overline x}_i(t')=
\cases
x_i(t')&\text{ for }s+1\le i\le l\\
\overline x_i(t')&\text{ otherwise}
\endcases
$$
$$
\tilde{\overline y}_i(t')=
\cases
\tilde{\overline x}_i(t')&\text{ for }s+1\le i\le l\\
\overline y_i(t')&\text{ otherwise}.
\endcases
$$
We are then in the form
of the conclusions of the Theorem.
\enddemo
\vskip .4truein
The proof of Theorem 3.10 also proves the following Theorem.
\proclaim{Theorem 3.11}
Let $n_{0,l}=m\left(k(U''(0))[[\overline y_1,\ldots,\overline y_l]]\right)$ in the assumptions of Theorem 3.10.
\item{1)} If $\Omega\in n_{0,l}^N$
 in the assumptions of Theorem 3.10, then a sequence of MTSs of type M2) and a MTS of type M1) (so that the CRUTS
along $\nu$ is in the
first $l$ variables)
are sufficient to transform $x_{l+1}$ into the form of the conclusions of Theorem 3.10.
\item{2)} Suppose that
$$
g=\overline y_1^{d_1}\cdots\overline y_s^{d_s}u(\overline y_1,\ldots,\overline y_l)+\Omega
$$
where $u$ is a unit power series with coefficients in $k(U''(0))$  and 
$\Omega\in n_{0,l}^N$
 with $N\nu(n_{0,l})>\nu(\overline y_1^{d_1}\cdots\overline y_s^{d_s})$.
 Then a sequence of MTSs of type M2) and a MTS of type M1) (so that the CRUTS along $\nu$ is in the first $l$ variables)
are sufficient to transform $g$ into the form
$$
g=\overline y_1(t')^{d_1(t')}\cdots\overline y_s(t')^{d_s(t')}\overline u(\overline y_1(t'),\ldots,\overline y_l(t'))
$$
where $\overline u$ is a unit power series.
\item{3)} Suppose that
$$
g=\overline y_1^{d_1}\cdots\overline y_s^{d_s}\Sigma(\overline y_1,\ldots,\overline y_l)+\Omega
$$
where $\nu(\overline y_1^{d_1}\cdots\overline y_s^{d_s})>A$  and $\Omega\in n_{0,l}^N$ 
with $N\nu(n_{0,l})>\nu(\overline y_1^{d_1}\cdots\overline y_s^{d_s})$.
 Then a sequence of MTSs of type M2) and a MTS of type M1) (so that the CRUTS along $\nu$ is in the first $l$ variables)
are sufficient to transform $g$ into the form
$$
g=\overline y_1(t')^{d_1(t')}\cdots\overline y_s(t')^{d_s(t')}\overline \Sigma(\overline y_1(t'),\ldots,\overline y_l(t'))
$$
where  $\nu(\overline y_1(t')^{d_1(t')}\cdots\overline y_s(t')^{d_s(t')})>A$.
\item{4)} Suppose that
$$
g=\overline y_1^{d_1}\cdots\overline y_s^{d_s}u(\overline y_1,\ldots,\overline y_l)+\Omega
$$
where $u$ is a unit power series with coefficients in $k(U''(0))$  and 
$\Omega\in m(U''(0))^N$
 with $N\nu(m(U''(0)))>\nu(\overline y_1^{d_1}\cdots\overline y_s^{d_s})$.
 Then there exists A CRUTS along $\nu$ as in the conclusions of Theorem 3.10
such that
$$
g=\overline y_1(t')^{d_1(t')}\cdots\overline y_s(t')^{d_s(t')}\overline u(\overline y_1(t'),\ldots,\overline y_l(t'))
$$
where $\overline u$ is a unit power series.
\endproclaim

\proclaim{Theorem 3.12}
Suppose that $ T''(0)\subset \hat R$ is a regular local ring essentially of finite type over $R$
such that the quotient field of $T''(0)$ is finite over $J$, $ U''(0)\subset \hat S$ is a
regular local ring essentially of finite type over $S$
such that the quotient field of $ U''(0)$ is finite over $K$, $ T''(0)\subset U''(0)$,
$T''(0)$ contains a subfield isomorphic to $k(c_0)$ for some $c_0\in k(T''(0))$ and $U''(0)$ contains a subfield isomorphic to
$k(U''(0))$. 
Suppose that $R$ has regular parameters $(x_1,\ldots, x_n)$,
$S$ has regular parameters $(y_1,\ldots,y_n)$,
$ T''(0)$ has regular parameters $(\overline x_1,\ldots,\overline x_n)$
 and $ U''(0)$ has regular parameters $(\overline y_1,\ldots, \overline y_n)$
such that
$$
\align
\overline x_1 &= \overline y_1^{c_{11}}\cdots \overline y_s^{c_{1s}}\phi_1\\
\vdots&\\
 \overline x_s &= \overline y_1^{c_{s1}}\cdots \overline y_s^{c_{ss}}\phi_s\\
\overline  x_{s+1} &= \overline y_{s+1}\\
\vdots&\\
\overline  x_l&=\overline y_l
\endalign
$$
where $\phi_1,\ldots,\phi_s\in  k(U''(0))$, $\nu(\overline x_1),\ldots,\nu(\overline x_s)$
are rationally independent, $\text{det}(c_{ij})\ne 0$. Suppose that there exists a regular local ring
$\tilde R\subset R$ such that $(x_1,\ldots,x_l)$ are regular parameters in $\tilde R$,  $k(\tilde R)\cong k(c_0)$ 
and 
$$
x_i=\cases
\gamma_i\overline x_i&1\le i\le l\\
\overline x_i&l<i\le n
\endcases
$$
 with 
$\gamma_i\in k(c_0)[[x_1,\ldots,x_l]]\cap T''(0)$ for $1\le i\le l$
and $\gamma_i\equiv 1\text{ mod }( x_1,\ldots,x_l)$, there exist
$\gamma_i^y\in U''(0)$ such that $y_i=\gamma_i^y\overline y_i$, $\gamma_i^y\equiv 1\text{ mod }m(U''(0))$
for $1\le i\le n$.

Then there exists a CRUTS along $\nu$ $(R, T''(t), T(t))$ and $(S,U''(t),U(t))$ with associated MTSs
$$
\matrix
S & \rightarrow & S(t)\\
\uparrow &&\uparrow\\
R &\rightarrow & R(t)
\endmatrix
$$
such that $T''(t)$ contains a subfield isomorphic to $k(c_0,\ldots,c_t)$, $U''(t)$ contains a 
subfield isomorphic to $k(U(t))$,
$R(t')$ has regular parameters $(x_1(t),\ldots,x_n(t))$,
$S(t')$ has regular parameters $(y_1(t'),\ldots,y_n(t'))$,
 $T''(t)$ has regular parameters
$(\tilde{\overline x}_1(t),\ldots \tilde{\overline x}_n(t))$,
 $U''(t)$ has regular parameters $(\tilde{\overline y}_1(t),\ldots, \tilde{\overline y}_n(t))$
where
$$
\align
\tilde{\overline x}_1(t) &= \tilde{\overline y}_1(t)^{c_{11}(t)}\cdots \tilde{\overline y}_s(t)^{c_{1s}(t)}
\phi_1(t)\\
\vdots&\\
 \tilde{\overline x}_s(t) &= \tilde{\overline y}_1(t)^{c_{s1}(t)}\cdots \tilde{\overline y}_s(t)^{c_{ss}(t)}
\phi_s(t)\\
\tilde{\overline x}_{s+1}(t) &= \tilde{\overline y}_{s+1}(t)\\
\vdots&\\
\tilde{\overline x}_{l+1}(t)&=\tilde{\overline y}_{l+1}(t)\\
\endalign
$$
such that $\phi_1(t),\ldots,\phi_s(t)\in k(U(t))$, $\nu(\tilde{\overline x}_1(t)),\ldots,\nu(\tilde{\overline x}_s(t))$
are rationally independent,\linebreak
 $\text{det}(c_{ij}(t))\ne 0$. 
$$
x_i(t)=\tilde{\overline x}_i(t)\text{ for }1\le i\le n.
$$
 For $1\le i\le n$ there exist $\gamma_i^y(t)\in U''(t)$ such that 
$y_i(t)=\gamma_i^y(t)\tilde{\overline y}_i(t)$, $\gamma_i^y(t)\equiv 1\text{ mod }m(U''(t))$.
\endproclaim

\demo{Proof}
We will construct a CRUTS $(R,T''(t),T(t))$ and $(S,U''(t),U(t))$ along $\nu$ with associated MTS
$$
\matrix
S=&S(0)&\rightarrow&S(t)\\
&\uparrow&&\uparrow\\
R=&R(0)&\rightarrow&R(t).
\endmatrix
$$
We will say that $\text{CN}(\beta)$ holds  (with $0\le\beta\le t$) if 
$T''(\beta)$ contains a subfield isomorphic to $k(c_0,\ldots,c_{\beta})$, $U''(\beta)$ contains a subfield isomorphic to
$k(U(\beta))$,
$R(\beta)$ has regular parameters $(x_1(\beta),\ldots,x_n(\beta))$, $T''(\beta)$ has regular parameters
$(\tilde{\overline x}_1(\beta),\ldots,\tilde{\overline x}_n(\beta))$, 
$U''(\beta)$ has regular parameters
$(\tilde{\overline y}_1(\beta),\ldots,\tilde{\overline y}_n(\beta))$, 
such that 
$$
x_i(\beta)=\cases
\gamma_i(\beta)\tilde{\overline x}_i(\beta)&1\le i\le l\\
\tilde{\overline x}_i(\beta)&l<i\le n
\endcases
$$
 with 
$$
\gamma_i(\beta)\in k(c_0,\ldots,c_{\beta})[[x_1(\beta),\ldots,x_l(\beta)]]\cap T''(\beta)
,
$$
$\gamma_i(\beta)\equiv 1\text{ mod }(x_1(\beta),\ldots,x_{l}(\beta))$ for $1\le i\le l$
and for $1\le i\le n$ there exist $\gamma_i^y(\beta)\in U''(\beta)$ such that
$y_i(\beta)=\gamma_i^y(\beta)\tilde{\overline y}_i(\beta)$, $\gamma_i^y(\beta)\equiv 1\text{ mod }m(U''(\beta))$.
We must further have
$$
\align
\tilde{\overline x}_1(\beta) &= \tilde{\overline y}_1(\beta)^{c_{11}(\beta)}\cdots \tilde{\overline y}_s(\beta)^{c_{1s}(\beta)}
\phi_1(\beta)\\
\vdots&\\
 \tilde{\overline x}_s(\beta) &= \tilde{\overline y}_1(\beta)^{c_{s1}(\beta)}\cdots \tilde{\overline y}_s(\beta)^{c_{ss}(\beta)}
\phi_s(\beta)\\
\tilde{\overline x}_{s+1}(\beta) &= \tilde{\overline y}_{s+1}(\beta)\\
\vdots&\\
\tilde{\overline x}_{l}(\beta)&=\tilde{\overline y}_{l}(\beta)\\
\endalign
$$
 with $\phi_i(\beta)\in k(S(\beta))$, $\nu(\tilde{\overline x}_1(\beta)),\ldots,
\nu(\tilde{\overline x}_s(\beta))$ are rationally independent, $\text{det}(c_{ij}(\beta))\ne0$,
and there exists a regular local ring $\tilde R(\beta)\subset R(\beta)$ such that
$(x_1(\beta),\ldots,x_l(\beta))$  are regular parameters in $\tilde R(\beta)$ and $k(\tilde R(\beta))\cong k(c_0,\ldots,c_{\beta})$.

By Theorem 1.12  and Theorems 3.9 (with $f=x_{l+1}$) and 3.10 we may assume that
$$
x_{l+1}= \overline x_{l+1} = P(\overline y_1,\ldots,\overline y_l)+
   \overline y_1^{d_1}\cdots\overline y_s^{d_s}\Sigma_0.
$$
where $\Sigma_0$ is a series with $\text{mult}(\Sigma_0(0,\ldots,0,\overline y_{l+1},0,\ldots,0)=1$.
Suppose that $\nu(P)=\infty$ (this includes the case $P=0$).
 Then by Theorems 3.9 and 3.11,  we have a CRUTS along $\nu$ in the first $l$ variables  such that
$$
x_{l+1}= \tilde{\overline y}_1(t)^{d_1(t)}\cdots\tilde{\overline y}_s(t)^{d_s(t)}
\Sigma_t(\tilde{\overline y}_1(t),\ldots,\tilde{\overline y}_{n}(t))\tag 3.84
$$
$\text{ mult }(\Sigma_t(0,\ldots,0,\tilde{\overline y}_{l+1}(t),0\ldots,0)=1$, and $\text{CN}(t)$ holds.
 We can thus make  changes of variables, replacing 
$\tilde{\overline x}_i(t)$ with $x_i(t)$ (for $1\le i\le n$) and $\tilde{\overline y}_i(t)$ with
$\tilde{\overline x}_i(t)$ (for $s+1\le i\le l$) so that (3.84) holds and $\text{CN}(t)$ holds
with $\gamma_i(t)=1$ for $1\le i\le l$, and $y_i(t)=x_i(t)$ for $1\le i\le l$.

Suppose that $\nu(P)<\infty$.
Set $d =\text{Det}(c_{ij})$. Let $(f_{ij})$ be the adjoint matrix of $(c_{ij})$, so that $(f_{ij}/d)$ is
the inverse of $(c_{ij})$. Let $\omega$ be a primitive $d^{th}$ root of unity in and algebraic closure $\overline k$
of $k$. Set
$$
g'=\prod_{i_1,\ldots,i_s=1}^{d}\left(x_{l+1}-
P(\omega^{i_1}\overline y_1,\ldots,\omega^{i_s}\overline y_s,\overline y_{s+1},\ldots,\overline y_l)\right).
$$
$$
g'\in k(U''(0))[[\overline y_1^{d},\ldots,\overline y_s^{d},\overline y_{s+1},\ldots,\overline y_{l}]][x_{l+1}].
$$
Let $G$ be the Galois group of a Galois closure of $k(U''(0))$ over $k(c_0)$.
Since $x_{l+1}$ is analytically independent of $\overline y_1^d,\ldots,\overline y_s^d,\overline y_{s+1},\ldots,\overline y_l$
(by Theorem 1.12)
we can define $g=\prod_{\tau\in G}\tau(g')$ where $\tau$ acts on the coefficients of $g'$.
$$
g\in k(c_0)[[\overline y_1^{d},\ldots,\overline y_s^{d},\overline y_{s+1},\ldots,\overline y_{l}]][x_{l+1}].
$$
Since  $\overline y_i^d=\overline x_1^{f_{i1}}\cdots\overline x_s^{f_{is}}\phi_1^{-f_{i1}}\cdots
\phi_s^{-f_{is}}$  for $1\le i\le s$,
by Lemma 3.2 we can perform a MTS of type M1) to get
$g\in k(c_0)[[\overline x_1(1),\ldots, \overline x_l(1)]][x_{l+1}]$. 
$\Sigma_0$ is irreducible in $U(1)$ (since $\text{mult}(\Sigma_0(0,\ldots,0,\overline y_{l+1},0,\ldots,0))=1$) and
$\Sigma_0\mid g$ in $U(1)$. 
$$
g = \sum_{i=1}^d a_i \overline x_{l+1}(1)^{f_i} + \sum_j a_j \overline x_{l+1}(1)^{f_j}
+ \sum_k a_k \overline x_{l+1}(1)^{f_k}+\overline x_{l+1}(1)^{\tau}\Gamma\tag 3.85
$$
where $a_i,a_j,a_k$ are series in $\overline x_1(1),\ldots,\overline x_l(1)$ with coefficients in $k(c_0)$,
$m>0$ and the first sum consists of the terms of (finite) minimum value $\rho$.
$$
\rho<\infty\text{ since }\text{mult}(g(0,\ldots,\overline x_{l+1}(1)))<\infty.
$$
$\tau\nu(\overline x_{l+1}(1))>\rho$,
the second sum consists of  (finitely many) remaining terms of finite value,
the third sum consists of  (finitely many) terms of infinite value. 
Let $r = \text{mult}(g(0,\cdots,0,\overline x_{l+1}(1))$. $1\le r <\infty$.

Suppose that $r>1$.  By Theorems 3.9 and 3.11, we can construct a 
CRUTS in the first $l$ variables, with associated MTS 
$$
\matrix
S(1) & \rightarrow & S(\alpha)\\
\uparrow &&\uparrow\\
R(1) &\rightarrow & R(\alpha)
\endmatrix
$$
such that $\text{CN}(\alpha)$ holds,
to get
$$
P = \tilde{\overline y}_1(\alpha)^{g_1(\alpha)}\cdots \tilde{\overline y}_s(\alpha)^{g_s(\alpha)}
\overline P(\tilde{\overline y}_1(\alpha),\ldots,\tilde{\overline y}_l(\alpha))
$$
where $\overline P$ is a unit, and
$$
a_{\zeta} = \tilde{\overline x}_1(\alpha)^{e_{\zeta}^1}\cdots\tilde{\overline x}_s(\alpha)^{e_{\zeta}^s}
\overline a_{\zeta}(\overline x_1(\alpha),\ldots,\overline x_l(\alpha))\tag 3.86
$$
for  $\zeta=i,j,k$ where $\overline a_i, \overline a_j$ are units and 
$\nu(\tilde{\overline x}_1(\alpha)^{e_{k}^1}\cdots\tilde{\overline x}_s(\alpha)^{e_{k}^s})>\rho$  for all $i,j,k$.
We have an expression 
$$
x_{l+1}(\alpha)=x_{l+1}= P + \tilde{\overline y}_1(\alpha)^{d_1(\alpha)}\cdots
 \tilde{\overline y}_s(\alpha)^{d_s(\alpha)}\Sigma_{\alpha}(\tilde{\overline y}_1(\alpha),\ldots,\tilde{\overline y}_{n}(\alpha)).
$$
where $\text{mult}(\Sigma_{\alpha}(0,\ldots,0,\tilde{\overline y}_{l+1}(\alpha),0\ldots,0)=1$.
If $\nu(P)> \nu(\tilde{\overline y}_1(\alpha)^{d_1(\alpha)}\cdots \tilde{\overline y}_1(\alpha)^{d_s(\alpha)})$,
then after possibly performing a CRUTS of type M1), 
so that 
$$
 \tilde{\overline y}_1(\alpha)^{d_1(\alpha)}\cdots \tilde{\overline y}_s(\alpha)^{d_s(\alpha)}\mid
\tilde{\overline y}_1(\alpha)^{g_1(\alpha)}\cdots \tilde{\overline y}_s(\alpha)^{g_s(\alpha)}
$$
in $U(\alpha+1)$,  we can set
$$
\Sigma_{\alpha+1}=\Sigma_{\alpha} + \tilde{\overline y}_1(\alpha)^{g_1(\alpha)-d_1(\alpha)}\cdots
\tilde{\overline y}_s(\alpha)^{g_s(\alpha)-d_s(\alpha)}\overline P
$$
to get 
$$
x_{l+1} = \tilde{\overline y}_1(\alpha+1)^{d_1(\alpha+1)}\cdots \tilde{\overline y}_{s}(\alpha+1)^{d_s(\alpha+1)}
\Sigma_{\alpha+1}.
$$

Now suppose that 
$\nu(P)\le  \nu(\tilde{\overline y}_1(\alpha)^{d_1(\alpha)}\cdots \tilde{\overline y}_s(\alpha)^{d_s(\alpha)})$.
After possibly performing a CRUTS of type M1), 
so that 
$$
 \tilde{\overline y}_1(\alpha)^{g_1(\alpha)}\cdots \tilde{\overline y}_s(\alpha)^{g_1(\alpha)}\mid
\tilde{\overline y}_1(\alpha)^{d_1(\alpha)}\cdots \tilde{\overline y}_s(\alpha)^{d_s(\alpha)}
$$
in $U(\alpha+1)$, we have
$$
\align
x_{l+1}(\alpha+1)=x_{l+1}&=
\tilde{\overline y}_1(\alpha+1)^{d_1(\alpha+1)}\cdots \tilde{\overline y}_s(\alpha+1)^{d_s(\alpha+1)}(
\overline P\\
&+ \tilde{\overline y}_1(\alpha+1)^{\epsilon_1(\alpha+1)}\cdots
\tilde{\overline y}_s(\alpha+1)^{\epsilon_s(\alpha+1)}
\overline \Sigma_{\alpha+1}(\tilde{\overline y}_1(\alpha+1),\ldots,\tilde{\overline y}_{n}(\alpha+1))).
\endalign
$$
where 
$$\text{mult}(\overline \Sigma_{\alpha+1}(0,\ldots,0,\tilde{\overline y}_{l+1}(\alpha+1),0,\ldots,0))=1,
$$
$\overline P+ \tilde{\overline y}_1(\alpha+1)^{\epsilon_1(\alpha+1)}\cdots
\tilde{\overline y}_s(\alpha+1)^{\epsilon_s(\alpha+1)}\Sigma_{\alpha+1}$ is a unit and $\text{CN}(\alpha+1)$ holds.

Set $(e_{ij})=(c_{ij}(\alpha+1))^{-1}$, $d=\text{det}(c_{ij}(\alpha+1))$. We can replace $\tilde{\overline y}_i(\alpha+1)$ with
$\tilde{\overline y}_i(\alpha+1)\gamma_1(\alpha+1)^{e_{i1}}\cdots\gamma_s(\alpha+1)^{e_{is}}$ for $1\le i\le s$, replace
$\tilde{\overline y}_i(\alpha+1)$ with $\tilde{\overline y}_i(\alpha+1)\gamma_i(\alpha+1)$ for $s+1\le i\le l$ and 
replace $U''(\alpha+1)$ with
$U''(\alpha+1)[\gamma_1(\alpha+1)^{\frac{1}{d}},\ldots,\gamma_s(\alpha+1)^{\frac{1}{d}}]_q$
where $q=U(\alpha+1)\cap \left(U''(\alpha+1)[\gamma_1(\alpha+1)^{\frac{1}{d}},\ldots,\gamma_s(\alpha+1)^{\frac{1}{d}}]\right)$
to get
$$
\align
x_1(\alpha+1)&=\tilde {\overline y}_1(\alpha+1)^{c_{11}(\alpha+1)}\cdots\tilde{ \overline y}_s(\alpha+1)^{c_{1s}(\alpha+1)}
\phi_1(\alpha+1)\\
&\vdots\\
x_s(\alpha+1)&=\tilde {\overline y}_1(\alpha+1)^{c_{s1}(\alpha+1)}\cdots\tilde {\overline y}_s(\alpha+1)^{c_{ss}(\alpha+1)}
\phi_s(\alpha+1)\\
x_{s+1}(\alpha+1)&=\tilde {\overline y}_{s+1}(\alpha+1)\\
&\vdots\\
x_l(\alpha+1)&=\tilde{ \overline y}_l(\alpha+1)
\endalign
$$
and  have an expression in $U(\alpha+1)\cong k(U(\alpha+1))
[[\tilde{\overline y}_1(\alpha+1),\ldots,\tilde{\overline y}_n(\alpha+1)]]$
$$
\align
x_{l+1} &=
\tilde{\overline y}_1(\alpha+1)^{d_1(\alpha+1)}\cdots \tilde{\overline y}_s(\alpha+1)^{d_s(\alpha+1)}
\tau[(
\overline P_{\alpha+1}(\tilde{\overline y}_1(\alpha+1),\ldots,\tilde{\overline y}_l(\alpha+1))\\
&+\tilde{\overline y}_1(\alpha+1)^{\epsilon_1(\alpha+1)}\cdots
\tilde{\overline y}_s(\alpha+1)^{\epsilon_s(\alpha+1)}
\Sigma_{\alpha+1}(\tilde{\overline y}_1(\alpha+1),\ldots,\tilde{\overline y}_{n}(\alpha+1))]).
\endalign
$$
where $\text{mult}(\Sigma_{\alpha+1}(0,\ldots,0,\tilde{\overline y}_{l+1}(\alpha+1),0\ldots,0)=1$, $\overline P_{\alpha+1}$ is
a power series with constant term 1 and $\tau\in k(U(\alpha+1))$.
We have $\overline x_{l+1}(\alpha+1) = x_{l+1}$ and $\overline y_{l+1}(\alpha+1) = \overline y_{l+1}$.
By Lemma 3.6 we can construct  MTSs $R(\alpha+1)\rightarrow R(\alpha+2)$ of type $\text{II}_{l+1}$ and
$S(\alpha+1)\rightarrow S(\alpha+2)$ of type I such that $R(\alpha+2)\,\,\hat{}\,$ 
has regular parameters $(x_1^*(\alpha+2),\ldots,
x_n^*(\alpha+2))$ 
$$
\align
 x_1(\alpha+1) &= x_1^*(\alpha+2)^{a_{11}(\alpha+2)}\cdots  x_s^*(\alpha+2)^{a_{1s}(\alpha+2)}
(x_{l+1}^*(\alpha+2)+c_{\alpha+2})^{a_{1,s+1}(\alpha+2)}\tag 3.87\\
\vdots&\\
 x_s(\alpha+1) &= x_1^*(\alpha+2)^{a_{s1}(\alpha+2)}\cdots  x_s^*(\alpha+2)^{a_{ss}(\alpha+2)}
(x_{l+1}^*(\alpha+2)+c_{\alpha+2})^{a_{s,s+1}(\alpha+2)}\\
 x_{l+1}(\alpha+1) 
&=  x_1^*(\alpha+2)^{a_{s+1,1}(\alpha+2)}\cdots  x_s^*(\alpha+2)^{a_{s+1,s}(\alpha+2)}
(x_{l+1}^*(\alpha+2)+c_{\alpha+2})^{a_{s+1,s+1}(\alpha+2)}
\endalign
$$
$S(\alpha+2)\,\,\hat{}\,$ has regular parameters $(\overline y_1(\alpha+2),\ldots,\overline y_n(\alpha+2))$ 
$$
\align
\tilde{\overline y}_1(\alpha+1) &=\overline y_1(\alpha+2)^{b_{11}(\alpha+2)}\cdots \overline y_s(\alpha+2)^{b_{1s}(\alpha+2)}\\
\vdots&\\
\tilde{\overline y}_s(\alpha+1) &=\overline y_1(\alpha+2)^{b_{s1}(\alpha+2)}\cdots \overline y_s(\alpha+2)^{b_{ss}(\alpha+2)}
\endalign
$$
such that $R(\alpha+2)\subset S(\alpha+2)$.

Set
$$
\gamma = \overline P_{\alpha+1}+\tilde{\overline y}_{1}(\alpha+1)^{\epsilon_1(\alpha+1)}\cdots
\tilde{\overline y}_{s}(\alpha+1)^{\epsilon_s(\alpha+1)}\Sigma_{\alpha+1}\tag 3.88
$$
so that
$$
x_{l+1} = \tilde{\overline y}_{1}(\alpha+1)^{d_1(\alpha+1)}\cdots
\tilde{\overline y}_{s}(\alpha+1)^{d_s(\alpha+1)}\Sigma_{\alpha+1}\tau\gamma.
$$
This shows that $\gamma\in U''(\alpha+1)$.
Set $(e_{ij})=(a_{ij}(\alpha+2))^{-1}$. 
By construction there are positive integers $f_{ij}$
such that
$$
\align
x_1^*(\alpha+2)&=x_1(\alpha+1)^{e_{11}}\cdots x_s(\alpha+1)^{e_{1s}}x_{l+1}(\alpha+1)^{e_{1,s+1}}\\
&=\overline y_1(\alpha+2)^{f_{11}}\cdots\overline y_s(\alpha+2)^{f_{1s}}\gamma^{e_{1,s+1}}
\tau^{e_{1,s+1}}\phi_1(\alpha+1)^{e_{11}}\cdots\phi_s(\alpha+1)^{e_{1s}}\tag 3.89\\
&\vdots\\
x_s^*(\alpha+2)&=x_1(\alpha+1)^{e_{s1}}\cdots x_s(\alpha+1)^{e_{ss}}x_{l+1}(\alpha+1)^{e_{s,s+1}}\\
&=\overline y_1(\alpha+2)^{f_{s1}}\cdots\overline y_s(\alpha+2)^{f_{ss}}\gamma^{e_{s,s+1}}
\tau^{e_{s,s+1}}\phi_1(\alpha+1)^{e_{s1}}\cdots\phi_s(\alpha+1)^{e_{ss}}\\
x_{l+1}^*(\alpha+2)+c_{\alpha+2}&=x_1(\alpha+1)^{e_{s+1,1}}\cdots x_s(\alpha+1)^{e_{s+1,s}}x_{l+1}(\alpha+1)^{e_{s+1,s+1}}\\
&=\overline y_1(\alpha+2)^{f_{s+1,1}}\cdots\overline y_s(\alpha+2)^{f_{s+1,s}}\gamma^{e_{s+1,s+1}}
\tau^{e_{s+1,s+1}}\\
&\,\,\,\,\,\,\,\,\,\phi_1(\alpha+1)^{e_{s+1,1}}\cdots\phi_s(\alpha+1)^{e_{s+1,s}}
\endalign
$$
in $S(\alpha+2)\,\,\hat{}\,$.
$\nu(x_{s+1}^*(\alpha+2)+c_{\alpha+2})=0$ implies $f_{s+1,1}=\cdots=f_{s+1,s}=0$.
Substituting (3.88), we have
$$
\align
x_{l+1}^*(\alpha+2)+c_{\alpha+2} &=
Q_0(\overline y_1(\alpha+2),\ldots,\overline y_l(\alpha+2))\\
&+\overline y_1(\alpha+2)^{\epsilon_1(\alpha+2)}\cdots
\overline y_s(\alpha+2)^{\epsilon_s(\alpha+2)}\Lambda_0(
\overline y_1(\alpha+2),\ldots,\overline y_n(\alpha+2))
\endalign
$$
where $Q_0$ is a unit and
$\text{mult}\Lambda_{0}(0,\ldots,0,\overline y_{l+1}(\alpha+2),0\ldots,0)=1$. 
 Define $\alpha_i\in \bold Q$ by
$$
\left(\matrix \alpha_1\\ \vdots \\ \alpha_s\endmatrix\right) =
\left(\matrix f_{11}&\cdots&f_{1s}\\
\vdots&&\vdots\\
f_{s1}&\cdots&f_{ss}
\endmatrix\right)^{-1}
\left(\matrix -e_{1,s+1}\\ \vdots \\ -e_{s,s+1}\endmatrix\right)
$$
and set 
$$
\hat y_i(\alpha+2)=\cases
\gamma^{-\alpha_i}\overline y_i(\alpha+2)& 1\le i\le s\\
\overline y_i(\alpha+2)&s<i
\endcases
$$
 to get
$$
\align
x_1^*(\alpha+2) &=\hat y_1(\alpha+2)^{c_{11}(\alpha+2)}\ldots\hat y_s(\alpha+2)^{c_{1s}(\alpha+2)}\psi_1\tag 3.90\\
\vdots\\
x_s^*(\alpha+2) &=\hat y_1(\alpha+2)^{c_{s1}(\alpha+2)}\ldots\hat y_s(\alpha+2)^{c_{ss}(\alpha+2)}\psi_s.
\endalign
$$
$(\hat y_1(\alpha+2),\ldots,\hat y_n(\alpha+2))$ are regular parameters in $S(\alpha+2)\,\,\hat{}\,$, 
$\psi_1,\ldots,\psi_s\in k(S(\alpha+2))$.

There are unit power series $Q_i$
and power series $\Lambda_i$
 such that 
$$
\align
\gamma^{\alpha_i}&= Q_i(\overline y_1(\alpha+2),\ldots,\overline y_l(\alpha+2))\\
&+\overline y_1(\alpha+2)^{\epsilon_1(\alpha+2)}\cdots\overline y_s(\alpha+2)^{\epsilon_s(\alpha+2)}
\Lambda_i(\overline y_1(\alpha+2),\ldots,\overline y_{n}(\alpha+2))
\endalign
$$
in $S(\alpha+2)\,\,\hat{}\,$ for  $1\le i\le s$, where $\text{mult}(\Lambda_i(0,\ldots,0,\overline y_{l+1}(\alpha+2),0,\ldots,0))=1$,
$$
Q_i(0,\ldots,0)=1.
$$
$$
\align
\overline y_i(\alpha+2) &\equiv Q_i(\overline y_1(\alpha+2),\ldots,\overline y_l(\alpha+2))\hat y_i(\alpha+2)\\
&\text{ mod }(\hat y_1(\alpha+2)^{\epsilon_1(\alpha+2)}\cdots\hat y_s(\alpha+2)^{\epsilon_s(\alpha+2)}
\hat y_i(\alpha+2))
\endalign
$$
for $1\le i\le s$.
We will show that there exist unit power series $\Omega_i$ such that
$$
\align
\gamma^{\alpha_i}&\equiv \Omega_i(\hat y_1(\alpha+2),\ldots,\hat y_l(\alpha+2))\\
&+\overline y_1(\alpha+2)^{\epsilon_1(\alpha+2)}\cdots\overline y_s(\alpha+2)^{\epsilon_s(\alpha+2)}
\Lambda_i(\overline y_1(\alpha+2),\ldots,\overline y_{n}(\alpha+2))\\
&\text{ mod } \hat y_1(\alpha+2)^{\epsilon_1(\alpha+2)}\cdots\hat y_s(\alpha+2)^{\epsilon_s(\alpha+2)}(\hat y_1(\alpha+2),\ldots,
\hat y_l(\alpha+2)).
\endalign
$$
This follows from induction, since for any series
 $A(\overline y_1(\alpha+2),\ldots,\overline y_l(\alpha+2))$,
 there exist series $A_{i_1\cdots i_l}$  such that $\text{mult}(A_{i_1,\ldots,i_s})>0$ and
$$
\align
&A(\overline y_1(\alpha+2),\ldots,\overline y_l(\alpha+2))\equiv
A(\hat y_1(\alpha+2),\ldots,\hat y_{l}(\alpha+2))\\
&+\sum_{i_1+\cdots+i_l>0} A_{i_1\cdots i_l}(\overline y_1(\alpha+2),\ldots,\overline y_l(\alpha+2))
\hat y_1(\alpha+2)^{i_1}\cdots\hat y_l(\alpha+2)^{i_l}\\
&\text{ mod } \hat y_1(\alpha+2)^{\epsilon_1(\alpha+2)}\cdots\hat y_l(\alpha+2)^{\epsilon_s(\alpha+2)}
(\hat y_1(\alpha+2),\ldots,
\hat y_l(\alpha+2)).
\endalign
$$
Thus 
$$
\align
\gamma^{\alpha_i} &\equiv \Omega_i(\hat y_1(\alpha+2),\ldots,\hat y_l(\alpha+2))\\
&+\hat y_1(\alpha+2)^{\epsilon_1(\alpha+2)}\cdots \hat y_s(\alpha+2)^{\epsilon_s(\alpha+2)}\Phi_i(\hat y_{l+1}(\alpha+2),
\ldots,\hat y_n(\alpha+2))\\
&\text{ mod } \hat y_1(\alpha+2)^{\epsilon_1(\alpha+2)}\cdots\hat y_s(\alpha+2)^{\epsilon_s(\alpha+2)}(\hat y_1(\alpha+2),
\ldots,\hat y_l(\alpha+2))
\endalign
$$
with 
$$
\text{ mult }(\Phi_i(\hat y_{l+1}(\alpha+2),0,\ldots,0)))=1.
$$
$R(\alpha+2)$ has regular parameters $(x_1(\alpha+2),\ldots,x_n(\alpha+2))$ defined by
$$
x_i(\alpha+2)=\cases x_i^*(\alpha+2)& 1\ne l+1\\
\prod (x_{l+1}^*(\alpha+2)+c_{\alpha+2}-\sigma(c_{\alpha+2}))& i=l+1
\endcases\tag 3.91
$$
where the product is over the distinct congugates $\sigma(c_{\alpha+2})\in\overline k$ of $c_{\alpha+2}$ over $k$.

Set $\tilde{\overline x}_i(\alpha+2)=x_i(\alpha+2)$,  $\tilde{\overline y}_i(\alpha+2)=\hat y_i(\alpha+2)$
for $1\le i\le n$. 
Set $U''(\alpha+2)=U'(\alpha+2)[\gamma^{\alpha_1},\ldots,\gamma^{\alpha_s}]_q$
where 
$$
q=m(U(\alpha+2))\cap U'(\alpha+2)[\gamma^{\alpha_1},\ldots,\gamma^{\alpha_s}].
$$
Then $\text{CN}(\alpha+2)$ holds.
We have 
$x_{l+1}^*(\alpha+2)+c_{\alpha+2} = \gamma^{e_{s+1,s+1}}\lambda$
for some $\lambda\in k(S(\alpha+2))$ by (3.89). $e_{s+1,s+1}\ne 0$  by theorem 1.12.
$$
\align
&x_{l+1}^*(\alpha+2) = \tilde P_{\alpha+2}(\tilde{\overline  y}_1(\alpha+2),\ldots,\tilde{\overline y}_l(\alpha+2))\\
&+\tilde{\overline y}_1(\alpha+2)^{\epsilon_1(\alpha+2)}\cdots\tilde{\overline y}_s(\alpha+2)^{\epsilon_s(\alpha+2)}
\tilde \Sigma_{\alpha+2}(\tilde{\overline y}_1(\alpha+2),\ldots,\tilde{\overline y}_{n}(\alpha+2))
\endalign
$$
where $\text{mult}(\tilde \Sigma_{\alpha+2}(0,\ldots,0,\tilde{\overline y}_{l+1}(\alpha+2),0\ldots,0))=1$.
$$
\align
x_{l+1}(\alpha+2) &= 
\prod_{\sigma}\left(\tilde P_{\alpha+2}(\tilde{\overline  y}_1(\alpha+2),\ldots,\tilde{\overline y}_l(\alpha+2))
+(c_{\alpha+2}-\sigma(c_{\alpha+2}))\right)\\
&+\tilde{\overline y}_1(\alpha+2)^{\epsilon_1(\alpha+2)}\cdots\tilde{\overline y}_s(\alpha+2)^{\epsilon_s(\alpha+2)}\\
&[\sum_{\sigma}\left(\tilde P_{\alpha+2}(\tilde{\overline  y}_1(\alpha+2),\ldots,\tilde{\overline y}_l(\alpha+2))
+(c_{\alpha+2}-\sigma(c_{\alpha+2}))\right)\tilde \Sigma_{\alpha+2}\\
&+\tilde{\overline y}_1(\alpha+2)^{\epsilon_1(\alpha+2)}\cdots\tilde{\overline y}_s(\alpha+2)^{\epsilon_s(\alpha+2)}
\Omega]
\endalign
$$
in $S(\alpha+2)\,\,\hat{}\,\cong k(S(\alpha+2))[[\tilde{\overline y}_1(\alpha+2),\ldots,\tilde{\overline y}_n(\alpha+2)]]$
and the product and sum are over the distinct conjugates $\sigma(c_{\alpha+2})$ of $c_{\alpha+2}$ in $\overline k$ over $k$.
If $c_{\alpha+2}\not\in k$, we have 
$$
\sum_{\sigma}(c_{\alpha+2}-\sigma(c_{\alpha+2}))\ne 0
$$
since if this sum were 0, we would have $c_{\alpha+2}$ invariant under the automorphism group of a Galois closure of
$k(c_{\alpha+2})$ over $k$ which is impossible. We have an expression
$$
\align
x_{l+1}(\alpha+2)&=
P_{\alpha+2}(\tilde{\overline y}_1(\alpha+2),\ldots,\tilde{\overline y}_l(\alpha+2))\tag 3.92\\
&+\tilde{\overline y}_1(\alpha+2)^{\epsilon_1(\alpha+2)}\cdots\tilde{\overline y}_s(\alpha+2)^{\epsilon_s(\alpha+2)}
\Sigma_{\alpha+2}(\tilde{\overline y}_1(\alpha+2),\ldots,\tilde{\overline y}_{n}(\alpha+2))
\endalign
$$
where $P_{\alpha+2}$, $\Sigma_{\alpha+2}$ are power series with coefficients in $k(S(\alpha+2))$
and 
$$
\text{mult}(\Sigma_{\alpha+2}(0,\ldots,0,\tilde{\overline y}_{l+1},0,\ldots,0))=1.
$$
If $c_{\alpha+2}\in k$, $x_{l+2}(\alpha+2) = x_{l+2}^*(\alpha+2)$ and this expression is immediate. By (3.85) and (3.86),
$$
\align
g&= \sum_{i=1}^da_i'\tilde{\overline x}_1(\alpha+1)^{e_i^1(\alpha+1)}\cdots\tilde{\overline x}_s(\alpha+1)^{e_i^s(\alpha+1)}
\tilde{\overline x}_{l+1}(\alpha+1)^{f_i}\\
&+\sum_ja_j'\tilde{\overline x}_1(\alpha+1)^{e_j^1(\alpha+1)}\cdots\tilde{\overline x}_s(\alpha+1)^{e_j^s(\alpha+1)}
\tilde{\overline x}_{l+1}(\alpha+1)^{f_j}\\
&+\sum_ka_k'\tilde{\overline x}_1(\alpha+1)^{e_k^1(\alpha+1)}\cdots\tilde{\overline x}_s(\alpha+1)^{e_k^s(\alpha+1)}
\tilde{\overline x}_{l+1}(\alpha+1)^{f_k}\\
&+\tilde{\overline x}_{l+1}(\alpha+1)^{\tau}\Gamma
\endalign
$$
with $a_i',a_j',a_k'\in k(c_0,\ldots,c_{\alpha+1})[[\tilde{\overline x}_1(\alpha+1),\ldots,
\tilde{\overline x_l}(\alpha+1)]]$.
Since 
$$
\tilde{\overline x}_{l+1}(\alpha+1) = x_{l+1}(\alpha+1)=x_{l+1}
$$
 and $\text{CN}(\alpha+1)$ holds, so that
$$
k(c_0,\ldots,c_{l+1})[[x_1(\alpha+1),\ldots,x_l(\alpha+1)]] 
= k(c_0,\ldots,c_{l+1})[[\tilde{\overline x}_1(\alpha+1),\ldots,\tilde{\overline x}_l(\alpha+1)]],
$$
we have an expansion
$$
\align
g&=\sum_{i=1}^d \tilde a_i x_1(\alpha+1)^{e_i^1(\alpha+1)}\cdots x_s(\alpha+1)^{e_i^s(\alpha+1)}
 x_{l+1}(\alpha+1)^{f_i}\\
&+\sum_{j}\tilde a_j x_1(\alpha+1)^{e_j^1(\alpha+1)}\cdots x_s(\alpha+1)^{e_j^s(\alpha+1)}
 x_{l+1}(\alpha+1)^{f_j}\\
&+\sum_{k}\tilde a_k x_1(\alpha+1)^{e_k^1(\alpha+1)}\cdots x_s(\alpha+1)^{e_k^s(\alpha+1)}
 x_{l+1}(\alpha+1)^{f_k}\\
&+ x_{l+1}(\alpha+1)^{\tau}\Gamma
\endalign
$$
with
$\tilde a_i,\tilde a_j,\tilde a_k\in k(c_0,\ldots,c_{\alpha+1})[[ x_1(\alpha+1),\ldots, x_l(\alpha+1)]]$,
$\tilde a_i,\tilde a_j$ units for all $i,j$ and 
$$
\nu(x_1(\alpha+1)^{e_k^1(\alpha+1)}\cdots x_s(\alpha+1)^{e_k^s(\alpha+1)}
 x_{l+1}(\alpha+1)^{f_k})>\rho
$$
for all $k$. In $R(\alpha+2)\,\,\hat{}\,$, by (3.87)
$$
\align
g&=\sum_{i=1}^d\tilde a_i\prod_{a=1}^sx_a^*(\alpha+2)^{a_{1a}(\alpha+2)e_i^1(\alpha+1)+\cdots+a_{sa}(\alpha+2)e_i^s(\alpha+1)
+a_{s+1,a}(\alpha+2)f_i}\\
&(x_{l+1}^*(\alpha+2)+c_{\alpha+2})^
{a_{1,s+1}(\alpha+2)e_i^1(\alpha+1)+\cdots+a_{s,s+1}(\alpha+2)e_i^s(\alpha+1)
+a_{s+1,s+1}(\alpha+2)f_i}\\
&+\sum_{j}\tilde a_j\prod_{a=1}^sx_a^*(\alpha+2)^{a_{1a}(\alpha+2)e_j^1(\alpha+1)+\cdots+a_{sa}(\alpha+2)e_j^s(\alpha+1)
+a_{s+1,a}(\alpha+2)f_j}\\
&(x_{l+1}^*(\alpha+2)+c_{\alpha+2})^
{a_{1,s+1}(\alpha+2)e_j^1(\alpha+1)+\cdots+a_{s,s+1}(\alpha+2)e_j^s(\alpha+1)
+a_{s+1,s+1}(\alpha+2)f_j}\\
&+\sum_{k}\tilde a_k\prod_{a=1}^sx_a^*(\alpha+2)^{a_{1a}(\alpha+2)e_k^1(\alpha+1)+\cdots+a_{sa}(\alpha+2)e_k^s(\alpha+1)
+a_{s+1,a}(\alpha+2)f_k}\\
&(x_{l+1}^*(\alpha+2)+c_{\alpha+2})^
{a_{1,s+1}(\alpha+2)e_k^1(\alpha+1)+\cdots+a_{s,s+1}(\alpha+2)e_k^s(\alpha+1)
+a_{s+1,s+1}(\alpha+2)f_k}\\
&+\left[\left(\prod_{a=1}^sx_a^*(\alpha+2)^{a_{s+1,a}(\alpha+2)}\right)
(x_{l+1}^*(\alpha+2)+c_{\alpha+2})^{a_{s+1,s+1}(\alpha+2)}\right]^{\tau}\Gamma
\endalign
$$
with $\tilde a_i,\tilde a_j,\tilde a_k\in k(c_0,\ldots,c_{\alpha+2})[[x_1^*(\alpha+2),\ldots,x_{l+1}^*(\alpha+2)]]$,
 $\tilde a_i,\tilde a_j$ units and 
$$
\nu(\prod_{a=1}^sx_a^*(\alpha+2)^{a_{1a}(\alpha+2)e_k^1(\alpha+1)+\cdots+a_{sa}(\alpha+2)e_k^s(\alpha+1)
+a_{s+1,a}(\alpha+2)f_k})>\rho.
$$
We have that
$$
\sum_{a=1}^s(a_{1a}(\alpha+2)e_i^1(\alpha+1)+\cdots+a_{sa}(\alpha+2)e_i^s(\alpha+1)
+a_{s+1,a}(\alpha+2)f_i)\nu(x_a^*(\alpha+2))
$$
are equal for $1\le i\le d$ since the corresponding terms in $g$ have equal value $\rho$.
Set 
$$
\alpha_{\zeta} = (a_{1\zeta}(\alpha+2)e_1^1(\alpha+1)+\cdots+a_{s\zeta}(\alpha+2)e_1^s(\alpha+1)
+a_{s+1,\zeta}(\alpha+2)f_1
$$
for $1\le \zeta\le s$.
 Since
$\nu(x_1^*(\alpha+2)),\ldots,\nu(x_s^*(\alpha+2))$ are rationally independent, we have
$$
\alpha_\zeta
 =a_{1\zeta}(\alpha+2)e_i^1(\alpha+1)+\cdots+a_{s\zeta}(\alpha+2)e_i^s(\alpha+1)
+a_{s+1,\zeta}(\alpha+2)f_i
$$
for $1\le i\le d$ and $1\le \zeta\le s$.
Set
$$
m_{f_i}=
a_{1,s+1}(\alpha+2)e_i^1(\alpha+1)+\cdots+a_{s,s+1}(\alpha+2)e_i^s(\alpha+1)
+a_{s+1,s+1}(\alpha+2)f_i
$$
for $1\le i\le d$.
Set
$$
\epsilon=\frac
{\text{det}
\left(\matrix a_{11}(\alpha+2)&\cdots& a_{1s}(\alpha+2)\\
\vdots&&\vdots\\
a_{s1}(\alpha+2)&\cdots& a_{ss}(\alpha+2)
\endmatrix\right)}
{\text{det}
\left(\matrix a_{11}(\alpha+2)&\cdots& a_{1,s+1}(\alpha+2)\\
\vdots&&\vdots\\
a_{s+1,1}(\alpha+2)&\cdots& a_{s+1,s+1}(\alpha+2)
\endmatrix\right)}.
$$
$\epsilon$ is a nonzero integer. By Cramer's rule,
$$
f_i-f_1 = \epsilon(m_{f_i}-m_{f_1})
$$
for $1\le i\le d$.
Assume that $\epsilon>0$. 
We can perform a CRUTS of type M1) with associated MTSs 
 $R(\alpha+2)\rightarrow R(\alpha+3)$ and $S(\alpha+2)\rightarrow S(\alpha+3)$
where $R(\alpha+3)$ has regular parameters 
$(x_1(\alpha+3),\ldots, x_n(\alpha+3))$ and $T''(\alpha+3)$ has regular parameters 
$(\overline x_1(\alpha+3),\ldots, \overline x_n(\alpha+3))$ such that 
$$
\align
x_1(\alpha+2) =x_1^*(\alpha+2)&=x_1(\alpha+3)^{a_{11}(\alpha+3)}\cdots x_s(\alpha+3)^{a_{1s}(\alpha+3)}\tag 3.93\\
&\vdots\\
x_s(\alpha+2)=x_s^*(\alpha+2)&=x_1(\alpha+3)^{a_{s1}(\alpha+3)}\cdots x_s(\alpha+3)^{a_{ss}(\alpha+3)}
\endalign
$$
and 
$$
\overline x_i(\alpha+3)=\cases
x_i(\alpha+3)&i\le s\\
x_i^*(\alpha+2)& s<i
\endcases
$$
$x_i(\alpha+2)=x_i(\alpha+3)$ for $s<i$
to get
$$
g=\overline x_1(\alpha+3)^{a_1'}\cdots \overline x_s(\alpha+3)^{a_s'}((\overline x_{l+1}(\alpha+3)+c_{\alpha+2})^{m_{f_1}}\Phi
+\overline x_1(\alpha+3)\cdots \overline x_s(\alpha+3)G)
$$
in $k(c_0,\ldots,c_{\alpha+3})[[\overline x_1(\alpha+3),\ldots,\overline x_{l+1}(\alpha+3)]]$, with
$$
\Phi=
\tilde a_1+\tilde a_2(\overline x_{l+1}(\alpha+3)+c_{\alpha+2})^{\frac{f_2-f_1}{\epsilon}}+
\cdots+\tilde a_d(\overline x_{l+1}(\alpha+3)+c_{\alpha+2})^{\frac{f_d-f_1}{\epsilon}}.
$$
In the case $\epsilon< 0$, we must consider an expression
$$
g=\overline x_1(\alpha+3)^{a_1'}\cdots \overline x_s(\alpha+3)^{a_s'}((\overline x_{l+1}(\alpha+3)+c_{\alpha+2})^{m_{f_d}}\Phi'
+\overline x_1(\alpha+3)\cdots \overline x_s(\alpha+3)G')
$$
with
$$
\Phi'=
\tilde a_1(\overline x_{l+1}(\alpha+3)+c_{\alpha+2})^{\frac{f_1-f_d}{\epsilon}}+
\cdots+\tilde a_d.
$$
Again assume that $\epsilon>0$. The proof when $\epsilon<0$ is similar.
Let 
$$
r' =\text{mult}(\Phi(0,\ldots,0,\overline x_{l+1}(\alpha+3))).
$$
$\rho=\nu(a_i)+f_i\nu(\overline x_{l+1}(1))\le r\nu(\overline x_{l+1}(1))$ implies $f_d\le r$. Set
$\eta_{i}=\frac{f_i-f_1}{\epsilon}$. The residue of $\tilde a_i$ in
$$
T(\alpha+3)/(\overline x_1(\alpha+3),\ldots,\overline x_l(\alpha+3),\overline x_{l+2}(\alpha+3),\ldots,
\overline x_n(\alpha+3))\cong k(T(\alpha+3))[[\overline x_{l+1}(\alpha+3)]]
$$
is a nonzero constant $\hat a_i\in k(c_0,\ldots,c_{\alpha+2})$ for $1\le i\le d$. Set 
$$
\zeta(t) = \hat a_1+\hat a_2t^{\eta_2}+\cdots+\hat a_dt^{\eta_d},
$$
$r'=\text{mult}(\zeta(\overline x_{l+1}(\alpha+3)+c_{\alpha+2}))$.
$r'\le\eta_d=\frac{f_d-f_1}{\epsilon}\le r$. Suppose that $r'=r$. Then $f_1=0$, $f_d=r$, $\epsilon=1$ and
$\zeta(t)=\hat a_d(t-c_{\alpha+2})^r$. Thus there exist
nonzero $(\overline x_{l+1}(\alpha+3)+c_{\alpha+2})^r$ and  $(\overline x_{l+1}(\alpha+3)+c_{\alpha+2})^{r-1}$
terms in $\Phi(0,\ldots,0,\overline x_{l+1}(\alpha+3))$ and
$f_{d-1}=f_d-1=r-1$. Thus $a_d$ is a unit so that $a_d=\tilde a_d$ and
$\nu(a_{d-1}\overline x_{l+1}(1)^{r-1})=\nu(a_d\overline x_{l+1}(1)^r)$ implies
$\nu(\overline x_{l+1}(1))=\nu(a_{d-1})$. 
$$
\sum_{i=1}^d\hat a_i(\overline x_{l+1}(\alpha+3)+c_{\alpha+2})^{\frac{f_i-f_1}{\epsilon}}
=\hat a_d\overline x_{l+1}(\alpha+3)^r.\tag 3.94
$$
 Expanding out the LHS of (3.94), we have
$$
r\hat a_dc_{\alpha+2}+\hat a_{d-1}=0
$$
which implies
$$
\frac{c_{\alpha+2}}{\hat a_{d-1}}=-\frac{1}{r\hat a_{d}}.
$$
Let $\alpha=\hat a_d\in k(c_0)$ be the constant term of the power series
$a_d\in k(c_0)[[x_1(1),\ldots,x_l(1)]]$.
$$
\align
\frac{\overline x_{l+1}(1)}{a_{d-1}}&=\frac{a_d\overline x_{l+1}(1)^r}{a_da_{d-1}\overline x_{l+1}(1)^{r-1}}
=\frac{\tilde a_d(\overline x_{l+1}(\alpha+3)+c_{\alpha+2})^{m_{f_d}}}
{\tilde a_d\tilde a_{d-1}(\overline x_{l+1}(\alpha+3)+c_{\alpha+2})^{m_{f_{d-1}}}}\\
&=\frac{\tilde a_d(\overline x_{l+1}(\alpha+3)+c_{\alpha+2})^{\frac{f_d-f_1}{\epsilon}}}
{\tilde a_d\tilde a_{d-1}(\overline x_{l+1}(\alpha+3)+c_{\alpha+2})^{\frac{f_{d-1}-f_1}{\epsilon}}}\\
&=\frac{\tilde a_d(\overline x_{l+1}(\alpha+3)+c_{\alpha+2})}{\tilde a_d\tilde a_{d-1}}
\endalign
$$ 
which has residue $-\frac{1}{r\alpha}$ in $k(R(\alpha+3))\subset \Cal O_{\nu}/m_{\nu}$.
(Here $\Cal O_v$ is the valuation ring of our extension of $\nu$ to the quotient field of $S(\alpha+3)\,\,\hat{}\,$.)
 There exists $Q\in \tilde R(1)$
such that $Q$ is equivalent to $-\frac{1}{r\alpha}a_{d-1}$ modulo a sufficiently high power of the maximal ideal 
of $k(c_0)[[x_1(1),\ldots,x_l(1)]]$ (recall that $c_1=1$) so that
$\nu(x_{l+1}(1)-Q)>\nu(x_{l+1}(1))$ (Recall that $\overline x_{l+1}(1)=x_{l+1}(1)$).

Since $a_d\overline x_{l+1}(1)^r$ is a
minimum value term of $g$, we have $\nu(\overline x_{l+1}(1))\le \nu(g)$. Make a change of variables in $R(1)$ and $T''(1)$,
 replacing $x_{l+1}(1)$ and 
$\overline x_{l+1}(1)$ with 
$$
\overline x_{l+1}^{(1)}(1)=x_{l+1}^{(1)}(1)=x_{l+1}(1)-Q
$$
$\text{CN}(1)$ holds for these new variables. Further, in $S(1)\,\,\hat{}\,$, we have
$$
\overline x_{l+1}^{(1)}(1)=P^{(1)}(\overline y_1(1),\ldots,\overline y_l(1))+\overline y_1(1)^{d_1(1)}
\cdots\overline y_s(1)^{d_s(1)}\Sigma_0
$$
where $\text{mult}(\Sigma_0(0,\ldots,0,\overline y_{l+1}(1),0,\ldots,0))=1$.
Then repeat the above procedure with
this change of variable and our previous $g$. 
If $\nu(P^{(1)})>\nu(\tilde{\overline y}_1(\alpha)^{d_1(\alpha)}\cdots
\tilde{\overline y}_s(\alpha)^{d_s(\alpha)})$ the above algorithm produces an expression
$$
x_{l+1}^{(1)}(1)=\tilde{\overline y}_1(\alpha+1)^{d_1(\alpha+1)}\cdots\tilde{\overline y}_s(\alpha+1)^{d_s(\alpha+1)}
\Sigma_{\alpha}^{(1)}
$$
where $\text{mult}(\Sigma_{\alpha}^{(1)}(0,\ldots,0,\tilde{\overline y}_{l+1}(\alpha+1),0,\ldots,0))=1$.
So suppose that 
$$
\nu(P^{(1)})\le \nu(\tilde{\overline y}_1(\alpha)^{d_1(\alpha)}\cdots\tilde{\overline y}_s(\alpha)^{d_s(\alpha)}).
$$
If we do not get a reduction  $r_1<r$, we have
$$
\nu(\overline x_{l+1}(1))<\nu(\overline x_{l+1}^{(1)}(1))\le\nu(g).
$$
We can repeat this process. By Lemma 1.3, We eventually get a reduction $r'<r$, or
$\nu(g)=\infty$ and  we can construct (as in (3.31) in the proof of $A(\overline m)$ in Theorem 1)
$$
\phi = \text{Lim}_{i\rightarrow \infty}Q_i(\overline x_1(1),\ldots,\overline x_{l}(1))\in 
k(c_0)[[\overline x_1(1),\ldots,\overline x_l(1)]]
$$
such that 
$$
g = u(\overline x_{l+1}(1)-\phi)^r+h.\tag 3.95
$$
where $h\in a_lk(c_0)[[\overline x_1(1),\ldots,\overline x_{l+1}]]$ with
$$
a_l = \{f\in k(c_0)[[\overline x_1(1),\ldots,\overline x_l(1)]]\,|\,\nu(f)>\infty\}
$$ 
and
$u(\overline x_1(1),\ldots,\overline x_{l+1})$ is a unit power series.

Suppose that $r'<r$. In our construction,
$\overline y_{l+1}(\alpha+3)=\overline y_{l+1}$, so that if
$$
\Sigma_0 = \Sigma_0^{\alpha+3}(\overline y_1(\alpha+3),\ldots,\overline y_n(\alpha+3)),
$$
$\text{mult}(\Sigma_0^{\alpha+3}(0,\ldots,0,\overline y_{l+1}(\alpha+3),0,\ldots,0)) = 1$.
Set 
$$
g_1=(\overline x_{l+1}(\alpha+3)+c_{\alpha+2})^{m_{f_1}}\Phi
+\overline x_1(\alpha+3)\cdots \overline x_s(\alpha+3)G.
$$
$x_{l+1}(\alpha+2)=\eta x_{l+1}^*(\alpha+2)$ where $\eta\in k(c_{\alpha+2})[x^*_{l+1}(\alpha+2)]$
is a unit
which implies 
$$
x_{l+1}(\alpha+3)=\eta\overline x_{l+1}(\alpha+3).
$$
 Thus
$$
\align
g_1 &\in 
k(c_0,\ldots,c_{\alpha+3})[[\overline x_1(\alpha+3),\ldots,\overline x_{l+1}(\alpha+3)]]\\
&=k(c_0,\ldots,c_{\alpha+3})[[ x_1(\alpha+3),\ldots, x_{l+1}(\alpha+3)]]\subset
R(\alpha+3)\,\,\hat{}\,
\endalign
$$
$\Sigma_0\mid g_1$ so that
$$
1\le \text{mult}(g_1(0,\ldots,0, x_{l+1}(\alpha+3))=
\text{mult}(g_1(0,\ldots,0,\overline x_{l+1}(\alpha+3))=r_1<r.
$$
By (3.92), there is an expression in $S(\alpha+3)\,\,\hat{}\,$
$$
\align
x_{l+1}(\alpha+3) &= x_{l+1}(\alpha+2)\\
&=P_{\alpha+3}(\tilde{\overline y}_1(\alpha+3),\ldots,\tilde{\overline y}_l(\alpha+3))\\
&\,\,\,+\tilde{\overline y}_1(\alpha+3)^{d_1(\alpha+3)}\cdots\tilde{\overline y}_s(\alpha+3)^{d_s(\alpha+3)}
\Sigma_{\alpha+3}(\tilde{\overline y}_1(\alpha+3),\ldots,\tilde{\overline y}_n(\alpha+3))
\endalign
$$
where
$\text{mult}(\Sigma_{\alpha+3}(0,\ldots,0,\tilde{\overline y}_{l+1}(\alpha+3),0,\ldots,0)) = 1$.
Now set $\tilde{\overline x}_i(\alpha+3)=x_i(\alpha+3)$ for $1\le i\le n$. By (3.90), (3.91) and (3.93) $\text{CN}(\alpha+3)$
holds (with $\gamma_i(\alpha+3)=1$ for $1\le i\le l$).

By induction on r we can now repeat the procedure following (3.85),
with $R(1),\tilde R(1), S(1)$ replaced with $R(\alpha+3), \tilde R(\alpha+3), S(\alpha+3)$ respectively,
$c_0$ with a primitive element of $k(c_0,\ldots,c_{\alpha+3)}$ over $k$, $g$ with $g_1$,
 to eventually  get $t$ such that $\text{CN}(t)$ holds with $\tilde{\overline x}_i(t)=x_i(t)$ for $1\le i\le n$ and
$$
x_{l+1}(t)= \tilde{\overline y}_1(t)^{d_1(t)}\cdots\tilde{\overline y}_s(t)^{d_s(t)}
\Sigma_t(\tilde{\overline y}_{1}(t),\ldots,\tilde{\overline y}_{n}(t))
$$
where $\text{mult}(\Sigma_t(0,\ldots,0,\tilde{\overline y}_{l+1}(t),0,\ldots,0)=1$ or we have 
$$
x_{l+1}(t) = P_t(\tilde{\overline y}_1(t),\cdots,\tilde{\overline y}_l(t))+\tilde{\overline y}_1(t)^{d_1(t)}\cdots
\tilde{\overline y}_s(t)^{d_s(t)}\Sigma_t(\tilde{\overline y}_1(t),\ldots,
\tilde{\overline y}_{n}(t))
$$
where $\text{mult}(\Sigma_t(0,\ldots,0,\tilde{\overline y}_{l+1}(t),0\ldots,0)=1$,
$P_t,\Sigma_t$ are series with coefficients in $k(S(t))$ and 
$$
g=u( x_1(t),\ldots, x_{l+1}(t))
 x_1(t)^{d_1}\cdots, x_{s}(t)^{d_s}[ x_{l+1}(t)+
\Phi( x_{1}(t),\ldots, x_{l}(t))]^a\tag 3.96
$$
for some $a>0$ where $u,\Phi$ are series with coefficients in $k(c_0,\ldots,c_t)$,
$u$ is a unit and $\Sigma_0\mid g$. Suppose that (3.96) holds. Define
$$
\Sigma_0^t(\tilde{\overline y}_1(t),\ldots,\tilde{\overline y}_n(t))=\Sigma_0.
$$
$\text{mult}(\Sigma_0^t(0,\ldots,0,\tilde{\overline y}_{l+1}(t),0,\ldots,0)=1$.
We have regular parameters $(\hat y_1(t),\ldots,\hat y_n(t))$ in $S(t)\,\,\hat{}\,$ defined by
$$
\hat y_i(t) = \cases
\tilde{\overline y}_i(t)&i\ne l+1\\
\Sigma_0&i=l+1
\endcases
$$
There is an expression
$$
x_{l+1}(t) = \overline P_t(\hat y_1(t),\ldots,\hat y_l(t))+\hat y_1(t)^{d_1(t)}\cdots\hat y_s^{d_s(t)}
\hat \Sigma_t(\hat y_1(t),\ldots,\hat y_n(t))
$$
with 
$\text{mult}(\hat \Sigma_t(0,\ldots,0,\hat y_{l+1}(t),0,\ldots,0)=1$.
Thus
$$
\hat y_{l+1}(t)\mid x_{l+1}(t)+
\Phi( x_{1}(t),\ldots, x_{l}(t))
$$
in $S(t)\,\,\hat{}\,$. Since $P_t+\Phi\in k(S(t))[[\tilde{\overline y}_1(t),\ldots,\tilde{\overline y}_l(t)]]$,
$$
P_t+\Phi=\Omega\tilde{\overline y}_1(t)^{d_1(t)}\cdots
\tilde{\overline y}_s(t)^{d_s(t)}
$$
with $\Omega\in k(S(t))[[\tilde{\overline y}_1(t),\ldots,\tilde{\overline y}_{l}(t)]]$ and
$$
\Sigma_0\mid \Omega+\Sigma_t.
$$
Set $m_{t,l}=m\left(k(c_0,\ldots,c_t)[[x_1(t),\ldots,x_l(t)]]\right)$.
Choose $N$ so that 
$$
N\nu(m_{t,l})> \nu(\tilde{\overline y}_1(t)^{d_1(t)}\cdots\tilde{\overline y}_s(t)^{d_s(t)}).
$$

 There exists $\Phi'\in k(c_0,\ldots,c_t)[x_1(t),\ldots,x_l(t)]$ such that $\Phi'\equiv\Phi\text{ mod }m_{t,l}^N$.
Make a change of variables, replacing $x_{l+1}(t)$ with $x_{l+1}(t)+\Phi'$ to get
$$
x_{l+1}(t)= \tilde{\overline y}_1(t)^{d_1(t)}\cdots\tilde{\overline y}_s(t)^{d_s(t)}(\Omega+\Sigma_t)+(\Phi'-\Phi)
$$
By Theorem 3.11,
we can perform a CRUTS along $\nu$ in the first $l$ variables, with associated MTSs
 $R(t)\rightarrow R(t')$, $S(t)\rightarrow S(t')$   to get
$$
\tilde{\overline x}_{l+1}(t')=x_{l+1}(t')= 
\tilde{\overline y}_1(t')^{d_1(t')}\cdots\tilde{\overline y}_s(t')^{d_s(t')}\Sigma_{t'},
$$
where $\text{mult}(\Sigma_{t'}(0,\ldots,0,\tilde{\overline y}_{l+1}(t'),0,\ldots,0))=1$
and such that $\text{CN}(t')$ holds.

We thus have regular parameters $(x_1(t'),\ldots,x_n(t'))$ in $R(t')$ and
$(y_1(t'),\ldots,y_n(t'))$ in $S(t')$, units $\tau_1(t'),\ldots,\tau_s(t')\in S(t')$ such that
$$
\align
x_1(t')&= y_1(t')^{c_{11}(t')}\cdots y_s(t')^{c_{1s}(t')}\tau_1(t')\\
&\vdots\\
x_s(t')&= y_1(t')^{c_{s1}(t')}\cdots y_s(t')^{c_{ss}(t')}\tau_s(t')\\
x_{s+1}(t')&=y_{s+1}(t')\\
&\vdots\\
x_l(t')&=y_l(t')\\
x_{l+1}(t')&=y_1(t')^{d_1(t')}\cdots y_s(t')^{d_s(t')}y_{l+1}(t').
\endalign
$$
Let $\phi_i(t')$ be the residue of $\tau_i(t')$ in $k(S(t'))$,
$\overline \tau_i = \frac{\tau_i(t')}{\phi_i(t')}$. Let $(e_{ij})=(c_{ij}(t'))^{-1}$. Define 
$$
\overline y_i(t')=\cases
\overline \tau_1^{e_{i1}}\cdots  \overline \tau_s^{e_{is}}y_i(t')&1\le i\le s\\
y_i(t')&s<i, i\ne l+1\\
\overline\tau_1^{-e_{11}d_1(t')-\cdots-e_{s1}d_s(t')}\cdots
\overline\tau_s^{-e_{1s}d_1(t')-\cdots-e_{ss}d_s(t')}y_{l+1}(t')&i=l+1
\endcases
$$
We have
$$
\align
x_1(t')&= \overline y_1(t')^{c_{11}(t')}\cdots \overline y_s(t')^{c_{1s}(t')}\phi_1(t')\\
&\vdots\\
x_s(t')&= \overline y_1(t')^{c_{s1}(t')}\cdots \overline y_s(t')^{c_{ss}(t')}\phi_s(t')\\
x_{s+1}(t')&=\overline y_{s+1}(t')\\
&\vdots\\
x_l(t')&=\overline y_l(t')\\
x_{l+1}(t')&=\overline y_1(t')^{d_1(t')}\cdots\overline  y_s(t')^{d_s(t')}\overline y_{l+1}(t')
\endalign
$$
in 
$$
U''(t')=S(t')[\overline \tau_1^{e_{11}},\ldots,\overline\tau_s^{e_{ss}}]_{(\overline y_1(t'),\ldots,\overline y_n(t'))}.
$$

By Lemma 3.5, we can perform a MTS of type $\text{II}_{l+1}$ $R(t')\rightarrow R(t'+1)$
$$
\align
\overline x_1(t') &=\overline x_1(t'+1)^{a_{11}(t'+1)}\cdots \overline x_s(t'+1)^{a_{1s}(t'+1)}
c_{t'+1}^{a_{1,s+1}(t'+1)}\\
\vdots&\\
\overline x_s(t') &=\overline x_1(t'+1)^{a_{s1}(t'+1)}\cdots \overline x_s(t'+1)^{a_{ss}(t'+1)}
c_{t'+1}^{a_{s,s+1}(t'+1)}\\
  x_{l+1}(t') 
&= \overline x_1(t'+1)^{a_{s+1,1}(t'+1)}\cdots \overline x_s(t'+1)^{a_{s+1,s}(t'+1)}(\overline x_{l+1}(t'+1)+1)
c_{t'+1}^{a_{s+1,s+1}(t'+1)}
\endalign
$$
and a MTS of type $\text{II}_{l+1}$ (possibly followed by a transformation of type I) $S(t')\rightarrow S(t'+1))$
$$
\align
\overline y_1(t') &=\overline y_1(t'+1)^{b_{11}(t'+1)}\cdots \overline y_s(t'+1)^{b_{1s}(t'+1)}
d_{t'+1}^{b_{1,s+1}(t'+1)}\\
\vdots&\\
 \overline y_s(t') &=\overline y_1(t'+1)^{b_{s1}(t'+1)}\cdots \overline y_s(t'+1)^{b_{ss}(t'+1)}
d_{t'+1}^{b_{s,s+1}(t'+1)}\\
 \overline y_{l+1}(t') 
&= \overline y_1(t'+1)^{b_{s+1,1}(t'+1)}\cdots \overline y_s(t'+1)^{b_{s+1,s}(t'+1)}(\overline y_{l+1}(t'+1)+1)
d_{t'+1}^{b_{s+1,s+1}(t'+1)}
\endalign
$$
such that $R(t'+1)\subset S(t'+1)$, and 
$\overline x_{l+1}(t'+1) =\overline y_{l+1}(t'+1)$.
By adding an appropriate series $\Omega$ to $\overline x_{l+1}(t'+1)$, we will have regular parameters
in $R(t'+1)\rightarrow S(t'+1)$ as desired.
\enddemo
\heading
4. Monomialization
\endheading

\proclaim{Theorem 4.1}
Suppose that $R\subset S$ are excellent regular local rings such that $\text{dim}(R)=\text{dim}(S)$, containing a 
field $k$ of characteristic 0 such that the quotient field $K$ of $S$ is a finite extension of the quotient field $J$ of $R$.
 Suppose that $\nu$ is a valuation of $K$ with 
valuation ring $V$ such that $V$ dominates $S$.
Suppose that $\nu$ has rank 1 and rational rank $s$. Then there exist sequences of
monodial transforms $R\rightarrow R'$ and $S\rightarrow S'$ along $\nu$ such that 
$\text{dim}(R')=\text{dim}(S')$, $S'$ dominates $R'$, $\nu$ dominates $S'$ and there are regular parameters
$(x_1',\ldots,x_n')$ in $R'$, $(y_1',\ldots,y_n')$ in $S'$, units $\delta_1,\ldots,
\delta_s\in S'$ and a matrix $(c_{ij})$ of nonnegative integers such that $\text{det}(c_{ij})\ne 0$ and
$$
\align
x_1'&=(y_1')^{c_{11}}\cdots (y_s')^{c_{1s}}\delta_1\\
&\vdots\\
x_s'&=(y_1')^{c_{s1}}\cdots (y_s')^{c_{ss}}\delta_s\\
x_{s+1}'&=y_{s+1}'\\
&\vdots\\
x_n'&=y_n'.
\endalign
$$
\enddemo

\demo{Proof}
By Theorem 1.7, applied to the lift to $V$ of a transcendence basis of $V/m_{\nu}$ over $R/m$
(which is finite by Theorem 1 [Ab2] or Appendix 2 [ZS]), there exists a MTS along 
$\nu$, $R\rightarrow R_1$, such
that $\text{dim}_{R_1}(\nu)=0$. Let $m_1$ be the maximal ideal of $R_1$. 
$\text{trdeg}_{R/m}(R_1/m_1)=\text{dim}_R(\nu)$ and $\text{dim}(R_1)=\text{dim}(R)-\text{dim}_R(\nu)$
by the dimension formula (Theorem 15.6 [M]).
By Theorem 1.6, there exists a MTS $S\rightarrow S_1$ along $\nu$  such that
$S_1$ dominates $R_1$. Let $n_1$ be the maximal ideal of $S_1$. $S$ is essentially of finite type (a spot) over $R$
by Theorem 1.11, since $\text{dim}(R)=\text{dim}(S)$. 
Hence $S_1$ is a spot over $R_1$.
By the dimension formula, 
$$
\text{dim}(R_1)=\text{dim}(S_1)+\text{trdeg}_{R_1/m_1}(S_1/n_1)=\text{dim}(S_1),
$$
since $\text{trdeg}_{R_1/m_1}(S_1/n_1)\le\text{dim}_{R_1}(\nu)=0$. We may thus assume that $\text{dim}_R(\nu)=0$.

Let $\{\overline t_i\}$ be a transcendence basis of 
$R/m$ over $k$. Let $t_i$ be lifts of $\overline t_i$ to $R$. Then the field $L$ obtained by adjoining all of the
$t_i$ to $k$ is contained in $R$, and $\nu$ is trivial on $L-\{0\}$. hence we can replace $k$ by $L$.
We may thus assume that assumptions 1) and 2) of Chapter 3 hold.

There exist $f_1,\ldots,f_s\in J$ such that $\nu(f_1),\ldots,\nu(f_s)$ are positive and rationally independent.
 By Theorem 1.7, there 
exists a MTS $R\rightarrow R'$ along $\nu$, such that $f_1,\ldots,f_s\in R'$. By Theorem 1.5, there exists a MTS $R'\rightarrow R''$
along $\nu$
such that $f_1\cdots f_s$ is a SNC divisor in $R''$. Then $R''$ has regular parameters $(x_1'',\ldots,x_n'')$ such that
 $\nu(x_1''),\ldots,\nu(x_s'')$ are rationally independent.
By Theorem 1.6, there exists a MTS $S\rightarrow S'$ 
along $\nu$, such that $R''\subset S'$. We may thus assume that 
 there exist regular parameters  $(x_1,\ldots,x_n)$ in $R$
such that $\nu(x_1),\ldots,\nu(x_s)$ are rationally independent.

 By Theorem 1.5, after replacing $S$ with a MTS along $\nu$ we may assume that  $x_1\cdots x_n$ has SNCs in $S$.
Thus there are regular parameters $(y_1,\ldots,y_n)$ in $S$ and units $\psi_i$ such that 
$$
x_i=y_1^{c_{i1}}\cdots y_n^{c_{in}}\psi_i
$$
for $1\le i\le s$.
Thus $\nu(y_1),\ldots,\nu(y_n)$ span a rational vector space of dimension $s$. After possibly reindexing the
$y_i$, we may assume that $\nu(y_1),\ldots,\nu(y_s)$ are rationally independent. By 1) of Theorem 3.9 with
$R=S$, $f=x_1\ldots x_s$ and 
 Theorem 3.11,
 we can
replace $S$ with a MTS along $\nu$ to get 
$$
x_i=y_1^{c_{i1}}\cdots y_s^{c_{is}}\psi_i
$$
for $1\le i\le s$,
where $\psi_i$ are units and $\text{det}(c_{ij})\ne 0$.

Let $\phi_i$ be the residue of $\psi_i$ in $S/n$. For $1\le j\le s$ set
$$
\epsilon_i=\prod_{i=1}^s\left(\frac{\psi_j}{\phi_j}\right)^{e_{ij}}
$$
where $(e_{ij})= (c_{ij})^{-1}$, a matrix with rational coefficients.
$\epsilon_j\in \hat S$ for $1\le j\le s$.

 Set $ T''(0)= R$, $\overline x_i=x_i$
for $1\le i\le n$.
Set $U''(0)=S[d_0,\epsilon_1,\ldots,\epsilon_s]_q$ where 
$d_0\in \hat S$ is such that $k(d_0)\cong k(S)$, $q=m(\hat S)\cap S[d_0,\epsilon_1,\ldots,\epsilon_s]$.
$ U''(0)$ has regular parameters
$$
\overline y_j=\cases
\epsilon_j y_j&1\le j\le s\\
y_j& s<j.
\endcases
$$ 
$$
\overline x_i = \overline y_1^{c_{i1}}\cdots \overline y_s^{c_{is}}\phi_i
$$
for $1\le i\le s$. Set $\tilde R(0) = k[x_1,\ldots,x_s]_q$, 
$q=m(R)\cap k[x_1,\ldots, x_s]$, $c_0=1$. Thus the assumptions of Theorem 3.12 are satisfied 
with $l=s$ and by the conclusions of Theorem 3.12
applied $n-s$ times consecutively,
we can construct  MTSs $R\rightarrow R'$, $S\rightarrow S'$ such that $V$ dominates $S'$, $S'$ dominates $R'$
and $R'$ has regular parameters $(x_1',\ldots,x_n')$, $S'$ has regular parameters
$(y_1',\ldots, y_n')$ satisfying the conclusions of the Theorem.
\enddemo

\proclaim{Corollary 4.2}
Suppose that $R\subset S$ are excellent regular local rings such that $\text{dim}(R)=\text{dim}(S)$, containing a 
field $k$ of characteristic 0 and with a common quotient field $K$. Suppose that $\nu$ is a valuation of $K$ with 
valuation ring $V$ such that $V$ dominates $S$.
Suppose that $\nu$ has rank 1 and rational rank $s$. Then there exist sequences of
monodial transforms $R\rightarrow R'$ and $S\rightarrow S'$ along $\nu$ such that 
$\text{dim}(R')=\text{dim}(S')$, $S'$ dominates $R'$, $\nu$ dominates $S'$ and there are regular parameters
$(x_1',\ldots,x_n')$ in $R'$, $(\tilde y_1,\ldots,\tilde y_n)$ in $S'$ such that
$$
\align
x_1'&=\tilde y_1^{c_{11}}\cdots \tilde y_s^{c_{1s}}\\
&\vdots\\
x_s'&=\tilde y_1^{c_{s1}}\cdots \tilde y_s^{c_{ss}}\\
x_{s+1}'&=\tilde y_{s+1}\\
&\vdots\\
x_n'&=\tilde y_n
\endalign
$$
where $\text{det}(c_{ij})=\pm 1$ and $k(R')\cong k(S')$.
\endproclaim

\demo{Proof} We can construct MTSs along $\nu$ $R\rightarrow R'$, $S\rightarrow S'$ such that the
conclusions of Theorem 4.1 hold.
To finish the proof, we must show that $\text{det}(c_{ij})=\pm 1$ and $k(R')\cong k(S')$.
We will analyze $(c_{ij})$ by constructing MTSs which may not be dominated by $\nu$.
Since interchanging the variables $x_i'$ will only change the sign of $\text{det}(c_{ij})$, we may assume that $c_{11}\ne0$.

\subheading{Case 1} Suppose that $c_{11}<c_{21}$. Then we can perform a MTS $S'\rightarrow S(1)$ where
$S(1)$ has regular parameters $(y_1(1),\ldots ,y_n(1))$ such that
$$
y_i'=\cases
y_1(1)y_2(1)^m\cdots y_s(1)^m& i=1\\
y_i(1)& i\ne 1
\endcases
$$
Then for $m>>0$ the monoidal transform $R'\rightarrow R(1)$ factors through $S(1)$, where $R(1)$ 
has regular parameters $(x_1(1),\ldots,x_n(1))$ defined by
$$
x_i'=\cases
x_1(1)x_2(1)&i=2\\
x_i(1)& i\ne 2
\endcases.
$$
Then
$$
\align
x_1(1)&=y_1(1)^{c_{11}}\cdots\\
x_2(1)&=y_1(1)^{c_{21}-c_{11}}\cdots\\
x_3(1)&=y_1(1)^{c_{31}}\cdots\\
&\vdots
\endalign
$$
\subheading{Case 2} Suppose that $c_{21}<c_{11}$. As in Case 1, we can perform MTSs to get
$$
\align
x_1(1)&=y_1(1)^{c_{11}-c_{21}}\cdots\\
x_2(1)&=y_1(1)^{c_{21}}\cdots\\
x_3(1)&=y_1(1)^{c_{31}}\cdots\\
&\vdots
\endalign
$$
\subheading{Case 3} Suppose that $c_{11}=c_{21}$ and $c_{1j}<c_{2j}$ for some $j$. Perform a MTS $S'\rightarrow S(1)$
where $S(1)$ has regular parameters $(y_1(1),\ldots ,y_n(1))$ such that
$$
y_i'=\cases
y_j(1)y_2(1)^m\cdots y_{j-1}(1)^my_{j+1}(1)^m\cdots y_s(1)^m& i=j\\
y_i(1)& i\ne j
\endcases
$$
Then for $m>>0$ the monoidal tranform $R'\rightarrow R(1)$ factors through $S(1)$, where $R(1)$ 
has regular parameters $(x_1(1),\ldots,x_n(1))$ defined by
$$
x_i'=\cases
x_1(1)x_2(1)&i=2\\
x_i(1)& i\ne 2
\endcases.
$$
Then
$$
\align
x_1(1)&=y_1(1)^{c_{11}}\cdots\\
x_2(1)&=y_2(1)^{c_{22}-c_{12}+m(c_{2j}-c_{1j})}\cdots\\
x_3(1)&=y_1(1)^{c_{31}}\cdots\\
&\vdots
\endalign
$$
\subheading{Case 4} In the remaining case $c_{11}=c_{21}$ and $c_{1j}\ge c_{2j}$ for all $j$. Then 
the monoidal tranform $R'\rightarrow R(1)$ factors through $S'$, where $R(1)$ 
has regular paramaters $(x_1(1),\ldots,x_n(1))$ defined by
$$
x_i'=\cases
x_1(1)x_2(1)&i=1\\
x_i(1)& i\ne 1
\endcases.
$$
Then
$$
\align
x_1(1)&=(y_2')^{c_{12}-c_{22}}\cdots \\
x_2(1)&=(y_1')^{c_{21}}(y_2')^{c_{22}}\cdots\\
x_3(1)&=(y_1')^{c_{31}}\cdots\\
&\vdots
\endalign
$$

By continuing to apply these four cases, we can construct $R'\rightarrow R(\alpha)$ and $S'\rightarrow S(\alpha)$ such that
$S(\alpha)$ dominates $R(\alpha)$,   
$$
\align
x_1(\alpha)&=y_1(\alpha)^{c_{11}(\alpha)}\cdots y_s^{c_{1s}(\alpha)}\phi_1\tag 4.1\\
&\vdots\\
x_s(\alpha)&=y_1(\alpha)^{c_{s1}(\alpha)}\cdots y_s^{c_{ss}(\alpha)}\phi_s\\
x_{s+1}(\alpha)&=y_{s+1}(\alpha)\\
&\vdots\\
x_n(\alpha)&=y_n(\alpha)
\endalign
$$
with $\phi_i$ units in $S(\alpha)$
and $c_{21}(\alpha)=0$.
By repeating the above procedure on successive rows we can construct a MTS (4.1) with
$$
c_{21}(\alpha)=\cdots=c_{s1}(\alpha)=0.
$$
Then the algorithm can be applied to the matrix obtained by removing the first row and column from $(c_{ij})$
 to construct (4.1) such that $(c_{ij}(\alpha))$ is a
upper triangular matrix. 

Set $q=(y_1(\alpha))$, $p=R(\alpha)\cap (y_1(\alpha))$. Our assumption that $(c_{ij})$ is upper triangular implies
$$
qS(\alpha)\,\,\hat{}\,\cap R(\alpha)\,\,\hat{}\,= x_1(\alpha)R(\alpha)\,\,\hat{}\,,
$$
so that $p=(x_1(\alpha))$ and $\text{dim }R(\alpha)_p=\text{dim }S(\alpha)_q$. By the dimension formula,
$A = (S(\alpha)/q)_q$ is finite over $(R(\alpha)/p)_p$. $S(\alpha)_q/pS(\alpha)_q\cong
A[y_1(\alpha)]/(y_1(\alpha)^{c_{11(\alpha)}})$ is then finite over $(R(\alpha)/p)_p$, so that
$R(\alpha)_p=S(\alpha)_q$ and $c_{11}(\alpha)=1$ by Theorem 1.10.

Now perform the MTS $S(\alpha)\rightarrow S(\alpha+1)$ where $S(\alpha+1)$ has regular parameters $(y_1(\alpha+1),\ldots,
y_n(\alpha+1))$ such that
$$
y_i(\alpha)=\cases
y_1(\alpha+1)y_2(\alpha+1)^{m_1}y_3(\alpha+1)^{m_2}\cdots y_s(\alpha+1)^{m_2}&i=1\\
y_i(\alpha+1)&i>1
\endcases
$$
where $m_1$ is chosen so that 
$$
m_1+c_{12}(\alpha)=rc_{22}(\alpha)
$$
 for some $r>0$ and $m_2$ is sufficiently large
that 
$$
m_2+c_{1j}(\alpha)>rc_{2j}(\alpha)
$$
for $3\le j\le s$. Then the MTS $R(\alpha)\rightarrow R(\alpha+1)$
factors through $S(\alpha+1)$, where
$$
x_i(\alpha)=\cases
x_1(\alpha+1)x_2(\alpha+1)^r&i=1\\
x_i(\alpha+1)&i>1
\endcases
$$
and
$$
\align
x_1(\alpha+1)&=y_1(1)y_3(1)^{c_{31}(\alpha+1)}\cdots\\
x_2(\alpha+1)&=y_2(1)^{c_{22}(\alpha+1)}\cdots\\
x_3(\alpha+1)&=y_3(1)^{c_{33}(\alpha+1)}\cdots\\
&\vdots
\endalign
$$
Now perform a series of similar MTSs to get (4.1) with $(c_{ij}(\alpha))$ an upper triangular matrix with
$$
c_{1j}(\alpha)=\cases
1&j=1\\
0&j>1.
\endcases
$$
Set $q=(y_1(\alpha),y_2(\alpha))$, $p=R(\alpha)\cap q=(x_1(\alpha),x_2(\alpha))$.
$$
S(\alpha)_q/pS(\alpha)_q\cong (S(\alpha)/q)_q[y_2(\alpha)]/(y_2(\alpha)^{c_{22}(\alpha)})
$$
is finite over $(R(\alpha)/p)_p$. By Theorem 1.10, $R(\alpha)_p=S(\alpha)_q$ and $c_{22}(\alpha)=1$.
We can repeat the above procedure to get (4.1) where $(c_{ij}(\alpha))$ is the identity matrix and $R(\alpha)=S(\alpha)$.

Thus $\text{det}(c_{ij})=\pm1$. Furthermore,
$$
k(R')\cong k(R(\alpha))\cong k(S(\alpha))\cong k(S').
$$
Set $(e_{ij})=(c_{ij})^{-1}$, a matrix with integral coefficients. Set
$$
\tilde y_i =\cases
  \delta_1^{e_{i1}}\cdots \delta_s^{e_{is}}  y_i' & 1\le i\le s\\
y_i' & s<i
\endcases
$$
then $(\tilde y_1,\ldots,\tilde y_n)$ are regular parameters in $S'$ satisfying the conclusions of the Theorem.
\enddemo
\vskip .3truein

Suppose that $R\subset S$ are excellent regular local rings such that $\text{dim}(R)=\text{dim}(S)=n$,
containing a field $k$ of characteristic 0, such that the quotient field $K$ of $S$ is a finite extension of the
quotient field $J$ of $R$.  Suppose that $\nu$ is a valuation of $K$
with valuation ring $V$ such that $V$ dominates $S$ and $\nu$ has rank $r$.
The primes of $V$ are a chain
$$
0=p_0\subset \cdots \subset p_r\subset V.
$$
We will begin by reviewing basic facts on specialization and composition of valuations 
(c.f. sections 8,9,10 of [Ab3] and chapter VI, section 10 of [ZS]).
Suppose that $\Gamma_{\nu}$ is the value group of $\nu$. The isolated subgroups of $\Gamma_{\nu}$ are a chain
$$
0=\Gamma_{r}\subset\cdots\subset \Gamma_0=\Gamma_{\nu}.
$$
Set 
$$
U_i=\{\nu(a)\mid a\in p_i\}
$$
Then the isolated subgroup $\Gamma_i$ of $p_i$ is defined to be the complement of $U_i$ and $-U_i$ in $\Gamma_{\nu}$.

For $i<j$ $\nu$ induces a valuation on the field $(V/p_i)_{p_i}$ with  valuation ring $(V/p_i)_{p_j}$ and value group
$\Gamma_i/\Gamma_j$. If $j=i+1$, $\Gamma_i/\Gamma_j$ has rank 1.

For all $i$, $V_{p_i}$ is a valuation ring of $K$ dominating $R_{p_i\cap R}$. Thus
$$
\text{trdeg}_{(R/P_i\cap R)_{p_i\cap R}}(V/p_i)_{p_i}<\infty
$$
by Theorem 1 [Ab2] or Appendix 2 [ZS]. We can lift transcendence bases of $(V/p_i)_{p_i}$ over 
$(R/p_i\cap R)_{p_i\cap R}$ for $1\le i\le r$ to $t_1,\ldots,t_a\in V$.
After possibly replacing the $t_i$ with $\frac{1}{t_i}$,
 we have $\nu(t_i)\ge0$ for all $t_i$. By Theorem 1.7, there
exists a MTS $R\rightarrow R'$ along $\nu$ such that $t_i\in R'$ for all $i$.
Let $p'_i=R'\cap p_i$. Then 
$$
\text{trdeg}_{(R'/p_i')_{p_i'}}(V/p_i)_{p_i}=0
$$
for $1\le i\le r$.
 By Theorem 1.6, there exists a MTS $S\rightarrow S'$ along $\nu$ such that $S'$ dominates $R'$.
Replacing $R$ by $R'$ and $S$ by $S'$, we may thus assume that 
$$
\text{trdeg}_{(R/p_i\cap R)_{p_i\cap R}}(V/p_i)_{p_i}=0
$$
for $1\le i\le r$. Then
$$
\text{trdeg}_{(R/p_i\cap R)_{p_i\cap R}}(S/p_i\cap S)_{p_i\cap S}=0
$$
for all $i$. By the dimension formula (c.f. Theorem 15.6 [M])
$$
\text{ dim } R/p_i\cap R = \text{ dim } S/p_i\cap S
$$
for $0\le i\le r$.

\proclaim{Theorem 4.3}
Suppose that $R\subset S$ are excellent regular local rings such that $\text{dim}(R)=\text{dim}(S)=n$,
containing a field $k$ of characteristic 0 such that the quotient field $K$ of $S$ is
a finite extension of the quotient field $J$ of $R$. Suppose that $\nu$ is a valuation of $K$
with valuation ring $V$ such that $V$ dominates $S$ and $\nu$ has rank $r$. Suppose that the segments of $\Gamma_{\nu}$ are
$$
0=\Gamma_r\subset\cdots\subset \Gamma_0=\Gamma_{\nu}
$$
with associated primes 
$$
0=p_0\subset\cdots\subset p_r\subset V.
$$
 Suppose that $\Gamma_{i-1}/\Gamma_{i}$ has rational rank $s_{i}$ for $1\le i\le r$ and
$$
\text{trdeg}_{(R/p_i\cap R)_{p_i\cap R}}(V/p_i)_{p_i}=0
$$
for $1\le i\le r$.
Set $t_i=\text{dim}(R/p_{i-1}\cap R)_{p_{i}\cap R}$ for $1\le i\le r$, so that $n=t_1+\cdots+t_r$. Then there exist MTSs
$R\rightarrow R'$ and $S\rightarrow S'$ along $\nu$ such that $S'$ dominates $R'$, $R'$ has regular parameters
$(z_1,\ldots, z_n)$, $S'$ has regular parameters $(w_1,\ldots, w_n)$ such that
$$
p_{i}\cap R' = (z_1,\ldots,z_{t_1+\cdots+t_i})
$$
$$
p_i\cap S' = (w_1,\ldots,w_{t_1+\cdots+t_i})
$$
for $1\le i\le r$ and
\vfill\eject
$$
\align
z_1&=w_1^{g_{11}(1)}\cdots w_{s_1}^{g_{1s_1}(1)}w_{t_1+1}^{h_{1,t_1+1}(1)}\cdots w_n^{h_{1n}(1)}\delta_{11}\\
&\vdots\\
z_{s_1}&=w_1^{g_{s_11}(1)}\cdots w_{s_1}^{g_{s_1s_1}(1)}w_{t_1+1}^{h_{s_1,t_1+1}(1)}\cdots w_n^{h_{s_1n}(1)}\delta_{1s_1}\\
z_{s_1+1}&=w_{s_1+1}w_{t_1+1}^{h_{s_1+1,t_1+1}(1)}\cdots w_n^{h_{s_1+1,n}(1)}\delta_{1s_1+1}\\
&\vdots\\
z_{t_1}&=w_{t_1}w_{t_1+1}^{h_{t_1,t_1+1}(1)}\cdots w_n^{h_{t_1n}(1)}\delta_{1t_1}\\
z_{t_1+1}&=w_{t_1+1}^{g_{11}(2)}\cdots w_{t_1+s_2}^{g_{1s_2}(2)}
w_{t_1+t_2+1}^{h_{1,t_1+t_2+1}(2)}\cdots w_n^{h_{1n}(2)}\delta_{21}\\
&\vdots\\
z_{t_1+s_2}&=w_{t_1+1}^{g_{s_21}(2)}\cdots w_{t_1+s_2}^{g_{s_2s_2}(2)}
w_{t_1+t_2+1}^{h_{s_2,t_1+t_2+1}(2)}\cdots w_n^{h_{s_2n}(2)}\delta_{2s_2}\\
z_{t_1+s_2+1}&=w_{t_1+s_2+1}w_{t_1+t_2+1}^{h_{s_2+1,t_1+t_2+1}(2)}\cdots w_n^{h_{s_2+1,n}(2)}\delta_{2s_2+1}\\
&\vdots\\
z_{t_1+t_2}&=w_{t_1+t_2}w_{t_1+t_2+1}^{h_{t_2,t_1+t_2+1}(2)}\cdots w_n^{h_{t_2n}(2)}\delta_{2t_2}\\
&\vdots\\
z_{t_1+\cdots+t_{r-1}+1}&=w_{t_1+\cdots+t_{r-1}+1}^{g_{11}(r)}\cdots w_{t_1+\cdots+t_{r-1}+s_r}^{g_{1s_r}(r)}\delta_{r1}\\
&\vdots\\
z_{t_1+\cdots+t_{r-1}+s_r}&=w_{t_1+\cdots+t_{r-1}+1}^{g_{s_r1}(r)}\cdots w_{t_1+\cdots+t_{r-1}+s_r}^{g_{s_rs_r}(r)}
\delta_{rs_r}\\
z_{t_1+\cdots+t_{r-1}+s_r+1}&=w_{t_1+\cdots+t_{r-1}+s_r+1}
\delta_{rs_r+1}\\
&\vdots\\
z_{t_1+\cdots+t_r}&=w_{t_1+\cdots+t_r}\delta_{rt_r}
\endalign
$$
where
$$
\text{det}\left(\matrix g_{11}(i)&\cdots&g_{1s_i}(i)\\
\vdots&&\vdots\\
 g_{s_i1}(i)&\cdots&g_{s_is_i}(i)
\endmatrix\right) \ne0,
$$
$\delta_{ij}$ are units in $S'$
for $1\le i\le r$, and $h_{jk}(i)$ are nonnegative integers.
\endproclaim

\demo{Proof} The proof is by induction on $r$. $r=1$ is immediate from Theorem 4.1. Suppose that the Theorem is true for
rank $r-1$. Set
$p_i(0)=p_i\cap R$, $q_i(0)=p_i\cap S$. Then there exist MTSs $R_{p_{r-1}(0)}\rightarrow T_1$ and $S_{q_{r-1}(0)}\rightarrow U_1$
such that $V_{p_{r-1}}$ dominates $U_1$, $U_1$ dominates $T_1$ and the conclusions of the Theorem hold for $T_1\subset U_1$.
 By Theorems 1.9 and 1.6 there exist MTSs along $\nu$
$R\rightarrow R(1)$ and $S\rightarrow S(1)$ such that $V$ dominates $S(1)$, $S(1)$ dominates $R(1)$ and if
$p_i(1)=p_i\cap R(1)$, $q_i(1)=p_i\cap S(1)$, $R(1)_{p_{r-1}(1)}\cong T_1$, $S(1)_{q_{r-1}(1)}\cong U_1$ and
$R(1)/p_i(1)$, $S(1)/q_i(1)$ are regular local rings for $1\le i\le r$.

 By assumption,
$R(1)_{p_{r-1}(1)}$ has regular parameters $(\tilde x_1,\ldots,\tilde x_{t_1+\cdots+t_{r-1}})$ and\linebreak
$S(1)_{q_{r-1}(1)}$ has regular parameters $(\tilde y_1,\ldots,\tilde y_{t_1+\cdots+t_{r-1}})$
 satisfying the conclusions of the 
Theorem. Set $\lambda = t_1+\cdots+t_{r-1}$. $R(1)$ has regular parameters $(x_1(1),\ldots, x_n(1))$ such that 
$p_{r-1}(1)=(x_1(1),\ldots,x_{\lambda}(1))$. Let
$\pi(1):R(1)\rightarrow R(1)/p_{r-1}(1)$. There exist $\overline a_j^i\in R(1)_{p_{r-1}(1)}$, $1\le i\le \lambda$,
$1\le j\le \lambda$ such that
$$
\tilde x_i=\overline a_1^ix_1(1)+\cdots+\overline a_{\lambda}^ix_{\lambda}(1)
$$
and $\text{det}(\overline a_j^i)\not\in p_{r-1}(1)_{p_{r-1}(1)}$. There exists $u(1)\in R(1)-p_{r-1}(1)$ such that
$u(1)\tilde x_i\in R(1)$ for $1\le i\le \lambda$ and if we define $\tilde x_i(1)=u(1)\tilde x_i$
$$
\tilde x_i(1)=a_1^ix_1(1)+\cdots+a_{\lambda}^ix_{\lambda}(1)
$$
for $1\le i\le \lambda$ where $a_j^i\in R(1)$ for all $i,j$ and $\text{det}(a_j^i)\not\in p_{r-1}(1)$. 
After reindexing the $\tilde x_i(1)$,
 we may assume that $a_1^1\not\in p_{r-1}(1)$. Let $b_1^1=\pi(1)(a_1^1)$.

$V/p_{r-1}$ is a rank 1, rational rank $s_r$ valuation ring. The quotient field of $V/p_{r-1}$ is algebraic over the quotient field
of $R(1)/p_{r-1}(1)$ so that if $L$ is the quotient field of $R(1)/p_{r-1}(1)$, then
$L\cap V/p_{r-1}$ is a rank 1, rational rank $s_r$ valuation ring. Let $\overline \nu$
denote the valuation induced by $\nu$ on $L$.

 By Theorems 4.1, 3.9 and 3.11  (with $R=S=R(1)/p_{r-1}(1)$)  there exists a
MTS 
$$
\tilde R(1)=R(1)/p_{r-1}(1)\rightarrow \tilde R(2)\rightarrow \cdots\rightarrow \tilde R(m)
$$
where each $\tilde R(i)\rightarrow\tilde R(i+1)$ is a monoidal transform and $\tilde R(m)$ has regular parameters
$(\overline y_{\lambda+1}(m),\ldots,\overline y_n(m))$ such that
$\overline \nu(\overline y_{\lambda+1}(m)),\ldots,\overline \nu(\overline y_{\lambda+s_r}(m))$ are rationally independent
and 
$b_1^1 = \overline y_{\lambda+1}(m)^{a_{\lambda+1}}\cdots\overline y_{\lambda+s_r}(m)^{a_{\lambda+s_r}}\overline u$ where
$\overline u\in \tilde R(m)$ is a unit. There exist regular parameters $(\overline y_{\lambda+1}(1),\ldots,\overline y_n(1))$
in $\tilde R(1)$ and $a\le n$ such that 
$$
\tilde R(2) = \tilde R(1)\left[\frac{\overline y_{\lambda+2}(1)}{\overline y_{\lambda+1}(1)},\cdots,
\frac{\overline y_{a}(1)}{\overline y_{\lambda+1}(1)}\right]_{\overline Q(2)}
$$
where $\overline Q(2)$ is a maximal ideal. Let $y_i(1)$ be lifts of $\overline y_i(1)$ to $R(1)$ for $\lambda+1\le i\le n$.
Then $(x_1(1),\ldots,x_{\lambda}(1),y_{\lambda+1}(1),\ldots,y_n(1))$ are regular parameters in $R(1)$.
We have a surjection
$$
\Phi_1:R(1)\left[\frac{ y_{\lambda+2}(1)}{ y_{\lambda+1}(1)},\cdots,
\frac{ y_{a}(1)}{ y_{\lambda+1}(1)}\right]\rightarrow 
\tilde R(1)\left[\frac{\overline y_{\lambda+2}(1)}{\overline y_{\lambda+1}(1)},\cdots,
\frac{\overline y_{a}(1)}{\overline y_{\lambda+1}(1)}\right].
$$
Let $Q_2=\Phi_1^{-1}(\overline Q(2))$. Set 
$$
R_2 = R(1)\left[\frac{ y_{\lambda+2}(1)}{ y_{\lambda+1}(1)},\cdots,
\frac{ y_{a}(1)}{ y_{\lambda+1}(1)}\right]_{Q_2}.
$$
$R(1)\rightarrow R_2$ is a monoidal transform along $\nu$ and $p_{r-1}\cap R_2 = (x_1(1),\ldots, x_{\lambda}(1))$,
$$
(R_2)_{p_{r-1}\cap R_2}\cong R(1)_{p_{r-1}(1)}\text{ and }
R_2/p_{r-1}\cap R_2\cong \tilde R(2).
$$
 We can inductively construct a MTS along $\nu$
$$
R(1)\rightarrow R_2\rightarrow\cdots\rightarrow R_m=R(2)\tag 4.2
$$
such that $R(2)_{p_{r-1}(2)}\cong R(1)_{p_{r-1}(1)}$, $R(2)/p_{r-1}(2)\cong \tilde R(m)$
 with $p_{r-1}(2)=p_{r-1}\cap R(2)$
 and 
$R(2)$ has regular parameters $(x_1(2),\ldots,x_n(2))$ such that
$$
x_i(2) = \cases x_i(1)&1\le i\le \lambda\\
y_i(2) & \lambda+1\le \lambda\le n
\endcases
$$
where $y_i(2)$ are lifts of $\overline y_i(m)$ to $R(2)$. Thus
$$
a_1^1 = x_{\lambda+1}(2)^{a_{\lambda+1}}\cdots x_{\lambda+s_r}(2)^{a_{\lambda+s_r}}u+b_1x_1(2)+\cdots+b_{\lambda}x_{\lambda}(2)
$$
where $u,b_1,\ldots,b_{\lambda}\in R(2)$ and $u$ is a unit. Thus
$$
\tilde x_1(1) = x_{\lambda+1}(2)^{a_{\lambda+1}}\cdots x_{\lambda+s_r}(2)^{a_{\lambda+s_r}}ux_1(2)+ax_1(2)^2+
\sum_{i=2}^{\lambda}a_ix_i(2)
$$
with $a_i,a\in R(2)$. Now perform a MTS $R(2)\rightarrow R(3)$ along $\nu$
$$
x_i(2)=\cases
x_{\lambda+1}(3)^{a_{\lambda+1}+1}\cdots x_{\lambda+s_r}(3)^{a_{\lambda+s_r}+1}x_1(3)&i=1\\
x_{\lambda+1}(3)^{2a_{\lambda+1}+2}\cdots x_{\lambda+s_r}(3)^{2a_{\lambda+s_r}+2}x_i(3)& 2\le i\le \lambda\\
x_i(3)&\lambda+1\le i\le n
\endcases
$$
Thus $x_i(3)\in R(2)_{p_{r-1}(2)}$ for $1\le i\le n$.
 Set $p_{r-1}(3) = p_{r-1}\cap R(3)$. Then 
$R(2)_{p_{r-1}(2)} = R(3)_{p_{r-1}(3)}$.
$$
\tilde x_1(1) = x_{\lambda+1}(3)^{2a_{\lambda+1}+1}\cdots x_{\lambda+s_r}(3)^{2a_{\lambda+s_r}+1}
(x_1(3)u+x_{\lambda+1}(3)\cdots x_{\lambda+s_r}(3)c)
$$
for some $c \in p_{r-1}(3)$. Set
$$
\hat x_i(3) = \cases 
x_1(3)u+x_{\lambda+1}(3)\cdots x_{\lambda+s_r}(3)c&i=1\\
x_i(3)& 2\le i\le n
\endcases
$$
Then $(\hat x_1(3),\ldots, \hat x_n(3))$ are regular parameters in $R(3)$ with $p_{r-1}(3) = (\hat x_1(3),\ldots,\hat x_{\lambda}(3))$.
$p_{r-1}(3)R(3)_{p_{r-1}(3)} = (\tilde x_1,\ldots,\tilde x_{\lambda})R(3)_{p_{r-1}(3)}$ implies there exists
$\overline a_j^i(3)\in R(3)_{p_{r-1}(3)}$ such that
$$
\tilde x_i = \cases
\overline a_1^1(3)\hat x_1(3)&i=1\\
\overline a_1^i(3)\hat x_1(3)+\cdots +\overline a_{\lambda}^i(3)\hat x_{\lambda}(3)& 2\le i\le \lambda
\endcases
$$
and 
$$
\text{det}\left(\matrix \overline a_2^2(3)&\cdots&\overline a_{\lambda}^2(3)\\
\vdots&&\vdots\\
\overline a_2^{\lambda}(3)&\cdots&\overline a_{\lambda}^{\lambda}(3)
\endmatrix
\right)
\not\in p_{r-1}(3)_{p_{r-1}(3)}.
$$
We can repeat the above argument to construct a MTS $R(1)\rightarrow R''$ along $\nu$ such that if $p_{r-1}''=p_{r-1}\cap R''$, 
$R(1)_{p_{r-1}(1)}=R''_{p_{r-1}''}$ and there exists a regular system of parameters $(x_1'',\ldots,x_n'')$ in $R''$ and
$u_1,\ldots,u_{\lambda}\in R''_{p_{r-1}''}-(p_{r-1}'')_{p_{r-1}''}$ such that $\tilde x_i=u_ix_i''$ for $1\le i\le \lambda$.
By Theorem 1.6 and the above argument, there exists a MTS $S(1)\rightarrow S''$ along $\nu$
such that if $q_{r-1}'' = p_{r-1}\cap S''$, $S(1)_{q_{r-1}(1)}=S''_{q_{r-1}''}$, $S''$ dominates $R''$
 and there are regular parameters
$(y_1'',\ldots,y_n'')$ in $S''$ and $v_1,\ldots,v_{\lambda}\in S''_{q_{r-1}''}-(q_{r-1}'')_{q_{r-1}''}$ such that
$\tilde y_i=v_iy_i''$ for $1\le i\le\lambda$. Thus we have
$$
\align
x_1'' &=\psi_1(y_1'')^{g_{11}(1)}\cdots (y_{s_1}'')^{g_{1s_1}(1)}(y''_{t_1+1})^{h_{1t_1+1}(1)}\cdots
(y_{\lambda}'')^{h_{1\lambda}(1)}\tag 4.3\\
&\vdots\\
x_{\lambda}''&=\psi_{\lambda}y_{\lambda}''
\endalign
$$
with $\psi_1,\ldots,\psi_{\lambda}\in S''_{q_{r-1}''}-(q_{r-1}'')_{q_{r-1}''}$. 
$\psi_i=\frac{f_i}{g_i}$ with $f_i,g_i\in S''-q_{r-1}''$, $f_i,g_i$ relatively prime in $S''$.
There are nonnegative integers $d_1^i$ such that

$g_ix_i'' = f_i(y_1'')^{d_1^i}\cdots(y_{\lambda}'')^{d_s^i}$ for $1\le i\le \lambda$ so that 
$g_i\mid f_i$ in $S''$ and $\psi_i\in S''-q_{r-1}''$ for $1\le i\le \lambda$.

Let $\pi':R''\rightarrow R''/p_{r-1}''$ and $\pi'':S''\rightarrow S''/q_{r-1}''$. Let $\overline x_i''=\pi'(x_i'')$ and 
$\overline y_i''=\pi''(y_i'')$. 
$\nu$ induces a rank 1 rational rank $s_r$ valuation on $\overline K = (S''/q_{r-1}'')_{q_{r-1}''}$. 
 By Theorem 4.1, there exist MTSs
$$
R''/p_{r-1}''=\tilde R\rightarrow \tilde R(1)\rightarrow\cdots\rightarrow \tilde R(m)=T
$$
and
 $$
S''/q_{r-1}''=\tilde S \rightarrow \tilde S(1)\rightarrow\cdots\rightarrow \tilde S(m)=U
$$
such that the valuation ring $(V/p_{r-1})\cap \overline K$ dominates $U$, $U$ dominates $T$, $T$ has regular parameters
$(\overline x_{\lambda+1},\ldots,\overline x_n)$, $U$ has regular parameters $(\overline y_{\lambda+1},\ldots,\overline y_n)$
such that\linebreak
  $\overline \nu(\overline y_{\lambda+1}),\ldots,\overline \nu(\overline y_{\lambda+s_r})$
are rationally independent, where $\overline \nu$ is the valuation induced by $\nu$ on the quotient field of $U$ and
$$
\align
\overline x_{\lambda+1} &= \overline y_{\lambda+1}^{g_{11}(r)}\cdots \overline y_{\lambda+s_r}^{g_{rs_r}(r)}
\overline \delta_{\lambda+1}\tag 4.4\\
&\vdots\\
\overline x_{\lambda+s_r} &= \overline y_{\lambda+1}^{g_{s_r1}(r)}\cdots \overline y_{\lambda+s_r}^{g_{s_rs_r}(r)}
\overline \delta_{\lambda+s_r}\\
\overline x_{\lambda+s_r+1} &= \overline y_{\lambda+s_r+1}\overline\delta_{\lambda+s_r+1}\\
&\vdots\\
\overline x_n&=\overline y_n\overline\delta_n
\endalign
$$
where $\overline \delta_i$ are units in $U$. 

Each $\tilde R(i)\rightarrow \tilde R(i+1)$ is a monoidal transform centered at a prime $\overline a_i$. By Theorems 3.9 and 3.11
and Lemma 3.2, there exist MTSs along $\nu$
$$
T = \tilde R(m)\rightarrow \cdots\rightarrow \tilde R(m')=T'
$$
and
$$
U = \tilde S(m)\rightarrow \cdots\rightarrow \tilde S(m')=U'
$$
such that $U'$ dominates $T'$, $T'$ has regular parameters $(\overline x_{\lambda+1},\ldots,\overline x_n)$,
$U'$ has regular parameters $(\overline y_{\lambda+1},\ldots,\overline y_n)$ such that (4.4) holds, and
$$
\overline a_iU' = (\overline y_{\lambda+1}^{d_1^i}\cdots\overline y_{\lambda+s_r}^{d_{s_r}^i})
$$
 for some
nonnegative integers $d_1^i,\ldots,d_{s_r}^i$  for $1\le i\le m$,
and
$$
\pi''(\psi_i)=\overline y_{\lambda+1}^{a_{i,\lambda+1}}\cdots\overline y_{\lambda+s_r}^{a_i{\lambda+s_r}}\overline u_i
$$
where $\overline u_i$ are units, $a_{ij}$ are positive integers for $1\le i\le \lambda$.

For $m\le i\le m'-1$
each $\tilde R(i)\rightarrow \tilde R(i+1)$ is a monoidal transform centered at a height 2 prime $\overline a_i$
(c.f Remark 3.1) such that $\overline a_iU' = (\overline y_{\lambda+1}^{d_1^i}\cdots \overline y_{\lambda+s_r}^{d_{s_r}^i})$
for some nonnegative integers $d_1^i,\ldots,d_{s_r}^i$. Consider the MTSs along $\nu$
$$
R''\rightarrow R(1)\rightarrow\cdots\rightarrow R(m')
$$
 and
$$
S''\rightarrow S(1)\rightarrow\cdots\rightarrow S(m')=\overline S
$$
constructed as in (4.2),
so that for $1\le i\le m'$, $R(i)_{p_{r-1}(i)}\cong R''_{p_{r-1}''}$,
$R(i)/p_{r-1}(i)\cong \tilde R(i)$, $S(i)_{q_{r-1}(i)}\cong S''_{q_{r-1}''}$,
$S(i)/q_{r-1}(i)\cong \tilde S(i)$ where $p_{r-1}(i)=p_{r-1}\cap R(i)$, $q_{r-1}=p_{r-1}\cap S(i)$
and $\overline S$ has regular parameters 
$(\tilde y_1,\ldots,\tilde y_n)$ such that $\tilde y_i$ has residue $\overline y_i$ in $U'$ for $\lambda+1\le i\le n$
and $\tilde y_i=y_i''$ for $1\le i\le \lambda$.
For $0\le i\le m'-1$, 
$R(i)\rightarrow R(i+1)$ is the blowup of $a_i\subset R(i)$ such that 
$a_i\tilde R(i)=\overline a_i$. Thus
$a_i\overline S/\overline q_{r-1} = (\overline y_{\lambda+1}^{d_1^i}\cdots \overline y_{\lambda+s_r}^{d_{s_r}^i})$
where $\overline q_{r-1} = p_{r-1}\cap \overline S$.
Set $\Phi_i= \tilde y_{\lambda+1}^{d_1^i}\cdots \tilde y_{\lambda+s_r}^{d_{s_r}^i}$.  Then
$$
a_i\overline S = (\Phi_i+\tilde y_1b_1^1(i)+\cdots +\tilde y_{\lambda}b_{\lambda}^1(i),\cdots,
\Phi_i+\tilde y_1b_1^t(i)+\cdots +\tilde y_{\lambda}b_{\lambda}^t(i))
$$
for some $t$, $b_j^k(i)\in\overline S$, $1\le i\le m'-1$.

Perform a MTS along $\nu$
$$
\overline S=\overline S(0)\rightarrow \overline S(1)\rightarrow\cdots\overline S(m')
$$
where $\overline S(j)$ has regular parameters
$(\overline y_1(j),\ldots,\overline y_n(j))$ defined by $\overline y_i(0)=\tilde y_i$ for $1\le i\le n$,
$$
\overline y_i(j)=\cases
\Phi_{j+1}\overline y_i(j+1)& 1\le i\le \lambda\\
\overline y_i(j+1)&i<\lambda\le n
\endcases
$$
for $0\le j\le m'-1$. Then we have 
 $a_i\overline S(m') = (\Phi_i)$ for $1\le i\le m'-1$. $R(m')\subset \overline S(m')$ (by Theorem 1.6) and
$\overline S(m')/\overline q(m')_{r-1}\cong U'$.

Let $z_i$ be  lifts of $\overline x_i$ to $R(m')$ for $\lambda+1\le i\le n$. Define regular parameters 
$(x_1(m'),\ldots,x_n(m'))$ in $R(m')$ by
$$
x_i(m')=\cases
x_i''&1\le i\le \lambda\\
z_i&\lambda+1\le i\le n
\endcases
$$
There exists a matrix of nonnegative integers $(e_{ij})$ such that 
$$
\align
x_1(m')&=\overline y_1(m')^{g_{11}(1)}\cdots\overline y_{s1}(m')^{g_{1s_1}(1)}
\overline y_{t_1+1}(m')^{h_{1,t_1+1}(1)}\cdots\overline y_{\lambda}(m')^{h_{1\lambda}(1)}\tag 4.5\\
&\,\,\,\,\,\overline y_{\lambda+1}(m')^{e_{1,\lambda+1}}\cdots\overline y_n(m')^{e_{1n}}\psi_1\\
&\vdots\\
x_{\lambda}(m')&=\overline y_{\lambda}(m')\overline y_{\lambda+1}(m')^{e_{\lambda,\lambda+1}}
\cdots\overline y_n(m')^{e_{\lambda,n}}\psi_{\lambda}\\
x_{\lambda+1}(m')&=\overline y_{\lambda+1}(m')^{g_{11}(r)}\cdots\overline y_{\lambda+s_r}(m')^{g_{s_rs_r}(r)}\delta_{\lambda+1}\\
&+f_1^{\lambda+1}\overline y_1(m')+\cdots+f_{\lambda}^{\lambda+1}\overline y_{\lambda}(m')\\
&\vdots\\
x_n(m')&=\overline y_n(m')\delta_n+f_1^n\overline y_1(m')+\cdots+f_{\lambda}^n\overline  y_{\lambda}(m')
\endalign
$$
where $\delta_i$ are lifts of $\overline \delta_i$ to $\overline S(m')$, $f_i^j\in \overline S(m')$.
For $1\le i\le \lambda$,
$$
\psi_i=u_i'\overline y_{\lambda+1}(m')^{a_{i,\lambda+1}}\cdots\overline y_{\lambda+s_r}^{a_i\lambda+s_r}+
h_1^i\overline y_1(m')+\cdots+h_{\lambda}^i\overline y_{\lambda}(m')
$$
where $u_i'$ are lifts of $\overline u_i$ to $\overline S(m')$, the $u_i'$ and $\delta_i$ are units in $\overline S(m')$. Choose 
$$
t>\text{max}\{a_{ij},g_{ij}(1)\}.
$$

Now perform a MTS $\overline S(m')\rightarrow \overline S(m'+1)$ along $\nu$ where $\overline S(m'+1)$ has regular
parameters $(\overline y_1(m'+1),\ldots,\overline y_n(m'+1))$ defined by
$$
\overline y_i(m')=\cases
\overline y_{\lambda+1}(m'+1)^t\cdots\overline y_n(m'+1)^t\overline y_i(m'+1)&1\le i\le \lambda\\
\overline y_i(m+1)& \lambda+1\le i\le n
\endcases
$$
to get
$$
\psi_i = u_i\overline  y_{\lambda+1}(m'+1)^{a_{i\lambda+1}}\cdots \overline  y_{\lambda+s_r}(m'+1)^{a_{\lambda+s_r}}
$$
for some units $u_i\in \overline S(m'+1)$, $1\le i\le \lambda$.
$\overline S(m'+1)/\overline q(m'+1)_{r-1}\cong U'$ and there is a matrix of
nonnegative integers $(b_{ij})$, units $u_{\lambda+1},\ldots,u_n\in \overline S(m+1)$ such that 
$$
\align
x_1(m')&=\overline y_1(m'+1)^{g_{11}(1)}\cdots \overline y_{s_1}(m'+1)^{g_{1s_1}(1)}
\overline y_{t'+1}(m'+1)^{h_{1,t_1+1}(1)}\cdots\overline y_{\lambda}(m'+1)^{h_{1\lambda}(1)}\tag 4.6\\
&\,\,\,\,\,\overline y_{\lambda+1}(m'+1)^{b_{1,\lambda+1}}\cdots \overline y_n(m'+1)^{b_{1n}}u_1\\
&\vdots\\
x_{\lambda}(m')&=\overline y_{\lambda}(m'+1)\overline y_{\lambda+1}(m'+1)^{b_{\lambda,\lambda+1}}
\cdots \overline y_n(m'+1)^{b_{\lambda n}}u_{\lambda}\\
x_{\lambda+1}(m')&=\overline y_{\lambda+1}(m'+1)^{g_{11}(r)}\cdots 
\overline y_{\lambda+s_r}(m'+1)^{g_{1s_r}(r)}u_{\lambda+1}\\ 
&\vdots\\
x_n(m')&=\overline y_n(m')u_n
\endalign
$$
\enddemo
\vskip .4truein
Theorem A is immediate from Theorem 4.3

\proclaim{Theorem 4.4}
Suppose that $R\subset S$ are excellent regular local rings such that $\text{dim}(R)=\text{dim}(S)=n$,
containing a field $k$ of characteristic 0 and with a common quotient field $K$. Suppose that $\nu$ is a valuation of $K$
with valuation ring $V$ such that $V$ dominates $S$ and $\nu$ has rank $r$. Suppose that the segments of $\Gamma_{\nu}$ are
$$
0=\Gamma_r\subset\cdots\subset \Gamma_0=\Gamma_{\nu}
$$
with associated primes 
$$
0=p_0\subset\cdots\subset p_r\subset V.
$$
 Suppose that $\Gamma_{i-1}/\Gamma_{i}$ has rational rank $s_{i}$ for $1\le i\le r$ and
$$
\text{trdeg}_{(R/p_i\cap R)_{p_i\cap R}}(V/p_i)_{p_i}=0
$$
for $1\le i\le r$.
Set $t_i=\text{dim}(R/p_{i-1}\cap R)_{p_{i}\cap R}$ for $1\le i\le r$, so that $n=t_1+\cdots+t_r$. Then there exist MTSs
$R\rightarrow R'$ and $S\rightarrow S'$ along $\nu$ such that $S'$ dominates $R'$, $R'$ has regular parameters
$(z_1,\ldots, z_n)$, $S'$ has regular parameters $(w_1,\ldots, w_n)$ such that
$$
p_{i}\cap R' = (z_1,\ldots,z_{t_1+\cdots+t_i})
$$
$$
p_i\cap S' = (w_1,\ldots,w_{t_1+\cdots+t_i})
$$
for $1\le i\le r$ and
\vfill\eject

$$
\align
z_1&=w_1^{g_{11}(1)}\cdots w_{s_1}^{g_{1s_1}(1)}\\
&\vdots\\
z_{s_1}&=w_1^{g_{s_11}(1)}\cdots w_{s_1}^{g_{s_1s_1}(1)}\\
z_{s_1+1}&=w_{s_1+1}\\
&\vdots\\
z_{t_1}&=w_{t_1}\\
z_{t_1+1}&=w_{t_1+1}^{g_{11}(2)}\cdots w_{t_1+s_2}^{g_{1s_2}(2)}\\
&\vdots\\
z_{t_1+s_2}&=w_{t_1+1}^{g_{s_21}(2)}\cdots w_{t_1+s_2}^{g_{s_2s_2}(2)}\\
z_{t_1+s_2+1}&=w_{t_1+s_2+1}\\
&\vdots\\
z_{t_1+t_2}&=w_{t_1+t_2}\\
&\vdots\\
z_{t_1+\cdots+t_{r-1}+1}&=w_{t_1+\cdots+t_{r-1}+1}^{g_{11}(r)}\cdots w_{t_1+\cdots+t_{r-1}+s_r}^{g_{1s_r}(r)}\\
&\vdots\\
z_{t_1+\cdots+t_{r-1}+s_r}&=w_{t_1+\cdots+t_{r-1}+1}^{g_{s_r1}(r)}\cdots w_{t_1+\cdots+t_{r-1}+s_r}^{g_{s_rs_r}(r)}\\
z_{t_1+\cdots+t_{r-1}+s_r+1}&=w_{t_1+\cdots+t_{r-1}+s_r+1}\\
&\vdots\\
z_{t_1+\cdots+t_r}&=w_{t_1+\cdots+t_r}
\endalign
$$
where
$$
\text{Det}\left(\matrix g_{11}(i)&\cdots&g_{1s_i}(i)\\
\vdots&&\vdots\\
 g_{s_i1}(i)&\cdots&g_{s_is_i}(i)
\endmatrix\right) = \pm 1
$$
and
$(R'/p_i\cap R')_{p_i\cap R'}\cong  (S'/p_i\cap S')_{p_i\cap S'}$
for $1\le i\le r$.
\endproclaim

\demo{Proof} The proof is a refinement of that of Theorem 4.3. The stronger Corollary 4.2 is used instead of Theorem 4.1. Formula 
(4.3) then becomes
$$
\align
x_1'' &=\psi_1(y_1'')^{g_{11}(1)}\cdots (y_{s_1}'')^{g_{1s_1}(1)}\tag 4.7\\
&\vdots\\
x_{\lambda}''&=\psi_{\lambda}y_{\lambda}''
\endalign
$$
(4.5) becomes
$$
\align
x_1(m')&=\overline y_1(m')^{g_{11}(1)}\cdots\overline y_{s1}(m')^{g_{1s_1}(1)}
\overline y_{\lambda+1}(m')^{e_{1,\lambda+1}}\cdots\overline y_n(m')^{e_{1n}}\psi_1\tag 4.8\\
&\vdots\\
x_{\lambda}(m')&=\overline y_{\lambda}(m')\overline y_{\lambda+1}(m')^{e_{\lambda,\lambda+1}}
\cdots\overline y_n(m')^{e_{\lambda,n}}\psi_{\lambda}\\
x_{\lambda+1}(m')&=\overline y_{\lambda+1}(m')^{g_{11}(r)}\cdots\overline y_{\lambda+s_r}(m')^{g_{s_rs_r}(r)}\delta_{\lambda+1}\\
&+f_1^{\lambda+1}\overline y_1(m')+\cdots+f_{\lambda}^{\lambda+1}\overline y_{\lambda}(m')\\
&\vdots\\
x_n(m')&=\overline y_n(m')\delta_n+f_1^n\overline y_1(m')+\cdots+f_{\lambda}^n\overline  y_{\lambda}(m')
\endalign
$$
(4.6) becomes
$$
\align
x_1(m')&=\overline y_1(m'+1)^{g_{11}(1)}\cdots \overline y_{s_1}(m'+1)^{g_{1s_1}(1)}
\overline y_{\lambda+1}(m'+1)^{b_{1,\lambda+1}}\cdots \overline y_n(m'+1)^{b_{1n}}u_1\tag 4.9\\
&\vdots\\
x_{\lambda}(m')&=\overline y_{\lambda}(m'+1)\overline y_{\lambda+1}(m'+1)^{b_{\lambda,\lambda+1}}
\cdots \overline y_n(m'+1)^{b_{\lambda n}}u_{\lambda}\\
x_{\lambda+1}(m')&=\overline y_{\lambda+1}(m'+1)^{g_{11}(r)}\cdots 
\overline y_{\lambda+s_r}(m'+1)^{g_{1s_r}(r)}u_{\lambda+1}\\ 
&\vdots\\
x_n(m')&=\overline y_n(m')u_n.
\endalign
$$

The MTS $R(m')\rightarrow R(m'+1)$, where $R(m'+1)$ has regular parameters 
$$
(x_1(m'+1),\ldots,x_n(m'+1))
$$
 defined by
$$
x_i(m')=\cases
x_i(m'+1)x_{\lambda+s_r+1}(m'+1)^{b_{i,\lambda+s_r+1}}\cdots x_n(m'+1)^{b_{in}}&1\le i\le \lambda\\
x_i(m'+1)&\lambda+1\le i\le n
\endcases
$$
factors through $\overline S(m'+1)$, and
$$
\align
x_1(m'+1)&=\overline y_1(m'+1)^{g_{11}(1)}\cdots \overline y_{s_1}(m'+1)^{g_{1s_1}(1)}
\overline y_{\lambda+1}(m'+1)^{b_{1,\lambda+1}}\cdots \overline y_{\lambda+s_r}(m'+1)^{b_{1,\lambda+s_r}}u_1'\tag 4.10\\
&\vdots\\
x_{\lambda}(m'+1)&=\overline y_{\lambda}(m'+1)\overline y_{\lambda+1}(m'+1)^{b_{\lambda,\lambda+1}}
\cdots \overline y_{\lambda+s_r}(m'+1)^{b_{\lambda ,\lambda+s_r}}u_{\lambda}'\\
x_{\lambda+1}(m'+1)&=\overline y_{\lambda+1}(m'+1)^{g_{11}(r)}\cdots 
\overline y_{\lambda+s_r}(m'+1)^{g_{1s_r}(r)}u_{\lambda+1}\\ 
&\vdots\\
x_n(m')&=\overline y_n(m')u_n
\endalign
$$
for some units $u_i'\in \overline S(m'+1)$. Since $\text{det}(g_{ij}(l))=\pm 1$ for $1\le l\le r$, 
we can make a change of variables
in $\overline S(m'+1)$, replacing $\overline y_i(m'+1)$ with a unit times $\overline y_i(m'+1)$ for all $i$ to get that the
$u_i$ and $u_j'$ in (4.10) are 1 for all $i,j$.
Let
$$
(h_{ij})=
\left(\matrix
g_{11}(r) & \cdots & g_{s_11}(r)\\
\vdots&&\vdots\\
g_{s_r1}(r)&\cdots& g_{s_rs_r}(r)\endmatrix\right)^{-1},
$$
an integral matrix.
$$
\align
\overline y_{\lambda+1}(m'+1)&=x_{\lambda+1}(m'+1)^{h_{11}}\cdots x_{\lambda+s_r}(m'+1)^{h_{1,s_r}}\\
&\vdots\\
\overline y_{\lambda+s_r}(m'+1)&=x_{\lambda+1}(m'+1)^{h_{s_r,1}}\cdots x_{\lambda+s_r}(m'+1)^{h_{s_r,s_r}}.
\endalign
$$
$\nu(\overline y_i(m'+1))>0$ for $\lambda+1\le i\le \lambda+s_r$, so by Lemmas 3.2 and 3.3 there exist
MTSs $R(m'+1)\rightarrow R(m'+2)$ and $\overline S(m'+1)\rightarrow \overline S(m'+2)$ along $\nu$
such that $R(m'+2)$ has regular parameters $(x_1(m'+2),\ldots,x_n(m'+2))$, 
 $S(m'+2)$ has regular parameters $(\overline y_1(m'+2),\ldots,\overline y_n(m'+2))$ defined by
$$
x_i(m'+1) = \cases
x_i(m'+2)&1\le i\le \lambda,\\
&\,\,\,\,\, \lambda+s_r<i\le n\\
x_{\lambda+1}(m'+2)^{a_{i1}(m'+2)}\cdots x_{\lambda+s_r}(m'+2)^{a_{i,\lambda+s_r}(m'+2)}
&\lambda+1\le i\le \lambda+s_r
\endcases
$$ 
$$
\overline y_i(m'+1) = \cases
\overline y_i(m'+2)&1\le i\le \lambda,\\
&\,\,\,\,\, \lambda+s_r<i\le n\\
\overline y_{\lambda+1}(m'+2)^{b_{i1}(m'+2)}\cdots \overline y_{\lambda+s_r}(m'+2)^{b_{i,\lambda+s_r}(m+2)}
&\lambda+1\le i\le \lambda+s_r.
\endcases
$$ 
such that $R(m'+2)\subset \overline S(m'+2)$ and
$$
\overline y_{\lambda+i}(m'+1)=x_{\lambda+1}(m'+2)^{e_{i1}}\cdots x_{\lambda+s_r}(m'+2)^{e_{is_r}}
$$
for $1\le i\le s_r$, where $e_{ij}\ge 0$ for all $i,j$. Set 
$$
d_{ij}=e_{1j}b_{i,\lambda+1}+\cdots +e_{s_rj}b_{i,\lambda+s_r}
$$
for $1\le i\le \lambda$, $1\le j\le s_r$. Then the MTS $R(m'+2)\rightarrow R(m'+3)$ where $R(m'+3)$ has regular 
parameters $(x_1(m'+3),\ldots,x_n(m'+3))$ defined by
$$
x_i(m'+2) = \cases 
x_{\lambda+1}(m'+3)^{d_{i1}}\cdots x_{\lambda+s_r}(m'+3)^{d_{is_r}}x_i(m'+3)& 1\le i\le \lambda\\
x_i(m'+3)& \lambda < i\le n
\endcases
$$
factors through $\overline S(m'+2)$ and the conclusions of the Theorem hold for the variables $x_i(m'+3)$ and $\overline y_i(m'+2)$.
\enddemo
\vskip .5truein

\heading
Factorization 1
\endheading
In this chapter we prove Theorem D, which shows that it is possible to factor
a birational map  along a valuation
 by alternating sequences of blowing ups and blowing downs. Theorem 4.4 reduces this to
a question of monomial morphisms and valuations of maximal rational rank. This reduces the problem to a question in
combinatorics.
Christensen, in [Ch], using elementary linear algebra,
 gives a proof, that in dimension 3, factorization holds along a rational rank 3 valuation.
 His algorithm produces
a factorization with one series of blowups and one series of blowdowns. We generalize his methods to   
 give a proof of factorization of monomial mappings in the special case of valuations of maximal rational rank.
 Then Theorem D follows from Theorem 4.4.

\proclaim{Lemma 5.1}
Suppose that $M=(a_{ij})$ is an $n\times n$ matrix such that the $a_{ij}\ge 0$ for all $i,j$ and
$\text{Det}(a_{ij})=\pm1$. Suppose that $R$ is a regular local ring with regular parameters
$(x_1,\ldots,x_n)$. Then there exists a regular local ring $S$ in the quotient field of $R$
such that $S$ has regular parameters $(y_1,\ldots,y_n)$ satisfying (5.1).
$$
\align
x_1&= y_1^{a_{11}}y_2^{a_{12}}\cdots y_n^{a_{1n}}\tag 5.1\\
\vdots&\\
x_n&= y_1^{a_{n1}}y_2^{a_{n2}}\cdots y_n^{a_{nn}}
\endalign
$$

\endproclaim

\demo{Proof} Set $(b_{ij})=M^{-1}$. There exists monomials $f_i$
in $x_1,\ldots,x_n$ for $0\le i\le n$ such that
$x_1^{b_{i1}}\cdots x_n^{b_{in}}=\frac{f_i}{f_0}$ for $1\le i\le n$. In
$R[\frac{f_1}{f_0},\cdots,\frac{f_n}{f_0}]$ we have
$x_i=(\frac{f_1}{f_0})^{a_{i1}}\cdots(\frac{f_n}{f_0})^{a_{in}}$ for
$1\le i\le n$ so that the maximal ideal $m = (x_1,\ldots,x_n,\frac{f_1}{f_0},\ldots,\frac{f_n}{f_0})$
is generated by $\frac{f_1}{f_0},\ldots,\frac{f_n}{f_0}$. Set $S= R[\frac{f_1}{f_0},\cdots,\frac{f_n}{f_0}]_m$
and $y_i=\frac{f_i}{f_0}$ for $1\le i\le n$. Then $S$ is a regular local ring and (5.1) holds.
\endproclaim

Suppose $R\rightarrow S$ is as in (5.1). An inverse monoidal transform (IMT)
 $R\rightarrow S(1) \rightarrow S$ consists of a regular local ring
$S(1)$ such that $R\subset S(1)\subset S$ which has regular parameters $(y_1(1),\ldots,y_n(1))$ such that 
$y_r(1)=y_ry_s$ for some $r\ne s$ and $y_i(1)=y_i$ for $i\ne r$.

\proclaim{Lemma 5.2} Suppose that (5.1) holds for $R\rightarrow S$ and the coefficients of the $s^{th}$ column of
$M$ minus the  $r^{th}$ column of $M$ are nonnegative ($a_{is}-a_{ir}\ge 0$ for all $i$).  Then there exists an IMT
$R\rightarrow S(1) \rightarrow S$ such that 
$$
\align
x_1&=y_1(1)^{a_{aa}(1)}\cdots y_n(1)^{a_{1n}(1)}\tag 5.2\\
&\vdots\\
x_n&=y_1(1)^{a_{n1}(1)}\cdots y_n(1)^{a_{nn}(1)},
\endalign
$$
$M(1)=(a_{ij}(1))$ is $M$ with the $r^{th}$ column subtracted from the $s^{th}$ column.
The adjoint matrix $A(1)$ of $M(1)$ is obtained from the adjoint matrix $A$ of $M$ by adding the
$s^{th}$ row of $A$ to the $r^{th}$ row of $A$.

\demo{Proof}  This follows from  Lemma 5.1. 
\enddemo

 Let $A = (A_{ij})$ be the adjoint matrix of $M$ in (5.1).
Consider a monoidal transform along $\nu$ $S\rightarrow S'$,  where $S'$ has regular parameters
$(y_1',\ldots,y_n')$ defined by 
$$
y_i=\cases
y_s'y_r'&i=r\\
y_i'&i\ne r
\endcases
$$
Of course, this means that $\nu(y_r)>\nu(y_s)$.
Then the matrix $M'=(a_{ij}')$ 
where $x_i=(y_1')^{a_{i1}'}\cdots(y_n')^{a_{in}'}$ for $1\le i\le n$
is obtained from $M$ by adding the $r^{th}$ column
to the $s^{th}$ column. The adjoint matrix of $M'$, $A'=(A_{ij}')$ is obtained from $A$ by subtracting
 the $s^{th}$ row from the $r^{th}$ row.

\proclaim{Theorem 5.3} 
Suppose that $R\subset S$ are excellent regular local rings of dimension $n$, containing a field $k$ of chracteristic 0, with a 
common quotient field $K$. Suppose that $\nu$ is a valuation of $K$ which dominates $S$, with valuation ring $V$.
Suppose that
\item{1)} $V$ has rational rank $n$
\item{2)} $R$ has regular parameters $(x_1,\ldots,x_n)$, $S$ has regular parameters $(y_1,\ldots,y_n)$ such that 
$$
\align
x_1&=  y_1^{a_{11}} y_2^{a_{12}}\cdots  y_n^{a_{1n}}\\
\vdots&\\
x_n&=  y_1^{a_{n1}} y_2^{a_{n2}}\cdots  y_n^{a_{nn}}
\endalign
$$
where $\text{Det}(a_{ij})=\pm1$.

Then there  exists a MTS along $\nu$
$$
S\rightarrow S(1)\rightarrow\cdots\rightarrow S(k)\tag 5.3
$$
where $S(i)$ has regular parameters $(y_1(i),\ldots,y_n(i))$ for $0\le i\le k$ with
$$
\align
x_1&=  y_1(i)^{a_{11}(i)} y_2(i)^{a_{12}(i)}\cdots  y_n(i)^{a_{1n}(i)}\tag 5.4\\
\vdots&\\
x_n&=  y_1(i)^{a_{n1}(i)} y_2(i)^{a_{n2}(i)}\cdots  y_n(i)^{a_{nn}(i)}.
\endalign
$$
such that if $M(k)=(a_{ij}(k))$
 is the coefficient matrix of $R\rightarrow S(k)$, with adjoint matrix $A(k)$, then all but at most
two of $A_{11}(k),A_{12}(k),\ldots,A_{1n}(k)$ are zero.
\endproclaim

\demo{Proof}
Set 
$M= (a_{ij})$. Let $A$ be the adjoint matrix of $M$. In a sequence such as (5.3), define $M(i)=(a_{jk}(i))$
and $A(i)=(A_{jk}(i))$ to be the adjoint matrix of $M(i)$.

We will call a
 monoidal transform $S(l)\rightarrow S(l+1)$ along $\nu$ allowable if it is centered at $P(l)=P_{ij}=(y_i(l),y_j(l))$
where $A_{1i}(l),A_{1j}(l)$ are nonzero and have the same sign.  
If $T\subset \{1,2,\ldots,n\}$ is a subset containing $i$ and $j$, and $P(l)$ is allowable, then
$$
\text{max}\{|A_{1k}(l+1)|: k\in T\}\le \text{max}\{|A_{1k}(l)|: k\in T\}.
$$
 
Suppose that there exists an infinite sequence of allowable monoidal transforms
$$
S\rightarrow S(1)\rightarrow\cdots\rightarrow S(l)\rightarrow\tag 5.5
$$
where $S(l)\rightarrow S(l+1)$ is centered at $P(l)$. We will derive a contradiction.
The Theorem will then follow since at least 
three $A_{1i}(l)$ nonzero imply two of them must have the same sign, which implies that there exists an
allowable monoidal transform.

Set
$$
\align
U(l)&= \{i: A_{1i}(l)\ne0\}\\
\alpha(l) &= |U(l)|\\
T(l) &= \{i : i\text{ occurs as an index in a }P(k)\text{ for some }k\ge l\}\\
\gamma(l)&= |T(l)|\\
\beta(l)&=\text{max}\{|A_{1i}(l)|:i\in T(l)\}\\
W(l) &= \{j\in T(l):|A_{1j}(l)|=\beta(j)\}\\
\delta(l)&=|W(l)|
\endalign
$$

We have $\alpha(l+1)\le \alpha(l)$, $\beta(l+1)\le\beta(l)$, $\gamma(l+1)\le\gamma(l)$
and if $\beta(l+1)=\beta(l)$  then $\delta(l+1)\le \delta(l)$.
Hence in the lexicographic ordering,
$$
(\alpha(l+1),\beta(l+1),\gamma(l+1),\delta(l+1))\le
(\alpha(l),\beta(l),\gamma(l),\delta(l))
$$
for all $l$. 

It suffices to show that this invariant decreases after a finite number of steps, so we may assume that
$$
(\alpha(l),\beta(l),\gamma(l),\delta(l))=(\alpha,\beta,\gamma,\delta)
$$
in (5.5) for all $l$, and derive a contradiction. Set $U=U(l)$, $T=T(l)$, $W=W(l)$.

If there is some $l$ such that $P(l)=P_{rs}$ with $r,s\in W$ and $\nu(y_r(l))>\nu(y_s(l))$, then
$$
A_{1r}(l+1)=A_{1r}(l)-A_{1s}(l)=0,
$$
 and $\alpha(l+1)<\alpha(l)$. This kind of monoidal transform can thus not
occur in (5.5). 

If some $P(l)=P_{ir}$ with $i\in T-W$, $r\in W$ and $\nu(y_r(l))>\nu(y_i(l))$, then
$A_{1r}(l+1)=A_{1r}(l)-A_{1i}(l)$. Hence $\beta(l+1)<\beta(l)$ or
$\beta(l+1)=\beta(l)$ and $\delta(l+1)<\delta(l)$. Thus such a monoidal transform cannot occur in (5.5).

Since $\gamma(l)$ cannot decrease, we must have infinitely many $l$ such that $P(l)=P_{ir}$ with
$r\in W$, $i\in T-W$ and $P(l)=P_{is}$ with $i,s\in T-W$ for all other $l$.

We must thus have $y_j(l)=y_j$ for $j\in W$ and for all $l$. Furthermore, $\nu(y_i(l))\le \nu(y_i)$
for all $i$ and $l$.

At each step
where $P(l)=P_{ir}$ with $r\in W$ and $i\in T-W$ we have
$$
y_i(l+1)=\frac{y_i(l)}{y_r(l)}=\frac{y_i(l)}{y_r}
$$
 and $\nu(y_i(l+1))=\nu(y_i(l))-\nu(y_r)$. After a finite number of
steps we must have $\nu(y_i(l))<0$ for some $i\in T-W$, a contradiction.
\enddemo
\vskip .5truein
When $n=3$, Theorem 5.4 is proved by Christensen [Ch].

\proclaim{Theorem 5.4}
Suppose that $R\subset S$ are excellent regular local rings of dimension $n\ge 3$, containing a field $k$ of chracteristic 0, with a 
common quotient field $K$. Suppose that $\nu$ is a valuation of $K$ which dominates $S$, with valuation ring $V$.
Suppose that
\item{1)} $V$ has rational rank $n$
\item{2)} $R$ has regular parameters $(x_1,\ldots,x_n)$, $S$ has regular parameters $(y_1,\ldots,y_n)$ such that 
$$
\align
x_1&=  y_1^{a_{11}} y_2^{a_{12}}\cdots  y_n^{a_{1n}}\\
\vdots&\\
x_n&=  y_1^{a_{n1}} y_2^{a_{n2}}\cdots  y_n^{a_{nn}}
\endalign
$$
where $\text{Det}(a_{ij})=\pm1$.

Then there is a sequence of regular local rings contained in $K$
$$
\matrix
 &          &    R_1  &             &     &           &        &          &         &          & R_{n-2} &         &\\
 &  \nearrow &        &   \nwarrow  &     &  \nearrow & \cdots & \nwarrow &         & \nearrow &         &\nwarrow &\\
R&          &         &             & S_1 &           &        &          & S_{n-3} &          &         & &S_{n-2}=S\\
\endmatrix
$$
 such that each local ring is dominated by $V$
and each arrow is a sequence of  monoidal transforms (blow ups of regular primes).
Furthermore, we have inclusions
$R\subset S_i$ for all $i$.
\endproclaim

\demo{Proof}
The proof is by induction on $n$. For $n=2$ there is a direct factorization by a MTS. Suppose that $n\ge 3$ and
the theorem is true for smaller values of $n$. We will show that there is a MTS $S\rightarrow S'$ along
$\nu$ and a sequence of IMTs
 $R\rightarrow S''\rightarrow S'$ such that a column of the matrix $M''$ of $R\rightarrow S''$
consists of a single 1 and zeros in the remaining entries. 
Without loss of generality, the first column of $M''$ has this form. By Lemma 5.2, there is then a 
sequence of IMTs $R\rightarrow S_{n-3}\rightarrow S''$ such that the matrix $\overline M$ of $R\rightarrow S_{n-3}$ has 
the form
$$
\overline M =
\left( \matrix 
1&0&\cdots&0\\
0&\overline a_{22}&\cdots&\overline a_{2n}\\
\vdots&\vdots&&\vdots\\
0&\overline a_{n2}&\cdots&\overline a_{nn}
\endmatrix\right).
$$
By induction on $n$, there will then exist a factorization of the desired form.

By Theorem 5.3, there exists a MTS $S\rightarrow S'$ along $\nu$ such that, after possibly interchanging variables,
$A_{1j}=0$ for $j>2$ and
$$
a_{11}A_{11}+a_{12}A_{12}=1\tag 5.6
$$

\subheading{Case 1} Suppose that $A_{11}<0$ and $A_{12}>0$. (The case $A_{12}<0$ and $A_{11}>0$
is similar.) Then
$1 = -a_{11}(-A_{11})+a_{12}A_{12}$. Set
$m=[\frac{-A_{11}}{A_{12}}]$, $n=[\frac{-A_{12}}{A_{11}}]$. 

Suppose that $m>0$. Note that $m=0$ implies $n>0$. 
$$
a_{i1}A_{11}+a_{i2}A_{12}+\cdots+a_{in}A_{1n}=
\cases 1& \text{ if }i=1\\
0 & \text{ if }i\ne 1
\endcases
$$
$$
\align
a_{i2}-a_{i1}m&\ge a_{i2}-a_{i1}(\frac{-A_{11}}{A_{12}})\\
&=\frac{1}{A_{12}}(a_{i1}A_{11}+a_{i2}A_{12}+\cdots+a_{in}A_{1n}).
\endalign
$$
Hence $a_{12}-a_{11}m\ge 1$ and $a_{i2}-a_{i1}m\ge 0$ for $2\le i\le n$.
Let $M'$ be the matrix obtained from $M$ by performing the column operation of
subtracting $m$ times the first column from the second column. All of the coefficients of $M'$ are positive, so by Lemma
5.2 there is an IMT $R\rightarrow S'\rightarrow S$ such that $M'$ is the matrix of $R\rightarrow S'$. 
We have $A_{1j}'=A_{1j}$ if $j\ne 1$ and $A_{11}'=A_{11}+mA_{12}$ so that $A_{11}<A_{11}'\le 0$. If $A_{11}'\ne0$, then
$$
m'=[\frac{-A_{11}'}{A_{12}'}]=[\frac{-A_{11}-mA_{12}}{A_{12}}] = [\frac{-A_{11}}{A_{12}}-m]=0
$$
so that $n'>0$.

Now suppose that $n>0$.
$$
\align
a_{i1}-a_{i2}n&\ge a_{i1}-a_{i2}(\frac{A_{12}}{-A_{11}})\\
&=\frac{1}{A_{11}}(a_{i1}A_{11}+a_{i2}A_{12}+\cdots+a_{in}A_{1n}).
\endalign
$$
Thus we have $a_{i1}-a_{i2}n\ge0$ for $2\le i\le n$. Suppose that $A_{11}\ne -1$. Then
$a_{11}-a_{12}n> -1$, and since this is an integer, $a_{11}-a_{12}n\ge 0$. We can then construct an  
IMT $R\rightarrow S'\rightarrow S$ such that    the matrix $M'$ of $R\rightarrow S'$
is obtained from $M$ by subtracting $n$ times the second column from the first column.
 We have $A_{1j}'=A_{1j}$ if $j\ne 2$ and $A_{12}'=A_{12}+nA_{11}$ so that $A_{12}>A_{12}'\ge 0$.
 If $A_{12}'\ne0$, then
$n'=0$ and $m'>0$ so that we can repeat Case 1.

Suppose that $A_{11}=-1$. $1=-a_{11}+a_{12}A_{12}$ implies $a_{12}>0$. 
$$
a_{i1}-(A_{12}-1)a_{i2} = -(a_{i1}A_{11}+a_{i2}A_{12}+\cdots+a_{in}A_{1n})+a_{i2}
$$
so that $a_{i1}-(A_{12}-1)a_{i2}= a_{i2}\ge 0$
if $i>1$, $a_{11}-(A_{12}-1)a_{12}=-1+a_{12}\ge 0$.
We can then construct an  
IMT $R\rightarrow S'\rightarrow S$ such that    the matrix $M'$ of $R\rightarrow S'$
is obtained from $M$ by subtracting $(A_{12}-1)$ times the second column from the first column.
Now construct the IMT  $R\rightarrow S''\rightarrow S'$ where   the matrix $M''$ of $R\rightarrow S''$
is obtained from $M'$ by subtracting  the first column from the second column. The second column
of $M''$ consists of a 1 in the first row, and the remaining  rows are 0. 

After a finite number of iterations of Case 1 we either prove the induction step, or reach the case 
$A_{12}=0$ or $A_{11}=0$.

\subheading{Case 2} Suppose that $A_{11}=0$ or $A_{12}=0$ (and $A_{13}=\cdots = A_{1n}=0$).
Without loss of generality we may assume that $A_{12}=0$. 
$1=a_{11}A_{11}$ implies $a_{11}=A_{11}=1$. for $i>1$ we have
$a_{i1} = a_{i1}A_{11}+a_{i2}A_{12}+\cdots+a_{in}A_{1n}=0$ so that the first column of $M$ 
consists of a 1 in the first row, and the remaining  rows are 0.

\subheading{Case 3} 
Suppose that $A_{11}>0$ and $A_{12}>0$. Then $a_{11}=A_{11}=1$ and $a_{12}=0$,
or $a_{11}=0$ and $a_{12}=A_{12}=1$. Without loss of generality we have the first case.
For $i>0$ we have
$0\le a_{i1}\le a_{i1}A_{11}+a_{i2}A_{12}+\cdots a_{in}A_{1n}=0$
Hence the first column of $M$ consists of a 1 in the first row, and the remaining  rows are 0.

This completes the induction step for the proof of the Theorem, since the case $A_{11}<0, A_{12}<0$ is not possible.
\enddemo

\subheading{Proof of Theorem D}
We can perform MTSs $R\rightarrow R'$ and $S\rightarrow S'$ so that the conclusions of Theorem 4.4 hold.
We can further replace $R'$ by a MTS $R'\rightarrow R''$ such that $S'$ dominates $R'$, the conclusions of Theorem 4.4 hold,
and if $s_i=2$ for some $i$, then
$$
\align
z_{t_1+\cdots+t_{i-1}+1}&=w_{t_1+\cdots+t_{i-1}+1}\\
z_{t_1+\cdots+t_{i-1}+2}&=w_{t_1+\cdots+t_{i-1}+2}
\endalign
$$
since factorization is possible if $n=2$.
Let
$\lambda_1,\ldots,\lambda_a$ be the $\lambda_i$ such that $1\le \lambda_i\le r$ and $s_{\lambda_i}>2$. Set
$$
\align
x_1&=z_{t_1+\cdots+t_{\lambda_i-1}+1}\\
&\vdots\\
x_{s_{\lambda_i}}&=z_{t_1+\cdots+t_{\lambda_i-1}+s_{\lambda_i}}
\endalign
$$
\vskip .2truein
$$
\align
y_1&=w_{t_1+\cdots+t_{\lambda_i-1}+1}\\
&\vdots\\
y_{s_{\lambda_i}}&=w_{t_1+\cdots+t_{\lambda_i-1}+s_{\lambda_i}}.
\endalign
$$
Set $\overline R_{\lambda_i}=k[x_1,\ldots,x_{s_{\lambda_i}}]_{(x_1,\ldots,x_{s_{\lambda_i}})}$,
 $\overline S_{\lambda_i}=k[y_1,\ldots,y_{s_{\lambda_i}}]_{(y_1,\ldots,y_{s_{\lambda_i}})}$.
Let $\overline K_{\lambda_i}$ be the quotient field of $\overline S_{\lambda_i}$. Then
$\overline R_{\lambda_i}\subset \overline S_{\lambda_i}$ and $\overline V_{\lambda_i}=V\cap \overline K_{\lambda_i}$
is a rank 1, rational rank $s_{\lambda_i}$ valuation ring dominating $\overline S_{\lambda_i}$.
By Theorem 5.4, for all $\lambda_i$, there exist MTSs of regular local rings contained in  $\overline K_{\lambda_i}$,

$$
\matrix
 &          &   (\overline R_{\lambda_i})_1  &    &     &  &&    & (\overline R_{\lambda_i})_{s_{\lambda_i}-2} &         &\\
 &  \nearrow &        &   \nwarrow  &     &  \nearrow &\cdots &  \nearrow &         &\nwarrow &\\
\overline R_{\lambda_i}&   &&  & (\overline S_{\lambda_i})_1 & & 
& & & &(\overline S_{\lambda_i})_{s_{\lambda_i}-2}=\overline S_{\lambda_i}\\
\endmatrix
$$
 such that each local ring is dominated by $\overline V_{\lambda_i}$ and 
$\overline R_{\lambda_i}\subset (\overline S_{\lambda_i})_j$ for all $j$.

We can perform the corresponding sequences of MTSs along $\nu$ on $R$ to construct a sequence of MTSs
$$
\align
&\matrix
 &          &   ( R_{\lambda_1})_1  &    &    &   &&    & ( R_{\lambda_1})_{s_{\lambda_1}-2} &         &\\
 &  \nearrow &        &   \nwarrow  &     &  \nearrow & \cdots          & \nearrow &         &\nwarrow &\\
 R&   &&  & ( S_{\lambda_1})_1 &   &  
& & & &( S_{\lambda_1})_{s_{\lambda_1}-2}= S_{\lambda_1}=R_{\lambda_2}\\
&&&&&&\vdots\\
&          &   ( R_{\lambda_a})_1  &            &        &   &&    & ( R_{\lambda_a})_{s_{\lambda_a}-2} &         &\\
 &  \nearrow &        &   \nwarrow  &     &  \nearrow & \cdots &  \nearrow &         &\nwarrow &\\
 R_{\lambda_a}&   &&  & ( S_{\lambda_a})_1 &   & 
& & & &( S_{\lambda_a})_{s_{\lambda_a}-2}= S\\
\endmatrix
\endalign
$$
$(s_{\lambda_1}-2)+(s_{\lambda_2}-2)+\cdots+(s_{\lambda_a}-2)\le n-2$ since 
$s_{\lambda_1}+\cdots+s_{\lambda_a}\le n$. Thus the conclusions of the Theorem hold.

\heading
Factorization 2
\endheading
In the special case of a monomial mapping, local factorization by one sequence of blowups followed by
one sequence of blowdowns follows from Morelli's 
Theorem on factorization of   birational morphisms of toric varieties [Mo], [AMR].
Theorem 6.1 states this result precisely.

\proclaim{Theorem 6.1} 
 Suppose that  $R,S$ are excellent regular local rings of dimension $n$,
containing a field $k$ of characteristic zero, with a common quotient field $K$, such that $S$ dominates $R$. 
Suppose that $R$ has regular parameters $(x_1,\ldots,x_n)$, $S$ has regular parameters
$(y_1,\ldots,y_n)$ and there exists a matrix $(a_{ij})$ of natural numbers such that $\text{Det}(a_{ij})=\pm1$ and
$$
\align
 x_1&=y_1^{a_{11}}\cdots y_n^{a_{1n}}\tag 6.1\\
&\vdots\\
x_n&=y_1^{a_{n1}}\cdots y_n^{a_{1n}}.
\endalign
$$
Let $V$ be a valuation ring of
$K$ which dominates $S$. Then there exists a regular local ring $T$, with quotient field $K$,
such that $T$ dominates $S$, $V$ dominates  $T$, and the inclusions $R\rightarrow T$ and $S\rightarrow T$
can be factored by sequences of monoidal transforms (blowups of regular primes).
$$
\matrix
&&V&&\\
&&\uparrow&&\\
&&T&&\\
&\nearrow&&\nwarrow&\\
R&&\longrightarrow && S
\endmatrix
$$
\endproclaim
\demo{Proof} With the given assumptions
$$
\text{Spec}(k[y_1,\ldots,y_n])\rightarrow\text{Spec}(k[x_1,\ldots,x_n])\tag 6.2
$$
is a toric birational morphism of toric varieties. There exist projective toric varieties $X$ and $Y$ and
 a birational projective toric morphism $f:X\rightarrow Y$ extending (6.2). 
By the main result of [Mo], [AMR] (Strong factorization of birational toric morphisms)
there exists a factorization
$$
\matrix
&&Z&&\\
&\swarrow&&\searrow&\\
X&&\longrightarrow && Y
\endmatrix
$$
where $Z$ is a projective toric variety, $Z\rightarrow X$ and $Z\rightarrow Y$ are composities of blowups of orbit closures.
$Z\rightarrow X$ and $Z\rightarrow Y$ induce MTSs along $\nu$ $R\rightarrow T$ and $S\rightarrow T$.
\enddemo

\subheading{Proof of Theorem G}

By Theorem A, we can perform sequences of monoidal transforms  $R\rightarrow R_1$ and $S\rightarrow S_1$ 
so that $V$ dominates $S_1$, $S_1$ dominates  $R_1$, and 
$R_1$ and $S_1$ have regular parameters satisfying (6.1). The proof of Theorem G now follows from Theorem 6.1.

\subheading{Proof of Theorem H}
If $K$ is a field containing a ground field $k$, and $v$ is a valuation of $K$, trivial on $k$,
 then the transcendence degree of 
$\Cal O_v/m_v$ over $k$ is called the dimension of $v$ (dim$(v)$). We have
$$
\text{rank}(v)\le \text{rrank}(v)\le\text{trdeg}_kK
$$
(C.f. the Corollary and note at the end of Chapter VI, Section 10 [ZS].)

Suppose that $\nu$ is a valuation associated to $V$. 
By Theorem 1.7, applied to the lift to $V$ of a transcendence basis of $V/m_{\nu}$, there exists a MTS along 
$\nu$, $R\rightarrow R_1$, such
that $\text{dim}_{R_1}(\nu)=0$. 
By assumption, $R_1$ is a localization of $k[f_1,\ldots,f_m]$ for some $f_1,\ldots,f_m\in K$, such that
$\nu(f_i)\ge0$ for all $i$. By Theorem 1.7, there exists a MTS $S\rightarrow S_1$ along $\nu$ such that
$f_1,\ldots,f_n$ are in $S_1$. Hence $S_1$ dominates $R_1$. 
$$
\text{dim}(R_1)=\text{trdeg}_k(K)-\text{trdeg}_k(R_1/m_1)=n-\text{dim}(\nu)
$$
and $\text{dim}(S_1)=n-\text{dim}(\nu)$. Now  the Theorem follows from Theorem G.

\heading
References
\endheading
\item{[Ab1]} Abhyankar, S., Local uniformization on algebraic surfaces over ground fields of
characteristic $p\ne 0$, Annals of Math, 63 (1956), 491-526.
\item{[Ab2]} Abhyankar, S., On the valuations centered in a local domain, Amer. J. Math 78 (1956), 321-348.
\item{[Ab3]} Abhyankar, S., Ramification theoretic methods in algebraic geometry, Princeton 
University Press, 1959.
\item{[Ab4]} Abhyankar, S., Resolution of singularities of embedded algebraic surfaces, Academic Press, New York, 1966.
\item{[Ab5]} Abhyankar, S., Algebraic Geometry for Scientists and Engineers, American Mathematical Society, 1990.
\item{[Ab6]} Abhyankar, S., Simultaneous resolution for algebraic surfaces, Amer. J. Math 78 (1956), 761-790.
\item{[AMR]} Abramovich, D., Matsuki, K., Rashid, S., A note on the factorization theorem of
toric birational maps after Morelli and its toroidal extension, preprint.
\item{[AK]} Akbulut S and King, H, Topology of algebraic sets, MSRI publications 25, Springer-Verlag Berlin.
\item{[Ca]} Cano, F.T., Desingularization strategies for three-dimensional vector fields,
Lecture Notes in Math. 1259, Springer-Vrlag, Heidelberg (1987).
\item{[Ch]} Christensen, C., Strong domination/ weak factorization of three dimensional regular local rings,
Journal of the Indian Math Soc., 45 (1981), 21-47.
\item{[Co]} Cossart, V., Polyedre caracteristique d'une singularite, Thesis, Universite de Paris-Sud, 
Centre d'Orsay (1987).
\item{[Cr]} Crauder, B., Birational morphisms of smooth algebraic 3-folds collapsing 3 surfaces to a point,
Duke Math. J. 48 (1981), 589-632.
\item{[C]} Cutkosky, S.D.,
Local Factorization of Birational Maps, Advances in Math. 132, (1997), 167-315.
\item{[D]} Danilov, V.I., The decomposition of certain birational morphisms, Math USSR Izv. 44 (1980), 465-477.
\item{[D1]} Danilov, V.I., Birational geometry of toric varieties, Math  USSR Izv. 21 (1983), 269-280.
\item{[E]} Ewald, G., Blow-ups of smooth toric 3-varieties, Abh. Math. Sem. Univ. Hamburg 57 (1987), 193-201.
\item{[G}] Giraud, J., Condition de Jung pour les rev\^etements radiciels de hauteur un,
Proc. Algebraic Geometry, Tokyo/Kyoto 1982. Lecture Notes in Math. 1016,
Springer-Verlag (183) 313-333.
\item{[Ha]} Hartshorne, R., Algebraic Geometry, Springer, New York, Heidelberg, Berlin (1977).
\item{[H]} Hironaka, H., Resolution of singularities of an algebraic variety over a field of
characteristic zero, Annals of Math, 79 (1964), 109-326.
\item{[H1]} Hironaka, H., On the theory of birational blowing up, Thesis, Harvard (1960).
\item{[HHS]} Heinzer, W., Huneke, H., Sally J., A criterion for spots, J. Math Kyoto Univ. 26 (1986) 667-671.
\item{[L1]} Lipman, J., Appendix to chapter II of Zariski, O., Algebraic Surfaces, second edition,
Springer Verlag, Berlin, (1971).
\item{[L2]} Lipman, J., Desingularization of two-dimensional schemes, Ann. Math.  107 (1978), 151-207.
\item{[M]} Matsumura, H., Commutative Ring Theory, Cambridge University press, Cambridge (1986).
\item{[Mo]} Morelli, R., The birational geometry of toric varieties, J. Algebraic Geometry 5 (1996)
751-782.
\item{[N]} Nagata, M., Local Rings, Interscience tracts in Pure and Applied Math. 13, J. Wiley, New York (1962).
\item{[P]} Pinkham, H., Factorization of birational maps in dimension three,
Proceedings of Symposia in Pure Mathematics, Vol. 40, Part 2, America Mathematical Society, providence (1983),
343-372.
\item{[S]} Sally, J., Regular overrings of regular local rings, Trans. Amer. Math. Soc. 171 (1972) 291-300.
\item{[Sch]} Schilling, O., The theory of valuations, Amer. Math. Soc, 1950.
\item{[Se]} Seidenberg, A., Reduction of the singularities of the differential equation $Ady=Bdx$, Amer. J. Math. 90
(1968), 248-269.
\item{[Sh]} Shannon, D.L., Monodial transforms, Amer. J. Math, 45 (1973), 284-320.
\item{[W]} Wlodarczyk J., Decomposition of birational toric maps in blowups and blowdowns. Trans. Amer. Math. Soc. 349 (1997),
373-411.
\item{[Z1]} Zariski, O., The reduction of the singularities of an algebraic surface, Annals of Math., 40 (1939) 639-689.
\item{[Z2]} Zariski, O., Local uniformization of algebraic varieties, Annals of Math., 41, (1940), 852-896.
\item{[Z3]} Zariski, O., Reduction of the singularities of algebraic three dimensional varieties, Annals of Math., 45 
(1944) 472-542.
\item{[Z4]} Zariski, O., Introduction to the problem of minimal models in the theory of
algebraic surfaces, Publications of the math. Soc. of Japan, (1958).
\item{[ZS]} Zariski, O. and Samuel, P., Commutative Algebra II, Van Nostrand, Princeton (1960).
\vskip 2truein

Department of Mathematics

University of Missouri

Columbia, Missouri, 65211, USA

\bye